\documentclass[a4paper,12pt,reqno]{amsbook}

\usepackage{graphicx}
\usepackage{amsmath}
\usepackage{amssymb}
\usepackage{amsfonts}
\usepackage{epsfig}
\usepackage{epigraph}
\usepackage{footnpag}

\makeatletter

\newtheorem{theorem}{Theorem}[section]

\newcommand{\OperatorFi}[3]{\ensuremath{
    {#1} \stackrel{\varphi_{#2}}{\longmapsto} {#3} }}
\newcommand{\OperatorFiFj}[5]{\ensuremath{
    {#1} \stackrel{\varphi_{#2}}{\longmapsto} {#3}
               \stackrel{\varphi_{#4}}{\longmapsto} {#5} }}

\renewcommand{\Im}{\text{\rm Im}}

\@addtoreset{equation}{chapter}

\newtheorem{corollary}[theorem]{Corollary}
\newtheorem{remark}[theorem]{Remark}
\newtheorem{proposition}[theorem]{Proposition}
\newtheorem{lemma}[theorem]{Lemma}
\newtheorem{fact}[theorem]{Fact}
\newtheorem{definition}[theorem]{Definition}

\newtheorem{conjecture}[theorem]{Conjecture}
\def\PerfProof{{\it Proof.\ }}
\def\PerfProofW{{\it \indent Proof\ }}
\def\UnionZeroInf{$\bigcup\limits_{n=0}^\infty$}

\makeatother

\makeindex

\begin{document}

\title[Admissible and perfect elements in modular lattices]
{Admissible and perfect elements in modular lattices}
         \author{Rafael Stekolshchik}
         \thanks{{\large email: rs2@biu.013.net.il}}

\date{}

\begin{abstract}
For the modular lattice $D^4 = \{1+1+1+1\}$ associated with the
extended Dynkin diagram $\widetilde{D}_4$ (and also for $D^r$,
where $r > 4$), Gelfand and Ponomarev introduced the notion of
{\it admissible} and {\it perfect} lattice elements and classified
them. In this work, we classify the admissible and perfect
elements in the modular lattice $D^{2,2,2} = \{2+2+2\}$ associated
with the extended Dynkin diagram $\widetilde{E}_6$. Gelfand and
Ponomarev constructed admissible elements for $D^r$ recurrently in
the length of multi-indices, which they called {\it admissible
sequences}. Here we suggest a direct method for creating
admissible elements. Admissible sequences and admissible elements
for $D^{2,2,2}$ (resp. $D^4$) form 14 classes (resp. 11 classes)
and possess some periodicity.

If under all indecomposable representations of a modular lattice
the image of an element is either zero or the whole representation
space, the element is said to be {\it perfect}. Our classification
of perfect elements for $D^{2,2,2}$ is based on the description of
admissible elements. The constructed set $H^+$ of perfect elements
is the union of $64$-element distributive lattices $H^+(n)$, and
$H^+$ is the distributive lattice itself. The lattice of perfect
elements $B^+$ obtained by Gelfand and Ponomarev for $D^4$ can be
imbedded into the lattice of perfect elements $H^+$, associated
with $D^{2,2,2}$. \vspace{5mm}

\end{abstract}

\subjclass{16G20, 06C05, 06B15}

\keywords{Modular lattices, Perfect polynomials, Coxeter Functor}

\maketitle

\newpage
~ \vspace{52mm} \\
~ \\
\begin{center}
 \begin{large}
 \textbf {\it Dedicated to my wife Rosa}
 \end{large}
\end{center}
\newpage

\tableofcontents
 \listoffigures
 \listoftables
\newpage

\chapter{\sc\bf Introduction}
  \label{intro}

\epigraph{\dots These developments, and several others that I have
not mentioned, are a belated validation of Garrett Birkhoff's
vision, which we learned in three editions of his Lattice Theory,
and they betoken Professor Gelfand's oft-repeated prediction that
lattice theory will play a leading role in the mathematics of the
twenty-first century.}{G-C Rota, \cite[p.1445]{R97}.}

\section{Modular and linear lattices, perfect elements}
  \index{modular law(=Dedekind's law)}
  \index{sum}
  \index{intersection}
  \index{inclusion}
  \index{idempotent operation}
  \index{absorbtion law}
  \index{lattice!- modular}
  \index{lattice!- linear}
  \index{tame quivers}
  \index{posets and quivers of tame type}

  The objects considered in this work have diverse applications from quantum
  logic of quantum mechanics (see Remark \ref{Manin_Rota})
  to representation problems of posets and quivers of {\it tame} type
  (see \S\ref{biblio_notes}).
  The potential of these applications is,
  however, still dormant to a considerable extent, even the basic
  notions are not widely known. In Appendix \ref{on_lattices} we give
  definitions and examples of some notions including
  modular and linear lattices, representations of modular lattices and
  perfect elements. In this section we briefly remind a number
  of notions.

  A {\it lattice} is a set $L$ with two commutative and
  associative operations: a {\it sum} and an {\it intersection}.
  If $a, b \in L$, then we denote the intersection
  by $ab$ and the sum by $a + b$.
  Both operations are {\it idempotent}
 $$
    aa = a, \quad a + a = a
 $$
   and satisfy the
  {\it absorbtion law}
 $$
    a(a + b) = a, \quad a + ab = a.
 $$
  On the lattice $L$, an {\it inclusion} $\subseteq$ is defined:
\begin{equation*}
    a \subseteq b \Longleftrightarrow a{b} = a
    \text{ or } a + b = b,
\end{equation*}
see \S\ref{lattice_1}.

   A lattice is said to be {\it modular} if, for every $b, a, c \in L$,
\begin{equation}
  \label{modular_law_1}
       a \subseteq b \Longrightarrow b(a + c) = a + b{c},
\end{equation}
see \S\ref{modular_1}.

Good examples of modular lattices, bringing us some intuition, are
the lattices of subspaces of a given vector space, normal
subgroups of a given group, ideals of a given ring.

\begin{figure}[h]
\includegraphics{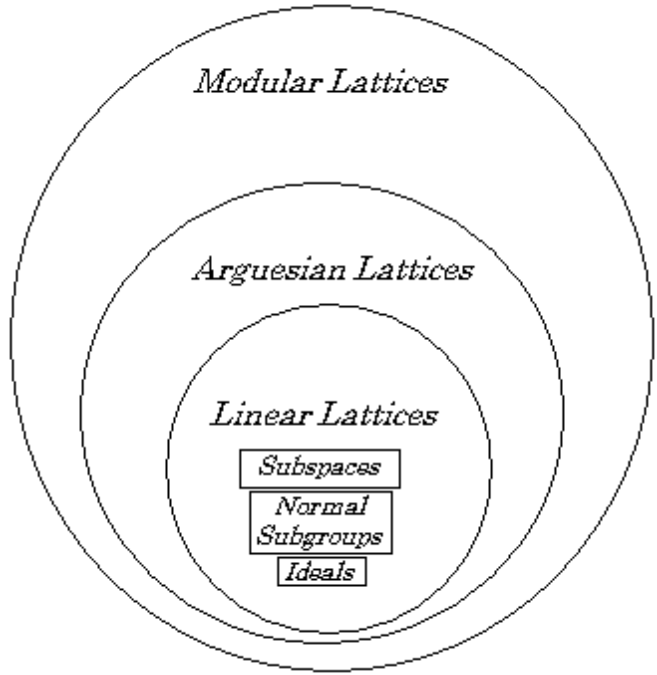}
\caption{\hspace{3mm}Modular, Arguesian and Linear Lattices}
\label{hierar_lattices}
\end{figure}

 \index{commuting equivalence relations}
 \begin{remark} {\rm
Actually, the above examples of modular lattices are more than
modular. In his famous work \cite{Hai85} M.~Haiman baptized these
lattice the {\it linear lattices}. The linear lattices are also
known as lattices of {\it commuting equivalence relations}, see
\S\ref{sect_lin_lat}. The term {\it linear} was suggested by
G.-C.Rota, see \cite[p.1]{Hai85}, \cite{FMR96}.

\index{Desargues's theorem}

 The linear lattices were studied by B.~J\'onsson in
\cite{Jo53}, \cite{Jo54}, where these lattices were named lattices
with {\it representations of type $1$}. In \cite{Jo53}
B.~J\'onsson also introduced the {\it Arguesian law} also studied
by M.~P.~Sch\"{u}tzenberger \cite{Sch45}. The Arguesian law, see
\cite[p.3]{MY2000} and \cite{Hai85}, is a lattice-theoretic form
of the Desargues theorem of projective geometry, see
\cite[Ch.4]{Gr98}. Lattices satisfying the Arguesian law are
called {\it Arguesian lattices}. For more details about linear and
Arguesian lattices, see \S\ref{sect_lin_lat}. On Fig.
\ref{hierar_lattices} the hierarchical
 relations between three classes of lattices: modular, Arguesian,
and linear are shown. }
 \end{remark}

   We denote by $D^r$ the free  modular lattice with $r$ generators
\begin{equation}
  \label{def_Dr}
     D^r = \{ e_1,\dots,e_r \},
\end{equation}
see \cite{GP74}, \cite{Fr80}. The word ``free'' means that no
relation exists between generators  $e_i, i = 1,\dots,r$.

 Denote by $D^{2,2,2}$ the modular lattice with $6$ generators $\{x_1, y_1,
x_2, y_2, x_3, y_3 \}$ satisfying the following relations:
\begin{equation}
  \label{def_D222}
     D^{2,2,2} =
       \{ x_1 \subseteq y_1, \quad
          x_2 \subseteq y_2, \quad
          x_3 \subseteq y_3 \}.
\end{equation}
Symbol $D$ in (\ref{def_Dr}) and $(\ref{def_D222})$ is used in
 honour of R.~Dedekind who derived the main properties of the
modular lattices, \cite{De1897}. The {\it modular law}
(\ref{modular_law_1}) is also known as {\it Dedekind's law} and
the modular lattices are sometimes called the {\it Dedekind
lattices}, see \S\ref{modular_1}.

Let $L$ be a modular lattice, $X$ a finite dimensional vector
space, $\mathcal{L}(X)$ the modular lattice of linear subspaces in
$L$. A morphism  $\rho: L \longrightarrow \mathcal{L}(X)$ is
called a {\it linear representation} of $L$ in the space $X$,
 see \S\ref{repres_lat}. By \S\ref{sect_lin_lat}
 (Proposition \ref{subsp_lin_lat}) the modular lattice $\mathcal{L}(X)$
 is a linear lattice. Thus, we consider representations of the modular
 lattice $L$ in the linear lattices $\mathcal{L}(X)$.

 Following \cite{GP74} we introduce now the notions of {\it perfect
elements} and {\it linear equivalence relation}.
 The element $a \in L$ is said to be {\it perfect} if
 $\rho(a) = 0$ or $\rho(a) = X$ for each indecomposable
representation $\rho$ of the modular lattice $L$, where $X$ is the
representation space of $\rho$, see \S\ref{repres_lat}.

 Two elements $a, b \in L$ are called {\it linearly
equivalent} if $\rho(a) = \rho(b)$ for all indecomposable
representations $\rho: L \longrightarrow \mathcal{L}(X)$ and we
write
\begin{equation}
  \label{lin_equiv_1}
   a \equiv b\mod\theta.
\end{equation}
The relation $\theta$ is called the {\it linear equivalence
relation}, see \S\ref{repres_lat}.

 \index{modular question}
 \index{quantum tautologies}
 \index{quantum logic}
 \index{Haiman's proof theory}
 \index{Birkhoff-von Neumann's logic of quantum mechanics}

\begin{remark}
 \label{Manin_Rota}
 {\rm
Manin \cite[p.98]{Man79} pointed out
 at a connection between {\it perfect elements} and
 {\it quantum tautologies} in the quantum logic  introduced
 by Birkhoff and von Neumann in \cite{BvN36}. In quantum
 logic, a {\it perfect element} is called a
 {\it modular question} \cite{Man79}.
In \cite{KS67}, Kochen and Specker considered related topics in
the foundations of the quantum mechanics and the quantum logic.

G.-C.Rota writes in \cite[p.1444]{R97} that Haiman's proof theory
\cite{Hai85} for linear lattices  brings to fruition the program
that was set forth in \cite{BvN36}: ``The
 authors\footnote{G.~Birkhoff~G. and J.~von Neuman \cite{BvN36}.}
did not know that the modular lattices of quantum mechanics are
linear lattices. In light of Haiman's proof theory, we may now
confidently assert that Birkhoff and von Neumann's logic of
quantum mechanics is indeed the long-awaited new ``logic''\dots''.

Finberg, Mainetti and Rota in \cite{FMR96} gave a simplified
 presentation of the work of Haiman \cite{Hai85}.

For more details about {\it linear lattices}, see
\S\ref{sect_lin_lat}. }
\end{remark}

\section{Outline of admissible elements}
First, in this section, we outline the idea of admissible elements
due to Gelfand and Ponomarev \cite{GP74}. In \cite{GP74},
\cite{GP76}, these elements were used to obtain the distributive
sublattice of perfect polynomials . We think today, that apart
from being helpful in construction of perfect elements, the
admissible elements are interesting in themselves.

 Further, we outline properties of
 admissible elements obtained in this work:
 {\it finite classification}, {\it $\varphi-$homomorphism},
 {\it reduction to atomic elements} and {\it periodicity}, see
 \ref{adm_steps}.

\subsection{The idea of admissible elements of Gelfand and Ponomarev }
  \label{idea_adm}
  \index{Gelfand-Ponomarev!- elementary maps $\varphi_i$}
  \index{Gelfand-Ponomarev!- admissible elements}

Here, we describe, omitting some details, the idea of constructing
the elementary maps and admissible elements of Gelfand-Ponomarev.
As we will see below, the admissible elements {\it grow}, in a
sense, from generators of the lattice or from the lattice's unity.

For every representation $\rho$ of the modular lattice $L$, a new
representation $\tilde\rho$ is uniquely constructed by means of
the Coxeter functor $\Phi^+$, introduced by Bernstein, Gelfand and
 Ponomarev in \cite{BGP73} to study representations of graphs, see
\S\ref{correctness}. Let $X$ (resp. $\tilde{X}$) be a space of the
representation $\rho$ (resp. $\tilde\rho$), and $\mathcal{L}(X)$
(resp. $\mathcal{L}(\tilde{X})$) be the set of subspaces of $X$
(resp. $\tilde{X}$):
\begin{equation}
  \rho\colon{L} \longrightarrow \mathcal{L}(X), \quad
  \tilde\rho\colon{L} \longrightarrow \mathcal{L}(\tilde{X}).
\end{equation}

{\it Construction of the elementary map $\varphi$}.
 Let there exist a map $\varphi$
 mapping every subspace
 $A \in \mathcal{L}(\tilde{X})$ to some subspace
 $B \in \mathcal{L}(X)$:
\begin{equation}
  \label{phi_def_0}
  \begin{array}{c}
    \varphi: \mathcal{L}(\tilde{X}) \longrightarrow \mathcal{L}(X), \vspace{2mm}\\
    \varphi{A} = B, \quad \text{ or } \quad
     A\stackrel{\varphi}{\longmapsto}B.
  \end{array}
\end{equation}
 It turns out that for many elements $a \in L$ there
 exists an element $b \in L$ (at least in $L = D^4$ and $D^{2,2,2}$) such that
\begin{equation}
 \label{phi_rho_ab}
   \varphi\tilde\rho(a) = \rho(b)
\end{equation}
for every pair of representations $(\rho, \tilde{\rho})$, where
$\tilde{\rho} = \Phi^+\rho$. The map $\varphi$ in
(\ref{phi_rho_ab}) is constructed in such a way, that $\varphi$ do
not depend on $\rho$. Thus, we can write
\begin{equation}
  \label{rel_elem_map_0}
   a \stackrel{\varphi}{\longmapsto} b.
\end{equation}
In particular, eq. (\ref{phi_rho_ab}) and (\ref{rel_elem_map_0})
are true for admissible elements, whose definition will be given
later, see \S\ref{adm_seq_polynom}. Eq. (\ref{phi_rho_ab}) is the
main property characterizing the admissible elements.

For $D^4$, Gelfand and Ponomarev constructed $4$ maps of form
(\ref{phi_def_0}): $\varphi_1, \varphi_2, \varphi_3, \varphi_4$
and called them {\it elementary maps}. In this work, for
$D^{2,2,2}$, we construct $3$ maps $\varphi_i$, one per chain $x_i
\subseteq y_i$, see (\S\ref{psi_homom_L6} and
Ch.\ref{section_Coxeter}). Given a sequence of relations
 \index{admissible elements! - growing from generators}
\begin{equation}
  \label{rel_elem_map_1}
   a_{i_1} \stackrel{\varphi_{i_1}}{\longmapsto}
   a_{i_2} \stackrel{\varphi_{i_2}}{\longmapsto} a_{i_3}
               \stackrel{\varphi_{i_3}}{\longmapsto}
   \dots
               \stackrel{\varphi_{i_{n-2}}}{\longmapsto}
   a_{i_{n-1}} \stackrel{\varphi_{i_{n-1}}}{\longmapsto} a_n,
\end{equation}
We say that elements
 $a_2,a_3, \dots, a_{i_{n-2}}, a_{i_{n-1}}, a_{i_n}$
 {\it grow} from the element $a_{i_1}$.
 By abuse of language, we say that elements growing from
 generators $e_i$ or unity $I$ are {\it admissible elements}, and
 the corresponding sequence of indices $\{i_{n}i_{n-1}i_{n-2}\dots{i}_2{i}_1\}$
 is said to be the {\it admissible sequence}.
 For details in the cases $D^{2,2,2}$ (resp. $D^4$)
 see \S\ref{subsec_seq_polynom}, (resp. \S\ref{adm_D4}).

\subsection{Reduction of the admissible elements to atomic elements}
  \label{adm_steps}
Here, we briefly give steps of creation of the admissible elements
in this work.

 \index{key property of admissible sequences}
 \index{admissible sequences! - key property}

 \subsubsection{Finite classification} The admissible elements (and admissible sequences) can be
 reduced to a finite number of classes, see \S\ref{adm_seq_polynom}.
 Admissible sequences for $D^{2,2,2}$ (resp. $D^4$) are depicted
 on Fig. \ref{diagram_123} (resp. Fig. \ref{pyramid_D4}).
 The key property of admissible sequences allowing to do this
 classification is the following relation between maps $\varphi_i$.

  \underline{For $D^{2,2,2}$} (see Proposition \ref{motiv_admis}):
\begin{equation*}
  \begin{split}
    iji & = iki, \\
    \varphi_i\varphi_j\varphi_i & = \varphi_i\varphi_k\varphi_i,
      \text{ where } \{i,j,k\} = \{1,2,3\}.
  \end{split}
\end{equation*}

 \underline{For $D^4$} (see Proposition \ref{motiv_admis_D4}):
\begin{equation*}
  \begin{split}
    ikj & = ilj, \\
    \varphi_i\varphi_k\varphi_j & = \varphi_i\varphi_l\varphi_j,
    \text{ where } \{i,j,k,l\} = \{1,2,3,4\}.
  \end{split}
\end{equation*}

 \subsubsection{$\varphi_i-$homomorphism} In a sense, the elementary maps $\varphi_i$
 are homomorphic with respect to admissible elements. More exactly,
 introduce a notion of $\varphi_i-$homomorphic elements. Let
\begin{equation}
  \label{phi_homom_0}
   a \stackrel{\varphi_i}{\longmapsto} \tilde{a}
    \quad \text{ and }\quad
   p \stackrel{\varphi_i}{\longmapsto} \tilde{p}.
\end{equation}
The element $a \in L$ is said to be {\it $\varphi_i-$homomorphic},
if
\begin{equation}
  \label{phi_homom_1}
   ap \stackrel{\varphi_i}{\longmapsto} \tilde{a}\tilde{p}
\end{equation}
for all $p \in L$, see \S\ref{psi_homom_L6}.

 \index{atomic elements}

  \subsubsection{Reduction to atomic elements} We permanently apply the following mechanism of creation of admissible elements
from more simple elements. Let $\varphi_i$ be any elementary map,
and
\begin{equation}
  \label{adm_chains}
    a \stackrel{\varphi_i}{\longmapsto} \tilde{a}, \quad
    b \stackrel{\varphi_i}{\longmapsto} \tilde{b}, \quad
    c \stackrel{\varphi_i}{\longmapsto} \tilde{c}, \quad
    p \stackrel{\varphi_i}{\longmapsto} \tilde{p}.
\end{equation}
Suppose $abcp$ is any admissible element and
 $a, b, c$ are $\varphi_i-$homomorphic.
By means of relation (\ref{phi_homom_1}) we construct new
admissible element
\begin{equation*}
  \tilde{a}\tilde{b}\tilde{c}\tilde{p}, \quad
\end{equation*}
The elements $a, b, c$ are called {\it atomic}, see
\S\ref{atomic_polynoms}.

 \index{periodicity of admissible elements}
 \index{admissible elements! - periodicity}

\subsubsection{Periodicity} 
\begin{figure}[h]
\includegraphics{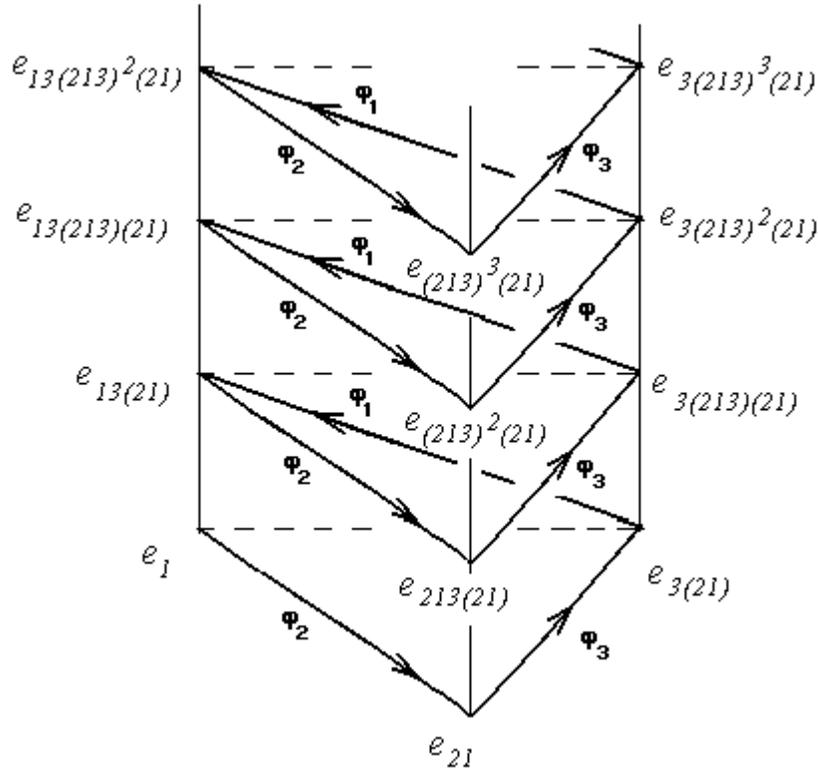}
\caption{\hspace{3mm}Periodicity of admissible elements}
\label{periodicity}
\end{figure}

Admissible sequences obtained by reducing in heading (1)
possess some {\it periodicity}. The corresponding admissible
elements are also periodic. On Fig. \ref{periodicity} we see an
example of admissible elements for $D^{2,2,2}$:
$$
  e_{(213)^p(21)}, \quad e_{3(213)^p(21)}, \quad
  e_{13(213)^p(21)}.
$$
Three vertical lines on Fig. \ref{periodicity} correspond to three
series of inclusions:
\begin{equation}
 \label{period_inclusion}
  \begin{split}
  & \dots \subseteq e_{(213)^{p+1}(21)}
         \subseteq e_{(213)^p(21)}
         \subseteq \dots
         \subseteq e_{(213)(21)}
         \subseteq e_{(21)}, \vspace{2mm} \\
 & \dots \subseteq e_{3(213)^{p+1}(21)}
         \subseteq e_{3(213)^p(21)}
         \subseteq \dots
         \subseteq e_{3(213)(21)}
         \subseteq e_{3(21)}, \vspace{2mm} \\
 & \dots \subseteq e_{13(213)^{p+1}(21)}
         \subseteq e_{13(213)^p(21)}
         \subseteq \dots
         \subseteq e_{13(213)(21)}
         \subseteq e_{13(21)}
         \subseteq e_1, \vspace{2mm} \\
 \end{split}
\end{equation}
Relation (\ref{period_inclusion}) is easily obtained from Table
 \ref{table_alpha}. From the same table we can see periodicity of
indices of atomic elements entering in the decomposition of
admissible elements.
  The triangular helix on Fig. \ref{periodicity} corresponds
to relations
\begin{equation}
  \label{helix_adm}
  \begin{split}
   & e_{(213)^p(21)} \stackrel{\varphi_3}{\longmapsto} e_{3(213)^p(21)}, \\
   & e_{3(213)^p(21)} \stackrel{\varphi_1}{\longmapsto} e_{13(213)^p(21)}, \\
   & e_{13(213)^p(21)} \stackrel{\varphi_2}{\longmapsto} e_{(213)^{p+1}(21)}.
  \end{split}
\end{equation}
Relation (\ref{helix_adm}) is a particular case of the main
theorem on admissible elements (Theorem \ref{th_adm_classes}).

\section{Bibliographic notes: word problems, tame problems}
 \label{biblio_notes}
  \index{perfect element}
  \index{modular lattice!- $D^r$}
  \index{tame problem}
  \index{word problem in $D^r$}
  \index{word problem}
  \index{Gelfand's conjecture}

Admissible elements for $D^r$ in \cite{GP74}, \cite{GP76} are
built recurrently in the length of multi-indices named {\it
admissible sequences}. In this work we suggest a direct method for
creating admissible elements. For $D^{2,2,2}$ (resp. $D^4$), the
admissible sequences and admissible elements form 14 classes
(resp. 11 classes) and possess some periodicity properties, see
\S\ref{subsec_seq_polynom}, Tables \ref{table_admissible},
\ref{table_adm_elems} (resp. \S\ref{adm_D4}, Tables
 \ref{table_admissible_ExtD4}, \ref{table_adm_elem_D4}).

 For the definition of admissible
polynomials due to Gelfand and Ponomarev and examples obtained
from this definition, see \S\ref{coinc_GP_D4}.
 Dlab and Ringel in \cite{DR80} and Cylke in \cite{Cyl82} proved
 that perfect elements\footnote{In the context of this article
 {\it perfect element} is the same as {\it perfect polynomial},
 similarly {\it admissible element} is the same that
 {\it admissible polynomial}} constructed
 by Gelfand and Ponomarev coincide with all possible perfect
 elements modulo linear equivalence, see \S\ref{repres_lat}.

 The study of $D^r$, where $r \geq 4$
 might be a rather hard task. We just point to two works on
 the word problem for the free modular lattices.
 Herrmann proved in \cite{H83} that the word problem in $D^4$ is
 unsolvable, and Freese in \cite{Fr80} proved that the word problem
 in $D^r$ is unsolvable for $r \geq 5$. Furthermore, Herrmann proved
 that
   {\it the modular lattice freely generated by a partially
   ordered set $P$ containing one of lattices
   $$
      1 + 1 + 1 + 1, \quad 2 + 2 + 2, \quad 1 + 2 + 3,
   $$
   has an unsolvable word problem, \cite[p.526]{H83}.}
Thus, the modular lattice $D^{2,2,2}$, which is of interest to us,
 has also an unsolvable word problem. We mention also the following
fact proved by Wille \cite{W73}: {\it a modular lattice is finite
if and only if it does not contain}
   $$
      1 + 1 + 1 + 1 \quad \text{ and } \quad 1 + 2 + 2.
   $$

Gelfand conjectured that the free linear lattice in four
generators has the solvable word problem. If this conjecture is
true, it would distinguish the class of linear lattices from that
of modular lattices, see \cite[p.4]{KY2003}.

 \index{indecomposable quadruples}
 \index{quadruples}
 \index{representation problem of tame type}

 The problem of classifying all indecomposable
representations of $D^4$ ({\it indecomposable quadruples}) of
subspaces of arbitrary finite dimensional spaces is equivalent to
that of determining all indecomposable representations of the
graph
 $$ \begin{array}{ccccc} \\
    \circ & & & & \circ \\
    \vspace{0.5mm} &\searrow & &\swarrow& \\
    \vspace{0.5mm} & & \circ & & \\
    \vspace{0.5mm} &\nearrow & &\nwarrow& \\
    \vspace{0.5mm} \circ & & & & \circ \vspace{1mm} \\
    \end{array}
$$
which corresponds to the extended Dynkin diagram $\tilde D_4$.
This is a representation problem of {\it tame} type solved by
Nazarova \cite{Naz67} and by Gelfand and Ponomarev \cite{GP70} by
different means.

The problem of classifying all indecomposable representations of
$D^{2,2,2}$ is equivalent to that of determining all
indecomposable representations of the graph
$$
 \begin{array}{rcccccccl}
     \circ & \rightarrow
   & \circ & \rightarrow
   & \circ & \leftarrow
   & \circ & \leftarrow  \circ \\
   & & & & \uparrow & & &  \\
   & & & & \circ  \\
   & & & & \uparrow & & &   \\
   & & & & \circ \\
 \end{array}
$$
which corresponds to the extended Dynkin diagram $\tilde{E}_6$.
This is also a representation problem of {\it tame} type solved by
Nazarova \cite{Naz73} and by Dlab and Ringel \cite{DR74},
\cite{DR76}.

 \index{poset of finite type}
 \index{poset of tame type of finite growth}

 Ponomarev considered in detail perfect
 elements for lattices related to the Dynkin diagram $A_n$ \cite{Pon90}.
 For every Dynkin diagram, he described {\it perfect elements} by means of
 so-called {\it hereditary subsets} in the corresponding
 Auslander-Reiten diagram, \cite[Theorem 1.6]{Pon90}.

 Cylke [Cyl93] has constructed a  sublattice $B^+(S) \cup B^-(S)$
 in the modular lattice $L(S)$ freely generated by a poset
 $S$ of finite type, or tame type of finite growth. The sublattice
 $B^+(S) \cup B^-(S)$ is a natural generalization
 of the sublattices $B^+ \cup B^-$ of perfect elements constructed by
 Gelfand-Ponomarev [GP74]. The sublattice
 $B^+(S) \cup B^-(S)$ is built by means of {\it poset differentiation}.
 Nazarova and Roiter used this technique to construct
 indecomposable representations of posets
 and quivers \cite{NR72}, \cite{Naz73}, \cite{Ro85}, \cite{GR92}.

\section{Outline of this work}

  \index{modular lattice!- $D^{2,2,2}$}

 In this work we construct
 {\it admissible} and {\it perfect} elements for the modular lattice
 $D^{2,2,2}$, associated with the extended Dynkin diagram
 $\widetilde{E}_6$, see
 Ch.\ref{section_Coxeter} -- Ch.\ref{section_Perfect_Union},
 and {\it admissible} elements for the modular lattice
 $D^4$, associated with the extended Dynkin diagram
 $\widetilde{D}_4$, see Ch.\ref{sect_adm_seq_D4}.
 For $D^4$, the admissible elements constructed in this work
 and admissible elements constructed by Gelfand and Ponomarev [GP74]
 coincide at least modulo linear equivalence, see \S\ref{coinc_GP_D4}.

 The diagrams $\widetilde{E}_6$
and $\widetilde{D}_4$ are {\it tame quivers}, i.e., their
representations can be classified \cite{Naz73}, \cite{DR76}.

  During the construction of perfect elements in the lattice
 $D^{2,2,2}$ various lattice polynomials in $D^{2,2,2}$ naturally
 appear: {\it atomic}, {\it $\varphi-$homomorphic},
 {\it admissible}, {\it cumulative} and {\it perfect}.
 {\it Atomic polynomials} form a basis for constructing
 {\it admissible elements}. The {\it cumulative polynomials} are constructed
 by means of admissible elements,
 and {\it perfect} elements are constructed
 from the {\it cumulative} elements:
 \begin{equation*}
 {\bf Atomic \hspace{3mm}\Longrightarrow\hspace{3mm} Admissible
             \hspace{3mm}\Longrightarrow\hspace{3mm} Cumulative
             \hspace{3mm}\Longrightarrow\hspace{3mm} Perfect  } \vspace{3mm}
 \end{equation*}

In the sequel to this paper we intend to show that admissible
sequences play an important role in the study of the preprojective
algebras, \cite{GP79}, \cite{DR79}, \cite{Rin96}.

The results of the work were partially published in \cite{St89}.
\\


In Ch.\ref{section_Coxeter} we introduce {\it atomic}, {\it
admissible} and {\it cumulative} lattice polynomials and consider
basic properties of these polynomials. Table \ref{table_adm_elems}
gives a full list of admissible polynomials for $D^{2,2,2}$, see
also \S\ref{adm_seq_polynom}.

 \index{modular lattice!- $D^{2,2,2}$}
 \index{distributive lattice!- $H^+(n)$}
 \index{distributive lattice!- $H^+$}

In Ch.\ref{section_Perfect_Union} we obtain sublattices of perfect
elements $H^+(n)$ in $D^{2,2,2}$ and show that the lattice
$$
  H^+ = \bigcup\limits_{n=0}^\infty{H}^+(n)
$$
is a distributive lattice of perfect elements. The union $H^+$ has
an interesting architecture and contains distributive lattice
$B^+$ of perfect elements of $D^4$, see
 Fig. \ref{cube_comparison}, \S\ref{sublattice_Hn}.

In Ch.\ref{sect_adm_seq_D4}, by repeating the technique developed
for $D^{2,2,2}$ in Ch.\ref{section_Coxeter}, we construct
admissible elements in $D^4$.

 For the convenience of the reader, we recall in
Appendix \ref{on_lattices} some properties of the modular
lattices, their representations, and give a brief introduction to
linear lattices.

Proofs of a number of theorems and verification of several
properties are moved to Appendices \ref{sect_proof_adm} and
 \ref{sect_inclusion_theorem}.

 In this work we consider the
 diagram $\widetilde{E}_6$ only with central orientation in which all arrows are directed
 to the vertex $x_0$.
 The modular lattice $D^{2,2,2}$ is generated by partially ordered set $2+2+2$.
 The generators $x_i, y_i, (i = 1,2,3)$ satisfy the inclusions
 $x_i \subseteq y_i$. For convenience, we include unity $I$
 in the lattice $D^{2,2,2}$:
 \vspace{2mm}
\begin{equation}
\label{Dynkin}
\begin{array}{rcccccccccccccccccl}
     x_1 & \rightarrow
   & y_1 & \rightarrow
   & x_0 & \leftarrow
   & y_2 & \leftarrow
   & x_2 & \hspace{5mm}
   & x_1 & \subseteq
   & y_1 & \subseteq
   & I   & \supseteq
   & y_2 & \supseteq
   & x_2 \\
   & & & & \uparrow & & & & & &
   & & & & \cup\shortmid & & & & \\
   & & \widetilde{E}_6    & & y_3 & & & & & &
   & & D^{2,2,2} & & y_3 & & & & \\
   & & & & \uparrow & & & & & &
   & & & & \cup\shortmid & & & & \\
   & & & & x_3 & & & & & &
   & & & & x_3 & & & & \vspace{2mm}
\end{array}
\end{equation}
 Following Bernstein-Gelfand-Ponomarev \cite{BGP73} and
Gelfand-Ponomarev \cite{GP74} we use the Coxeter functor $\Phi^+$,
and the elementary linear maps $\varphi_i$. \vspace{2mm}

\section{Atomic polynomials}
 \label{atomic_polynoms}

 {\it Atomic} lattice polynomials have a simple lattice
 definition. They are called {\it atomic} because {\it admissible},
 {\it cumulative} and {\it perfect} lattice polynomials
 (as well as {\it invariant} in \cite{St92}) are defined by means of
 these polynomials.

  \index{atomic elements(=polynomials) $a^n_{ij}$, $A^n_{ij}$}

 For $D^{2,2,2}$,
 the definition of $a_n^{ij}$, $A_n^{ij}$,
 $n \in  \{ 0,1,2,\dots \}$,
 is cyclic through the indices $i,j,k$, where the triple $\{{i,j,k}\}$
 is a permutation of $\{1,2,3\}$. We set
\begin{equation}
 \begin{array}{cc}
  \label{eq_atomic}
  & a_n^{ij} =
 \begin{cases}
   a_n^{ij} = I  \qquad \qquad  \hspace{5mm} \text{ for } n = 0,\\
   a_n^{ij} = x_i + y_j{a}_{n-1}^{jk} \hspace{3.5mm}  \text{ for } n \geq 1,
 \end{cases} \vspace{3mm} \\
  & A_n^{ij} =
 \begin{cases}
   A_n^{ij} = I \qquad \qquad  \hspace{4.7mm} \text{ for } n = 0,\\
   A_n^{ij} = y_i + x_j{A}_{n-1}^{ki} \hspace{2.3mm}  \text{ for } n \geq 1,
 \end{cases}
 \end{array}
\end{equation}

The elements $a_n^{ij}$, $A_n^{ij}$ are said to be {\it atomic
elements} in $D^{2,2,2}$: \vspace{2mm}

 {\it Examples of atomic elements in $D^{2,2,2}$.}
 \begin{equation*}
  \begin{array}{lll}
  &  a_1^{12} = x_1 + y_2,
  &  A_1^{12} = y_1 + x_2,  \\
  &  a_2^{12} = x_1 + y_2(x_2 + y_3),
  &  A_2^{12} = y_1 + x_2(y_3 + x_1), \\
  &  a_3^{12} = x_1 + y_2(x_2 + y_3(x_3 + y_1)),
  &  A_3^{12} = y_1 + x_2(y_3 + x_1(y_2 + x_3)), \\
  &  a_4^{12} = x_1 + y_2(x_2 + y_3(y_3 + y_1(x_1 + y_2))),
  &  A_4^{12} = y_1 + x_2(y_3 + x_1(y_2 + x_3(y_1 + x_2))), \\
  &  \qquad \qquad \qquad \qquad \cdots
  &  \qquad \qquad \qquad \qquad \cdots
   \end{array}
 \end{equation*}

For $D^4$, we define {\it atomic} lattice polynomials $a_n^{ij}$,
where $i,j \in \{1,2,3,4\}$, $n \in \mathbb{Z}_{+}$, as follows
\begin{equation}
  \label{eq_atomic_D4}
  a_n^{ij} =
 \begin{cases}
   a_n^{ij} = I \qquad \qquad \qquad \qquad  \hspace{11.9mm} \text{ for } n = 0,\\
   a_n^{ij} = e_i + e_j{a}_{n-1}^{kl} = e_i + e_j{a}_{n-1}^{lk}
       \hspace{2.5mm} \text{ for } n \geq 1,
 \end{cases}
\end{equation}
where $\{i,j,k,l\}$ is the permutation of the quadruple
$\{1,2,3,4\}$, and $\{e_1, e_2, e_3, e_4 \}$ are generators in
$D^4$, see \S\ref{atomic_D4}.

 {\it Examples of the atomic elements in $D^4$}.
 \begin{equation*}
  \begin{array}{lll}
  &  a_1^{12} = e_1 + e_2, \\
  &  a_2^{12} = e_1 + e_2(e_3 + e_4), \\
  &  a_3^{12} = e_1 + e_2(e_3 + e_4(e_1 + e_2)), \\
  &  a_4^{12} = e_1 + e_2(e_2 + e_3(e_1 + e_2(e_3 + e_4))),\\
  &  \qquad \qquad \qquad \qquad \cdots
   \end{array}
 \end{equation*}

\section{The elementary maps $\varphi_i$ of Gelfand-Ponomarev}
 \label{psi_homom_L6}
 \index{elementary maps $\varphi_i$}
 \index{Gelfand-Ponomarev!- elementary maps $\varphi_i$}

Throughout the work the maps $\varphi_i$ play a central role. For
exact definition of the maps $\varphi_i$, see \S\ref{sec_fi}.
Gelfand and Ponomarev introduced the maps $\varphi_i$ in
\cite[p.27]{GP74} for the modular lattices $D^r$ and called them
{\it elementary maps}.

\index{Coxeter functor $\Phi^+$}
\begin{remark}[An explanation of indices]
  \label{rem_index}
{\rm Let $X_0$ be the representation space of the representation
$\rho$, $X^n_0$ be the representation space of the representation
$\rho^{n} = (\Phi^+)^n\rho$, where $n \geq 1$ and $\Phi^+$ is the
Coxeter transformation. The lower index $0$ means that $X^n_0$ is
the image of the generator $x_0$ of the lattice $D^{2,2,2}$ under
the representation $(\Phi^+)^{n}\rho$, see (\ref{Dynkin}):
\begin{equation}
 \begin{split}
      & (\Phi^+)^{n}\rho(x_0) = X^n_0 \text{ for } n \geq 1,
       \text{ where } \rho(x_0) = X_0. \\
      & (\Phi^+)^{n}\rho(y_i) = \hspace{1mm} Y^n_i \text{ for } n \geq 1,
       \text{ where } \rho(y_i) = Y_i,
        (i = 1,2,3). \\
      & (\Phi^+)^{n}\rho(x_i) = X^n_i \text{ for } n \geq 1,
       \text{ where } \rho(x_i) = X_i,
       (i = 1,2,3).
 \end{split}
\end{equation}
The indices $i=1,2,3$ of $\varphi_i$ mean that $\varphi_i$
operates to the subspace $Y_i \subseteq X_0$, where $Y_i$ is the
image of the generator $y_i$ under the representation
$\Phi^+\rho$, see (\ref{repr_rho}): }
\end{remark}

\begin{equation}
\label{repr_rho}
\begin{array}{rcccccccccccccccccl}
     X_1 & \subseteq
   & Y_1 & \subseteq
   & X_0 & \supseteq
   & Y_2 & \supseteq
   & X_2 & \hspace{3mm}
   & X^1_1 & \subseteq
   & Y^1_1 & \subseteq
   & X^1_0 & \supseteq
   & Y^1_2 & \supseteq
   & X^1_2 \\
   & & & & \cup\shortmid & & & & & &
   & & & & \cup\shortmid & & & & \\
   & & \rho      & & Y_3 & & & & & &
   & & \Phi^+\rho& & Y^1_3 & & & & \\
   & & & & \cup\shortmid & & & & & &
   & & & & \cup\shortmid & & & & \\
   & & & & X_3 & & & & & &
   & & & & X^1_3 & & & &
\end{array}
\end{equation}

We introduce now $\varphi_i-$homomorphic polynomials playing an
important role in our further considerations.

\index{$\varphi_i-$homomorphic polynomial}

   An element $a \subseteq D^{2,2,2}$ is said to be
    {\it $\varphi_i-$homomorphic}, if
    \begin{equation}
       \varphi_i\Phi^+\rho(ap) =
          \varphi_i\Phi^+\rho(a)\varphi_i\Phi^+\rho(p)
            \text{ for all } p \subseteq D^{2,2,2}.
    \end{equation}

    An element $a \subseteq D^{2,2,2}$ is said to be
     {\it $(\varphi_i, y_k)-$homomorphic}, if
    \begin{equation}
        \varphi_i\Phi^+\rho(ap) =
          \varphi_i\Phi^+\rho(y_k{a})\varphi_i\Phi^+\rho(p)
            \text{ for all } p \subseteq y_k.
    \end{equation}

Here $\Phi^+\rho(a)$, $\Phi^+\rho(p)$ are the images of $a, p \in
D^{2,2,2}$ under representation $\Phi^+\rho$ (see Remark
\ref{rem_index}). For the detailed definition of $\varphi_i$ and
related notions, see \S\ref{seq_assoc}--\S\ref{sec_fi}.

The notion of $\varphi_i-$homomorphic elements in $D^4$ is
similarly introduced in \S\ref{subs_homom_D4}.

\vspace{3mm} All {\it atomic} polynomials are
$\varphi_i-$homomorphic. More exactly, for $D^{2,2,2}$, we have
 (see Theorem \ref{th_homomorhism}, \S\ref{sect_Homomorphic}):

  1) The polynomials $a^{ij}_n$ are $\varphi_i-$ and $\varphi_j-$homomorphic.

  2) The polynomials $A^{ij}_n$ are $\varphi_i-$homomorphic.

  3) The polynomials $A^{ij}_n$ are ($\varphi_j,y_k)-$homomorphic.
\vspace{2mm}

 For $D^4$, we have (Theorem
\ref{th_homomorhism_D4}, \S\ref{subs_homom_D4}):

  1) The polynomials $a^{ij}_n$ are $\varphi_i-$homomorphic.

  2) The polynomials $a^{ij}_n$ ($\varphi_j,e_k)-$homomorphic.

\section{Admissible sequences and admissible elements}
 \index{admissible sequences!- for $D^{2,2,2}$}
 \label{adm_seq_polynom}

\subsection{Admissible sequences for $D^{2,2,2}$}
  \label{subsec_seq_polynom}

The {\it admissible} lattice polynomials are introduced in
\S\ref{sect_adm_classes}. They appear when we use different maps
$\varphi_i$, where $i=1,2,3$. The admissible polynomials can be
indexed by means of a finite number of types of index sequences,
called {\it admissible sequences}. Consider a finite sequence of
indices $s = i_n\dots{i}_1$, where $i_p \in \{1,2,3\}$. The
sequence $s$ is said to be {\it admissible} if
\begin{equation*}
 \begin{split}
 & \text{(a) Adjacent indices are distinct $(i_p \neq i_{p+1})$ }.
 \\
 & \text{(b) In each subsequence $iji$, we can replace index
 $j$ by $k$. In other words: } \\
 & \dots{iji}\dots = \dots{iki}\dots,
   \text{ where all indices $i,j,k$ are distinct }.
  \end{split}
\end{equation*}

The admissible sequence with $i_1 = 1$ for $D^{2,2,2}$ may be
transformed to one of the seven types (Proposition \ref{list_adm},
see \S\ref{list_admis}),

\begin{equation}
 \label{starting_1}
 \begin{split}
 & \qquad 1) ~(213)^m(21)^n ,
   \qquad 2) ~3(213)^m(21)^n ,
   \qquad 3) ~13(213)^m(21)^n , \\
 & \qquad 4) ~(312)^m(31)^n ,
   \qquad 5) ~2(312)^m(31)^n ,
   \qquad 6) ~12(312)^m(31)^n , \\
 & \qquad 7) ~1(21)^n = 1(31)^n .
  \end{split}
\end{equation}
similarly for $i_1 = 2,3$, for example, for $i_1 = 2$, the
admissible sequences are:
\begin{equation}
 \label{starting_2}
 \begin{split}
 & \qquad 1) ~(123)^m(12)^n ,
   \qquad 2) ~3(123)^m(12)^n ,
   \qquad 3) ~23(123)^m(12)^n , \\
 & \qquad 4) ~(321)^m(32)^n ,
   \qquad 5) ~1(321)^m(32)^n ,
   \qquad 6) ~21(321)^m(32)^n , \\
 & \qquad 7) ~2(12)^n = 2(32)^n .
  \end{split}
\end{equation}

 \index{admissible elements! - $f_\alpha$, $e_\alpha$, $g_{{\alpha}0}$}

\begin{figure}[h]
\includegraphics{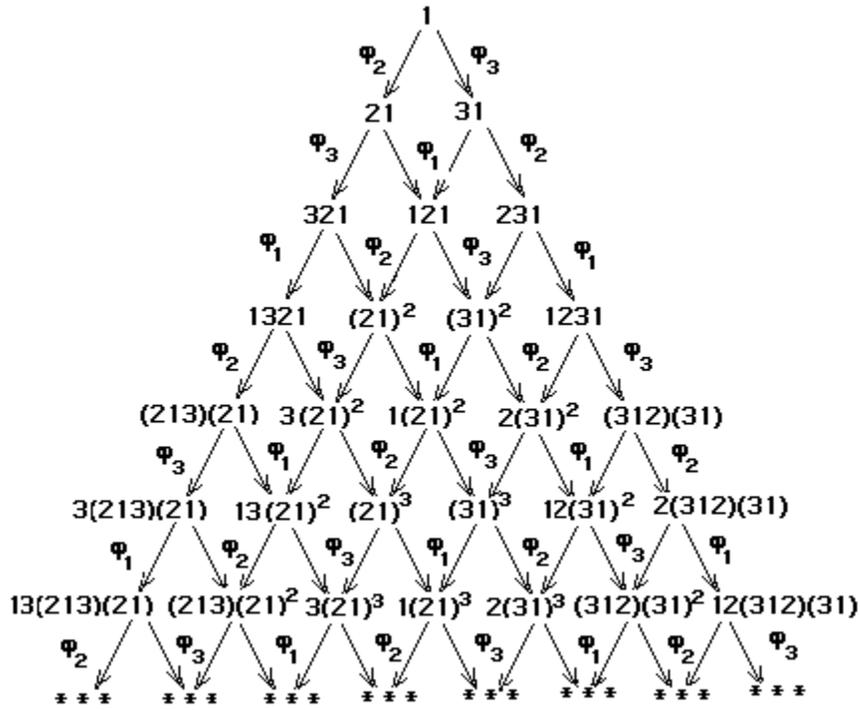}
\caption{\hspace{3mm}How maps $\varphi_i$ act on the admissible sequences}
\label{diagram_123}
\end{figure}

A description of all {\it admissible sequences} for $D^{2,2,2}$ is
given in Table \ref{table_admissible}, and a description of all
{\it admissible elements} is given in Table \ref{table_adm_elems}.
For every
$$
  z_{\alpha} = e_{\alpha}, f_{\alpha}, g_{{\alpha}0}
$$
from Table \ref{table_adm_elems}, where $\alpha$ is an admissible
sequence, the following statement takes place:

Let $\alpha$ be an admissible sequence $i_n\dots{i}_1$ and $i \neq
i_n$. Then $i\alpha$ is also an admissible sequence and
$$
   \varphi_i\Phi^+\rho(z_{\alpha}) = \rho(z_{i\alpha})
   \text{ for each representation $\rho$ }
$$
 (Theorem \ref{th_adm_classes}, \S\ref{sect_adm_classes}; see Fig. \ref{diagram_123}).

\subsection{Admissible sequences for $D^4$}
 \label{adm_D4}
 \index{modular lattice!- $D^4$}
 \index{admissible sequences!- for $D^4$}

The admissible sequence for $D^4$ is defined as follows. Consider
a finite sequence of indices $s = i_n\dots{i}_1$, where $i_p \in
\{1,2,3\}$. The sequence $s$ is said to be {\it admissible} if
\begin{equation*}
 \begin{split}
  & \text{ (a) Adjacent indices are distinct ($i_p \neq i_{p+1}$)}. \\
  & \text{ (b) In each subsequence $ijl$, we can replace index
    $j$ by $k$. In other words: } \\
  & \dots{ijl}\dots = \dots{ikl}\dots , \text{ where all indices } i,j,k,l
 \text{ are distinct}.
 \end{split}
\end{equation*}

Any admissible sequence with $i=1$ for $D^4$ may be transformed to
one the next $11$ types (Proposition \ref{full_adm_seq_D4}, Table
 \ref{table_admissible_ExtD4}):
\begin{equation}
\begin{split}
 & 1) ~(21)^t(41)^r(31)^s = (21)^t(31)^s(41)^r, \\
 & 2) ~(31)^t(41)^r(21)^s = (31)^t(21)^s(41)^r, \\
 & 3) ~(41)^t(31)^r(21)^s = (41)^t(21)^s(31)^r, \\
 & 4) ~1(41)^t(31)^r(21)^s = 1(31)^t(41)^s(21)^r =
 1(21)^t(31)^s(41)^r,\\
 & 5) ~2(41)^r(31)^s(21)^t = 2(31)^{s+1}(41)^{r-1}(21)^t, \\
 & 6) ~3(41)^r(21)^s(31)^t = 3(21)^{s+1}(41)^{r-1}(31)^t, \\
 & 7) ~4(21)^r(31)^s(41)^t = 3(31)^{s+1}(21)^{r-1}(41)^t, \\
 & 8) ~(14)^r(31)^s(21)^t = (14)^r(21)^{t+1}(31)^{s-1} =
       (13)^s(41)^r(21)^t = \\
 &     ~(13)^s(21)^{t+1}(41)^{r-1} = (12)^{t+1}(41)^r(31)^{s-1} =
       (12)^{t+1}(31)^{s-1}(41)^r, \\
 & 9) ~2(14)^r(31)^s(21)^t, \\
 & 10) ~3(14)^r(31)^s(21)^t, \\
 & 11) ~4(14)^r(31)^s(21)^t. \
 \end{split}
 \end{equation}

\begin{figure}[h]
\includegraphics{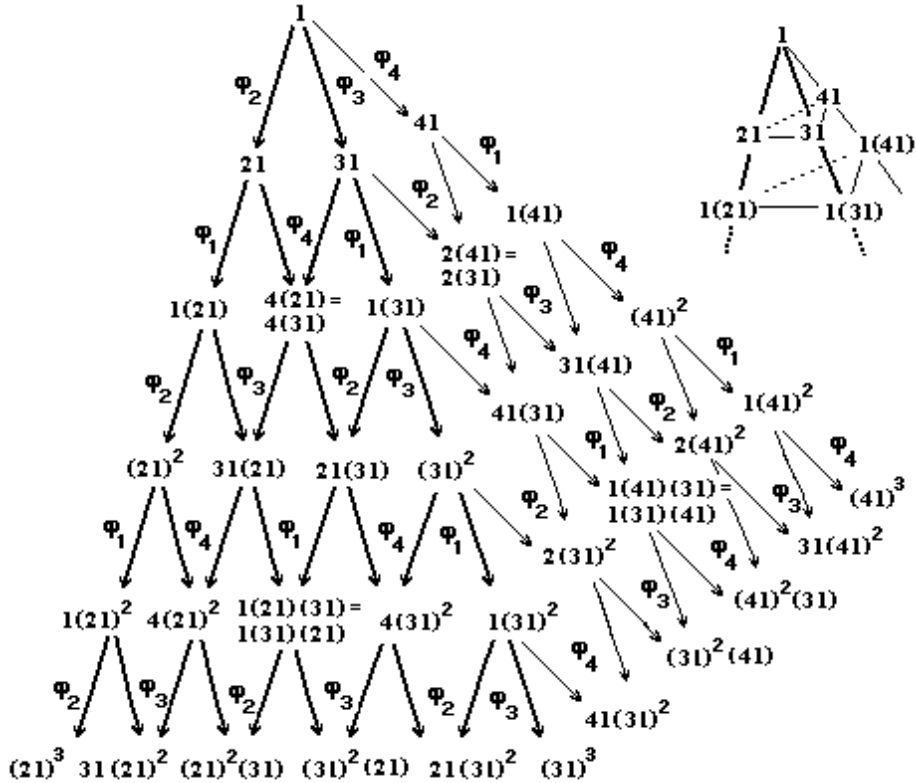}
\caption{\hspace{3mm}Action of maps $\varphi_i$ on the admissible sequences
for ${\it D}^4$}
\label{pyramid_D4}
\end{figure}

Let $\alpha$ be an admissible sequence $i_n\dots{i}_1$ and $i \neq
 i_n$, where $i_s \in \{1,2,3,4\}$. Let $z_\alpha = e_\alpha,
f_{{\alpha}0}$ be admissible elements for $D^4$ from Table
\ref{table_adm_elem_D4}. Then $i\alpha$ is also an admissible
sequence and
\begin{equation}
 \label{main_prop_adm_D4}
  \varphi_i\Phi^+\rho(z_\alpha) = \rho(z_{i\alpha})
  \text{ for each representation } \rho,
\end{equation}
see Theorem \ref{th_adm_classes_D4}, Fig. \ref{pyramid_D4}. The
admissible elements for $D^4$ form a finite family and are
directly constructed, see Table \ref{table_adm_elem_D4}.

Admissible polynomials given by Table \ref{table_adm_elem_D4}
coincide$\mod\theta$\footnote{For definition of$\mod\theta$, see
(\ref{lin_equiv}) from \S\ref{repres_lat}} with polynomials
constructed by Gelfand and Ponomarev. For small lengths of the
admissible sequences, we prove that these polynomials coincide
without restriction$\mod\theta$, see Propositions
\ref{coincidence_GP}, \ref{coincidence_E}, \ref{coincidence_F}
from \S\ref{coinc_GP_D4}.

\begin{remark}
 \index{Grassmann-Cayley algebra}
 \index{Howe-Huang invariants}
{\rm  In \S\ref{GC_algebra} we mention a connection between some
Howe-Huang projective invariants of quadruples \cite{HH96},
 the Grassmann-Cayley algebra \cite{BBR85} and {\it admissible
elements} in $D^4$}.
\end{remark}

\section{Inclusion of admissible elements in the modular lattice}
One can easily check, with definitions from Table
 \ref{table_adm_elems}, that
\begin{equation} \label{incl}
    f_\alpha \subseteq e_\alpha , \quad
    g_{\alpha0} \subseteq e_\alpha.
\end{equation}
It is more difficult to prove the following theorem (see Theorem
\ref{inclusion}, \S\ref{sect_inclusion},
 Appendix \ref{sect_inclusion_theorem}):

\vspace{2mm}
 {\bf Theorem.}
  {\it For every admissible sequence
  ${\alpha}$ and  ${\alpha}i$, where $i = 1,2,3$
  from Table \ref{table_admissible}, we have}
  \begin{equation*}
      e_{{\alpha}i} \subseteq g_{\alpha0}, \quad
         i=1,2,3.  \vspace{3mm}
  \end{equation*}

A proof of the theorem is given in Appendix
\ref{sect_inclusion_theorem}. This proof uses case by case
consideration of {\it admissible sequences} \hspace{1mm}$\alpha$
from Table \ref{table_adm_elems}.

\section{Examples of the admissible elements}

\subsection{The lattice $D^{2,2,2}$}
 \label{exampl_adm}

Let $n$ be the length of an admissible sequence $\alpha$, let
$k$ and $p$ be indices defined by Table \ref{table_adm_elems}. \\

  \underline{For $n=1$}: \indent $\alpha = 1$,
  (Table \ref{table_adm_elems}, Line 7, $k = 0$). We have
\begin{equation*}
  e_1 = y_1, \indent
  f_1 = x_1 \subseteq e_1, \indent
  g_{10} = y_1(y_2 + y_3). \\
\end{equation*}

  \underline{For $n=2$}: \indent $\alpha = 21$,
  (Table \ref{table_adm_elems}, Line 9, $k = 0, p = 0$). Here,
\begin{equation*}
  \begin{split}
  & e_{21} = y_2{a}_1^{13} = y_2(x_1 + y_3), \indent
    e_{12} = y_1(x_2 + y_3)
       \subseteq {g}_{10}, \indent \text{(permutation $(21)
       \rightarrow (12)$)}.\\
  & f_{21} = y_2{y}_3 \subseteq e_{21}, \indent
    g_{210} = e_{21}(x_2 + A_1^{13}) =
      y_2(x_1 + y_3)(x_2 + x_3 + y_1).
 \end{split}
\end{equation*}

  \underline{For $n=3$}: \indent $\alpha = 321$,
  (Table \ref{table_adm_elems}, Line 10, $k = 0, p = 0$). In this
  case:
\begin{equation*}
  \begin{split}
   & e_{321} = y_3{A}_1^{12}{A}_1^{23} = y_3(y_1 + x_2)(y_2 + x_3),\\
   & e_{213} = y_2(y_3 + x_1)(x_1 + x_2) \subseteq g_{210},
           \text{(permutation $(321) \rightarrow (213)$)}, \\
   & f_{321} = x_3{A}_1^{12} = x_3(y_1 + x_2) \subseteq e_{321}, \\
   & g_{3210} = e_{321}(y_1{y}_3 + a_1^{12}) =
           y_3(y_1 + x_2)(y_2 + x_3)(x_1 + y_2 + y_1{y}_3).
 \end{split}
\end{equation*}

\underline{For $n=4$}: \indent $\alpha = 1321$,
  (Table \ref{table_adm_elems}, Line 11, $k = 0 , p = 0$ ). Here, we have
\begin{equation*}
  \begin{split}
  & e_{1321} = y_1{a}_1^{13}{a}_2^{32} =
           y_1(x_1 + y_3)(x_3 + y_2(x_2 + y_1)), \\
  & e_{3213} = y_3(x_3 + y_2)(x_2 + y_1(x_1 + y_3)),
         \text{ (permutation $(1321) \rightarrow (3213)$)}.
 \end{split}
\end{equation*}
By the modular law (\ref{modular_law}) of $D^{2,2,2}$ and since
    $x_1 \subseteq y_1$, we have
\begin{equation*}
\begin{split}
  &  y_1(x_1 + y_3) = x_1 + y_1{y}_3 \hspace{3mm}   \text{ and }  \\
  &  e_{3213} = y_3(x_3 + y_2)(x_2 + x_1 + y_1{y}_3)
         \subseteq {g}_{3210}.
\end{split}
\end{equation*}
Further,
\begin{equation*}
 f_{1321} = y_1{y}_2{a}_1^{13} =
         {y}_1{y}_2({x}_1 + {y}_3)  \subseteq {e}_{1321}
\end{equation*}
and
\begin{equation*}
\begin{split}
 & g_{13210} = e_{1321}(x_1 + a_1^{32}{A}_1^{32}) =
         e_{1321}(x_1 + (x_3 + y_2)(y_3 + x_2)) = \\
 & y_1(x_1 + y_3)(x_3 + y_2(x_2 + y_1))(x_1 + (x_3 + y_2)(y_3 + x_2)).
\end{split}
\end{equation*}
\vspace{2mm}

\subsection{The lattice $D^4$}
  \label{examples_D4}
Admissible sequences and admissible elements are taken from Table
 \ref{table_adm_elem_D4}.

  \underline{For $n=1$}: \indent $\alpha = 1$,
  (Table \ref{table_adm_elem_D4}, Line $G11$, $r = 0, s = 0, t = 0$). We have
\begin{equation*}
 \begin{split}
 & e_1 = e_1 \quad \text{(admissible element $e_1$ coincides with
  generator $e_1$)}, \\
 & f_{10} = e_1(e_2 + e_3 + e_4) \subseteq e_1.
 \end{split}
\end{equation*}

  \underline{For $n=2$}: \indent $\alpha = 21$,
  (Table \ref{table_adm_elem_D4}, Line $G21$, $r = 0,s = 0, t = 1$). We have
\begin{equation*}
 \begin{split}
  e_{21} = &  e_2(e_3 + e_4), \\
  f_{210} = & e_2(e_3 + e_4)(e_4 + e_3(e_1 + e_2) + e_1) = \\
 &         e_2(e_3 + e_4)(e_4 + e_2(e_1 + e_3) + e_1) =  \\
 &         e_2(e_3 + e_4)(e_2(e_1 + e_4) + e_2(e_1 + e_3)) \subseteq e_{21}.
 \end{split}
\end{equation*}

  \underline{For $n=3$}: \indent Consider two admissible sequences:
  $\alpha = 121$ and  $\alpha = 321 = 341$. \vspace{2mm}

  1) $\alpha = 121$,
  (Table \ref{table_adm_elem_D4}, Line $G11$, $r = 0,s = 0, t = 1$). We have
\begin{equation*}
 \begin{split}
   e_{121} = & e_1{a}^{34}_2 = e_1(e_3 + e_4(e_1 + e_2)) =  \\
   & e_1(e_3(e_1 + e_2) + e_4(e_1 + e_2)), \\
  f_{1210} = & e_{121}(e_1{a}^{32}_1 + a^{24}_1{a}^{34}_1) = \\
   & e_{121}(e_1(e_2 + e_3) + (e_2 + e_4)(e_3 + e_4)) =  \\
   & e_1(e_3 + e_4(e_1 + e_2))(e_1(e_2 + e_3) + (e_2 + e_4)(e_3 + e_4))
    \subseteq e_{121}.
 \end{split}
\end{equation*}

  2) $\alpha = 321 = 341$,
  (Table \ref{table_adm_elem_D4}, Line $G31$, $t = 0, s = 0, r = 1$). We have
\begin{equation*}
 \begin{split}
   e_{321} = e_{341} = & e_3{a}^{21}_1{a}^{14}_1 =
     e_3(e_2 + e_1)(e_4 + e_1),  \\
  f_{3210} = f_{3410} = & e_{321}({a}^{14}_2 + e_3a^{24}_1) = \\
   & e_{321}(e_1 + e_4(e_3 + e_2) + e_3(e_2 + e_4)) =  \\
   & e_{321}(e_1 + (e_3 + e_2)(e_4 + e_2)(e_3 + e_4)) = \\
   & e_3(e_2 + e_1)(e_4 + e_1)(e_1 + (e_3 + e_2)(e_4 + e_2)(e_3 + e_4))
    \subseteq e_{321}.
 \end{split}
\end{equation*}

  \underline{For $n=4$}: \indent $\alpha = 2341 = 2321 = 2141$,
  (Table \ref{table_adm_elem_D4}, Line $F21$, $s = 0, t = 1, r = 1$). We have
\begin{equation*}
 \begin{split}
  e_{2141} = & e_2{a}^{41}_2{a}^{34}_1 =
     e_2(e_4 + e_1(e_3 + e_2))(e_3 + e_4), \\
  f_{21410} = & e_{2141}(a^{34}_2 + e_1{a}^{24}_2) = \\
  & e_{2141}(e_2(e_4 + e_1(e_3 + e_2)) + e_1(e_2 + e_4(e_1 +  e_3)))
    \subseteq e_{21}.
 \end{split}
\end{equation*}

\section{Cumulative polynomials}
  \label{def_cumul}

\index{cumulative polynomials $x_t(n)$, $y_t(n)$, $x_0(n)$}

The cumulative polynomials $x_t(n), y_t(n)$, where $t=1,2,3$, and
$x_0(n)$ (all of length $n$) are sums of all admissible elements
of the same length $n$, where $n$ is the length of the
multi-index. In other words, cumulative polynomials are as follows
\begin{equation*}
\begin{split}
& x_t(n) = \sum{f}_{i_{n}\dots{i}_{2}t}, \hspace{5mm} t = 1,2,3, \\
& y_t(n) = \sum{e}_{i_{n}\dots{i}_{2}t}, \hspace{5mm} t = 1,2,3, \\
& x_0(n) = \sum{g}_{i_{n}\dots{i}_{2}0}.
\end{split}
\end{equation*}
From (\ref{incl}) and Theorem \ref{inclusion} we deduce that:
\begin{equation*}
   x_t(n) \subseteq  y_t(n) \subseteq x_0(n), \hspace{5mm} t = 1,2,3.
\end{equation*}
The {\it cumulative polynomials} satisfy the same inclusions as
the corresponding generators in $D^{2,2,2}$.

\section{Examples of the cumulative polynomials}
 \label{examp_cumul}

\underline{For $n=1$}:
\begin{equation*}
  \begin{split}
 & x_1(1) = f_1 = x_1, \\
 & y_1(1) = e_1 = y_1, \\
 & x_0(1) = g_0 = I.
 \end{split}
\end{equation*}
  i.e., the {\it cumulative polynomials} for $n=1$ coincide with
the generators of $D^{2,2,2}$. \vspace{2mm}

\underline{For $n=2$}: \quad
 By the modular law (\ref{modular_law}) and since
        $y_3(x_1 + y_2) \subseteq  x_1 + y_3$,
        we have
\begin{equation*}
 \begin{split}
  x_1(2) = & f_{21} + f_{31} = y_2{y}_3, \indent
        (f_{21} = f_{31} = y_2{y}_3),\\
  y_1(2) = & e_{21} + e_{31} =
       y_2(x_1 + y_3) + y_3(x_1 + y_2) = \\
 & (x_1 + y_3)(y_2 + y_3(x_1 + y_2)) =
   (x_1 + y_3)(y_2 + y_3)(x_1 + y_2), \\
  x_0(2) = & g_{10} +  g_{20} + g_{30} =
    y_1(y_2 + y_3) +  y_2(y_1 + y_3) +  y_3(y_2 + y_3) = \\
 &  (y_2 + y_3)(y_1 + y_3)(y_1 + y_2). \vspace{2mm}
 \end{split}
\end{equation*}

\underline{For $n=3$}: \quad (see Fig. \ref{diagram_123})
\begin{equation*}
\begin{split}
  & x_1(3) = f_{321} + f_{121} + f_{231}, \\
  & y_1(3) = e_{321} + e_{121} + e_{231}.
\end{split}
\end{equation*}
 From Table \ref{table_adm_elems} (Line 8, $k = 0$)
        $f_{121} =  y_1(x_2 + x_3)$
        and  by \S\ref{exampl_adm} it follows that
\begin{equation*}
   x_1(3) = x_3(y_1 + x_2) + y_1(x_2 + x_3) + x_2(y_1 + x_3) =
        (x_2 + x_3)(y_1 + x_2)(y_1 + x_3).
\end{equation*}
 Again, from Table \ref{table_adm_elems} (Line 8, $k = 0$)
        $e_{121} = y_1(x_2 + y_3)(y_2 + x_3)$
        and  by \S\ref{exampl_adm}, we have
\begin{equation*}
 y_1(3) =
        y_3(y_1 + x_2)(y_2 + x_3) +
        y_1(x_2 + y_3)(y_2 + x_3) +
        y_2(y_1 + x_3)(y_3 + x_2). \\
\end{equation*}
Finally,
\begin{equation*}
\begin{split}
 x_0(3) = & g_{210} + g_{310} + g_{320} +
                     g_{120} +  g_{230} + g_{130} = \\
 & \sum{y}_i(x_j + y_k)(y_j + x_i + x_k),
\end{split}
\end{equation*}
where the sum runs over all permutations $\{i, j, k\}$ of
$\{1,2,3\}$.

\section{The sublattice $H^+(n)$ of the perfect elements}
 \label{sublattice_Hn} The perfect elements (see
 \S\ref{repres_lat}) are constructed in the following way:
 \index{perfect element}

\begin{equation}
\label{def_gen_abc}
 \begin{split}
  & \underline{\text{For $n = 0$:}} \\
  & \qquad a_i(0) = y_j + y_k, \\
  & \qquad b_i(0) = x_i + y_j + y_k, \\
  & \qquad c_1(0) = c_2(0) = c_3(0) = \sum{y}_i. \\
  & \underline{\text{For $n \geq 1$:}} \\
  & \qquad a_i(n) = x_j(n) + x_k(n) + y_j(n+1) + y_k(n+1), \\
  & \qquad b_i(n) = a_i(n) + x_i(n+1) =
     x_j(n) + x_k(n) + x_i(n+1) + y_j(n+1) + y_k(n+1),\\
  & \qquad c_i(n) = a_i(n) + y_i(n+1) = x_j(n) + x_k(n) +
       y_i(n+1) + y_j(n+1) + y_k(n+1).
 \end{split}
\end{equation}
\begin{figure}[h]
\centering
\includegraphics{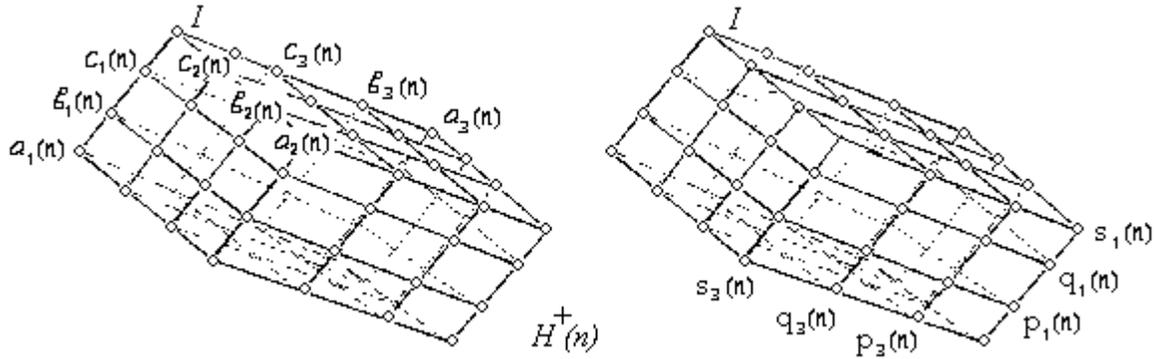}
\caption[\hspace{3mm}The $64$-element distributive lattice
         $H^+(n)$ of perfect elements]
{\hspace{3mm}The $64$-element distributive lattice $H^+(n)$ of
perfect elements.
       Two generator systems:
 $\{ a_i \subseteq b_i \subseteq c_i \}$,
       \S\ref{sublattice_Hn} and
 $\{ p_i \subseteq q_i \subseteq s_i \}$,
       \S\ref{another_gen}.}
\label{cubic64}
\end{figure}
{\bf Proposition} (see Proposition \ref{cumul_polyn},
\S\ref{sect_adm_classes}). {\it If $z_t(n)$ is one of the {\it
cumulative polynomials} $x_t(n)$, $y_t(n)$, where $t = 1,2,3$, or
 $x_0$, then}
\begin{equation*}
  \sum_\text{$i=1,2,3$}\varphi_i\Phi^+\rho(z_t(n)) = \rho(z_t(n+1)).
\end{equation*}

We will also show (Proposition \ref{new_perfect}) that the
property of element $z$ to be {\it perfect} follows from the same
property of element $u$ if
\begin{equation*}
  \sum_\text{$i=1,2,3$}\varphi_1\Phi^+\rho(z) = \rho(u).
\end{equation*}
The elements $z_t(n)$ are not {\it perfect}, but {\it perfect}
elements are expressed as sums of {\it cumulative} polynomials.
\vspace{2mm} \\
{\bf Proposition} (see Proposition \ref{perfect_abc},
\S\ref{chains}). {\it The elements $a_i(n)$, $b_i(n)$, $c_i(n)$
are perfect for all $n$}.

The elements $a_i(n) \subseteq b_i(n) \subseteq c_i(n)$, where $i
= 1,2,3$, generate the $64$-element distributive sublattice
$H^+(n)$ of {\it perfect elements} for $n \geq 1$ (see
 Fig. \ref{cubic64}) and the $27$-element distributive sublattice
$H^+(n)$ of {\it perfect elements} for $n = 0$. The sublattice
$H^+(0)$ contains only $27$ elements because
 $c_1(0) = c_2(0) = c_3(0)$.

\begin{figure}[b]
\includegraphics{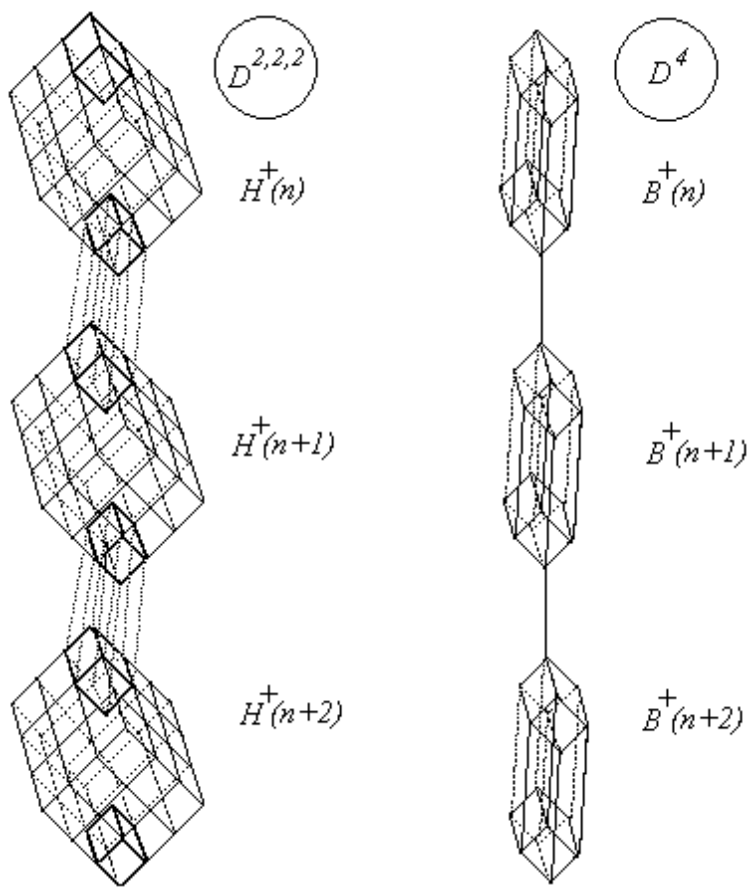}
\caption{\hspace{3mm}Comparison between perfect sublattices
                     in $D^4$ and $D^{2,2,2}$}
\label{cube_comparison}
\end{figure}

\section{The union of sublattices $H^+(n)$}
  \label{union_sublattices}

\subsection{Perfect cubes in the free modular lattice $D^r$}
  \label{cubes}
  \index{modular lattice!- $D^r$}
  \index{distributive lattice!- $B^+$}
  \index{distributive lattice!- $B^+(n)$}
  \index{Boolean cube(=perfect Boolean cube)}

In \cite{GP74}, Gelfand and Ponomarev constructed the sublattice
$B$ of perfect elements for the free modular lattice $D^r$ with
$r$ generators:
\begin{equation*}
   B = B^+ \bigcup B^+, \hspace{3mm} \text{ where }
   B^+ = \bigcup\limits_{n=1}^\infty B^+(n), \hspace{5mm}
   B^- = \bigcup\limits_{n=1}^\infty B^-(n).
\end{equation*}
They proved that every sublattice $B^+(n)$ (resp. $B^-(n)$) is
$2^r$-element Boolean algebra, so-called {\it Boolean cube} (which
can be also named {\it perfect Boolean cube}) and these {\it cubes
} are ordered in the following way. Every element of the cube
$B^+(n)$ is included in every element of the cube $B^+(n+1)$,
i.e.,
\begin{equation*}
 \left .
 \begin{array}{ll}
    v^+(n) \in B^+(n) \\
    v^+(n+1) \in B^+(n+1)
 \end{array}
 \right \}
 \Longrightarrow
 v^+(n+1) \subseteq v^+(n).
\end{equation*}
By analogy, the dual relation holds:
\begin{equation*}
 \left .
 \begin{array}{ll}
    v^-(n) \in B^-(n) \\
    v^-(n+1) \in B^-(n+1)
 \end{array}
 \right \}
 \Longrightarrow
 v^-(n) \subseteq v^-(n+1).
\end{equation*}

\begin{figure}[h]
\centering
\includegraphics{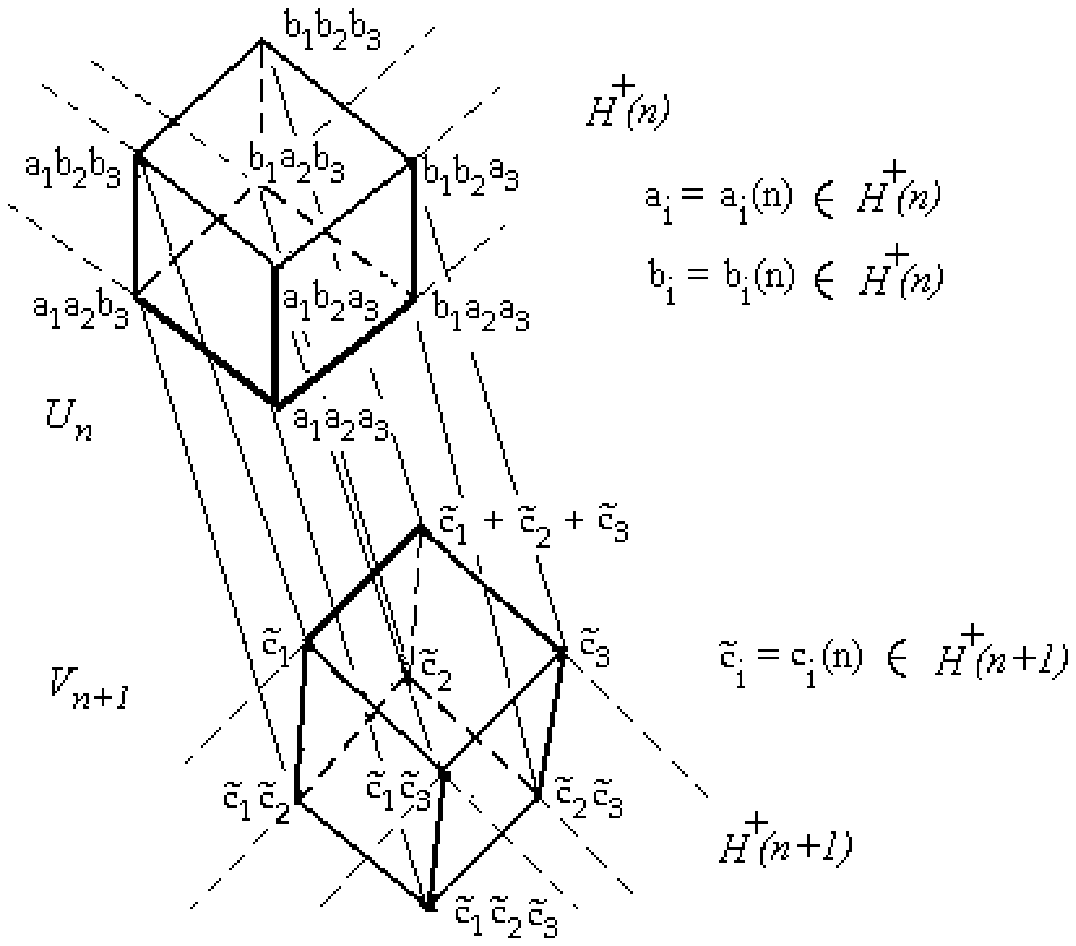}
\caption{\hspace{3mm}The $16$-element Boolean algebra
        $U_n \bigcup V_{n+1}$}
\label{boolean16}
\end{figure}

\subsection{Perfect sublattices in the lattice $D^{2,2,2}$}
  \label{union_L6}
Let us return to the lattice $D^{2,2,2}$. Consider the modular
lattice $H^+$ generated by the union \UnionZeroInf$H^+(n)$ and the
modular lattice $H^-$ generated by the union
 \UnionZeroInf$H^-(n)$, where $H^-(n)$ is the dual modular lattice for
$H^+(n)$. The union \UnionZeroInf$H^+(n)$ is not organized as
simple as \UnionZeroInf$B^+(n)$ in the case of $D^r$. Denote by
$U_n$ the lower $8$-element cube of $H^+(n)$ and by $V_{n+1}$ the
upper $8$-element cube of $H^+(n+1)$, see Fig. \ref{boolean16}.

 {\bf Proposition} (see Proposition \ref{UV_16bool_alg},
\S\ref{sect_bool_alg_16}). {\it The union $U_n \bigcup V_{n+1}$ is
the $16$-element Boolean algebra}.

It is necessary to draw $8$ additional edges in the places where
$H^+(n)$ joins $H^+(n+1)$.

 {\bf Theorem} (see Theorem \ref{main_perf_th},
 \S\ref{sect_merge_constr}).
 {\it The union \UnionZeroInf$H^+(n)$ is a distributive
 lattice$\mod\theta$. The diagram $H^+$ is obtained by uniting the
diagrams $H^+(n)$ and joining the cubes $U_n$ and $V_{n+1}$ for
all $n \geq 0$, i.e., it is necessary to draw $8$ additional edges
for all $n \geq 0$, see Fig. \ref{boolean16}}.

Denote by $H^-$ the lattice {\it dual} to the lattice $H^+$.
\begin{remark}
  {\rm For a comparison between sublattices of the perfect elements
$$
  H^+ \subseteq D^{2,2,2} \text{ and } B^+ \subseteq D^4,
$$
  compare the Hasse diagrams of $H^+ $ and $B^+$ on Fig. \ref{cube_comparison}.
  From Fig. \ref{cube_comparison} we see that the
  perfect sublattice $B^+$ can be injectively mapped (imbedded) into
  the perfect sublattice $H^+\mod\theta$. Indeed,
$$
  p_1(n) + p_2(n) + p_3(n) =  b_1(n)b_2(n)b_3(n)
           \subseteq c_1(n)c_2(n)c_3(n),
$$
see \S\ref{another_gen}, Proposition \ref{two_set_gen}, Table
\ref{two_set_gen}.
 }
\end{remark}

\begin{conjecture} {\rm
  \label{conj_1} The lattice $H^+ \bigcup H^-$ contains all perfect
  elements of $D^{2,2,2}\mod\theta$. }
\end{conjecture}

\chapter[Elementary maps $\varphi_i$, atomic, admissible,cumulative polynomials]
        {\sc\bf Elementary maps $\varphi_i$ and atomic,
        admissible, cumulative polynomials}
  \label{section_Coxeter}

\section{Construction of the Coxeter functor $\Phi^+$
  by means of the lattice operations}   \label{cox_functor}
 \index{Gelfand-Ponomarev!- Coxeter functor $\Phi^+$}
 \index{Coxeter functor $\Phi^+$}

Following Gelfand and Ponomarev \cite{GP74, GP76} we construct the
representation $\Phi^+\rho$ by means of the lattice operations.
Our construction is different from the classic definition of the
Coxeter functor of Bernstein-Gelfand-Ponomarev \cite{BGP73}.

The representation $\Phi^+\rho$ is constructed in the space
$$
   X^1_0 = \{ (\eta_1, \eta_2, \eta_3) \mid \eta_i \in Y_{i},
   \quad \sum\eta_i = 0 \},
$$
where $i$ runs over the set $\{1,2,3\}$. Set
$$
    R = \bigoplus\limits_\text{$i = 1,2,3$}Y_i,
$$
i.e.,
$$
  R = \{ (\eta_1, \eta_2, \eta_3) \mid \eta_i \in Y_i, \quad i = 1,2,3 \}.
$$
Then we see that $X^1_0 \subseteq R$. Hereafter $\eta_i$ denotes a
vector from the space $Y_i$ and $\xi_i$ denotes a vector from the
space $X_i$.

Introduce the spaces $G_i, H_i, G_{i}^{'}, H_{i}^{'} \subseteq R$,
where $i = 1,2,3$:

\begin{equation}
 \label{spaces_G_H}
\begin{array}{lll}
& G_1 = \{(\eta_1, 0, 0) \mid \eta_1 \in Y_1\},
  & H_1 = \{(\xi_1, 0, 0) \mid \xi_1 \in X_1\}, \\
& G_2 = \{(0, \eta_2, 0) \mid \eta_2 \in Y_2\},
  & H_2 = \{(0, \xi_2, 0) \mid \xi_2 \in X_2\}, \\
& G_3 = \{(0, 0, \eta_3) \mid \eta_3 \in Y_3\},
  & H_3 = \{(0, 0, \xi_3) \mid \xi_3 \in X_3\},  \\
& G_{1}^{'} = \{(\xi_1, \eta_2, \eta_3) \mid
           \xi_1 \in X_1,  \eta_i \in Y_i\},
  & H_{1}^{'} = \{(0, \eta_2, \eta_3) \mid
            \eta_i \in Y_i\}, \\
& G_{2}^{'} = \{(\eta_1, \xi_2, \eta_3) \mid
           \xi_2 \in X_2,  \eta_i \in Y_i\},
  & H_{2}^{'} = \{(\eta_1, 0, \eta_3) \mid
            \eta_i \in Y_i\}, \\
& G_{3}^{'} = \{(\eta_1, \eta_2, \xi_3) \mid
           \xi_3 \in X_3,  \eta_i \in Y_i\},
  & H_{3}^{'} = \{(\eta_1, \eta_3, 0) \mid
            \eta_i \in Y_i\}.
\end{array}
\end{equation}
\indent Now, let
\begin{equation}  \label{def_X}
         Y_{i}^1  = G_{i}^{'}X_0^{1},
\end{equation}
\begin{equation}  \label{def_Y}
         X_{i}^1  = H_{i}^{'}X_0^{1},
\end{equation}
\begin{equation}  \label{def_Ro}
         \rho^1 = \Phi^+\rho  =
     [ X_i^1 \subseteq Y_i^1 \subseteq X_0^1 \mid i=1,2,3 ].
\end{equation}
This construction of the Coxeter functor $\Phi^+$
is well-defined as follows from  \S\ref{correctness}.

\section{The associated representations $\nu^0, \nu^1$}
 \label{seq_assoc}
\indent With representations $\rho$ and $\rho^1$
 we associate the following representations $\nu^0, \nu^1$ in $R$:

\begin{equation}  \label{nu0}
         \nu^0(y_i) = X_0^{1} + G_i, \indent
         \nu^0(x_i) = X_0^{1} + H_i, \indent i = 1,2,3,
\end{equation}
\begin{equation} \label{nu1}
         \nu^1(y_i) = X_0^{1}G_i^{'},  \indent
         \nu^1(x_i) = X_0^{1}H_i^{'},  \indent i = 1,2,3.
\end{equation}
    Let $\mu : X_0^{1} \longrightarrow R$ be an inclusion. Then
\begin{equation} \label{nu_muro}
        \nu^1(a) = \mu\rho^1(a)
        \hspace{2mm}\text{ for every } a \in D^{2,2,2}.
\end{equation}
    Let $\nabla : R \longrightarrow X_0$ be a projection,
    $\nabla : (\eta_1, \eta_2, \eta_3)  \longrightarrow
    \Sigma\eta_i$. Then
    $\ker\nabla = X_0^{1}$.
\begin{lemma}
 \label{preserv_intersection}
   For every  $a \in D^{2,2,2}$
   and for every subspace $B \subseteq R$, we have
\begin{equation}
   \label{pr_inters}
    \nabla(\nu^0(a)B) =
    \nabla(\nu^0(a))\nabla(B).
\end{equation}
\end{lemma}
\PerfProof The inclusion
\begin{equation*}
    \nabla(\nu^0(a)B) \subseteq
    \nabla(\nu^0(a))\nabla(B)
\end{equation*}
is clear, because $\nu^0(a)B \subseteq \nu^0(a)$ and $\nu^0(a)B
\subseteq B$. Conversely, let $z \in \nabla(\nu^0(a))\nabla(B)$,
then there exist vectors
\begin{equation}
  \label{eq_nables}
   u \in \nu^0(a) \text{ and } v \in B
   \text{ such that }  \nabla(u) = \nabla(v).
\end{equation}
  According to (\ref{nu0})
\begin{equation*}
  \nu^0(a) = X^1_0 + A
\end{equation*}
 for some subspace $A \subseteq R$.
 From (\ref{eq_nables}), we have
\begin{equation*}
\begin{split}
 & v - u = w \in \ker\nabla = X^1_0,
  \text{ i.e., } \\
 & v = u + w \in  X^1_0 + A = \nu^0(a),
  \text{ and hence }  v \in \nu^0(a)B.
\end{split}
\end{equation*}
Thus, $z = \nabla(v) \in \nabla(\nu^0(a)B)$ and (\ref{pr_inters})
is proved. \qedsymbol \vspace{2mm}

    From (\ref{pr_inters}) it follows that the map
    $\nabla : R \longrightarrow R/X^1_0$ gives a representation
    by the formula
$$
    a \longmapsto \nabla(\nu^0(a)).
$$
We have
\begin{equation}
 \label{nabla_y}
    \nabla(\nu^0(y_i)) = \nabla{G}_i = Y_i = \rho(y_i),
    \indent i = 1,2,3,
\end{equation}
\begin{equation}
 \label{nabla_x}
    \nabla(\nu^0(x_i)) = \nabla{H}_i = X_i = \rho(x_i),
    \indent i = 1,2,3,
\end{equation}
\begin{equation} \label{nablanu_ro}
    \nabla(\nu^0(a))  = \rho(a)
    \text{ for every } a \in D^{2,2,2}.
\end{equation}

\section{The elementary maps $\varphi_i$}
 \label{sec_fi}

 Following \cite{GP74}, we introduce the linear maps
 $$
 \varphi_i : X_0^{1}  \longrightarrow X_0,
  (\eta_1, \eta_2, \eta_3) \longmapsto \eta_i.
 $$
 Let $\pi_i$ be a projection $R$ on $G_i$. Then
\begin{equation} \label{def_fi}
  \varphi_{i} = \nabla\pi_{i}\mu.
\end{equation}
    It follows from the definition that
\begin{equation}  \label{sum_fi}
        \varphi_1  + \varphi_2 + \varphi_3  = 0.
\end{equation}
    From (\ref{sum_fi}) we immediately deduce that
\begin{equation}  \label{f12}
        \varphi_1{B} + \varphi_2{B} =
        \varphi_1{B} + \varphi_3{B} =
        \varphi_1{B} + \varphi_2{B}
\end{equation}
for every subspace $B \subseteq X_0^{1}$.

\section{Basic relations for the maps $\psi_i$ and $\varphi_i$}

We introduce the maps $\psi_i\colon{D}^{2,2,2} \longrightarrow
\mathcal{L}(R)$ useful in the study of maps $\varphi_i$. Set
\begin{equation}  \label{psi}
    \psi_{i}(a) = X_0^1 + G_i(H'_i + \nu^1(a)).
\end{equation}
The map $\psi_i$ is said to be a {\it joint} map.
\begin{proposition}
 \label{basic_eq}
The joint maps $\psi_i$ satisfy the following basic relations:
\begin{enumerate}
 \item $\psi_i(x_i) = X_0^1$, \vspace{2mm}
 \item $\psi_i(x_j) = \psi_k(x_j) = \nu^0(y_i{y}_k)$, \vspace{2mm}
 \item $\psi_i(y_i) =
      \nu^0(x_i(y_j + y_k)) $, \vspace{2mm}
 \item $\psi_i(y_j) = \nu^0(y_i(x_j + y_k))$, \vspace{2mm}
 \item $\psi_i(I) = \nu^0(y_i(y_j + y_k))$, \vspace{2mm}
 \item $\psi_i(y_j{y}_k) = \nu^0(y_i(x_j + x_k))$. \vspace{2mm}
\end{enumerate}
\end{proposition}
\PerfProof 1) From (\ref{nu1}) and (\ref{psi}), we have
$\psi_i(x_i) = X_0^1 + G_i{H'_i}$. Since
       $G_i{H'_i} = 0$, we have $\psi_i(x_i) = X_0^1$.
       \vspace{2mm}

 2) From (\ref{nu1}) and (\ref{psi}), we have $\psi_i(x_j) =
       X_0^1 + G_i(H'_i + X_0^1{H'_j})$.
       Since $H'_j = G_i + G_k$
       (where $i,j,k$ are distinct), it follows that
$$
  \psi_i(x_j) = X_0^1 + G_i(G_j + G_k + X_0^1(G_i + G_k)).
$$
       By the permutation property (\ref{permutation1})
$$
       \psi_i(x_j) = X_0^1 + G_i(X_0^1 + (G_j + G_k)(G_i + G_k)) =
              X_0^1 + G_i(X_0^1 + G_k + G_j(G_i + G_k)).
$$
       Since $G_j(G_i + G_k) = 0$, it follows that
$$
       \psi_i(x_j) = X_0^1 + G_i(X_0^1 + G_k) =
       (X_0^1 + G_i)(X_0^1 + G_k) =
               \nu^0(y_i{y}_k). \vspace{2mm}
$$

 3) As above, $\psi_i(y_i) =
       X_0^1 + G_i(H'_i + X_0^1{G'_i})$.
       Since $H_i^{'} \subseteq G_i^{'}$, and
       by the modular law (\ref{modular_law}), we have
       $\psi_i(y_i) =  X_0^1 + G_i{G'_i}(H'_i + X_0^1)$.
       Since $G_i{G'_i} = H_i$ and
       $H'_j = G_i + G_k$ (where $i,j,k$ are distinct),
       we have
$$
       \psi_i(y_i) = X_0^1 + H'_i(G_j + G_k + X_0^1) =
       (X_0^1 + H'_i)(G_j + G_k + X_0^1) = \nu^0(x_i(y_j + y_k)).
       \vspace{2mm}
$$

 4)  Again, $\psi_i(y_j) = X_0^1 + G_i(H'_i + X_0^1{G'_j})$.
       From $G_i \subseteq G'_j$ for $i \neq j$ and by
       the permutation property (\ref{permutation1}), we have
$$
       \psi_i(y_j) = X_0^1 + G_i(H'_i{G'_j} + X_0^1).
$$
       Since $H'_iG'_j =  H_j + G_k$, we have
$$
       \psi_i(y_j) = X_0^1 + G_i(X_0^1 + H_j + G_k) =
       (X_0^1 + G_i)(X_0^1 + H_j + G_k) = \nu^0(y_i(x_j + y_k)).
       \vspace{2mm}
$$

     5) Similarly,
\begin{equation*}
\begin{split}
   &    \psi_i(I) = X_0^1 + G_i(H'_i + X_0^1) =
       X_0^1 + G_i(G_j + G_k + X_0^1) = \\
   &       (X_0^1 + G_i)(G_j + G_k + X_0^1) = \nu^0(y_i(y_j + y_k)).
   \vspace{2mm}
\end{split}
\end{equation*}

    6) As above,
$$
   \psi_i(y_j{y}_k) = X_0^1 + G_i(H'_i + X_0^1{G'_j}{G'_k}).
$$
       From $G'_j{G'_k} = H_j + H_k + G_i$
       and by the permutation (\ref{permutation1}), we have
\begin{equation*}
\begin{split}
    &   \psi_i(y_j{y}_k) =
       X_0^1 + G_i((G_j + G_k)(H_j + H_k + G_i) + X_0^1) = \\
    &   X_0^1 + G_i(H_j + H_k + G_i(G_j + G_k) + X_0^1).
\end{split}
\end{equation*}
       The intersection vanishes: $G_i(G_j + G_k) = 0$,
       and hence
$$
     \psi_i(y_j{y}_k) =
       X_0^1 + G_i(H_j + H_k + X_0^1) = \nu^0(y_i(x_j + x_k)).
       \qed \vspace{2mm}
$$

\begin{remark}[The lattice description of the projection, \cite{GP74}]
  {\rm If $\pi\colon{R} \longrightarrow R$ is a projection,
  and $A$ is a subspace of $R$,
  then
\begin{equation} \label{projection}
    \pi{A} = \Im\pi(\ker\pi + A).
\end{equation}
Indeed,
$$
   v \in \Im\pi(\ker\pi + A) \Longrightarrow
   v = \pi{v} \in \pi(\ker\pi + A), \indent
   \text{ i.e., }  v \in \pi{A}.
$$
Conversely,
\begin{equation*}
    v \in \pi{A} \Longrightarrow {v} \in \Im\pi, \text{ besides }
    v = \pi{v} \Longrightarrow \pi{v} \in \pi{A} \Longrightarrow  v \in \ker\pi + A.
    \qed \vspace{2mm}
\end{equation*}
}
\end{remark}

\begin{proposition}[A relation between the $\varphi_i$ and $\psi_i$]
    The following properties establish a relation
    between the $\varphi_i$ and $\psi_i$:
    \label{phi_and_psi}
    \begin{enumerate}
    \item
    Let $a,b \subseteq D^{2,2,2}$.
    If $\psi_i(a) = \nu^0(b)$, then
    $\varphi_i\Phi^+\rho(a) = \rho(b)$. \vspace{2mm}
    \item
    Let $c \subseteq D^{2,2,2}$. Then
     $\nabla\psi_i(c) = \varphi_i\Phi^+\rho(c)$. \vspace{2mm}
    \item
    Let $a,b,c \subseteq D^{2,2,2}$.
    If $\psi_i(a) = \nu^0(b)$ and
    $\psi_i(ac) = \psi_i(a)\psi_i(c)$, then \vspace{2mm}
\begin{equation}
    \varphi_i\Phi^+\rho(ac) = \varphi_i\Phi^+\rho(a)\varphi_i\Phi^+\rho(c).
\end{equation}
\end{enumerate}
\item
\end{proposition}
\PerfProof
     (1) From (\ref{def_fi}) we deduce that
       $\varphi_i\Phi^+\rho(a) = \nabla\pi_i\mu\Phi^+\rho(a)$.
    By (\ref{nu_muro}), we have
       $\varphi_i\Phi^+\rho(a) = \nabla\pi_i\nu^1(a)$.
    By (\ref{projection}) for
    $\Im\pi_i = G_i$ and $\ker\pi_i = H_i^{'}$, we have
$$
    \pi_i\nu^1(a) = G_i(H_i^{'} + \nu^1(a)),
$$
    and
$$
     \nabla\pi_i\nu^1(a) =
       \nabla(G_i(H_i^{'} + \nu^1(a))) =
       \nabla(X_0^1 + G_i(H_i^{'} + \nu^1(a))) =
       \nabla(\psi_i(a)) = \nabla(\nu^0(b)).
$$
       From (\ref{nablanu_ro}) $\nabla(\psi_i(a) = \rho(b)$.
       \vspace{2mm}

 (2) By definition (\ref{psi}) of $\psi_i$, we have
$$
      \nabla\psi_i(c) =
         \nabla(X_0^1 + G_i(H_i{'} + \nu^1(c))).
$$
      Since $X_0^1 \subseteq \ker\nabla$, we have
      $\nabla\psi_i(c) = \nabla(G_i(H_i{'} + \nu^1(c))$.
      By lattice description of the projection (\ref{projection}), we have
      $\nabla\psi_i(c) = \nabla(\pi_i\nu^1(c))$.
      By (\ref{nu_muro})
         $\nabla\psi_i(c) = \nabla(\pi_i\mu\Phi^+\rho(c))$.
      Finally, by definition (\ref{def_fi}) of $\varphi_i$,
      we have
        $\nabla\psi_i(c) = \varphi_i\Phi^+\rho(c)$.
      \vspace{2mm}

 (3) By the hypothesis
      $\psi_i(ac) = \nu^0(a)\psi_i(c)$.
      By Lemma \ref{preserv_intersection} and
      eq.(\ref{pr_inters}), we have
$$
        \nabla\psi_i(ac) =
             \nabla\psi_i(a)\nabla\psi_i(c).
$$
      The proof follows from (2).
\qedsymbol \vspace{2mm} \\

From Propositions \ref{basic_eq}, \ref{phi_and_psi} we have

\begin{corollary}
 \label{cor_psi}
 For the elementary map $\varphi_i$, the following basic relations hold:
\begin{enumerate}
 \item $\varphi_i\Phi^+\rho(x_i) = 0$, \vspace{2mm}
 \item $\varphi_i\Phi^+\rho(x_j) = \varphi_k\Phi^+\rho(x_j) =
               \rho(y_i{y}_k)$, \vspace{2mm}
 \item $\varphi_i\Phi^+\rho(y_i) = \rho(x_i(y_j + y_k))$, \vspace{2mm}
 \item $\varphi_i\Phi^+\rho(y_j) = \rho(y_i(x_j + y_k))$, \vspace{2mm}
 \item $\varphi_i\Phi^+\rho(I) = \rho(y_i(y_j + y_k))$, \vspace{2mm}
 \item $\varphi_i\Phi^+\rho(y_j{y}_k) = \rho(y_i(x_j + x_k))$. \vspace{2mm}
 \end{enumerate}
\end{corollary}
\PerfProof
  It is necessary to explain only (1). Since the $i${th}
   coordinate in the space $H_i^{'}$ is $0$, it follows that
$$
  \varphi_i\Phi^+\rho(x_i) =  \nabla\pi_i\mu\Phi^+\rho(x_i) =
  \nabla\pi_i\nu^1(x_i) =
  \nabla\pi_i(X_0^1{H}_i^{'}) = 0.
  \qed
$$

\section{Additivity and multiplicativity of the joint maps $\psi_i$}
\begin{proposition}
  \label{additivity}
   The map $\psi_i$ is additive and
   quasimultiplicative\footnote{The map $\psi_i$ is multiplicative if
      $\psi_i(a)\psi_i(b) = \psi_i(ab)$. We will say that
      $\psi_i$ is {\it quasimultiplicative} if element $ab$
      in this relation is a little deformed, or more exactly, if $ab$
      is deformed to $(a + x_i)(b + x_j{x}_k)$ or to
      $a(b + x_i + x_j{x}_k)$ as in (2) and (3) of Proposition
      \ref{additivity}.}
   with respect to the
   lattice operations + and $\cap$. The following relations hold
   for every $a,b \in D^{2,2,2}$:
\begin{enumerate}
   \item $\psi_i(a) + \psi_i(b) = \psi_i(a+b)$, \vspace{2mm}
   \item $\psi_i(a)\psi_i(b)$ =
      $\psi_i((a + x_i)(b + x_j{x}_k))$, \vspace{2mm}
   \item $\psi_i(a)\psi_i(b)$ =
      $\psi_i(a(b + x_i + x_j{x}_k))$. \vspace{2mm}
\end{enumerate}
\end{proposition}
\PerfProof
   1) By the modular law (\ref{modular_law})
\begin{equation*}
\begin{split}
 &  \psi_i(a) + \psi_i(b) =
   X_0^1 + G_i(H'_i + \nu^1(a)) +  G_i(H'_i + \nu^1(b)) = \\
 & X_0^1 + G_i
   \left (H'_i + \nu^1(a) + G_i(H'_i + \nu^1(b)) \right ) =
   X_0^1 + G_i \left ((H'_i + \nu^1(b))(H'_i + G_i) + \nu^1(a) \right ).
\end{split}
\end{equation*}
   Since $H_i^{'} + G_i = R$, it follows that
$$
   \psi_i(a) + \psi_i(b) =
       X_0^1 + G_i \left (H'_i + \nu^1(b) + \nu^1(a) \right ) =
       X_0^1 + G_i(H'_i + \nu^1(b + a) ).
     \vspace{2mm}
$$

2) By definition
   (\ref{nu1}) $\nu^1(b) \subseteq X_0^1$, and
   by the permutation property (\ref{permutation2}) we have
$$
   X_0^1 + G_i(H'_i + \nu^1(a)) =
   X_0^1 + H'_i(G_i + \nu^1(a)).
$$
   By the modular law (\ref{modular_law}) and by (\ref{permutation2})
   we have
\begin{equation*}
\begin{split}
   & \psi_i(a)\psi_i(b) = \\
   & \left (X_0^1 + G_i(H'_i + \nu^1(a)) \right )
     \left (X_0^1 + G_i(H'_i + \nu^1(b)) \right ) = \\
   & X_0^1 + G_i(H'_i + \nu^1(a))
     \left (X_0^1 + H'_i(G_i + \nu^1(b)) \right ).
\end{split}
\end{equation*}
   Further,
$$
   \psi_i(a)\psi_i(b) =
      X_0^1 + G_i \left (X_0^1(H'_i + \nu^1(a)) +
         H'_i(G_i + \nu^1(b)) \right ).
$$
   Since
$$
    X_0^1(H'_i + \nu^1(a)) = X_0^1{H'_i} + \nu^1(a)
    \indent \text{ and } \indent X_0^1{H'_i} = \nu^1(x_i),
$$
   we have
 \begin{equation} \label{mul}
     \psi_i(a)\psi_i(b) =
        X_0^1 + G_i \left (\nu^1(x_i) + \nu^1(a) +
          H'_i(G_i + \nu^1(b)) \right ).
 \end{equation}
    By the permutation property (\ref{permutation1}) and by (\ref{mul}) we have
 \begin{equation} \label{mul1}
     \psi_i(a)\psi_i(b) =
       X_0^1 + G_i \left (H'_i +
          (\nu^1(x_i) + \nu^1(a))(G_i + \nu^1(b)) \right ).
 \end{equation}
    Since
$$
   G_i = H'_j{H'_k} \indent  \text{ and } \indent
      \nu^1(x_i) + \nu^1(a) = \nu^1(x_i + a) =
      X_0^1(\nu^1(x_i + a)),
$$
   it follows that
\begin{equation*}
\begin{split}
      & \psi_i(a)\psi_i(b) =
          X_0^1 + G_i \left (H'_i +
          \nu^1(x_i + a)(X_0^1{H'_j}{H'_k} + \nu^1(b)) \right ) = \\
      & X_0^1 + G_i \left (H'_i +
          \nu^1(x_i + a)(\nu^1(x_j{x}_k) + \nu^1(b)) \right ) = \\
      & X_0^1 + G_i \left (H'_i +
          \nu^1(x_i + a)\nu^1(x_j{x}_k + b) \right ) =
          \psi_i((a + x_i)(b + x_j{x}_k)). \vspace{2mm}
\end{split}
\end{equation*}

3) From (\ref{mul}) and since
  $\nu^1(x_i) = X_0^1{H'_i} \subseteq H'_i$, we have
\begin{equation}
 \label{eq_mul}
     \psi_i(a)\psi_i(b) = X_0^1 + G_i \left (\nu^1(a) +
      H'_i(G_i + \nu^1(b) + \nu^1(x_i)) \right ).
 \end{equation}
 Again, by (\ref{permutation1}) we have
\begin{equation*}
\begin{split}
   & \psi_i(a)\psi_i(b) = X_0^1 + G_i \left (H'_i +
          \nu^1(a)(G_i + \nu^1(b) + \nu^1(x_i)) \right ) = \\
   &  X_0^1 + G_i \left ( H'_i + \nu^1(a)(X_0^1{H'_j}{H'_k}
             + \nu^1(b) + \nu^1(x_i)) \right ).
\end{split}
\end{equation*}
 Thus,  $\psi_i(a)\psi_i(b) = \psi_i(a(b + x_i + x_j{x}_k))$.
       \qedsymbol \vspace{2mm}
\begin{corollary}[Atomic multiplicativity]
  \label{cor_mul}
  1) Let one of the following inclusions hold:
\begin{equation*}
\begin{split}
   & {\rm(i)} \indent x_i + x_j{x}_k \subseteq a, \\
   & {\rm (ii)} \indent x_i + x_j{x}_k \subseteq b,  \\
   & {\rm (iii)} \indent x_j \subseteq a, \indent
              x_j{x}_k \subseteq b, \\
   & {\rm (iv)} \indent x_i \subseteq b, \indent
              x_j{x}_k \subseteq a. \\
\end{split}
\end{equation*}
 Then the joint map $\psi_i$ operates as a homomorphism on the elements
 $a$ and $b$ with respect to the
 lattice operations $+$ and $\cap$, i.e.,
$$
   \psi_i(a) + \psi_i(b) = \psi_i(a+b), \hspace{3mm}
    \psi_i(a)\psi_i(b) = \psi_i(a)\psi_i(b).
$$
  2) The joint map ${\psi_i}$ applied to the following atomic elements
    is the intersection preserving map, i.e., multiplicative with respect
    to the operation $\cap$:
\begin{equation}
     \label{homo1}
     \psi_i(ba_n^{ij}) =
             \psi_i(b)\psi_i(a_n^{ij})
             \text{ for every } b \subseteq D^{2,2,2},
\end{equation}
\begin{equation}
     \label{homo2}
     \psi_j(ba_n^{ij}) =
             \psi_j(b)\psi_j(a_n^{ij})
             \text{ for every } b \subseteq D^{2,2,2},
\end{equation}
\begin{equation}
     \label{homo2A}
     \psi_i(bA_n^{ij}) = \psi_i(b)\psi_i(A_n^{ij})
             \text{ for every } b \subseteq D^{2,2,2}.
\end{equation}
   3) The element $x_i$ can be inserted in any
      additive expression under the action of the joint map $\psi_i$:
\begin{equation}
     \label{homo3}
            \psi_i(b + x_i) = \psi_i(b)
            \text{ for every } b \subseteq D^{2,2,2}.
\end{equation}
\end{corollary}
\PerfProof
    1) Follows from Proposition \ref{additivity}.
        \vspace{2mm}

    2) Let us consider $a_n^{ij}$. For $n \geq 2$, we have
 $$
    a_n^{ij} = x_i + y_j{a}_{n-1}^{jk} =
    x_i + y_j(x_j + y_k{a}_{n-2}^{ki}) = \\
    x_i + x_j + y_j{y}_k{a}_{n-2}^{ki},
$$
    thus
    $a_n^{ij} \supseteq x_i + x_j$.
    Now, we consider atomic element $A_n^{ij}$.
    For $n \geq 2$, we have
$$
    A_n^{ij} = y_i + x_j{A}_{n-1}^{ki} =
    y_i + x_j(y_k + x_i{A}_{n-2}^{ji}) \supseteq x_i + x_j{x}_k.
$$
    For $n = 1$, the relation is proved directly.
    \vspace{2mm}

     3) Follows from Corollary \ref{cor_psi}, (1) since we have
$$
     \psi_i(b + x_i) = \psi_i(b) + \psi_i(x_i) \text{ and }
     \psi_i(x_i) = X_0^1 \subseteq \psi_i(b)
$$
     for every $b \in D^{2,2,2}$.
     \qedsymbol \vspace{2mm}

\section{Basic properties of atomic elements $a_n^{ij}$ and $A_n^{ij}$}
 \label{basic_atomic_L6}

The {\it atomic elements} are defined in \S\ref{atomic_polynoms}.
These polynomials are building bricks in the further
constructions.
\begin{table}[h]
  \renewcommand{\arraystretch}{1.7}
  \begin{tabular} {|c|c|c|c|}
  \hline
    1.1 & $a_n^{ij} = x_i + y_j{a}_{n-1}^{jk}$  (definition)
     &
    1.2 & $A_n^{ij} = y_i + x_j{A}_{n-1}^{ki}$  (definition) \\
  \hline
    2.1 &
      $a_0^{ij} \supseteq a_1^{ij} \supseteq \cdots
              a_n^{ij} \supseteq \cdots \supseteq x_i + x_j$
      &
    2.2 &
      $A_0^{ij} \supseteq A_1^{ij} \supseteq \cdots
              A_n^{ij} \supseteq \cdots$ \\
   \hline
    3.1 &
      $x_k{A}_n^{ij} = x_k{a}_n^{ji}$
      &
    3.2 &
      $A_n^{ij} = y_i + x_j{a}_{n-1}^{ik}$ \\
   \hline
    4.1 &
      $y_i{y}_j{A}_n^{ki} = y_i{y}_j{a}_n^{ik}$
      &
    4.2  & ${y}_i{a}_{n+1}^{ij} = y_i(x_i + y_j{A}_n^{kj})$     \\
   \hline
    5.1 &
      $a_n^{ij}a_m^{jk} =
          A_n^{ji}{a}_m^{jk}$  for $n \leq m + 1$
      &
    5.2 &
      $y_i(x_j + x_k){a}_n^{ij}$ =
            $y_i(x_j + x_k){A}_n^{ji}$ \\
   \hline
    6.1 &
       $ A_n^{ij}{a}_m^{ki} = y_i{a}_{m-1}^{ij}  + x_k{A}_n^{ij}$
      &
    6.2 &
        $A_n^{ij}{a}_m^{jk} = y_i{a}_m^{jk} + x_j{A}_{n-1}^{ki}$  \\
   \hline
    7.1 &
      $y_i(x_j + x_k) +
                     y_i{y}_j{a}_{n-2}^{ik} =
                        y_i{a}_n^{kj}$
      &
    7.2 &
      $y_i{y}_j +
                 y_i(x_j + x_k)A_{n-2}^{ki} = y_i{A}_n^{jk}$   \\
   \hline
    8.1 &
      $x_i(y_j + y_k) + y_i{y}_j{A}_{n-1}^{kj}$ =
      &
    8.2 &
      $y_i{y}_j + x_j(y_i + y_k){A}_{n-1}^{ki}$ = \\
      &
        $y_i(y_j + y_k){a}_n^{ij}$
      &  &
        $y_j(y_i + y_k){A}_n^{ij}$  \\
   \hline
  \end{tabular}
\vspace{2mm} \caption{\hspace{3mm}The basic properties of atomic
elements
 $a_n^{ij}$ and $A_n^{ij}$ in $D^{2,2,2}$}
  \label{table_atomic}
\end{table}

       The basic properties of atomic elements are given
       in Table \ref{table_atomic}.
       These properties, which
       will be used in further considerations.
       For convenience, we also give definition of $a_n^{ij},
       A_n^{ij}$ in relations (1.1)--(1.2) of Table \ref{table_atomic}.
       Indices $i,j,k$ in
       Table \ref{table_atomic} are distinct and $\{ i,j,k \} = \{1,2,3\}$.
       Properties 3.1, 3.2, 4.1, 4.2, 5.1, 5.2 from
       Table \ref{table_atomic} mean that
       atomic elements $a_n^{ij}$ and $A_n^{ji}$ coincide
       for many lattice polynomials and can
       substitute one another in many cases.
\vspace{2mm}

 {\it Proof of Properties 2.1 -- 8.2 (Table \ref{table_atomic}).}

    2.1) First, $a_{n-1}^{ij} = x_i + y_j{a}_{n-2}^{jk}$,
    $a_{n}^{ij} = x_i + y_j{a}_{n-1}^{jk}$.
    By induction, if $a_{n-2}^{jk} \supseteq a_{n-1}^{jk}$,
    then $a_{n-1}^{ij} \supseteq a_{n}^{ij}$. For $n = 1$,
    we have $a_0^{ij} = I \supseteq a_1^{ij}$. Further,
$$
    a_n^{ij} =
       x_i + y_j(x_j + y_k{a}_{n-2}^{ki}) =
       x_i + x_j + y_j{y}_k{a}_{n-2}^{ki} \supseteq
       x_i + x_j.
     \qed  \vspace{2mm}
$$

  2.2) Similarly,
$$
     A_{n-1}^{ij} = y_i + x_j{A}_{n-2}^{ki}
       \supseteq  y_i + x_j{A}_{n-1}^{ki} = A_{n}^{ij}.
     \qed \vspace{2mm}
$$

  3.1) By the induction hypothesis
  $x_j{A}_{n-1}^{ki} = x_j{a}_{n-1}^{ik}$, then
$$
 x_k{A}_n^{ij} =
      x_k(y_i + x_j{A}_{n-1}^{ki}) = x_k(y_i + x_j{a}_{n-1}^{ik}).
$$
      By property 1.1 from Table \ref{table_atomic}
      we have $a_{n-1}^{ik} \supseteq x_k$,
      thus by the permutation property (\ref{permutation1}) we deduce
$$
   x_k{A}_n^{ij} =
        x_k(x_j + y_i{a}_{n-1}^{ik}) = x_k{a}_n^{ji}.
     \qed \vspace{2mm}
$$

   3.2)  By property 3.1 from Table \ref{table_atomic} we have
   $A_n^{ij} =
   y_i + x_j{A}_{n-1}^{ki} = y_i + x_j{a}_{n-1}^{ik}$.
    \qedsymbol \vspace{2mm}

   4.1) We have
$$
   y_i{y}_j{a}_n^{ik} =
        y_i{y}_j(x_i + y_k{a}_{n-1}^{kj}) =
        y_i{y}_j(x_i + y_k(x_k + y_j{a}_{n-2}^{ji})),
$$
i.e.,
\begin{equation} \label{prop4}
    y_i{y}_j{a}_n^{ik} =
        y_i{y}_j(x_i + x_k + y_k{y}_j{a}_{n-2}^{ji})).
\end{equation}
On the other hand,
$$
    y_i{y}_j{A}_n^{ki} = y_i{}_j(y_k + x_i{A}_{n-1}^{jk}) =
         y_i{y}_j(y_k + x_i(y_j + x_k{A}_{n-2}^{ij})).
$$
Since $A_{n-2}^{ij} \supseteq x_i$, we have
        by the permutation property (\ref{permutation1}):
$$
    y_i{y}_j{A}_n^{ki} = y_i{y}_j(y_k +
            x_i(x_k + y_j{A}_{n-2}^{ij})).
$$
Further, since $y_j{A}_{n-2}^{ij} \supseteq y_i{y}_j$,
        again, by (\ref{permutation1}), we have
$$
         y_i{y}_j{A}_n^{ki} =
          y_i{y}_j(x_i + y_k(x_k + y_j{A}_{n-2}^{ij})) =
           y_i{y}_j(x_i + x_k +  y_k{y}_j{A}_{n-2}^{ij}).
$$
By the induction hypothesis
$$
       y_i{y}_j{A}_n^{ki} = y_i{y}_j(x_i +
                x_k + y_k{y}_j{a}_{n-2}^{ji})
$$
and by (\ref{prop4}) we have
              $y_i{y}_j{A}_n^{ki} = y_i{y}_j{a}_n^{ik}$.
  \qedsymbol \vspace{2mm}

 4.2) By property 4.1 we have
$$
    y_i{a}_{n+1}^{ij} =
          y_i(x_i + y_j{a}_n^{jk}) = x_i + y_i{y}_j{a}_n^{jk} =
          x_i + y_i{y}_j{A}_n^{kj} = y_i(x_i + y_j{A}_n^{kj}).
          \qed \vspace{2mm}
$$

    5.1) By property 2.1 and (\ref{permutation1}) we have
$$
     a_n^{ij}{a}_m^{jk} =
         (x_i + y_j{a}_{n-1}^{jk}){a}_m^{jk} =
         (y_j + x_i{a}_{n-1}^{jk}){a}_m^{jk}.
$$
      By property 3.1
      $a_n^{ij}a_m^{jk} =
         (y_j + x_i{A}_{n-1}^{kj})a_m^{jk} =
         A_{n}^{ji}a_m^{jk}$.
   \qedsymbol \vspace{2mm}

    5.2) We have
$$
     y_i(x_j + x_k)a_n^{ij} =
     y_i(x_j + x_k)(x_i + y_j{a}_{n-1}^{jk}).
$$
      Since  $a_{n-1}^{jk} \supseteq  x_j + x_k$,
      it follows that
\begin{equation*}
\begin{split}
    & y_i(x_j + x_k)a_n^{ij} =
      y_i(x_j + x_k)(y_j + x_i{a}_{n-1}^{jk}) = \\
    & y_i(x_j + x_k)(y_j + x_i{A}_{n-1}^{kj}) =
            y_i(x_j + x_k)A_n^{ji}.
         \qed \vspace{2mm}
\end{split}
\end{equation*}

  6.1) Since   $y_i{a}_{m-1}^{ij} \subseteq A_{n}^{ij}$,
      it follows that
$$
   A_{n}^{ij}a_m^{ki} =
        A_{n}^{ij}(x_k + y_i{a}_{m-1}^{ij}) =
        x_k{A}_{n}^{ij} + y_i{a}_{m-1}^{ij}.
    \qed \vspace{2mm}
$$

  6.2) Since  $x_j{A}_{n-1}^{ki} \subseteq a_m^{jk}$,
      it follows that
$$
       A_{n}^{ij}a_m^{jk} = a_m^{jk}(y_i + x_j{A}_{n-1}^{ki}) =
           y_i{a}_{m}^{jk} + x_j{A}_{n-1}^{ki}.
       \qed \vspace{2mm}
$$

  7.1) Since
$$
    y_i(x_j + x_k) + y_i{y}_j{a}_{n-2}^{ik} =
      y_i(x_j + x_k + y_i{y}_j{a}_{n-2}^{ik}) =
       y_i(x_k + y_j(x_j + y_i{a}_{n-2}^{ik})),
$$
by definition (1.1, Table \ref{table_atomic}) we have
        $y_i(x_k + y_j(x_j + y_i{a}_{n-2}^{ik})) =
              y_i{a}_n^{kj}$.
\qedsymbol \vspace{2mm}

   7.2) As above, using the modular law
(\ref{modular_law}) and
    the permutation property (\ref{permutation1})  we have
\begin{equation*}
\begin{split}
   & y_i{y}_j + y_i(x_j + x_k)A_{n-2}^{ki} =
     y_i((x_j + x_k)A_{n-2}^{ki} + y_i{y}_j) =\\
   & y_i(x_k + x_j{A}_{n-2}^{ki} + y_i{y}_j) =
     y_i(x_k + y_j(y_i + x_j{A}_{n-2}^{ki})) =
     y_i(x_k + y_j{A}_{n-1}^{ij}) = \\
   & y_i(y_j + x_k{A}_{n-1}^{ij}) =
     y_i{A}_n^{jk}.
     \qed \vspace{2mm}
\end{split}
\end{equation*}

  8.1) Since
$$
    x_i(y_j + y_k) + y_i{y}_j{A}_{n-1}^{kj} =
      (y_j + y_k)(x_i + y_i{y}_j{A}_{n-1}^{kj}),
$$
    by property 4.1 we have
\begin{equation*}
\begin{split}
   & x_i(y_j + y_k) + y_i{y}_j{A}_{n-1}^{kj} =
       (y_j + y_k)(x_i + y_i{y}_j{a}_{n-1}^{jk}) =  \\
   & y_i(y_j + {y}_k)(x_i + y_j{a}_{n-1}^{jk}) =
       y_i(y_j + y_k){a}_n^{ij}.
       \qed \vspace{2mm}
\end{split}
\end{equation*}

    8.2) As above,
$$
     y_i{y}_j + x_j(y_i + y_k)A_{n-1}^{ki} =
      y_j(y_i + y_k)(y_i + x_j{A}_{n-1}^{ki}) =
       y_i(y_j + y_k)A_n^{ij}.
       \qed \vspace{2mm}
$$

\section{Action of the maps $\psi_i$ and $\varphi_i$
            on atomic elements}
\begin{proposition}[Action of joint maps on atomic elements]
  The joint maps $\psi_i$ applied to the
 atomic elements $a^{ij}_n$ and $A^{ij}_n$  $(n \geq 1)$
 satisfy the following relations:
   \label{action_psi}
 \begin{enumerate}
   \item  $\psi_i(a^{ij}_n) = \nu^0(y_i{A}^{kj}_n)$, \vspace{2mm}
   \item $\psi_j(a^{ij}_n) =
        \begin{cases}
          \nu^0(y_j{A}^{kj}_n) & \text{ for } n > 1, \\
          \nu^0(y_j(y_i + y_k){A}^{kj}_1) & \text{ for } n = 1,
        \end{cases}$  \vspace{2mm}
   \item $\psi_k(y_i{a}^{ij}_n) =
                         \nu^0(y_k{A}^{ji}_{n+1})$, \vspace{2mm}
   \item $\psi_i(A^{ij}_n) =
            \nu^0(y_i(y_j + y_k){a}^{ik}_n)$, \vspace{2mm}
   \item $\psi_j(y_k{A}^{ij}_n) =
        \nu^0(y_j(x_k + y_i){a}^{ik}_n)$, \vspace{2mm}
   \item $\psi_k(y_j{A}^{ij}_n) =
        \nu^0(y_k{a}^{ji}_{n+1})$.  \vspace{2mm}
  \end{enumerate}
\end{proposition}

\begin{remark} {\rm
     The {\it joint} map $\psi_k$ increases the lower index $n$
     of the atomic element and swaps indices $i$ and $j$ in
     headings (3) and (6). The new atomic elements appear in this way.
\begin{equation*}
      ij \Longrightarrow ji, \indent n  \Longrightarrow n+1 .
\end{equation*}
     In headings (1) and (2), the lower index $n$ does not grow
     and the upper pair $ij$ becomes $kj$.
     The \underline{first index} of the pair is changed and
     the atomic elements ``$a$'' are transformed
     to atomic elements ``$A$'':
\begin{equation*}
      ij \Longrightarrow kj, \indent a \Longrightarrow A.
\end{equation*}
     In headings (4) and (5), the lower index $n$ does not grow
     and the upper pair $ij$ becomes $ik$.
     The \underline{second index} of the pair is changed and
     the atomic elements ``$A$'' are transformed
     to atomic elements ``$a$'':
\begin{equation*}
      ij \Longrightarrow ik, \indent A \Longrightarrow a.
\end{equation*}
     Map $\psi_k$ transforms
     $a \Longrightarrow A$ in headings (1), (2), (3)
     and $A \Longrightarrow a$ in headings (4), (5), (6).
     }
\end{remark}
\begin{remark}
    {\rm  It is not so interesting to
     consider $\psi_j(A^{ij}_n)$ because }
$$
     \psi_j(A^{ij}_n) = \psi_j(y_i + x_j{A}^{ki}_n) =
     \psi_j(y_i).
$$
    {\rm By basic relations (Proposition \ref{basic_eq})}
     $\psi_j(y_i) = \nu^0(y_j(x_i + y_k))$.
\end{remark}

     For further moving of indices, see
     \S\ref{list_admis}-\S\ref{sect_adm_classes} devoted
     to the admissible elements. \\
\PerfProofW {\it Proposition \ref{action_psi}}
  uses a simultaneous induction on (1)--(6).
  For convenience, without loss of generality,
  we let $i = 1$, $j = 2$. \\
  \indent 1) For $n = 1$, by Proposition \ref{basic_eq} we have
$$
   \psi_1(a^{12}_1) =
    \psi_1(x_1 + y_2) =
     \psi_1(y_2) = \nu^0(y_1(x_2 + y_3)) =
     \nu^0(y_1{A}^{32}_1).
$$
   For $n \geq 2$, by Corollary \ref{cor_psi} we have
$$
    \psi_1(a^{12}_n) = \psi_1(x_1 + y_2{a}^{23}_{n-1}) =
     \psi_1(y_2{a}^{23}_{n-1}).
$$
     Since
$$
    y_2{a}^{23}_{n-1} = y_2(x_2 + y_3{a}^{31}_{n-2}) =
            x_2 + y_2{y}_3{}^{31}_{n-2}
$$
    and
            $a^{31}_{n-2} \supseteq x_1$,
    we see by multiplicativity
            (Corollary \ref{cor_mul}) that
$$
    \psi_1(a^{12}_n) = \psi_1(x_2) +
        \psi_1(y_2{y}_3)\psi_1(a^{31}_{n-2}).
$$
    By basic relations
     (Proposition \ref{basic_eq})
$$
    \psi_1(a^{12}_n) = \nu^0(y_1{y}_3) +
       \nu^0(y_1(x_2 + x_3))\psi_1(a^{31}_{n-2}).
$$
    By the induction hypothesis
     from (2) we have $\psi_1(a^{31}_{n-2}) =
     \nu^0(y_1{A}^{21}_{n-2})$ and
$$
     \psi_1(a^{12}_n) = \nu^0(y_1{y}_3 +
          y_1(x_2 + x_3){A}^{21}_{n-2}).
$$
By property 7.2 from Table \ref{table_atomic} we have
$$
      \psi_1(a^{12}_n) = \nu^0(y_1{A}^{32}_{n}).
      \qed \vspace{2mm}
$$

  2)  For $n = 1$, we have
\begin{equation*}
\begin{split}
 & \psi_2(a^{12}_1) =  \psi_2(x_1) + \psi_2(y_2) =
  \nu^0(y_2{y}_3 + x_2(y_1 + y_3))= \\
 & \nu^0(y_2(y_1 + y_3)(y_3 + x_2))=
  \nu^0(y_2(y_1 + y_3)A^{32}_1).
\end{split}
\end{equation*}
  For $n \geq 2$, we have
  $\psi_2(a^{12}_n) = \psi_2(x_1) +
                    \psi_2(y_2{a}^{23}_{n-1})$.
  By the multiplicativity property (Corollary \ref{cor_mul}) we have
$$
  \psi_2(a^{12}_n) = \nu^0(y_2{y}_3) +
              \psi_2(y_2)\psi_2(a^{23}_{n-1}).
$$
  By the induction hypothesis, from (1) we have
       $\psi_2(a^{23}_{n-1}) = \nu^0(y_2{A}^{13}_{n-1})$ and
$$
  \psi_2(a^{12}_n) = \nu_0(y_2{y}_3 +
       x_2(y_1 + y_3{A}^{13}_{n-1})).
$$
  Since $A^{13}_{n-1} \subseteq y_1 + y_3$, we see that
$$
  \psi_2(a^{12}_n) =
       \nu^0(y_2{y}_3 + x_2{A}^{13}_{n-1}) =
       \nu^0(y_2(y_3 + x_2{A}^{13}_{n-1})) = \nu^0(y_2{A}^{32}_n).
       \qed \vspace{2mm}
$$

   3) We have
$$
  \psi_3(y_1{a}^{12}_n) = \psi_3(y_1(x_1 + y_2{a}^{23}_{n-1})) =
      \psi_3(x_1 + y_1{y}_2{a}^{23}_{n-1}).
$$
Further, by multiplicativity (Corollary \ref{cor_mul}) we have
$$
   \psi_3(y_1{a}^{12}_n) = \psi_3(x_1) +
         \psi_3(y_1{y}_2)\psi_3(a^{23}_{n-1}).
$$
From basic relations (Proposition \ref{basic_eq}) it follows that
$$
  \psi_3(y_1{a}^{12}_n) =
         \nu^0(y_3{y}_2) +
            \nu^0(y_3(x_1 + x_2))\psi_3(a^{23}_{n-1}).
$$
By induction, from (2)
      $\psi_3(y_1{a}^{12}_n) =
           \nu^0(y_3{y}_2 + y_3(x_1 + x_2){A}^{13}_{n-1})$.
      Again, from Table \ref{table_atomic}, property 7.2 we have
         $\psi_3(y_1{a}^{12}_n) = \nu^0(y_3{A}^{21}_{n+1})$.
      In this case $\psi_k$ increases the lower index of the atomic element
      $a^{ij}_n$.
\qedsymbol \vspace{2mm}

  4) As above,
$$
\psi_1(A^{12}_n) = \psi_1(y_1 + x_2{A}^{31}_{n-1}) =
      \psi_1(y_1) + \psi_1(x_2{A}^{31}_{n-1}).
$$
     From Table \ref{table_atomic}, property 3.1 we have
$$
       \psi_1(A^{12}_n) =
        \psi_1(y_1) + \psi_1(x_2{a}^{13}_{n-1}) =
        \psi_1(y_1) + \psi_1(x_2)\psi_1(a^{13}_{n-1}).
$$
     By basic relations (Proposition \ref{basic_eq})
     and by induction hypothesis, from (1)
     we have
$$
       \psi_1(A^{12}_n) = \nu^0(x_1(y_2 + y_3) + y_1{y}_3{A}^{23}_{n-1}).
$$
     From Table \ref{table_atomic}, property 8.1 we have
$$
     \psi_1(A^{12}_n) = \nu^0(y_1(y_2 + y_3){a}^{13}_n).
     \qed \vspace{2mm}
$$

  5) Again,
    $\psi_2(y_3{A}^{12}_n) =
            \psi_2(y_3(y_1 + x_2{A}^{31}_{n-1}))$.
     Since $y_3 \subseteq A^{31}_{n-1}$, we have
     by the permutation property (\ref{permutation1}):
$$
    \psi_2(y_3{A}^{12}_n) =
             \psi_2(y_3(x_2 + y_1{A}^{31}_{n-1})).
$$
     By (\ref{homo1})-(\ref{homo3}) and
            Corollary \ref{cor_mul} we have
$$
    \psi_2(y_3{A}^{12}_n) = \psi_2(y_3)(\psi_2(x_2) +
                      \psi_2(y_1{A}^{31}_{n-1})) =
             \psi_2(y_3)\psi_2(y_1{A}^{31}_{n-1}).
$$
    By the induction hypothesis, from (6) it follows that
$$
        \psi_2(y_3{A}^{12}_n) = \nu^0(y_2(x_3 + y_1){a}^{13}_n))
$$
         because
          $\psi_2(y_1{A}^{31}_{n-1}) = \nu^0(y_2{a}^{13}_n)$.
\qedsymbol \vspace{2mm}

  6) We have
$$
   \psi_2(y_1{A}^{31}_n) =
          \psi_2(y_1(y_3 + x_1{A}^{23}_{n-1})) =
          \psi_2(y_1{y}_3 + x_1{a}^{32}_{n-1}) =
          \psi_2(y_1{y}_3) + \psi_2(x_1)\psi_2(a^{32}_{n-1}).
$$
    By the induction hypothesis, from (2) it follows that
$$
     \psi_2(y_1{A}^{31}_n) = \nu^0(y_2(x_1 + x_3) +
                 y_2{y}_3{A}^{12}_{n-1}).
$$
    From Table \ref{table_atomic}, property 4.1 we have
$$
    \psi_2(y_1{A}^{31}_n) = \nu^0(y_2(x_1 + x_3) +
                 y_2{y}_3{a}^{21}_{n-1}).
$$
    Finally, from Table \ref{table_atomic}, property 7.1
      $\psi_2(y_1{A}^{31}_n) =
            \nu^0(y_2{a}^{13}_{n+1})$.
    In this case $\psi_k$ increases the lower index of the
    atomic element  $A^{ij}_n$.
\qedsymbol \vspace{2mm}

 From Propositions (\ref{phi_and_psi}) and
(\ref{action_psi}) we have
\begin{corollary}[Action of the elementary maps]
 The elementary maps $\varphi_i$ applied to the
atomic elements $a^{ij}_n$ and $A^{ij}_n$ $(n \geq 1)$ satisfy the
following relations:
   \label{act_runner_map}
   \item \vspace{2mm}
   \begin{enumerate}
   \item 
     $\varphi_i\Phi^+\rho(a^{ij}_n) = \rho(y_i{A}^{kj}_n)$, \vspace{2mm}
   \item 
     $\varphi_j\Phi^+\rho(a^{ij}_n) =
        \begin{cases}
          \rho(y_j{A}^{kj}_n) & {\rm for \hspace{1mm}} n > 1, \\
          \rho(y_j(y_i + y_k)A^{kj}_1)
                                         & {\rm for \hspace{1mm}} n = 1,
        \end{cases}$ \vspace{2mm}
   \item 
     $\varphi_k\Phi^+\rho(y_i{a}^{ij}_n) =
                         \rho(y_k{A}^{ji}_{n+1})$, \vspace{2mm}
   \item 
     $\varphi_i\Phi^+\rho(A^{ij}_n) = \rho(y_i(y_j + y_k){a}^{ik}_n)$,
     \vspace{2mm}
   \item 
     $\varphi_j\Phi^+\rho(y_k{A}^{ij}_n) =
        \rho(y_j(x_k + y_i)a^{ik}_n)$, \vspace{2mm}
   \item 
     $\varphi_k\Phi^+\rho(y_j{A}^{ij}_n) =
        \rho(y_k{a}^{ji}_{n+1})$. \vspace{2mm}
   \end{enumerate}
\end{corollary}

\section{$\varphi_i-$homomorphic theorem}
 \label{sect_Homomorphic}
\begin{theorem}
  \label{th_homomorhism}
  1) The polynomials $a^{ij}_n$
            are $\varphi_i-$ and $\varphi_j-$homomorphic.
            \vspace{2mm}

  2) The polynomials $A^{ij}_n$ are $\varphi_i-$homomorphic.
            \vspace{2mm}

  3) The polynomials $A^{ij}_n$ are ($\varphi_j,y_k)-$homomorphic.
            \vspace{2mm}
\end{theorem}
\PerfProof
  1),2) Follow from multiplicativity (\ref{homo1}) -- (\ref{homo2A}),
  Proposition \ref{action_psi} and Proposition \ref{phi_and_psi}.
\qedsymbol \vspace{2mm}

  3) For convenience, without loss of
generality,
  we let $i=1, j=2, k=3$. We need to prove that
   \begin{equation} \label{need_prove1}
        \varphi_2\Phi^+\rho(A^{12}_n{p}) =
          \varphi_2\Phi^+\rho(y_3{A}^{12}_n)\varphi_2\Phi^+\rho(p)
            \text{ for every } p \subseteq y_3.
   \end{equation}
   Since $\psi_2(A^{12}_n{p}) =
           \psi_2((y_1 + x_2{A}^{31}_{n-1})p)$ and
   $p \subseteq y_3 \subseteq A^{31}_{n-1}$,
   we have by the permutation property
   (\ref{permutation1}):
$$
    \psi_2(A^{12}_n{p}) =
         \psi_2((x_2 + y_1{A}^{31}_{n-1})p) =
            \psi_2(x_2 + y_1{A}^{31}_{n-1})\psi_2(p).
$$
   The element $x_2$ can be dropped according to (\ref{homo3}),
   and hence,
$$
   \psi_2(A^{12}_n{p}) = \psi_2(y_1{A}^{31}_{n-1})\psi_2(p).
$$

   By property (6) from Proposition \ref{action_psi} we have
\begin{equation}
     \label{loc_psi1}
       \psi_2(A^{12}_n{p}) = \nu^0(y_2{a}^{13}_n)\psi_2(p).
\end{equation}
   Since $y_2{a}^{13}_n =
      y_2(y_1 + y_3){a}^{13}_n$, we have
   $\nu^0(y_2{a^{13}}_n) =
      \nu^0(y_2(y_1 + y_3))\nu^0(a^{13}_n)$
   and by (5) from Proposition \ref{basic_eq} we have
   $\nu^0(y_2{a}^{13}_n) = \psi_2(I)\nu^0(a^{13}_n)$.
   Further, $\psi_2(I) \supseteq \psi_2(p)$ and (\ref{loc_psi1})
   is equivalent to
   \begin{equation}  \label{loc_psi2}
        \psi_2(A^{12}_n{p}) = \nu^0(a^{13}_n)\psi_2(p).
   \end{equation}
   We have $\psi_2(p) \subseteq \psi_2(y_3)$ together
   with $p \subseteq y_3$, and therefore
   $\psi_2(p) = \psi_2(p)\psi_2(y_3)$. From (\ref{loc_psi2})
   and Proposition \ref{basic_eq} we have
\begin{equation*}
\begin{split}
   & \psi_2(A^{12}_n{p}) =
    \nu^0(a^{13}_n)\psi_2(p)\psi_2(y_3) = \\
   & \nu^0(a^{13}_n)\psi_2(p)\nu^0(y_2(x_1 + y_3)) =
    \nu^0(y_2(x_1 + y_3)a^{13}_n)\psi_2(p).
\end{split}
\end{equation*}
   By (5) from Proposition \ref{action_psi}
   $\psi_2(A^{12}_n{p}) =
        \psi_2(y_3{A}^{12}_n)\psi_2(p)$.
   Applying projection $\nabla$ from \S\ref{seq_assoc} we get
   \begin{equation}  \label{loc_psi3}
    \nabla\psi_2(A^{12}_np) =
                  \nabla\psi_2(y_3{A}^{12}_n)\nabla\psi_2(p)
   \end{equation}
   because  $\psi_2(y_3{A}^{12}_n) = \nu^0(y_2(x_1 + y_3){a}^{13}_n)$,
      see Lemma \ref{preserv_intersection}, relation (\ref{pr_inters}).

    By Proposition \ref{phi_and_psi}
    $\nabla\psi_i(c) =
    \varphi_i\Phi^+\rho(c) \text{ for every } c \subseteq D^{2,2,2}$,
    thus we get (\ref{need_prove1}), and therefore heading (3) of the
    $\varphi-$homomorphic Theorem \ref{th_homomorhism} is proven.
\qedsymbol \vspace{2mm}

The final proposition of this section shows when the lower indices
of atomic elements are increased by action some elementary
mappings $\varphi_i$:
\begin{corollary}
   The index growth takes place in the following cases:
\begin{enumerate}
  \item
     $\varphi_k\Phi^+\rho(y_i{a}^{ij}_n) =
                  \rho(y_k{A}^{ji}_{n+1})$; \vspace{2mm}
  \item
     $\varphi_k\Phi^+\rho(y_j{a}^{ij}_n) =
                  \rho(y_k{a}^{ji}_{n+1})$. \vspace{2mm}
\end{enumerate}
\end{corollary}

\section{The perfectness of $y_i + y_j$} To prove the {\it
perfectness} of the element $y_1 + y_2$, two simple lemmas are
necessary.
\begin{lemma}
  Let  $A, B \subseteq X, A, B$ be subspaces
  in the finite-dimensional vector space $X$.
  There exists a subspace $C \subseteq B$, such that
  \begin{equation} \label{decomp1}
        C + AB = B \text{ and } CA = 0.
  \end{equation}
\end{lemma}

  Indeed, take any direct complement $C$ of $AB$, i.e.,
  $C \oplus AB = B$. If $v \in CA$ is non-zero, then $v \in CA \subseteq BA$
  and the sum $C + AB$ is not direct.
 \qedsymbol \vspace{2mm}

\begin{lemma}
   If $U,V,W \subseteq X$ and $(U + V)W = 0$, then
\begin{equation}
  \label{decomp2}
   U(V + W) = UV \text{ and } V(U + W) = UV.
\end{equation}
\end{lemma}
   Indeed,
$$
   U(V + W) \subseteq (U + V)(V + W) = (U + V)W + V = V,
$$
   i.e., $U(V + W) \subseteq UV$.
   The inverse inclusion is obvious. Similarly,
$$
   V(U + W) \subseteq (U + V)(U + W) = U + W(U + V) = U,
$$
   i.e., $V(U + W) \subseteq UV$
      and the inverse inclusion is obvious.
\qedsymbol \vspace{2mm}
\begin{proposition}  \label{first_perfect}
  1) The element $y_i + y_j$ where $i \neq j$, is perfect in $D^{2,2,2}$ .
    \vspace{2mm}

  2) If $\rho(y_i + y_j)$ = 0 for some
   indecomposable representation $\rho$, where $i \neq j$, then
   $\rho$ is one from the seven projective representations
   (Table \ref{RepOrder12}).
    \vspace{2mm}
\end{proposition}
\PerfProof  1) Consider $y_1 + y_2$. Let $B = X_3, ~A = Y_1 + Y_2$
in (\ref{decomp1}), then there exists a subspace $C \subseteq X_3$
such that
\begin{equation}  \label{decomp3}
    C + X_3(Y_1 + Y_2) = X_3, \indent  C(Y_1 + Y_2) = 0.
\end{equation}
Therefore
    $Y_3 \supseteq C + Y_3(Y_1 + Y_2)$.
Again use (\ref{decomp1}) with
$$
   B = Y_3, A = C + Y_3(Y_1 + Y_2).
$$
  Then, there exists a subspace $D \subseteq Y_3$ such that
\begin{equation} \label{decomp4}
     D + C + Y_3(Y_1 + Y_2) = Y_3, \indent
     D[C + Y_3(Y_1 + Y_2)] = 0.
\end{equation}
 Using (\ref{decomp2}) with $U = Y_3(Y_1 + Y_2)$,
 $V = C$ and $W = D$, we obtain from
 (\ref{decomp2}) and (\ref{decomp4}):
$$
   U(V + W) = Y_3(Y_1 + Y_2)(C + D) =
       C{Y}_3(Y_1 + Y_2).
$$
 By (\ref{decomp3}) we have
     $Y_3(Y_1 + Y_2)(C + D)$ = 0.
 Since $D + C \subseteq Y_3$, we see that
 $(Y_1 + Y_2)(C + D) = 0$. Therefore, $D + C$
 and $Y_1 + Y_2$ form a direct sum. We complement this sum to $X_0$:
 \begin{equation} \label{decomp5}
        (D + C) \oplus (Y_1 + Y_2) \oplus F = X_0.
 \end{equation}
 We will show that the decomposition (\ref{decomp5}) forms a decomposition of
 the representation $\rho$ in $X_0$.

 First,
 \begin{equation}  \label{decomp6}
       Z = Z(Y_1 + Y_2) + Z(D + C + F)
       \text{ for }
       Z = Y_1, Y_2, X_1, X_2
 \end{equation}
because $Z \subseteq Y_1 + Y_2$ and the sum (\ref{decomp5}) is
direct. By (\ref{decomp3}) we have
$$
   X_3 \subseteq X_3(Y_1 + Y_2) +  X_3(D + C + F) \subseteq X_3,
$$
i.e.,
\begin{equation}   \label{decomp7}
    X_3 = X_3(Y_1 + Y_2) + X_3(D + C + F).
\end{equation}
Similarly, by (\ref{decomp4}) we have
$$
    Y_3 \subseteq Y_3(Y_1 + Y_2)
        + Y_3(D + C + F) \subseteq Y_3,
$$
i.e.,
\begin{equation}
  \label{decomp8}
    Y_3 = Y_3(Y_1 + Y_2) + Y_3(D + C + F).
\end{equation}
So, by (\ref{decomp6})--(\ref{decomp8}) $Y_1 + Y_2$ and $D + C +
F$ form a decomposition of the representation $\rho$. Since
$\rho(y_1 + y_2) = Y_1 + Y_2$, we have either
 $\rho(y_1 + y_2) = X_0$ or $\rho(y_1 + y_2) = 0$,
i.e., $y_1 + y_2$ is a perfect element. \qedsymbol

  2) If $\rho(y_1 + y_2) = 0$, then
  $Y_1 = X_1 = Y_2 = X_2 = 0$ and $\rho$ is one of
   $\rho_{x_0}, \rho_{y_3}, \rho_{x_3}$,
   see Table \ref{RepOrder12}.
\qedsymbol

\begin{corollary}  \label{second_perfect}
   The elements $y_1 + y_2 + x_3$ and
   $y_1 + y_2 + y_3$ are perfect.
\end{corollary}
\PerfProof
   Let $Y_1 + Y_2 = X_0$. Then $Y_1 + Y_2 + X_3 = X_0$.
   If $Y_1 + Y_2 = 0$, then, in the indecomposable representation
   $\rho: D^{2,2,2} \longrightarrow \mathcal{L}(X_0)$, we have
   $X_0 = Y_3 = X_3$.
   The same for the element $y_1 + y_2 + y_3$.

 If $\rho(y_1 + y_2) = 0$, then
  $Y_1 = X_1 = Y_2 = X_2 = 0$ and  $\rho$ is one of
   $\rho_{x_0}, \rho_{y_3}, \rho_{x_3}$,
   see Table \ref{RepOrder12}.
\qedsymbol

\section{The elementary map $\varphi_i$: a way to construct new perfect elements}

In this section we list three fundamental properties of the
elementary map $\varphi_i$. Proposition \ref{new_perfect} gives a
way to construct new {\it perfect} elements from already existing
ones. Proposition \ref{motiv_admis} motivates the construction of
{\it admissible} sequences and {\it admissible} elements.
\begin{proposition} \label{new_perfect}
  Let z be a perfect element in $D^{2,2,2}$ and,
for every indecomposable representation $\rho$, let
\begin{equation}
     \varphi_i\Phi^+\rho(z) + \varphi_j\Phi^+\rho(z) =
     \rho(u),  \text{ where } i \neq j.
\end{equation}
Then u is also the perfect element.
\end{proposition}
\PerfProof
   If $\Phi^+\rho(z) = 0$, then $\rho(u) = 0$.
   If $\Phi^+\rho(z) = X^1_0$,
      then by Corollary \ref{cor_psi} we have
\begin{equation*}
\begin{split}
  & \rho(u) =
        \varphi_i{X}^1_0 + \varphi_j{X}^1_0 =
        \varphi_i\Phi^+\rho(I) + \varphi_j\Phi^+\rho(I) = \\
  &     Y_i(Y_j + Y_k) + Y_j(Y_i + Y_k) =
        (Y_i + Y_j)(Y_i + Y_k)(Y_j + Y_k)
\end{split}
\end{equation*}
        and by Proposition \ref{first_perfect}
      the element $u$ is also perfect.
\qedsymbol

\begin{proposition}
 \label{motiv_admis}
  For $\{i,j,k\} = \{1,2,3\}$, the following relations hold
\begin{equation}
     \varphi_i\varphi_j\varphi_i +  \varphi_i\varphi_k\varphi_i = 0, \vspace{2mm}
\end{equation}
\begin{equation}
     \varphi_i^3 = 0.
\end{equation}
\end{proposition}
\PerfProof
  For every vector $v \in X_0^1$, by definition of $\varphi_i$
  we have $(\varphi_i + \varphi_j + \varphi_k)(v) = 0$, see eq. (\ref{sum_fi}).
  In other words, $\varphi_i + \varphi_j + \varphi_k$ = 0. Therefore,
$$
    \varphi_i\varphi_j\varphi_i +  \varphi_i\varphi_k\varphi_i =
    \varphi_i(\varphi_j + \varphi_k)\varphi_i = -\varphi_i^3.
$$
  So, it suffices to prove that $\varphi_i^3$ = 0.
  For every
   $z \subseteq D^{2,2,2}$, by headings (1),(3) and (5)
  of Corollary \ref{cor_psi} we have
\begin{equation*}
\begin{split}
  & \varphi_i^3((\Phi^+)^3\rho)(z) \subseteq
      \varphi_i^2[\varphi_i((\Phi^+)^3\rho)(I)] =  \\
  & \varphi_i^2[((\Phi^+)^2\rho)(y_i(y_j + y_k))]
      \subseteq \varphi_i^2[((\Phi^+)^2\rho)(y_i)] =  \\ \vspace{2mm}
  &  \varphi_i[(\Phi^+\rho)(x_i(y_j + y_k))]
      \subseteq  \varphi_i[(\Phi^+\rho)(x_i)] =
  \varphi_i\Phi^+\rho(x_i) = 0.
  \qed
\end{split}
\end{equation*}

 \begin{corollary}
  \label{basic_rel}
     The relation
 \begin{equation}
      \varphi_i\varphi_j\varphi_i(B)
       = \varphi_i\varphi_k\varphi_i(B)
 \end{equation}
 takes place for every subspace $B \subseteq X_0^3$,
 where $X_0^3$ is the representation space
 of $(\Phi_i^+)^3\rho$.
 \end{corollary}

\section{List of the admissible sequences}
  \label{list_admis}
Recall that the {\it admissible sequences} are introduced in
\S\ref{adm_seq_polynom}. Let us construct new admissible sequences
acting by the {\it elementary} map $\varphi_i$.

  Let the action of $\varphi_i$ on an admissible sequence $\alpha$
  be defined so that
  \begin{enumerate}
  \item  The index $i$ is added in front of the
  sequence $\alpha = i_1{i}_2\dots{i}_n$, i.e.,
  \begin{equation}
        \varphi_i(i_1{i}_2\dots{i}_n) = i{i}_1{i}_2\dots{i}_n.
  \end{equation}
  \item New sequence $\varphi_i(\alpha)$
  should also be an admissible sequence, in other words, $i \neq i_1$.
  \end{enumerate}

The rule $i \neq i_1$ in the definition of admissible sequence is
motivated by the property $\varphi_i = -(\varphi_j + \varphi_k)$,
so that $\varphi_i$ can be always replaced by
 $\varphi_j + \varphi_k$.

\begin{proposition} \label{list_adm}
The admissible sequence starting at $i_1 = 1$ may be transformed
to one of the seven types  given in Table \ref{table_admissible},
column 1 (similarly for $i_1 = 2$ and $i_1 = 3$, see
(\ref{starting_2})). Here, $n \in \{1,2,3,\dots \}, m \in
\{0,1,2,3,\dots \}$. In case 7, $n \geq 0$. In cases 1 and 4, $m >
0$.

For $m = 0$ cases 1 and 4 are not given in the table: in case 1
(resp. case 4), action $\varphi_1$ transforms $(21)^n$ to
$1(21)^n$ (resp. $(31)^n$ to $1(31)^n$).
\end{proposition}

\begin{table}[h]
  \renewcommand{\arraystretch}{1.35}
  \begin{tabular} {|c||c|c|c|c|}
  \hline \hline
     & Admissible & Action  & Action & Action  \cr
     $N$  & Sequence & $\varphi_1$ & $\varphi_2$ & $\varphi_3$
     \\
  \hline  \hline  
    $1$  & $(213)^m(21)^n$
     & $13(213)^{m-1}(21)^{n+1}$
     & --
     & $3(213)^m(21)^n$  \\
  \hline         
    $2$  & $3(213)^m(21)^n$
     & $13(213)^m(21)^n$
     & $(213)^m(21)^{n+1}$
     & --                   \\
  \hline        
    $3$  & $13(213)^m(21)^n$
     & --
     & $(213)^{m+1}(21)^n$
     & $3(213)^m(21)^{n+1}$   \\
  \hline \hline   
   $4$  & $(312)^m(31)^n$
     & $12(312)^{m-1}(31)^{n+1}$
     & $2(312)^m(31)^n$
     & -- \\
  \hline         
   $5$  & $2(312)^m(31)^n$
     & $12(312)^m(31)^n$
     & --
     & $(312)^m(31)^{n+1}$ \\
  \hline         
   $6$  & $12(312)^m(31)^n$
     & --
     & $2(312)^m(31)^{n+1}$
     & $(312)^{m+1}(31)^n$ \\
  \hline \hline   
   $7$  & $1(21)^n = 1(31)^n$
     & --
     & $(21)^{n+1}$
     & $(31)^{n+1}$ \\
  \hline
  \end{tabular}
  \vspace{2mm}
  \caption{\hspace{3mm}For $D^{2,2,2}$: the admissible sequences starting at $i_1 = 1$}
  \label{table_admissible}
\end{table}

\PerfProof
  It suffices to prove that every action of $\varphi_i$ on any admissible
  sequence from the list (Table \ref{table_admissible}, column 1) gives us again
  an admissible sequence from the same list.

  \underline{Line 1, action $\varphi_1$}.
  Applying $\varphi_1$ to the sequence $(213)^m(21)^n$ we get
$$
  1(213)^m(21)^n = 1(213)(213)^{m-1}(21)^n.
$$
  Since
  1213(2...) = 1321(2...), we see that
  $1(213)^m(21)^n = 132[1(213)^{m-1}(21)^n$]. By the induction hypothesis
$$
  1(213)^{m-1}(21)^n = 13(213)^{m-2}(21)^{n+1},
$$
  so
$$
  1(213)^m(21)^n = 132[13(213)^{m-2}(21)^{n+1}] =
  13(213)^{m-1}(21)^{n+1},
$$
i.e., we get a sequence from Line 3. \qedsymbol

  \underline{Line 1, action $\varphi_3$}. Applying
$\varphi_3$ to the sequence
  $(213)^m(21)^n$ we just get a sequence from Line 2.
\qedsymbol

  \underline{Line 2, action $\varphi_1$}. Action $\varphi_1$ on the
$3(213)^m(21)^n$ leads to Line 3. \qedsymbol

  \underline{Line 2, action $\varphi_2$}. Applying $\varphi_2$
  to the sequence $3(213)^m(21)^n$, we get
$$
  23(213)^m(21)^n = 23(213)(213)^{m-1}(21)^n.
$$
  Since 23213... = 21323..., we have
  $23(213)^m(21)^n =   213[23(213)^{m-1}(21)^n]$.
  Again, by the induction hypothesis
$$
  23(213)^{m-1}(21)^n = (213)^{m-1}(21)^{n+1},
$$
  so
$$
   23(213)^m(21)^n = 213[(213)^{m-1}(21)^{n+1}] =
  (213)^m(21)^{n+1},
$$
i.e., a sequence from Line 1. \qedsymbol

 \underline{Line 3, action $\varphi_2$}. Applying
$\varphi_2$ we immediately get Line 1. \qedsymbol

 \underline{Line 3, action $\varphi_3$}. Apply $\varphi_3$ to the sequence
  $13(213)^m(21)^n$. Again,
$$
  313(213)^m(21)^n =
  313(213)(213)^{m-1}(21)^n = 321[313(213)^{m-1}(21)^n].
$$
  By the induction hypothesis
$$
  313(213)^m(21)^n =
  321[3(213)^{m-1}(21)^{n+1}] = 3(213)^m(21)^{n+1}.
  \qed
$$

  \underline{Lines 4, 5, 6}. Similarly as for Lines 1, 2, 3. \qedsymbol

  \underline{Line 7}. Applying $\varphi_2$ and $\varphi_3$
we get Line 1 and Line 4 with $m = 0$, respectively: $(21)^{n+1}$
and $(31)^{n+1}$. \qedsymbol

\section{The theorem on the admissible element classes}
 \label{sect_adm_classes}
 The lattice polynomials indexed by admissible
sequences are also said to be {\it admissible}. We define
admissible elements $f_\alpha, e_\alpha, g_{\alpha0}$ in Table
 \ref{table_adm_elems}. By definition $g_0 = I$.

 \index{admissible elements! -
   $f_\alpha$, $e_\alpha$, $g_{{\alpha}0}$}

{
\begin{table}[h]
 \renewcommand{\arraystretch}{1.5}
 \begin{tabular} {|c||c|c|c|c|}
  \hline \hline
    $N$ & $\alpha$
            & $f_\alpha$ & $e_\alpha$ & $g_{\alpha0}$
     \\
  \hline  \hline
    1 & $\gamma$
      & $y_1y_2a^{13}_qA^{32}_{k-1}$
      & $y_2A^{32}_{k-1}a^{21}_qA^{13}_ka^{13}_q$
      & $e_\alpha(x_1+a^{32}_qA^{32}_k)$    \\
  \hline
    2 & $3\gamma$
      & $y_3(x_1+x_2)A^{23}_qa^{31}_{k-1}$
      & $y_3a^{31}_{k-1}A^{12}_{q+1}a^{12}_kA^{23}_q$
      & $e_\alpha(y_2y_3+A^{12}_qa^{31}_k)$ \\
  \hline
    3 & $13\gamma$
      & $y_3y_1a^{32}_{q+1}A^{21}_{k-1}$
      & $y_1A^{21}_{k-1}a^{13}_{q+1}A^{32}_ka^{32}_{q+1}$
      & $e_\alpha(x_3+a^{21}_{q+1}A^{21}_k)$    \\
  \hline
    4 & $(213)\gamma$
      & $y_2(x_3+x_1)A^{12}_{q+1}a^{23}_{k-1}$
      & $y_2a^{23}_{k-1}A^{31}_{q+2}a^{31}_kA^{12}_{q+1}$
      & $e_\alpha(y_1y_2+A^{31}_{q+1}a^{23}_k)$ \\
  \hline
    5 & $3(213)\gamma$
      & $y_2y_3a^{21}_{q+2}A^{13}_{k-1}$
      & $y_3A^{13}_{k-1}a^{32}_{q+2}A^{21}_ka^{21}_{q+2}$
      & $e_\alpha(x_2+a^{13}_{q+2}A^{13}_k)$    \\
  \hline
    6 & $13(213)\gamma$
      & $y_1(x_2+x_3)A^{31}_{q+2}a^{12}_{k-1}$
      & $y_1a^{12}_{k-1}A^{23}_{q+3}a^{23}_kA^{31}_{q+2}$
      & $e_\alpha(y_1y_3+A^{23}_{q+2}a^{12}_k)$ \\
  \hline \hline
    7 & $1(21)^{2k}$
      & $x_1A^{23}_ka^{23}_k$
      & $y_1A^{31}_ka^{12}_kA^{23}_ka^{23}_k$
      & $e_\alpha(y_2a^{21}_k+y_3a^{31}_k)$ \\
  \hline
    8 & $1(21)^{2k+1}$
      & $y_1(x_2+x_3)A^{31}_ka^{12}_k$
      & $y_1A^{31}_ka^{12}_kA^{23}_{k+1}a^{23}_{k+1}$
      & $e_\alpha(y_1a^{12}_{k+1}+y_1a^{13}_{k+1})$ \\
  \hline  \hline
    9 & $\beta$
      & $y_2y_3a^{21}_qA^{13}_k$
      & $y_2A^{32}_ka^{21}_qA^{13}_ka^{13}_{q+1}$
      & $e_\alpha(x_2+a^{13}_qA^{13}_{k+1})$    \\
  \hline
   10 & $3\beta$
      & $x_3A^{12}_{q+1}a^{12}_k$
      & $y_3a^{31}_kA^{12}_{q+1}a^{12}_kA^{23}_{q+1}$
      & $e_\alpha(y_1y_3+A^{23}_qa^{12}_{k+1})$ \\
  \hline
   11 & $13\beta$
      & $y_1y_2a^{13}_{q+1}A^{32}_k$
      & $y_1A^{21}_ka^{13}_{q+1}A^{32}_ka^{32}_{q+2}$
      & $e_\alpha(x_1+a^{32}_{q+1}A^{32}_{k+1})$    \\
  \hline
   12 & $(213)\beta$
      & $x_2A^{31}_{q+2}a^{31}_k$
      & $y_2a^{23}_kA^{31}_{q+2}a^{31}_kA^{12}_{q+2}$
      & $e_\alpha(y_2y_3+A^{12}_{q+1}a^{31}_{k+1})$ \\
  \hline
   13 & $3(213)\beta$
      & $y_3y_1a^{32}_{q+2}A^{21}_k$
      & $y_3A^{13}_ka^{32}_{q+2}A^{21}_ka^{21}_{q+3}$
      & $e_\alpha(x_3+a^{21}_{q+2}A^{21}_{k+1})$    \\
  \hline
    14 & $13(213)\beta$
      & $x_1A^{23}_{q+3}a^{23}_k$
      & $y_1a^{12}_kA^{23}_{q+3}a^{23}_kA^{31}_{q+3}$
      & $e_\alpha(y_1y_2+A^{31}_{q+2}a^{23}_{k+1})$ \\
  \hline \hline
  \end{tabular}
  \vspace{2mm}
  \caption{\hspace{3mm}For $D^{2,2,2}$: the admissible elements starting at $i_1=1$}
  \label{table_adm_elems}
\begin{equation*}
\begin{split}
 & \text{ In Lines 1--6: } \gamma = (213)^{2p}(21)^{2k}, ~k > 0, ~p \geq 0.
   \text{ In Lines 7--8: } k \geq 0. \\
 & \text{ In Lines 9--14: } \beta =  (213)^{2p}(21)^{2k+1}, ~k \geq 0, ~p \geq 0.
   \text{ In all lines, } q = k + 3p.
 \end{split}
 \end{equation*}
\end{table}
}

\begin{theorem}[On classes of admissible elements]
   \label{th_adm_classes}
   Let $\alpha = i_n{i}_{n-1}\dots{1}$
   be an admissible sequence and $i \neq i_n$. Then
   ${i}\alpha$ is admissible and, for
    $z_\alpha = f_\alpha, e_\alpha, g_{\alpha0}$ from
    Table \ref{table_adm_elems}, the following relation holds:
 \begin{equation} \label{adm_classes}
    \varphi_i\Phi^+\rho(z_\alpha) = \rho(z_{i\alpha}).
 \end{equation}
\end{theorem}

  For the proof of Theorem \ref{th_adm_classes},
  see Appendix \ref{sect_proof_adm}.

\section{The inclusion theorem and the cumulative polynomials}
 \label{sect_inclusion}
\begin{theorem} \label{inclusion}
  For every admissible sequences\hspace{0.3mm}
  ${\alpha}$, ${\alpha}i$, $(i = 1,2,3)$ from
  Table \ref{table_admissible}, the following inclusion
  holds\footnote{Please do not confuse polynomials $e_{{\alpha}i}$
  considered here with polynomials $e_{i{\alpha}}$ considered,
  for example, in Theorem \ref{th_adm_classes}.}
  \begin{equation*}
      e_{{\alpha}i} \subseteq g_{\alpha0}, \quad i=1,2,3.  \vspace{3mm}
  \end{equation*}
\end{theorem}

  For the proof of Theorem \ref{inclusion},
  see \S\ref{proof_incl_th}.

It is easy to check by the definitions from Table
 \ref{table_adm_elems} that
\begin{equation} \label{incl_1}
    f_\alpha \subseteq e_\alpha,
    \quad g_{\alpha0} \subseteq e_\alpha.
\end{equation}
   So, from Theorem  \ref{inclusion} ({\it Inclusion Theorem})
   and inclusions \ref{incl_1} we get
\begin{equation} \label{incl_2}
    f_{{\alpha}i} \subseteq e_{{\alpha}i} \subseteq g_{\alpha0}.
\end{equation}
 \index{admissible element(=polynomial)}

  Recall the definition of {\it cumulative polynomials} from
  \S\ref{def_cumul}.
  {\it Cumulative polynomials} of length $n$ are constructed as sums of
all admissible elements of the same length $n$, where $n$ is the
length of the multi-index:
\begin{equation*}  \label{def_cumu1}
\begin{split}
& x_t(n) = \sum{f}_{i_{n}\dots{i}_{2}t}, \indent t = 1,2,3, \\
& y_t(n) = \sum{e}_{i_{n}\dots{i}_{2}t}, \indent t = 1,2,3, \\
& x_0(n) = \sum{g}_{i_{n}\dots{i}_{2}0}.
\end{split}
\end{equation*}
From (\ref{incl_2}) and Theorem \ref{th_adm_classes} the next
inclusions take place:
\begin{corollary} \label{corollary_cumulative}
\begin{equation}
   x_t(n) \subseteq  y_t(n) \subseteq x_0(n),
   \indent t = 1,2,3.
\end{equation}
\end{corollary}
Thus, {\it cumulative} polynomials satisfy the same inclusions as
the corresponding generators in $D^{2,2,2}$.

\begin{proposition}[On the cumulative polynomials]
  \label{cumul_polyn}
If $z_t(n)$ is one of the {\it cumulative polynomials}
 $x_t(n), y_t(n) (t = 1,2,3)$ or $x_0(n)$,
then
\begin{equation}
 \label{zt_n}
\sum_\text{$i = 1,2,3$}\varphi_i\Phi^+\rho(z_t(n)) =
\rho(z_t(n+1)).
\end{equation}
\end{proposition}
\PerfProof  By Theorem \ref{th_adm_classes}
\begin{equation*}
\varphi_i\Phi^+\rho(z_\alpha) = \rho(z_{i\alpha})
 \text{ for every }
 z_\alpha = f_\alpha, e_\alpha, g_{\alpha0}.
\end{equation*}
The latter relation is true only under the condition that
$i\alpha$ is {\it admissible}, i.e., $i \neq i_n$ for
 $\alpha = i_n{i}_{n-1}\dots{i}_1$. By (\ref{f12}) in the sum of (\ref{zt_n})
we can always exclude one of the elementary maps $\varphi_1,
\varphi_2, \varphi_3$ such that $i \neq i_n$. \qedsymbol

\chapter{\sc\bf Perfect polynomials in $D^{2,2,2}$}
    \label{section_Perfect_Union}

\section{The chains $a_i(n) \subseteq b_i(n) \subseteq c_i(n)$.
         Perfectness and distributivity of $H^+(n)$.}
  \label{chains}
Recall the definition of the elements $a_i(n), b_i(n), c_i(n)$.
The perfect elements (see \S\ref{repres_lat}) are constructed in
the following way, see Definition \ref{def_gen_abc} from
\S\ref{sublattice_Hn}:
\begin{equation}
\begin{split}
  & \underline{\text{For } n = 0}: \\
  & \qquad a_i(0) = y_j + y_k, \\
  & \qquad b_i(0) = x_i + y_j + y_k, \\
  & \qquad c_1(0) = c_2(0) = c_3(0) = \sum{y}_i. \\
  & \underline{\text{For } n \geq 1}:  \\
  & \qquad a_i(n) = x_j(n) + x_k(n) + y_j(n+1) + y_k(n+1), \\
  & \qquad b_i(n) = a_i(n) + x_i(n+1) = x_j(n) + x_k(n) +
       x_i(n+1) + y_j(n+1) + y_k(n+1), \\
  & \qquad c_i(n) = a_i(n) + y_i(n+1) = x_j(n) + x_k(n) +
       y_i(n+1) + y_j(n+1) + y_k(n+1).
\end{split}
\end{equation}

By definition we have
\begin{equation}
    a_i(n) \subseteq b_i(n) \subseteq c_i(n).
\end{equation}
The perfectness of the elements $a_i(0), b_i(0), c_i(0)$ is proved
in Proposition \ref{first_perfect} and Corollary
\ref{second_perfect}. The perfectness of the elements
 $a_i(1), b_i(1), c_i(1)$ is proved in
\S\ref{ai_perfect}-\S\ref{ci_perfect}.

 \begin{proposition}
   \label{perfect_abc}
  1) For every element $v_i(n) = a_i(n)$, $b_i(n)$, $c_i(n),$
   the following relation holds:
 \begin{equation}  \label{vn_n}
    \rho(v_i(n)) = \sum_\text{$k=1,2,3$}\varphi_k\Phi^+\rho(v_i(n-1)).
 \end{equation}

  2) The elements $a_i(n), b_i(n), c_i(n)$ are perfect
   for every $n \geq 1$ and $i = 1,2,3$.
 \end{proposition}

 The proof of Proposition \ref{perfect_abc} is given in
 \S\ref{perfectness_abc_n}.

Let us consider sublattice $H^+(0)$ generated by $3$ chains
\begin{equation}
  s_i = \{ a_i(n) \subseteq b_i(n) \}, \indent i=1,2,3
\end{equation}
and sublattice $H^+(n)$ $(n \geq 1)$ generated by 3 chains
\begin{equation}
  s_i = \{ a_i(n) \subseteq b_i(n) \subseteq c_i(n) \}, \indent i=1,2,3.
\end{equation}

According to definition of perfect elements (see
 \S\ref{repres_lat}), $H^+(n)$ forms a sublattice of perfect
elements in $D^{2,2,2}$.
\begin{proposition}
 \label{distrib_Hn}
  $H^+(n)$ is a distributive sublattice for every $n \geq 0$.
\end{proposition}

 A proof of Proposition \ref{distrib_Hn} is given in
 \S\ref{sect_distrib_Hn}.

\section{The multiplicative form and the cardinality of $H^+(n)$}
\begin{proposition}
  Every element of the lattice $H^+(n)$, where $n \geq 0$, is of the form
  $v_1{v}_2{v}_3$, where ${v}_i \in \{ a_i(n), b_i(n), c_i(n), I_n \}$.
\end{proposition}
\PerfProof Let $v_1{v}_2{v}_3$ and $v_1'{v}_2'{v}_3'$ be two
elements of $H^+(n)$. Then
\begin{equation}
   v_1{v}_2{v}_3 + v_1'{v}_2'{v}_3' =
      \bigcap_{i=1}^3\sup(v_i, v_i').
\end{equation}
 Indeed, if $v_i \supseteq v_i'$ for every $i = 1,2,3$, then
 $v_1{v}_2{v}_3 + v_1'{v}_2'{v}_3' = v_1{v}_2{v}_3$.
 If $v_1 \supseteq v_1'$ and $v_2 \subseteq v_2'$,
 $v_3 \subseteq v_3'$, then by (\ref{distributive_4})
 we have
$$
  v_1{v}_2{v}_3 + v_1'{v}_2'{v}_3' =
  v_1{v}_2'{v}_3'(v_2{v}_3 + v_1') = v_1{v}_2'{v}_3'.
  \qed \vspace{2mm}
$$

\begin{corollary}
  1) The lattice $H^+(0)$ contains $\leq 27$ elements.
  \vspace{2mm}

  2) The lattice $H^+(n)$ for $n>0$ contains $\leq 64$ elements.
  \vspace{2mm}
\end{corollary}

\begin{proposition}
  \label{27_0_distinct_elem}
   The lattice $H^+(0)$ contains $27$ distinct elements.
\end{proposition}
\begin{proposition}
  \label{64_1_distinct_elem}
   The lattice $H^+(1)$ contains $64$ distinct elements.
\end{proposition}

For proof of Propositions \ref{27_0_distinct_elem}  and
\ref{64_1_distinct_elem}, see \S\ref{proofs_27_64}.

Proposition \ref{64_1_distinct_elem} is true for arbitrary $n$:

\begin{proposition}
  \label{contains_64}
   The lattice $H^+(n)$ for $n \in
\{1,2,3,\dots \}$  contains $64$ distinct elements.
\end{proposition}

\PerfProof Let $H^+(n)$ contain $64$ distinct elements. We will
prove that $H^+(n+1)$ also contains $64$ distinct elements. Let
$v_1(n){v}_2(n){v}_3(n)$ and
 $v_1'(n){v}_2'(n){v}_3'(n)$ be two distinct elements in
the $H^+(n)$ and $\rho$ is the preprojective representation
separating them. For $n=1$, the representation $\rho$ is given by
Tables \ref{table_char_fun1_gen}, \ref{table_char_fun1}. For
example,
\begin{equation} \label{separate_1}
\begin{split}
  & \rho(v_1(n){v}_2(n){v}_3(n)) = 0,    \\
  & \rho(v_1'(n){v}_2'(n){v}_3'(n)) = X_0^1.
\end{split}
\end{equation}
Therefore,
\begin{equation} \label{separate_2}
   \text{ there exists an } i \text{ such that } \rho(v_i(n)) = 0,
\end{equation}
and
\begin{equation} \label{separate_21}
   \rho(v_i'(n)) = X_0^1 \text{ for every } i.
\end{equation}
\indent For every preprojective representation $\rho$, there
exists preprojective representation $\tilde{\rho}$ such that
$\Phi^+\tilde{\rho}$ = $\rho$. (Naturally, $\tilde{\rho}$ =
$\Phi^-{\rho}$, see \S\ref{sect_proj_repr}, \cite{BGP73}). Then,
by (\ref{vn_n}) for $i$ from (\ref{separate_2}) we have
\begin{equation} \label{sum_vin}
     \tilde{\rho}(v_i(n+1)) =
     \sum_\text{$p = 1,2,3$}\varphi_p\Phi^+\tilde{\rho}(v_i(n)) =
     \sum_\text{$p = 1,2,3$}\varphi_p\rho(v_i(n)) = 0,
\end{equation}
i.e.,
\begin{equation} \label{separate_3}
  \tilde{\rho}(v_1(n+1){v}_2(n+1){v}_3(n+1)) = 0.
\end{equation}
On the other hand,
\begin{equation} \label{sum_x01}
     \tilde{\rho}(v_i'(n+1)) =
     \sum_\text{$p = 1,2,3$}\varphi_p\Phi^+\tilde{\rho}(v_i'(n)) =
     \sum_\text{$p = 1,2,3$}\varphi_p\rho(v_i'(n)) =
     \sum_\text{$p = 1,2,3$}\varphi_p{X}_0^1.
\end{equation}
By Corollary \ref{cor_psi} we have
\begin{equation} \label{sum_yi}
\begin{split}
&  \sum_\text{$p = 1,2,3$}\varphi_p{X}_0^1 =
   \sum_\text{$p = 1,2,3$}\varphi_p\Phi^+\tilde{\rho}(I) = \\
&  \tilde{\rho}(y_1(y_2 + y_3) +
     y_2(y_1 + y_3) + y_3(y_1 + y_2)) =
     \tilde{\rho}(\bigcap_{i \neq j}(y_i + y_j)).
\end{split}
\end{equation}
 Since $y_i + y_j$ is perfect and
 $\tilde{\rho} = \Phi^-\rho$ is non-projective, by
 Proposition \ref{first_perfect} we have
 \begin{equation} \label{yi_yj_x0}
     \tilde{\rho}(y_i + y_j) = X_0
     \text{ for every } \{i,j\}, i \neq j.
\end{equation}
 Thus,  by
(\ref{sum_x01}), (\ref{sum_yi})  and \ref{yi_yj_x0}
\begin{equation} \label{rho_v_is_X0}
    \tilde{\rho}(v_i'(n+1)) = X_0 \text{ for every } i,
\end{equation}
and therefore
\begin{equation} \label{separate_4}
 \tilde{\rho}(v_1'(n+1){v}_2'(n+1){v}_3'(n+1)) = X_0.
\end{equation}
Thus, by (\ref{separate_3}) and (\ref{separate_4}) the
representation $\tilde{\rho}$ separates $v_1(n+1)v_2(n+1)v_3(n+1)$
and $v_1'(n+1)v_2'(n+1)v_3'(n+1)$. \qedsymbol

\begin{corollary}
  $H^+(n)$ is a distributive $64$-element lattice for $n \geq 1$.
  It is the direct product $(Fig. \ref{cubic64})$ of three chains
  $\{ a_i(n) \subseteq b_i(n) \subseteq c_i(n) \}$.
\end{corollary}

\section{The additive form for elements of $H^+(n)$.
    Considerations$\mod\theta$}
 \label{additive_form}

\section{Sums and intersections}
\begin{proposition}
  \label {every_perfect}
   1) For every perfect elements $v_i, v_j, v_k$,
   the following relation holds
\begin{equation}  \label{sum_cap}
        \sum_\text{$p = 1,2,3$}\varphi_p\Phi^+\rho(\bigcap_\text{$t=i,j,k$}v_t) =
        \bigcap_\text{$t=i,j,k$}(\sum_\text{$p = 1,2,3$}\varphi_p\Phi^+\rho(v_t)).
\end{equation}
   2) For every perfect element v and
  every u, the following relation holds
\begin{equation}  \label{sum_cap_uv}
        \sum_\text{$p = 1,2,3$}\varphi_p\Phi^+\rho(vu) =
        \sum_\text{$p = 1,2,3$}\varphi_p\Phi^+\rho(v)
        \sum_\text{$p = 1,2,3$}\varphi_p\Phi^+\rho(u).
\end{equation}
\end{proposition}
\PerfProof \indent 1) If $\Phi^+\rho(v_t) = 0$ for some $t$, then
\begin{equation*}
   \Phi^+\rho(\bigcap_\text{$t=i,j,k$}v_t) = 0
    \text{ and }
   \sum_\text{$p = 1,2,3$}\varphi_p\Phi^+\rho(v_t) = 0
\end{equation*}
and (\ref{sum_cap}) is true. If $\Phi^+\rho(v_t) = X_0^1$ for
every $t$, then by Corollary \ref{cor_psi} and (\ref{sum_yi}) we
have
\begin{equation}  \label{sum_cap_1}
\begin{split}
  & \Phi^+\rho(\bigcap_\text{$t=i,j,k$}v_t) = X_0^1 =
          \Phi^+\rho(I), \\
  & \sum_\text{$p = 1,2,3$}\varphi_p\Phi^+\rho(v_t) =
   \rho(\bigcap_\text{$i \neq j$}(y_i + y_j)).
\end{split}
\end{equation}
Hereafter
$$
  \bigcap_\text{$i \neq j$}(y_i + y_j) \text{ means }
  (y_1 + y_2)(y_1 + y_3)(y_2 + y_3).
$$
From (\ref{sum_cap_1}) we have
\begin{equation}  \label{sum_cap_2}
\begin{split}
  & \sum_\text{$p = 1,2,3$}\varphi_p\Phi^+\rho(\bigcap_\text{$t=i,j,k$}v_t) =
         \sum_\text{$p = 1,2,3$}\varphi_p\Phi^+\rho(I) =
   \rho(\bigcap_\text{$i \neq j$}(y_i + y_j)), \\
  & \bigcap_\text{$t=i,j,k$}(\sum_\text{$p = 1,2,3$}\varphi_p\Phi^+\rho(v_t)) =
   \rho(\bigcap_\text{$i \neq j$}(y_i + y_j)),
\end{split}
\end{equation}
and again (\ref{sum_cap}) is true. Note that by Proposition
\ref{first_perfect} in (\ref{sum_cap_2})
 $\rho(\bigcap(y_i + y_j)) = X_0$ except for the case
where $\rho$ is a {\it projective} representation (Table
\ref{proj_repr_E6}), for which
 $\rho(\bigcap(y_i + y_j))= 0$. \qedsymbol

 2) As in 1), if $\Phi^+\rho(v) = 0$, then
\begin{equation*}
   \Phi^+\rho(vu) = 0
    \text{ and }
   \sum_\text{$p = 1,2,3$}\varphi_p\Phi^+\rho(v) = 0
\end{equation*}
and (\ref{sum_cap_uv}) is true.

If $\rho$ is {\it projective}, then $\Phi^+\rho = 0$ and
(\ref{sum_cap_uv}) is true as well.

Consider the case $\Phi^+\rho \neq 0$, i.e., $\rho$ is not
 a projective representation. If $\Phi^+\rho(v) = X^1_0$, then
 $\Phi^+\rho(vu) = \Phi^+\rho(u)X^1_0 = \Phi^+\rho(u)$ and
\begin{equation}  \label {sum_cap_uv_1}
 \sum_\text{$p = 1,2,3$}\varphi_p\Phi^+\rho(vu) =
 \sum_\text{$p = 1,2,3$}\varphi_p\Phi^+\rho(u).
\end{equation}

On the other hand,
\begin{equation}  \label {sum_cap_uv_2}
  \sum_\text{$p = 1,2,3$}\varphi_p\Phi^+\rho(v) =
  \sum_\text{$p = 1,2,3$}\varphi_p\Phi^+\rho(I) =
   \rho(\bigcap_\text{$i \neq j$}(y_i + y_j))
\end{equation}
and by Proposition \ref{first_perfect} we have
$\sum\varphi_p\Phi^+\rho(v) = X_0$ and
\begin{equation}  \label {sum_cap_uv_3}
  \sum_\text{$p = 1,2,3$}\varphi_p\Phi^+\rho(v)
  \sum_\text{$p = 1,2,3$}\varphi_p\Phi^+\rho(u) =
  (\sum_\text{$p = 1,2,3$}\varphi_p\Phi^+\rho(u))X_0 =
  \sum_\text{$p = 1,2,3$}\varphi_p\Phi^+\rho(u).
\end{equation}
So, (\ref{sum_cap_uv}) is true. \qedsymbol

 All considerations given below up to the end of
 \S\ref{additive_form} are represented {\it modulo linear
equivalence}, see \S\ref{repres_lat}.
\begin{proposition}
  \label{prop_x_cap}
 The next relation takes place modulo linear equivalence
\begin{equation} \label{x_cap_a_0}
      x_0(n+2) = \bigcap_\text{$i=1,2,3$}a_i(n) \qquad \mod\theta.
\end{equation}
\end{proposition}
\PerfProof We will prove (\ref{x_cap_a_0}) for $n=0,1$ without
restriction$\mod\theta$.

 \underline{Case $n = 0$}. By
Table \ref{table_adm_elems} , Line 7, $k = 0$,  we have
\begin{equation} \label{x_cap_a}
  x_0(2) = g_{10} + g_{20} + g_{30} =
  \sum_\text{$r,s,t$}y_r(y_s + y_t) = \bigcap_\text{$r \neq s$}(y_r + y_s).
\end{equation}
The sum is taken over all permutations $\{r,s,t\}$ of $\{1,2,3\}$.
On the other hand, by Table \ref{table_char_fun0}
\begin{equation}
 \label{case_0a}
 \bigcap_\text{$i=1,2,3$}a_i(0) = \bigcap_{r \neq s}(y_r + y_s).
\end{equation}

 \underline{Case $n = 1$}. We have
\begin{equation} \label{x_cap_a_1}
  x_0(3) = g_{210} + g_{120} + g_{310} + g_{130} + g_{230} + g_{320}.
\end{equation}
By Table \ref{table_adm_elems} , Line 9, $k = 0, p = 0, q = 0$,
we have, for example,
$$
 g_{210} = y_2{a}^{13}_1(x_2 + A^{13}_1) = y_2(x_1 + y_3)(x_2 + x_3 + y_1).
$$
From here,
\begin{equation*}
\begin{split}
  & g_{210} + g_{310} = y_2(x_1 + y_3)(x_2 + x_3 + y_1) +
    y_3(x_1 + y_2)(x_2 + x_3 + y_1) = \\
  &  (x_1 + y_3)(x_1 + y_2)
     [y_2(x_2 + x_3 + y_1) + y_3(x_2 + x_3 + y_1)] = \\
  &  (x_1 + y_3)(x_1 + y_2)[x_2 + y_2(x_3 + y_1) +
      x_3 + y_3(x_2 + y_1)] = \\
  &  (x_1 + y_3)(x_1 + y_2)
     [(x_3 + y_2)(x_3 + y_1) + (x_2 + y_3)(x_2 + y_1)].
\end{split}
\end{equation*}
Thus,
\begin{equation}
\begin{split}
   & g_{210} + g_{310} =
     (x_1 + y_3)(x_1 + y_2)
     [(x_3 + y_2)(x_3 + y_1) +
      (x_2 + y_3)(x_2 + y_1)], \\
   & g_{120} + g_{320} =
     (x_2 + y_3)(x_2 + y_1)
     [(x_3 + y_1)(x_3 + y_2) +
      (x_1 + y_3)(x_1 + y_2)], \\
   & g_{130} + g_{230} =
     (x_3 + y_2)(x_3 + y_1)[(x_2 + y_1)(x_2 + y_3) +
      (x_1 + y_2)(x_1 + y_3)].
\end{split}
\end{equation}
On the other hand, by (\ref{designation_t}) we have
\begin{equation*}
\begin{split}
  & t_1(1) = x_1(1) + y_1(2) = x_1 + e_{21} + e_{31} = \\
  & x_1 + y_2(x_1 + y_3) + y_3(x_1 + y_2) =
     (x_1 + y_2)(x_1 + y_3),
\end{split}
\end{equation*}
 i.e.,
\begin{equation}
\begin{split}
   & g_{210} + g_{310} = t_1(1)(t_2(1) + t_3(1)), \\
   & g_{120} + g_{320} = t_2(1)(t_1(1) + t_3(1)), \\
   & g_{130} + g_{230} = t_3(1)(t_1(1) + t_2(1)),
\end{split}
\end{equation}
and therefore
\begin{equation}  \label{case_1a}
  x_0(3) = \bigcap_{i \neq j}(t_i(1) + t_j(1)) =
                 \bigcap_{k=1,2,3}a_k(1).
\end{equation}
\indent {\it Induction step}.
  As in (\ref{sum_vin}),
\begin{equation} \label{induct_step_1}
     \tilde{\rho}(x_0(n+1)) =
     \sum_\text{$p = 1,2,3$}\varphi_p\Phi^+\tilde{\rho}(x_0(n)) =
     \sum_\text{$p = 1,2,3$}\varphi_p\Phi^+\tilde{\rho}
              (\bigcap_{i=1,2,3}a_i(n-2)).
\end{equation}
By (\ref{sum_cap}), Proposition \ref{every_perfect}
\begin{equation} \label{induct_step_2}
     \tilde{\rho}(x_0(n+1)) =
     \bigcap_\text{$i=1,2,3$}(\sum_\text{$p = 1,2,3$}\varphi_p\Phi^ +
            \tilde{\rho}a_i(n-2)).
\end{equation}
From here, by (\ref{vn_n}), Proposition \ref{perfect_abc}
\begin{equation} \label{induct_step_3}
     \tilde{\rho}(x_0(n+1)) =
     \bigcap_\text{$i=1,2,3$}\tilde{\rho}(a_i(n-1)),
\end{equation}
which is to be proved. \qedsymbol

\begin{corollary}
 \label{x0_is_perfect}
  For $n \geq 2$, the elements $x_0(n)$ are perfect.
\end{corollary}
\PerfProof It follows from the fact that the intersection of
perfect elements is also perfect.
  \qedsymbol \vspace{4mm}

\begin{conjecture}
 \label{conj_2}
  Relation (\ref{x_cap_a_0}) takes place
  without restriction$\mod\theta$.
\end{conjecture}

For $n=0$, this conjecture was proven above in (\ref{x_cap_a}),
(\ref{case_0a}) and, for $n=1$, it was proven in (\ref{case_1a}).

\section{Other generators:
 $p_i(n) \subseteq q_i(n) \subseteq s_i(n)$} \label{another_gen}
Let us introduce the following polynomials:
\begin{equation}
 \label{another_gen_eq}
 \begin{split}
  &  p_i(n) = x_0(n+2) + x_i(n+1),     \\
  &  q_i(n) = x_0(n+2) + y_i(n+1),     \\
  &  s_i(n) = x_0(n+2) + y_i(n+1) + x_i(n).
 \end{split}
\end{equation}
Compare two generator systems on Fig. \ref{cubic64}. The
generators $a_i(n) \subseteq b_i(n) \subseteq c_i(n)$ are given by
the upper edges of $H^+(n)$. The generators $p_i(n) \subseteq
q_i(n) \subseteq s_i(n)$ are given by the lower edges of $H^+(n)$.
\begin{proposition} \label{two_set_gen}
  Table \ref{mult_add_forms} gives
  relations between the upper and lower generators.
  These relations are true$\mod\theta$.
\end{proposition}
\begin{table}[h]
 \renewcommand{\arraystretch}{1.35}
 \begin{tabular} {||c||c|c|c||}
  \hline \hline
    $N$     & The multiplicative
            & The additive
            & The number of \cr
            & form
            & form
            & Polynomials \\
  \hline  \hline
   1 & $a_i(n)$ & $s_j(n) + s_k(n)$ & 3 \\
  \hline
   2 & $b_i(n)$
     & $s_j(n) + s_k(n) + p_i(n)$ & 3 \\
  \hline
   3 & $c_i(n)$
     & $s_j(n) + s_k(n) + q_i(n)$ & 3 \\
  \hline \hline
   4 & $a_i(n)a_j(n)$
     & $s_k(n)$  & 3 \\
  \hline
   5 & $b_i(n)b_j(n)$
     & $s_k(n) + p_i(n) + p_j(n)$  & 3 \\
  \hline
   6 & $c_i(n)c_j(n)$
     & $s_k(n) + q_i(n) + q_j(n)$  & 3 \\
  \hline
   7 & $a_i(n)b_j(n)$
     & $s_k(n) + p_j(n)$  & 6 \\
  \hline
   8 & $a_i(n)c_j(n)$
     & $s_k(n) + p_j(n)$  & 6 \\
  \hline
   9 & $b_i(n)c_j(n)$
     & $s_k(n) + p_i(n)+ q_j(n)$  & 6 \\
  \hline \hline
   10 & $a_i(n)a_j(n)b_k(n)$
     & $p_k(n)$  & 3 \\
  \hline
   11 & $a_i(n)a_j(n)c_k(n)$
     & $q_k(n)$  & 3 \\
  \hline
   12 & $b_i(n)b_j(n)a_k(n)$
     & $p_i(n) + p_j(n)$  & 3 \\
  \hline
   13 & $b_i(n)b_j(n)c_k(n)$
     & $p_i(n) + p_j(n) + q_k(n)$  & 3 \\
  \hline
   14 & $c_i(n)c_j(n)a_k(n)$
     & $q_i(n) + q_j(n)$  & 3 \\
  \hline
   15 & $c_i(n)c_j(n)b_k(n)$
     & $q_i(n) + q_j(n) + p_k(n)$ & 3 \\
  \hline
   16 & $a_i(n)b_j(n)c_k(n)$
     & $p_i(n) + q_j(n)$ & 6 \\
  \hline
   17 & $a_i(n)a_j(n)a_k(n)$
     & $x_0(n+2)$ & 1 \\
  \hline
   18 & $b_i(n)b_j(n)b_k(n)$
     & $p_i(n) + p_j(n) + p_k(n)$ & 1 \\
  \hline
   19 & $c_i(n)c_j(n)c_k(n)$
     & $q_i(n) + q_j(n) + q_k(n)$ & 1 \\
  \hline \hline
   20 & $I_n = a_i(n) + a_j(n)$
      & $s_i(n) + s_j(n) + s_k(n)$ & 1  \cr
      &  \text{ for all } $\{i,j\}$, $i \neq j$
      &   & \\
  \hline \hline
  \end{tabular}
  \vspace{2mm}
  \caption{\hspace{3mm}The multiplicative and additive forms of $H^+(n)$}
  \label{mult_add_forms}
 \end{table}
 For a proof, see \S\ref{sect_upper_lower}.

\section{Coupling together $H^+(n)$ and $H^+(n+1)$.
        Considerations$\mod\theta$} \label{joint_hh}
We consider ordering of $H^+(n)$ and $H^+(n+1)$ relative to each
other. The way $H^+(n)$ is connected with $H^+(n+1)$ is not as
simply as free modular lattice ${\it D}^r$ in \cite{GP76},
\cite{GP77}, see \S\ref{cubes}, \S\ref{union_L6}.
 For the way $H^+(n)$ is connected with $H^+(n+1)$,
 see the Hasse diagram on Fig. \ref{cube_comparison}.

\section{Boolean algebras $U_n$ and $V_{n+1}$}

Consider the lower cube $U_n \subseteq H^+(n)$ and the upper cube
$V_{n+1} \subseteq H^+(n+1)$:
\begin{equation}
\begin{split}
& U_n = \{b_1b_2b_3, a_1b_2b_3, b_1a_2b_3, b_1b_2a_3,
              a_1a_2b_3, a_1b_2a_3, b_1a_2a_3, a_1a_2a_3\}, \\
& V_{n+1} = \{\tilde{c}_1, \tilde{c}_2, \tilde{c}_3,
      \tilde{c}_1\tilde{c}_2, \tilde{c}_2\tilde{c}_3, \tilde{c}_1\tilde{c}_3,
                    \tilde{c}_1\tilde{c}_2\tilde{c}_3,
                    \tilde{c}_1 + \tilde{c}_2 + \tilde{c}_3  \},
\end{split}
\end{equation}
where $a_i = a_i(n)$, $b_i = b_i(n)$,
      $\tilde{c}_i = c_i(n+1)$.

\begin{table}[h]
 \renewcommand{\arraystretch}{1.35}
 \begin{tabular} {||c||c|c|c||}
  \hline \hline
    $N$ & The multiplicative
            & The additive
            & The complement \cr
            & form of $z$
            & form of $z$
            & element $z'$ \\
  \hline  \hline
   1 & $I_{U_n} = b_1(n)b_2(n)b_3(n)$
     & $p_1(n) + p_2(n) + p_3(n)$
     & $a_1(n)a_2(n)a_3(n)$ \\
  \hline
   2 & $a_1(n)b_2(n)b_3(n)$
     & $p_2(n) + p_3(n)$
     & $b_1(n)a_2(n)a_3(n)$ \\
  \hline
   3 & $b_1(n)a_2(n)b_3(n)$
     & $p_1(n) + p_3(n)$
     & $a_1(n)b_2(n)a_3(n)$ \\
  \hline
   4 & $b_1(n)b_2(n)a_3(n)$
     & $p_1(n) + p_2(n)$
     & $a_1(n)a_2(n)b_3(n)$ \\
  \hline
   5 & $a_1(n)a_2(n)b_3(n)$
     & $p_3(n)$
     & $b_1(n)b_2(n)a_3(n)$ \\
  \hline
   6 & $a_1(n)b_2(n)a_3(n)$
     & $p_2(n)$
     & $b_1(n)a_2(n)a_3(n)$ \\
  \hline
   7 & $b_1(n)a_2(n)b_3(n)$
     & $p_1(n)$
     & $a_1(n)b_2(n)a_3(n)$ \\
  \hline
   8 & $O_{U_n} = a_1(n)a_2(n)a_3(n)$
     & $x_0(n+2)$
     & $b_1(n)b_2(n)b_3(n)$ \\
  \hline \hline
  \end{tabular}
  \vspace{2mm}
  \caption{
   \hspace{3mm}The $8$-element Boolean algebra $U_n$}
  \label{8elem_U}
 \end{table}

In Table \ref{8elem_U}, the elements $I_{U_n}$ and $O_{U_n}$ are,
 respectively, the maximal and minimal elements of $U_n$.
 The distributivity of $U_n$ follows from the distributivity of $H^+(n)$,
 Proposition \ref{distrib_Hn}. It is easy to see that $U_n$ is an $8$-element
Boolean lattice, where

\begin{equation}
\begin{split}
  &   O_{U_n} = a_1(n)a_2(n)a_3(n),
      \hspace{1cm}
      I_{U_n} = b_1(n)b_2(n)b_3(n), \\
  & ({a_ia_jb_k})' = b_ib_ja_k, \text{ i.e., }
      a_ia_jb_k + b_ib_ja_k = I_{U_n} \text{ and }
      a_ia_jb_k \cap b_ib_ja_k = O_{U_n}.
\end{split}
\end{equation}
\footnotesize
\begin{table}[h]
 \renewcommand{\arraystretch}{1.35}
 \begin{tabular} {||c||c|c|c||}
  \hline \hline
    $N$     & The multiplicative
            & The additive
            & The complement \cr
            & form of $z$
            & form of $z$
            & element $z'$ \\
  \hline  \hline
   1 & $I_{V_{n+1}}$ =
     & $s_1(n+1) + s_2(n+1) + s_3(n+1)$
     & $c_1(n+1)c_2(n+1)c_3(n+1)$ \cr
     & $c_i(n+1) + c_j(n+1), i \neq j$
     & & \\
  \hline
   2 & $c_1(n+1)$
     & $q_1(n+1) + s_2(n+1) + s_3(n+1)$
     & $c_2(n+1)c_3(n+1)$ \\
  \hline
   3 & $c_2(n+1)$
     & $s_1(n+1) + q_2(n+1) + s_3(n+1)$
     & $c_1(n+1)c_3(n+1)$ \\
  \hline
   4 & $c_3(n+1)$
     & $s_1(n+1) + s_2(n+1) + q_3(n+1)$
     & $c_2(n+1)c_3(n+1)$ \\
  \hline
   5 & $c_1(n+1)c_2(n+1)$
     & $q_1(n+1) + q_2(n+1) + s_3(n+1)$
     & $c_3(n+1)$ \\
  \hline
   6 & $c_1(n+1)c_3(n+1)$
     & $q_1(n+1) + s_2(n+1) + q_3(n+1)$
     & $c_2(n+1)$ \\
  \hline
   7 & $c_2(n+1)c_3(n+1)$
     & $s_1(n+1) + q_2(n+1) + q_3(n+1)$
     & $c_1(n+1)$ \\
  \hline
   8 & $O_{V_{n+1}}$ =
     & $q_1(n+1) + q_2(n+1) + q_3(n+1)$
     & $c_1(n+1) + c_2(n+1) + $ \cr
     & $c_1(n+1)c_2(n+1)c_3(n+1)$
     & &   $c_3(n+1)$    \\
  \hline \hline
  \end{tabular}
  \vspace{2mm}
  \caption{
   \hspace{3mm}The $8$-element Boolean algebra $V_{n+1}$}
  \label{8elem_V}
 \end{table}
\normalsize

In Table \ref{8elem_V}, the elements $I_{V_{n+1}}$ and
 $O_{V_{n+1}}$ are, respectively, the maximal and minimal elements of $V_{n+1}$.
The distributivity of $V_{n+1}$ follows from the distributivity of
$H^+(n+1)$. It is easy to see that $V_{n+1}$ is an $8$-element
Boolean lattice, where

\begin{equation}
\begin{split}
  &   O_{V_{n+1}} = c_1(n+1)c_2(n+1)c_3(n+1),
      \indent
      I_{V_{n+1}} = c_i(n+1) + c_j(n+1), \\
  & \tilde{c}_i' = \tilde{c}_j\tilde{c}_k, \text{ i.e., }
     \tilde{c}_i + \tilde{c}_i' =
     \tilde{c}_i + \tilde{c}_j\tilde{c}_k = I_{V_{n+1}},
     \text{ and }
     \tilde{c}_i \cap \tilde{c}_i' =
     \tilde{c}_i \cap \tilde{c}_j\tilde{c}_k = O_{V_{n+1}}.
\end{split}
\end{equation}
\begin{figure}[h]
\includegraphics{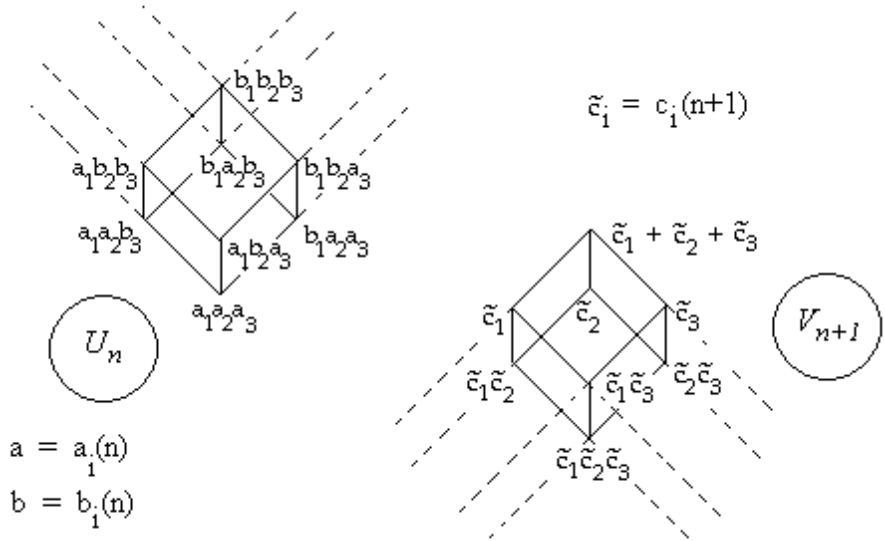}
\caption{\hspace{3mm}The Boolean algebras $U_n$ and $V_{n+1}$}
\label{8elem_U_V}
\end{figure}

\section{The cardinality of $U_n \bigcup V_{n+1}$}

\begin{proposition}
  The union $U_n \bigcup V_{n+1}$ contains $16$ distinct elements.
\end{proposition}

Consider case $n=1$. We will prove that $(\Phi^-)^2\rho_{x_0}$
(see Table \ref{RepOrder12}) separates $U_1$ and $V_2$. It is
sufficient to prove that
\begin{equation}
  \label{separate_U}
  \begin{split}
   & (\Phi^-)^2\rho_{x_0}(z) = X_0 \text{ for every }
    z \subseteq U_1, \\
   & (\Phi^-)^2\rho_{x_0}(z) = 0 \text { for every }
    z \subseteq V_2.
  \end{split}
\end{equation}

For the generators $a_i(1)$ and
 $b_i(1)$\footnote{The generators $a_i(n)$ and $b_i(n)$ are defined in
 (\ref{def_gen_abc}); for $n=1$, these generators are explicitly
 written in (\ref{generators_abc}), (\ref{generators_abc_1})},
 we have
$$
   (\Phi^-)^2\rho_{x_0}(a_i(1)) = X_j + X_k + ... \neq 0.
$$
 The generators $a_i(1), b_i(1)$ are perfect, i.e.,
$$
 (\Phi^-)^2\rho_{x_0}(a_i(1)) = X_0, \indent
 (\Phi^-)^2\rho_{x_0}(b_i(1)) = X_0,
$$
and therefore (\ref{separate_U}) is true.

On the other hand, since
\begin{equation}
  c_i(2) = x_j(2) + x_k(2) + \sum_\text{$q=1,2,3$}y_q(3)
\end{equation}
and $\Phi^+(\Phi^-)^2\rho_{x_0} = \Phi^-\rho_{x_0}$, we have
\begin{equation}
   (\Phi^-)^2\rho_{x_0}(c_i(2)) =
   \sum_\text{$p = 1,2,3$}\varphi_p\Phi^-\rho_{x_0}
      (x_j + x_k + \sum_\text{$q=1,2,3$}y_q(2)).
\end{equation}
According to Table \ref{RepOrder12} we have
 $\Phi^-\rho_{x_0}(x_j)= 0$ for every $j$, i.e.,
\begin{equation}
   (\Phi^-)^2\rho_{x_0}(c_i(2)) =
   \sum_\text{$p = 1,2,3$}\varphi_p\Phi^-\rho_{x_0}(\sum_\text{$q=1,2,3$}y_q(2)).
\end{equation}
By Table \ref{RepOrder12} we have $\rho_{x_0}(y_j) = 0$. Further,
$$
  \Phi^-\rho_{x_0}(y_j(2)) =
  \sum\varphi_q\Phi^+\Phi^-\rho_{x_0}(y_j) =
  \sum\varphi_q\rho_{x_0}(y_j) = 0.
$$
Thus, we have $(\Phi^-)^2\rho_{x_0}(c_i(2)) = 0$ and
(\ref{separate_U}) is true. The induction step is proved as in
Proposition \ref{contains_64} above. \qedsymbol

\section{The Boolean algebra $U_n \bigcup V_{n+1}$}

 \label{sect_bool_alg_16}
\begin{proposition}  \label{UV_16bool_alg}
$U_n \bigcup V_{n+1}$ is a $16$-element Boolean algebra.
\end{proposition}

\PerfProof
According to Proposition \ref{DC_lemma}
it is sufficient to prove relations
(\ref{c_less_d}), (\ref{di_less_ci_dj}), (\ref{dicj_less_ci})
for
\begin{equation*}
    c_i = c_i(n+1), \indent
    d_i = a_i(n){b}_j(n){b}_k(n) = p_j(n) + p_k(n)
\end{equation*}
(see Fig. \ref{8elem_U_V} and Table \ref{8elem_U}), i.e.,
\begin{equation} \label {c_less_d_uv}
      c_i(n+1) \subseteq p_j(n) + p_k(n),
\end{equation}
\begin{equation} \label {di_less_ci_dj_uv}
      p_j(n) + p_k(n) \subseteq c_i(n+1) + p_i(n) + p_k(n),
\end{equation}
\begin{equation} \label {dicj_less_ci_uv}
      (p_j(n) + p_k(n))c_j(n+1) \subseteq c_i(n+1),
\end{equation}
where $\{i,j,k\}$ is a permutation of $\{1,2,3\}$.

According to the definition of perfect elements from
 \S\ref{repres_lat}
\begin{equation*}
   c_i(n+1) = x_j(n+1) + x_k(n+1) + \sum_\text{$t=1,2,3$}y_t(n+2).
\end{equation*}
By definition (\ref{another_gen_eq})
\begin{equation*}
\begin{split}
  & p_j(n) = x_0(n+2) + x_j(n+1),  \\
  & p_j(n) + p_k(n) = x_0(n+2) + x_j(n+1) + x_k(n+1).
\end{split}
\end{equation*}

  1) Then the inclusion (\ref{c_less_d_uv}) follows from
\begin{equation} \label{sum_yt_x0}
   \sum_\text{$t=1,2,3$}y_t(n+2) \subseteq x_0(n+2),
\end{equation}
which, in turn, comes from the Corollary
\ref{corollary_cumulative} concerning {\it cumulative}
polynomials.

 2) Inclusion (\ref{di_less_ci_dj_uv}) follows from the relation
\begin{equation*}
   x_j(n+1) \subseteq c_i(n+1).
\end{equation*}

 3) The last inclusion (\ref{dicj_less_ci_uv}) is a little more
difficult. It is equivalent to
\begin{equation} \label{long_incl}
\begin{split}
  & (p_j(n) + p_k(n))
    (x_i(n+1) + x_k(n+1) + \sum_\text{$t=1,2,3$}y_t(n+2))
    \subseteq \\
  & x_j(n+1) + x_k(n+1) + \sum_\text{$t=1,2,3$}y_t(n+2).
\end{split}
\end{equation}
From (\ref{sum_yt_x0}) and since $x_k(n+1) \subseteq p_j(n)$, we
see that (\ref{long_incl}) is equivalent to
\begin{equation} \label{long_incl_1}
\begin{split}
  & x_k(n+1) + \sum_\text{$t=1,2,3$}y_t(n+2) + (p_j(n) + p_k(n))x_i(n+1)
   \subseteq \\
  & x_j(n+1) + x_k(n+1) + \sum_\text{$t=1,2,3$}y_t(n+2).
\end{split}
\end{equation}
 To prove (\ref{long_incl_1}), it suffices to check that
\begin{equation} \label{long_incl_2}
      (p_j(n) + p_k(n))x_i(n+1)
         \subseteq  y_i(n+2),
\end{equation}
i.e.,
\begin{equation} \label{long_incl_3}
      a_i(n){b}_j(n)b_k(n)x_i(n+1) \subseteq  y_i(n+2).
\end{equation}
 Since $x_i(n+1) \subseteq b_j(n)b_k(n)$
(by relation \ref{c_less_d}), we see that (\ref{long_incl_3}) is
equivalent to
\begin{equation} \label{long_incl_4}
      a_i(n)x_i(n+1) \subseteq  y_i(n+2).
\end{equation}
 We will prove (\ref{long_incl_3}) modulo linear equivalence
\begin{equation}
 \label{long_incl_5}
      a_i(n)x_i(n+1) \subseteq  y_i(n+2) \mod\theta.
\end{equation}
\begin{conjecture}
 \label{conj_3}
  Relation (\ref{long_incl_5})
  takes place without restriction$\mod\theta$.\\
 \end{conjecture}
 So, let us prove (\ref{long_incl_5}). Let $n=1, i=1$. We have to
prove that
\begin{equation} \label{short_incl}
         a_1(1)x_1(2) \subseteq y_1(3).
\end{equation}
According to (\ref{def_gen_abc}) or (\ref{generators_abc}),
(\ref{generators_abc_1}) we have
\begin{equation*}
\begin{split}
   a_1(1) = & x_2 + x_3 +  y_3(x_2 + y_1) + y_1(x_2 + y_3) +
   y_1(x_3 + y_2) + y_2(x_3 + y_1) = \\
  & (x_2 + y_3)(x_2 + y_1) + (x_3 + y_2)(x_3 + y_1).
\end{split}
\end{equation*}
Further, according to \S\ref{examp_cumul}, we have
$$
  x_1(2) = f_{21} + f_{31} = y_2{y}_3
$$
and
\begin{equation*}
\begin{split}
  y_1(3) = & y_3(y_1 + x_2)(y_2 + x_3) + y_1(x_2 + y_3)(x_3 + y_2) +
             y_2(y_1 + x_3)(y_3 + x_2) = \\
  & (x_2 + y_3)(x_3 + y_2)[y_3(y_1 + x_2)(y_2 + x_3) +
        y_1 + y_2(y_1 + x_3)(y_3 + x_2)] = \\
  & (x_2 + y_3)(x_3 + y_2)[(y_1 + x_2)(x_3 + y_2{y}_3) + y_1 +
    (y_1 + x_3)(x_2 + y_2{y}_3)] = \\
  & (x_2 + y_3)(x_3 + y_2)
    [(y_1 + x_2)(x_3 + y_1 + y_2{y}_3) +
     (y_1 + x_3)(x_2 + y_1 + y_2{y}_3)] = \\
  & (x_2 + y_3)(x_3 + y_2)(x_3 + y_1 + y_2{y}_3)
    (x_2 + y_1 + y_2{y}_3)(y_1 + x_2 + x_3).
\end{split}
\end{equation*}
Now, (\ref{short_incl}) is true since
$$
  x_1(2) = y_2{y}_3 \subseteq (x_2 + y_3)(x_3 + y_2)
      (x_3 + y_1 + y_2{y}_3)(x_2 + y_1 + y_2{y}_3)
$$
and
$$
  a_1(1) =
  (x_2 + y_3)(x_2 + y_1) + (x_3 + y_2)(x_3 + y_1)
  \subseteq y_1 + x_2 + x_3.
$$
Therefore, (\ref{short_incl}) is also true$\mod\theta$, i.e.,
(\ref{long_incl_5}) is true for $n = 1$.

{\it Induction step}. Let (\ref{long_incl_5}) be true for $n$,
i.e.,
\begin{equation} \label{long_incl_6}
      \Phi^+\rho(a_i(n){x}_i(n+1))  \subseteq \Phi^+\rho(y_i(n+2))
\end{equation}
for every indecomposable representation $\rho$. Applying
$\sum\varphi_i$ we have
\begin{equation} \label{long_incl_7}
  \sum_\text{$p = 1,2,3$}\varphi_p\Phi^+\rho(a_i(n)x_i(n+1))
  \subseteq
  \sum_\text{$p = 1,2,3$}\varphi_p\Phi^+\rho(y_i(n+2)).
\end{equation}
Since $a_i(n)$ is perfect, we have according to (\ref{sum_cap_uv})
and Proposition \ref{every_perfect}
\begin{equation} \label{long_incl_8}
     \sum_\text{$p = 1,2,3$}\varphi_i\Phi^+\rho(a_i(n))
     \sum_\text{$p = 1,2,3$}\varphi_i\Phi^+\rho(x_i(n+1))
     \subseteq
     \sum_\text{$p = 1,2,3$}\varphi_i\Phi^+\rho(y_i(n+2)).
\end{equation}
By Propositions \ref{cumul_polyn} and \ref{perfect_abc}
\begin{equation} \label{long_incl_9}
      \rho(a_i(n+1)x_i(n+2))  \subseteq \rho(y_i(n+3)).
\end{equation}
So, the induction step and Proposition \ref{UV_16bool_alg} are
proven. \qedsymbol

\begin{corollary} [Connection edges]
\label{connection_edges}
   The following $8$ inclusions hold:
\begin{equation}  \label{cor_incl_1}
  c_i(n+1) \subseteq a_i(n)b_j(n)b_k(n), \indent i = 1,2,3,
\end{equation}
\begin{equation}  \label{cor_incl_2}
  c_i(n+1)c_j(n+1) \subseteq a_i(n)a_j(n)b_k(n),
  \indent k = 1,2,3,
\end{equation}
\begin{equation}   \label{cor_incl_3}
  \bigcap_{i=1,2,3}c_i(n+1)
         \subseteq \bigcap_{i=1,2,3}a_i(n),
\end{equation}
\begin{equation}    \label{cor_incl_4}
  \sum_\text{$i=1,2,3$}c_i(n+1)
         \subseteq \bigcap_{i=1,2,3}b_i(n).
\end{equation}
  The $8$ inclusions (\ref{cor_incl_1}) -- (\ref{cor_incl_4})
  correspond to the $8$ edges shown on
  Fig. \ref{boolean16} or Fig. \ref{16elem_C_D}.
\end{corollary}
\PerfProof  Inclusions (\ref{cor_incl_1}) follow from
(\ref{c_less_d_uv}), inclusions (\ref{cor_incl_2}),
(\ref{cor_incl_3}) and (\ref{cor_incl_4}) follows from
(\ref{cor_incl_1}), see Table \ref{mult_add_forms}. \qedsymbol

\section{The theorem on perfect elements. Considerations$\mod\theta$}
 \label{sect_merge_constr}
We will prove the theorem describing connection construction of
sublattices $H^+(n)$.
\begin{proposition}
 \label{merge_1}
  1) For every $u \in U_n$ and every $v \in H^+(n+1)$, we have
    $u + v \in U_n$.

  2) For every $u \in H^+(n)$ and every $v \in V_{n+1}$, we have
    $uv \in V_{n+1}$.
\end{proposition}
\PerfProof  1) According to Table \ref{8elem_U}, the elements of
$U_n$ are
\begin{equation}  \label{generators_Un}
   p_i(n), \indent
   p_i(n) + p_j(n), \indent
   \sum{p}_i(n), \indent x_0(n+2).
\end{equation}
 It is sufficient to verify heading 1) for the generators
 $p_i(n+1), q_i(n+1) , s_i(n+1)$ of
 $H^+(n+1)$. According to definition (\ref{another_gen_eq})
 it is sufficient to verify that
 adding the elements
 $$
    x_0(n+3), \quad x_i(n+2), \quad y_i(n+2)
 $$
 does not lead out of the elements of type
(\ref{generators_Un}).
 According to (\ref{incl_1}) and Corollary (\ref{corollary_cumulative})
 we have
\begin{equation}  \label{addition_1}
   x_0(n+3) \subseteq y_i(n+2) \subseteq x_0(n+2)
\end{equation}
and adding the element $x_0(n+3)$ does not lead out of $U_n$. From
(\ref{addition_1}) it follows that
\begin{equation}  \label{addition_2}
  p_i(n) = x_0(n+2) + x_i(n+1)  \supseteq  y_i(n+2)
    \supseteq  x_i(n+2).
\end{equation}
Therefore adding $y_i(n+2)$ and $x_i(n+2)$ does not lead out of
 $U_n$. Finally,
\begin{equation}
  p_i(n)  \supseteq x_i(n+1), \quad
  p_j(n) + x_i(n+1) = p_i(n) + p_j(n),
\end{equation}
i.e., addition $x_i(n+1)$ does not go lead of the $U_n$ either.
\qedsymbol

 2) Every
 $u \in H^+(n)$ is of the form $u_1{u}_2{u}_3$, where  (Fig. \ref{cubic64})
\begin{equation*}
 u_i \in \{ a_i(n) \subseteq
                 b_i(n) \subseteq c_i(n) \subseteq I \}.
\end{equation*}
If $u_i \neq a_i(n)$ for every $i$, then $u \supseteq
\bigcap{b}_i(n) \supseteq \sum{c}_i(n+1) \supseteq v$ and $uv =
v$. So, it is sufficient to consider the elements
 $u \in H^+(n)$ of the type $a_i(n){u}'$,
 $a_i(n){a}_j(n){u}'$ and $a_i(n){a}_j(n)a_k(n)$,
 where $u' \supseteq v$, i.e., it is
sufficient to prove that
\begin{equation}
 \label {va_v}
  va_i(n) \in V_{n+1} \text{ for every }
   v \subseteq  V_{n+1}.
\end{equation}
By Corollary \ref{connection_edges} we have
 $c_i(n+1) \subseteq a_i(n)$ and
\begin{equation} \label{ciai_ci}
  c_i(n+1)a_i(n) = c_i(n+1),
\end{equation}
and so (\ref{va_v}) is proven for $v = c_i(n+1)$.

Now we consider the case $i \neq j$, for example,
$a_2(n)c_1(n+1)$. Further, we have
\begin{equation}
\begin{split}
 & c_1(n+1) = x_3(n+1) + x_2(n+1) + \sum_\text{$i=1,2,3$}y_i(n+2), \\
 & x_3(n+1) \subseteq  a_2(n), \\
 & \sum_\text{$i=1,2,3$}y_i(n+2) \subseteq x_0(n+2) =
     \bigcap_\text{$i=1,2,3$}a_i(n) \subseteq a_2(n).
\end{split}
\end{equation}
The last relation follows from the (\ref{x_cap_a_0}) and it is
true$\mod\theta$. Therefore
\begin{equation} \label{a2n_c1}
  c_1(n+1)a_2(n) =
  x_3(n+1) + \sum_\text{$i=1,2,3$}y_i(n+2) + x_2(n+1)a_2(n).
\end{equation}
According to (\ref{long_incl_5}) we have
 $x_2(n+1)a_2(n) \subseteq c_2(n+1)$ and
\begin{equation*}
\begin{split}
  & c_1(n+1)a_2(n) \subseteq
    x_3(n+1) + \sum_\text{$i=1,2,3$}y_i(n+2) +
    x_2(n+1)c_2(n+1) = \\
  & c_2(n+1)[x_3(n+1) + \sum_\text{$i=1,2,3$}y_i(n+2) +
    x_2(n+1)] =   \\
  & c_2(n+1)c_1(n+1) \subseteq a_2(n)c_1(n+1).
\end{split}
\end{equation*}
Thus,
\begin{equation}  \label{c1a2_c1c2}
  c_1(n+1)a_2(n) = c_1(n+1)c_2(n+1)
\end{equation}
and (\ref{va_v}) is proven for $v = c_j(n+1)$, $i \neq j$.
Finally, by (\ref{ciai_ci}), (\ref{c1a2_c1c2}) we have
\begin{equation*}
\begin{split}
  & a_j(n)\sum_\text{$i=1,2,3$}c_i(n+1) = a_j(n)[c_i(n+1) + c_j(n+1)] = \\
  & c_j(n+1) + c_i(n+1)a_j(n)  = c_j(n+1) + c_i(n+1)c_j(n+1),
\end{split}
\end{equation*}
i.e.,
\begin{equation}
   a_j(n)\sum_\text{$i=1,2,3$}c_i(n+1) = c_j(n+1).
 \qed \vspace{2mm}
\end{equation}

\begin{theorem} \label{main_perf_th}
 The union \UnionZeroInf$H^+(n)$ is a distributive
 lattice$\mod\theta$. The diagram $H^+$ is obtained by uniting the
diagrams of $H^+(n)$ for $n \geq 0$ and joining the cubes $U_n$
and $V_{n+1}$ for all $n \geq 0$, i.e., it is necessary to draw
$8$ additional edges for all $n \geq 0$. (Fig. \ref{boolean16}).
\end{theorem}

 \PerfProof
 We have to show that the {\it sum}
and the {\it intersection} do not lead out of $H^+$. The
distributivity follows from the absence of diamonds $M_3$, see
Proposition \ref{ch_distributive}, \cite{Bir48}. If
$$
  u \in H^+(n), \quad v \in H^+(m)
$$
and these lattices are not adjacent, i.e., $n \leq m-2$, then $v
\subseteq u$. Indeed, by Proposition \ref{UV_16bool_alg} (Fig.
\ref{boolean16}) we have
\begin{equation*}
  v \subseteq \sum_\text{$i=1,2,3$}c_i(m)
          \subseteq \bigcap_\text{$i=1,2,3$}b_i(m-1)
          \subseteq \bigcap_\text{$i=1,2,3$}c_i(m-1)
          \subseteq \bigcap_\text{$i=1,2,3$}a_i(m-2)
          \subseteq u .
\end{equation*}
So, we consider only the case $u \in H^+(n), v \in H^+(n+1)$.
Every element $u \in H^+(n)$ is of the form $u_1 + u_2$, where
 $u_1 \in H^+(n)$ and $u_2 \in U_n$. By Proposition
\ref{merge_1} we have
\begin{equation}  \label{merge_2}
   u + v = u_1 + u_2 + v =
   u_1 + (u_2 + v) \in H^+(n).
\end{equation}
Similarly, every element $v$ from  $H^+(n+1)$ is of the form
 $v_1{v}_2$, where $v_1 \in H^+(n+1)$ and
 $v_2 \in V_{n+1}$. By Proposition \ref{merge_1}
\begin{equation} \label{merge_3}
   vu = v_1{v}_2{u} = v_1(v_2{u}) \in H^+(n+1).
\end{equation}
So, Theorem \ref{main_perf_th} follows from (\ref{merge_2}) and
(\ref{merge_3}). This concludes the proof. \qedsymbol

\chapter{\sc\bf Atomic and admissible polynomials in $D^4$}
  \label{sect_adm_seq_D4}

\section{Admissible sequences in $D^4$}
  \label{subsect_adm_seq_D4}
For definition of admissible sequences in the case of the modular
lattice $D^4$, see
 \S\ref{adm_D4}. Essentially, the fundamental property of this definition is
\begin{equation}
  \label{main_D4}
   ijk = ilk \indent \text{ for all } \{i,j,k,l\} = \{1,2,3,4\}
\end{equation}
Relation (\ref{main_D4}) is our main tool in all further
calculations of admissible sequences of $D^4$.

 Without loss of generality only sequences starting at $1$ can
be considered. The following proposition will be used for the
classification of admissible sequences in $D^4$.

\begin{proposition} \label{relations_1_14}
  The following relations hold
 \item[1)]  $(31)^r(32)^s(31)^t = (32)^s(31)^{r+t},$
 \item[2)]  $(31)^r(21)^s(31)^t = (31)^{r+t}(21)^s,$
 \item[3)]  $(42)^r(41)^s = (41)^r(31)^s$, \indent $s \geq 1,$
 \item[4)]  $2(41)^r(31)^s = 2(31)^{s+1}(41)^{r-1},$
   \indent $r \geq 1,$
 \item[5)]  $(43)^r(42)^s(41)^t = (41)^r(21)^s(31)^t,$
 \item[6)]  $1(41)^r(21)^s = 1(21)^s(41)^r,$ \indent
            $1(i1)^r(j1)^s = 1(j1)^s(i1)^r,$ \indent
            $i,j \in \{2,3,4\}$,  $i \neq j,$
 \item[7)]  $(41)^r(21)^t(31)^s = (41)^r(31)^s(21)^t,$
 \item[8)]  $(13)^s(21)^r = (12)^r(31)^s,$
 \item[9)]  $12(41)^r(31)^s(21)^t =
   (14)^r(31)^{s+1}(21)^t = (14)^r(21)^{t+1}(31)^s,$
 \item[10)] $12(14)^r(31)^s(21)^t = (14)^r(31)^s(21)^{t+1},$
 \item[11)] $13(14)^r(31)^s(21)^t = (14)^r(31)^{s+2}(21)^{t-1},$
 \item[12)] $32(14)^r(31)^s(21)^t = (31)^s(21)^{t+1}(41)^r =
            34(14)^r(31)^s(21)^t,$
 \item[13)] $42(14)^r(31)^s(21)^t =
            (41)^s(21)^{t+1}(31)^r = 43(14)^r(31)^s(21)^t,$
 \item[14)] $23(14)^r(31)^s(21)^t = (21)^{t+1}(31)^s(41)^r =
            24(14)^r(31)^s(21)^t$.
\end{proposition}
\PerfProof 1) For $r=1$, we have
\begin{equation*}
 \begin{split}
   & (31)(32)^s(31)^t = 31(32)(32)...(32)(32)(31)^t =
                       31(41)(32)...(32)(32)(31)^t = \\
   & 31(41)(41)...(41)(41)(31)^t =
                       32(41)(41)...(41)(41)(31)^t =
                       32(32)(41)...(41)(41)(31)^t = \\
   & 32(32)(32)...(32)(31)(31)^t =
                       (32)^s(31)^{t+1}.
 \end{split}
\end{equation*}
 Applying induction on $r$ we get the relation
$$
   (31)^r(32)^s(31)^t = (32)^s(31)^{r+t}.
    \qed \vspace{2mm}
$$

2) For $t=1$, we have
\begin{equation*}
 \begin{split}
  & (31)^r(21)^s(31) =  (31)^r(21)(21)...(21)(21)31 =
                       (31)^r(21)(21)...(21)(34)31 = \\
  & (31)^r(21)(21)...(34)(34)31 =
                       (31)^r(34)(34)...(34)(34)21 =
                       (31)^r(34)(34)...(34)(21)21 = \\
  & (31)^r(34)(34)...(21)(21)21 =
                       (31)^r(31)(21)...(21)(21)21 =
                       (31)^{r+1}(21)^s.
 \end{split}
\end{equation*}
Applying induction on $t$ we get the relation
$$
  (31)^r(21)^s(31)^t = (31)^{r+t}(21)^s.
    \qed \vspace{2mm}
$$

3) For $s=1$, we have
\begin{equation*}
 \begin{split}
   & (42)^r41 =  42(42)(42)...(42)(42)41 =
              42(42)(42)...(42)(42)31 = \\
   & 42(42)(42)...(42)(31)31 =
        42(31)(31)...(31)(31)31 = \\
   & 41(31)(31)...(31)(31)31 =
              41(31)^r.
 \end{split}
\end{equation*}
Thus, by heading 2) we have
$$
  (42)^r(41)^s = 41(31)^r(41)^{s-1} = (41)^s(31)^r.
    \qed \vspace{2mm}
$$

4) For $s=0$, we have
$$
      2(41)^r = 2(41)(41)^{r-1} = 2(31)(41)^{r-1}
$$
and by heading 2):
$$
      2(41)^r(31)^s = 2(31)(41)^{r-1}(31)^s =
      2(31)^{s+1}(41)^{r-1}.
    \qed \vspace{2mm}
$$

5) Applying heading 3) we get
$$
   (43)^r(42)^s(41)^t = (43)^r(41)^s(31)^t.
$$
Again, applying heading 3) to  $(43)^r(41)^s$ we get
$$
   (43)^r(41)^s = (41)^r(21)^s
$$
and
$$
  (43)^r(42)^s(41)^t = (41)^r(21)^s(31)^t.
    \qed \vspace{2mm}
$$

6) For $s=1$, we have
\begin{equation*}
 \begin{split}
     & 1(41)^r(21) = 1(41)(41)...(41)(41)(21) =
       1(41)(41)...(41)(23)(21) = \\
     & 1(41)(41)...(23)(23)(21) = 1(23)(23)...(23)(23)(21) = \\
     & 1(23)(23)...(23)(23)(41) = 1(23)(23)...(23)(41)(41) = \\
     & 1(23)(41)...(41)(41)(41) = 1(21)(41)...(41)(41)(41) = \\
     & 1(21)(41)^r.
  \end{split}
\end{equation*}
   Thus, by heading 2) we have
$$
   1(41)^r(21)^s = 1(21)(41)^r(21)^{s-1} = 1(21)^s(41)^r.
    \qed \vspace{2mm}
$$

7) Follows from 6).
    \qedsymbol \vspace{2mm}

8) First,
\begin{equation}
  \label{1to1}
   1321 = 1421 = 1431 = 1231 = 1241 = 1341
\end{equation}

  For $r = 1$, by (\ref{1to1}) we have
\begin{equation*}
 \begin{split}
      & 13(21)^s = 13(21)(21)...(21)(21) =
                   12(31)(21)(21)...(21)(21)(21) = \\
      & 12(34)(21)(21)...(21)(21)(21) =
                   12(34)(34)(21)...(21)(21)(21) =  \\
      & 12(34)(34)(34)...(34)(34)(31) =
                   12(12)(12)(12)...(12)(34)(31) =  \\
      & 12(12)(12)(12)...(12)(12)(31) =
                   (12)^s(31).
 \end{split}
\end{equation*}
Suppose
$$
  (13)^r(21)^s = (12)^s(31)^r,
$$
then we have
\begin{equation*}
 \begin{split}
       & (13)^{r+1}(21)^s = 13(12)^s(31)^r =
               13(12)(12)...(12)(12)(31)^r = \\
       & 13(43)(43)...(43)(12)(31)^r =
               13(43)(43)...(43)(41)(31)^r = \\
       & 12(43)(43)...(43)(41)(31)^r =
               12(12)(43)...(43)(41)(31)^r = \\
       & 12(12)(12)...(12)(41)(31)^r =
               12(12)(12)...(12)(31)(31)^r =
               (12)^s(31)^{r+1}.
  \end{split}
\end{equation*}

Thus, by induction the following relation holds:
$$
  (13)^r(21)^s = (12)^s(31)^r
    \qed \vspace{2mm}
$$

9) First of all,
\begin{equation*}
 \begin{split}
      & 12(41)^r(31)^s(21)^t =
                12(41)(41)...(41)(41)(31)^s(21)^t = \\
      & 12(32)(41)...(41)(41)(31)^s(21)^t =
                12(32)(32)...(32)(31)(31)^s(21)^t =  \\
      & 14(14)(32)...(32)(31)(31)^s(21)^t =
                14(14)(14)...(14)(31)(31)^s(21)^t = \\
      & (14)^r(31)^{s+1}(21)^t.
  \end{split}
\end{equation*}
By heading 4) we have
$$
   4(31)^{s+1}(21)^t = 4(21)^{t+1}(31)^s.
$$
Thus,
$$
   (14)^r(31)^{s+1}(21)^t = (14)^r(21)^{t+1}(31)^s.
    \qed \vspace{2mm}
$$

10) First, we have
\begin{equation}
  \label{12_31}
    12(14)^r(31) = (14)^r(21)^2
\end{equation}
since
\begin{equation*}
 \begin{split}
     &  12(14)^r(31) =
        12(14)(14)...(14)(14)(31) =
        12(32)(14)...(14)(14)(31) = \\
     &  12(32)(32)...(32)(34)(31) =
        14(14)(32)...(32)(34)(31) = \\
     &  14(14)(14)...(14)(34)(21) =
        14(14)(14)...(14)(21)(21) = (14)^r(21)^2.
  \end{split}
\end{equation*}
By (\ref{12_31}) we have
$$
  12(14)^r(31)^s(21)^t = 12(14)^r(31)(31)^{s-1}(21)^t =
                (14)^r(21)^2(31)^{s-1}(21)^t.
$$
By headings 2) and 4):
$$
     (21)^2(31)^{s-1}(21)^t = (21)^{t+2}(31)^{s-1}
$$
and
$$
    12(14)^r(31)^s(21)^t = (14)^r(21)^{t+2}(31)^{s-1} =
    (14)^r(31)^s(21)^{t+1}.
    \qed \vspace{2mm}
$$

11) Applying permutation $2 \leftrightarrow 3$ to heading 10), we
get
$$
   13(14)^r(21)^t(31)^s = (14)^r(21)^t(31)^{s+1}.
$$
By heading 4)
$$
   13(14)^r(21)^t(31)^s = (14)^r(31)^{s+2}(21)^{t-1}.
    \qed \vspace{2mm}
$$

12) By heading 8)
$$
  (14)^r(31)^s = (13)^s(41)^r
$$
and
\begin{equation*}
\begin{split}
    & 32(14)^r(31)^s(21)^t = 32(13)^s(41)^r(21)^t =    \\
    & 32(13)(13)...(13)(13)(41)^r(21)^t =
      34(24)(13)...(13)(13)(41)^r(21)^t = \\
    & 34(24)(24)...(24)(21)(41)^r(21)^t =
      31(24)(24)...(24)(21)(41)^r(21)^t = \\
    & 31(31)(31)...(31)(21)(41)^r(21)^t =
         (31)^s(21)(41)^r(21)^t.
 \end{split}
\end{equation*}
By heading 2) we have
$$
  (31)^s(21)(41)^r(21)^t = (31)^s(21)^{t+1}(41)^r
$$
and
$$
  32(14)^r(31)^s(21)^t = (31)^s(21)^{t+1}(41)^r
    \qed \vspace{2mm}
$$

13) By heading 8):
$$
   (14)^r(31)^s = (13)^s(41)^r
$$
and
$$
   42(14)^r(31)^s(21)^t = 42(13)^s(41)^r(21)^t.
$$
Applying permutation $4 \leftrightarrow 3$ to heading 12) we get
$$
   42(13)^r(41)^s(21)^t = (41)^s(21)^{t+1}(31)^r.
    \qed \vspace{2mm}
$$

14) By heading 4) we have
$$
  23(14)^r(31)^s(21)^t = 23(14)^r(21)^{t+1}(31)^{s-1}.
$$
Applying permutation $2 \leftrightarrow 3$ to heading 12) we get
$$
     23(14)^r(21)^{t+1}(31)^{s-1} = (21)^{t+1}(31)^s(41)^r.
    \qed \vspace{2mm}
$$

\begin{table}[h]
 \renewcommand{\arraystretch}{1.9}
  \begin{tabular} {||c|c|c|c|c|c||}
  \hline \hline
     & Admissible Sequence
     & Action       & Action      & Action   & Action \cr
     &  & $\varphi_1$  & $\varphi_2$ & $\varphi_3$ & $\varphi_4$
     \\
  \hline         
    $F21$  & $(21)^t(41)^r(31)^s = (21)^t(31)^s(41)^r $
     & $G11$
     & -
     & $G31$
     & $G41$ \\
  \hline         
    $F31$  & $(31)^s(41)^r(21)^t = (31)^s(21)^t(41)^r $
     & $G11$
     & $G21$
     & -
     & $G41$ \\
  \hline      
    $F41$  & $(41)^r(31)^s(21)^t = (41)^r(21)^t(31)^s $
     & $G11$
     & $G21$
     & $G31$
     & - \\
  \hline \hline        
    $G11$  & $1(41)^r(31)^s(21)^t$ =
     & -
     & $F21$
     & $F31$
     & $F41$ \cr
     &    \hspace{3mm} $1(31)^s(41)^r(21)^t$ = & & & &  \cr
     &    \hspace{6mm} $1(21)^t(31)^s(41)^r$ & & & & \\
  \hline         
    $G21$  & $2(41)^r(31)^s(21)^t$ =
                $2(31)^{s+1}(41)^{r-1}(21)^t$
     & $H11$
     & -
     & $F31$
     & $F41$ \\
  \hline         
     $G31$  & $3(41)^r(21)^t(31)^s$ =
                $3(21)^{t+1}(41)^{r-1}(31)^s$
     & $H11$
     & $F21$
     & -
     & $F41$  \\
  \hline      
    $G41$  & $4(21)^t(31)^s(41)^r$ =
                $4(31)^{s+1}(21)^{t-1}(41)^r$
     & $H11$
     & $F21$
     & $F31$
     & -         \\
  \hline \hline     
    $H11$  & $(14)^r(31)^s(21)^t = (14)^r(21)^{t+1}(31)^{s-1}$ =
     & -
     & $H21$
     & $H31$
     & $H41$ \cr
     &   $(13)^s(41)^r(21)^t = (13)^s(21)^{t+1}(41)^{r-1}$ = & & & &  \cr
     &   $(12)^{t+1}(41)^r(31)^{s-1} = (12)^{t+1}(31)^{s-1}(41)^r$ & & & & \\
  \hline      
    $H21$  & $2(14)^r(31)^s(21)^t$
     & $H11$
     & -
     & $F31$
     & $F41$ \\
  \hline      
    $H31$  & $3(14)^r(31)^s(21)^t$
     & $H11$
     & $F21$
     & -
     & $F41$ \\
  \hline      
    $H41$  & $4(14)^r(31)^s(21)^t$
     & $H11$
     & $F21$
     & $F31$
     & -         \\
  \hline  \hline
  \end{tabular}
  \vspace{2mm}
  \caption{\hspace{3mm}Admissible sequences for the modular lattice $D^4$}
  \label{table_admissible_ExtD4}
\end{table}

\begin{proposition}
 \label{full_adm_seq_D4}
   Full list of admissible sequences starting at $1$ is given
   by Table \ref{table_admissible_ExtD4}.
\end{proposition}
  \underline{Note to Table \ref{table_admissible_ExtD4}}. Type
$Fij$ (resp. $Gij$, $Hij$) denotes the admissible sequence
starting at $j$ and ending at $i$. Sequences of type $Fij$ and
$H11$ contain an even number of symbols, sequences of type $Gij$
and $Hij, (i > 1)$ contain an odd number of symbols. For
differences in types $Fij, Gij, Hij$, see the table. \vspace{3mm}

\PerfProof It suffices to prove that maps $\varphi_i$, where
$i=1,2,3,4$, do not lead out of Table
 \ref{table_admissible_ExtD4}. The exponents $r,s,t$ may be any
non-negative integer number. The proof is based on the relations
from Proposition \ref{relations_1_14}. We refer only to number of
relation and drop reference to Proposition \ref{relations_1_14}
itself.

 \underline{Lines $F21$--$F41$}.
Consider, for example, $\varphi_3$. By heading 4) we get
\begin{equation}
 \label{act_psi_4_F21}
   \varphi_3((21)^t(41)^r(31)^s) = 3(41)^{r+1}(21)^{t-1}(31)^s,
\end{equation}
i.e., we get $G31$.
 \qedsymbol \vspace{2mm}

  \underline{Line $G11$}. Consider the action $\varphi_2$.
By heading 6) we get
$$
    \varphi_2(1(41)^r(31)^s(21)^t) = 2(1(41)^r(21)^t(31)^s) =
    21(41)^r(21)^t(31)^s.
$$
By heading 2) we have
\begin{equation}
  \label{act_psi_2_F21}
  \varphi_2(1(41)^r(31)^s(21)^t) =
    21(41)^r(21)^t(31)^s = (21)^{t+1}(41)^r(31)^s,
\end{equation}
i.e.,  we get $F21$.
 \qedsymbol \vspace{2mm}

 \underline{Line G21}.
By heading 9) applying $\varphi_1$ we get $H11$. For the action
$\varphi_3$, we have
\begin{equation}
  \label{act_psi_3_F31}
  32(41)^r(31)^s(21)^t = 31(41)^r(31)^s(21)^t =
  (31)^{s+1}(41)^r(21)^t,
\end{equation}
i.e., we get $F31$. For $\varphi_4$ from heading 4), we have
\begin{equation}
  \label{act_psi_4_F31}
  \begin{split}
  & 42(41)^r(31)^s(21)^t = 42(31)^{s+1}(41)^{r-1}(21)^t = \\
  & 41(31)^{s+1}(41)^{r-1}(21)^t = (41)^r(31)^{s+1}(21)^t,
  \end{split}
\end{equation}
i.e., we get $F41$. For $\varphi_1$ from heading 9), we have
\begin{equation}
  \label{act_psi_1_F31}
  \begin{split}
  & 12(41)^r(31)^s(21)^t = (14)^r(31)^{s+1}(21)^t = \\
  & (14)^r(21)^{t+1}(31)^s,
  \end{split}
\end{equation}
i.e., we get $H11$.
 \qedsymbol \vspace{2mm}

 \underline{Line G31}.
For the action $\varphi_1$ by relation $13(41)^t = 12(41)^t$ we
get $H11$. For the action $\varphi_2$ and by heading 6), we have
$$
  23(41)^r(31)^s(21)^t = 23(41)^r(21)^t(31)^s =
  21(41)^r(21)^t(31)^s = (21)^{t+1}(41)^r(31)^s
$$
i.e. we get $F21$.

For the action $\varphi_4$ from headings 6), 4), we have
$$
  43(41)^r(31)^s(21)^t = 43(41)^r(21)^t(31)^s =
  43(21)^{t+1}(41)^{r-1}(31)^s = 41(21)^{t+1}(41)^{r-1}(31)^s
$$
So, by heading 2) we have
$$
  43(41)^r(31)^s(21)^t =
  (41)^r(21)^{t+1}(31)^s,
$$
i.e., we get $F41$.
 \qedsymbol \vspace{2mm}

 \underline{Line G41}.
By heading 9) we get $H11$. For action $\varphi_2$ we have from 4)
$$
  24(21)^t(31)^s(41)^r = 24(31)^{s+1}(21)^{t-1}(41)^r  =
  21(31)^{s+1}(21)^{t-1}(41)^r = (21)^t(31)^{s+1}(41)^r
$$
i.e., we get $F21$. For the action $\varphi_3$ we have from 2)
$$
  34(21)^t(31)^s(41)^r = 31(21)^t(31)^s(41)^r =
  (31)^{s+1}(21)^t(41)^r
$$
i.e., we get $F31$
 \qedsymbol \vspace{2mm}

 \underline{Line H11}. Actions are rather trivial.
 \qedsymbol \vspace{2mm}

 \underline{Line H21}. For action $\varphi_1$ the relation follows from heading 10),
for $\varphi_3$ -- follows from heading 12), and for $\varphi_4$
-- follows from heading 13).
 \qedsymbol \vspace{2mm}

 \underline{Line H31}. Follows from headings 9), 14) and 13).
 \qedsymbol \vspace{2mm}

 \underline{Line H41}. Action $\varphi_1$ is trivial since $14(14)^r =
(14)^{r+1}$. Action $\varphi_2$ follows from heading 14) and
action $\varphi_3$ follows from 12).
 \qedsymbol \vspace{2mm}

Thus, all cases are considerd and the proposition is proved.
 \qedsymbol \vspace{2mm}

 \index{slice}

\begin{figure}[h]
\centering
\includegraphics{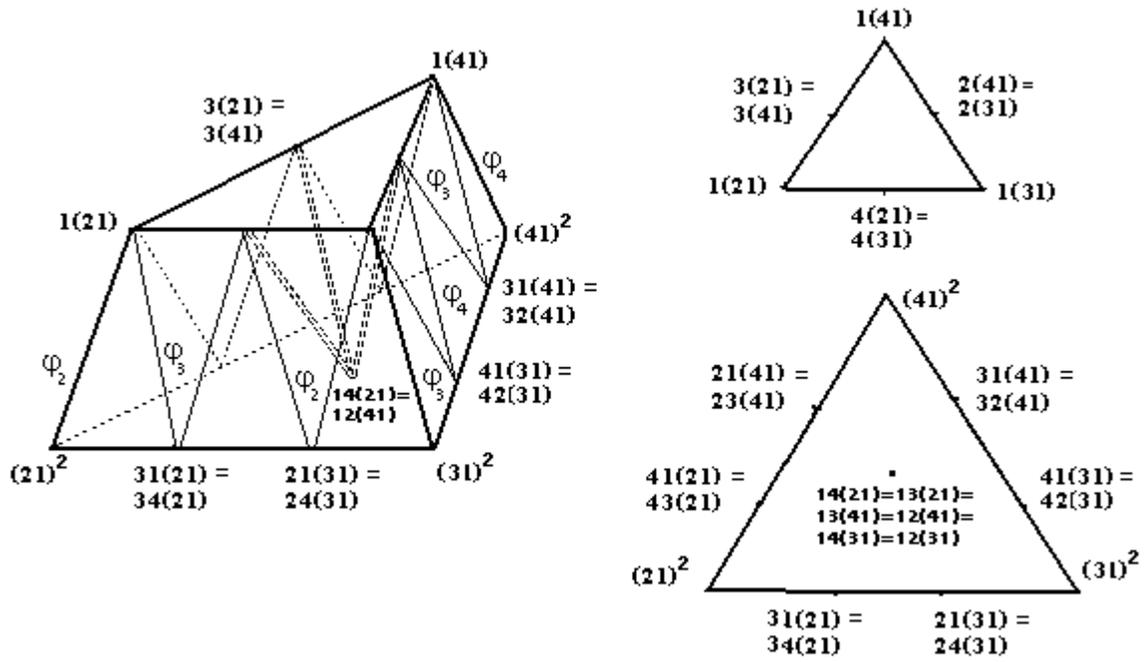}
\caption{\hspace{3mm}Slices of admissible sequences for $D^4$,
$l=3$ and $l=4$}
\label{section34}
\end{figure}

\begin{figure}[h]
\centering
\includegraphics{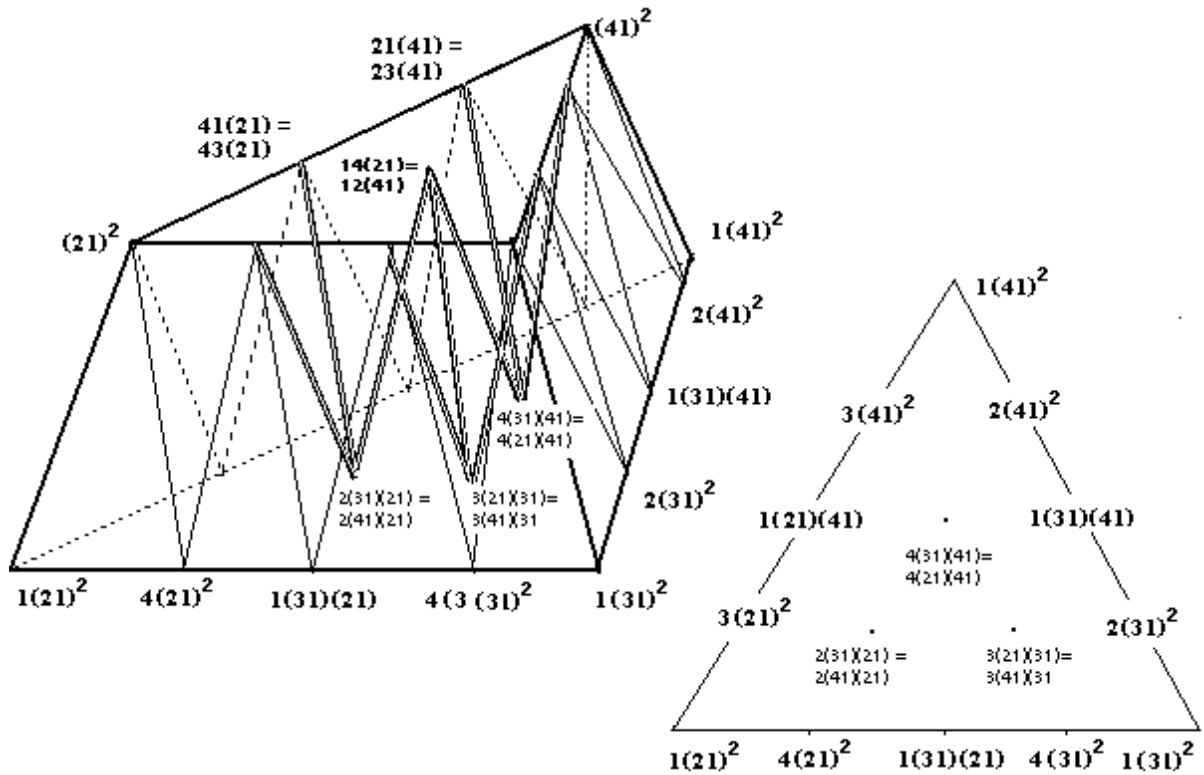}
\caption{\hspace{3mm}Slices of admissible sequences for $D^4$,
$l=4$ and $l=5$}
\label{section45}
\end{figure}

The pyramid on the Fig. \ref{pyramid_D4} has internal points. We
consider the slice $S(n)$ containing all sequences of the same
length $n$. The slices $S(3)$ and $S(4)$ are shown on the Fig.
\ref{section34}.
 The slice $S(4)$ contains only one internal point
\begin{equation}
  14(21) = 13(21) = 13(41) = 12(41) = 14(31) = 12(31)
\end{equation}
 The slices $S(4)$ and $S(5)$ are shown on the Fig. \ref{section45}.
 The slice $S(5)$ contains $3$ internal points
\begin{equation}
\begin{split}
  & 2(31)(21) = 2(41)(21), \\
  & 3(21)(31) = 3(41)(31), \\
  & 4(21)(41) = 4(31)(41).
\end{split}
\end{equation}

\begin{remark}
 {\rm The slice $S(n)$ contains $\displaystyle\frac{1}{2}{n(n+1)}$ different admissible
 sequences. Actions of $\varphi_i$, where $i = 1,2,3,4$ move
 every line in the triangle $S(n)$, which is parallel to
 some edge of the triangle, to the edge of $S(n+1)$. The edge
 containing $k$ points is moved to $k+1$ points in the $S(n+1)$. }
\end{remark}

\section{Atomic polynomials and elementary maps in the modular lattice $D^4$}
 \label{atomic_D4}

  \index{atomic elements $a^n_{ij}$ in $D^4$}

The free modular lattice $D^4$ is generated by $4$ generators:
$$
    D^4 = \{e_1,e_2,e_3,e_4\}.
$$

Recall, that {\it atomic} lattice polynomials $a_n^{ij}$, where
$i,j \in \{1,2,3,4\}$, $n \in \{0,1,2,3,\dots \}$, for the case of
$D^4$, are defined as follows
\begin{equation*}
  a_n^{ij} =
 \begin{cases}
   a_n^{ij} = I \qquad \qquad \qquad \qquad \qquad \hspace{3mm} \text{ for } n = 0,\\
   a_n^{ij} = e_i + e_j{a}_{n-1}^{kl} = e_i + e_j{a}_{n-1}^{lk}
        \hspace{1.7mm} \text{ for } n \geq 1,
 \end{cases}
\end{equation*}
where $\{i,j,k,l\}$ is the permutation of $\{1,2,3,4\}$, see
\S\ref{atomic_polynoms}.

\begin{proposition}
  1) The following property of the atomic elements take place
   \begin{equation}
     \label{main_atomic_D4}
      e_j{a}_n^{kl} = e_j{a}_n^{lk} \text{ for } n \geq 1,
      \text{ and distinct indices } j,k,l.
   \end{equation}

  2) The definition of the atomic elements $a_n^{ij}$ in (\ref{eq_atomic_D4}) is
   well-defined. \vspace{2mm}

  3) We have
   \begin{equation}
     \label{incl_atomic_D4}
        a^{ij}_n \subseteq a^{ij}_{n-1} \subseteq \dots \subseteq
        a^{ij}_2 \subseteq a^{ij}_1 \subseteq a^{ij}_0 = I
        \text{ for all } i \neq j.
   \end{equation}

  4) To equalize the lower indices of the admissible polynomials
  $f_{{\alpha}0}$
  (see Table \ref{table_adm_elem_D4} and Theorem \ref{th_adm_classes_D4})
  we will use the following relation:
   \begin{equation}
     \label{equalize_D4}
        e_j + e_i{a}^{jk}_{t+1}{a}^{lk}_{s-1} =
        e_j + e_k{a}^{ij}_{t}{a}^{il}_{s}
        \quad \text{ for all } \{i,j,k,l\} = \{1,2,3,4\}.
   \end{equation}
\end{proposition}
 \PerfProof 1) Suppose (\ref{main_atomic_D4}) is true for the index
 $n-1$:
   \begin{equation*}
      e_j{a}_{n-1}^{kl} = e_j{a}_{n-1}^{lk}.
   \end{equation*}
   Since $e_j \subseteq a_{n-1}^{ji}$,
   by the permutation property (\ref{permutation1}), we have
   \begin{equation}
      e_j{a}_n^{kl} = e_j(e_k + e_l{a}_{n-1}^{ij}) =
        e_j(e_k + e_l{a}_{n-1}^{ji}) =
          e_j(e_l + e_k{a}_{n-1}^{ji}) =
           e_j{a}_n^{lk}. \qed \vspace{2mm}
   \end{equation}

2) follows from 1). \qedsymbol \vspace{2mm}

3) By induction hypothesis we have $a^{kl}_{n-1} \subseteq
a^{kl}_{n-2}$, and therefore
$$
  a^{ij}_n = e_i + e_ja^{kl}_{n-1} \subseteq
  e_i + e_ja^{kl}_{n-2}= a^{ij}_{n-1}.
  \qed \vspace{2mm}
$$

4) Without loss of generality, we will show that
 \begin{equation}
     \label{equalize_D4_2}
        e_2 + e_1{a}^{24}_{t+1}{a}^{34}_{s-1} =
        e_2 + e_4{a}^{12}_{t}{a}^{13}_{s}.
 \end{equation}
 By permutation property (\ref{permutation1}) we have
 \begin{equation*}
  \begin{split}
       & e_2 + e_1{a}^{24}_{t+1}{a}^{34}_{s-1} =
         e_2 + e_1(e_4 + e_2a^{13}_{t}){a}^{43}_{s-1} = \\
       & e_2 + e_1(e_4 + e_2a^{13}_{t}{a}^{43}_{s-1}) =
         e_2 + e_4(e_1 + e_2a^{13}_{t}{a}^{43}_{s-1}) = \\
       & e_2 + e_4a^{13}_{t}(e_1 + e_2{a}^{43}_{s-1}) =
         e_2 + e_4a^{13}_{t}{a}^{12}_{s}.
  \qed \vspace{2mm}
  \end{split}
 \end{equation*}

Now we briefly recall definitions due to Gelfand and Ponomarev
\cite{GP74}, \cite{GP76}, \cite{GP77} of spaces $G_i$, $G_i{'}$,
representations $\nu^0$, $\nu^1$, joint maps $\psi_i$, and
elementary maps $\varphi_i$, where $i = 1,2,3,4$. To compare these
definitions with a case of the modular lattice $D^{2,2,2}$, see
Ch.\ref{section_Coxeter}. \vspace{2mm}

We denote by
$$
  \{ Y_1, Y_2, Y_3, Y_4 \mid Y_i \subseteq X_0, i=1,2,3,4 \}
$$
the representation $\rho$ of $D^4$ in the finite dimensional
vector space
 $X_0$, and by
$$
  \{ Y^1_1, Y^1_2, Y^1_3, Y^1_4 \mid Y^1_i \subseteq X^1_0, i=1,2,3,4 \}
$$
the representation $\Phi^+\rho$ of $D^4$. Here $Y_i$ (resp.
$Y^1_i$) is the image of the generator $e_i$
 under the representation $\rho$ (resp. $\Phi^+\rho$).
\begin{equation}
\label{repr_rho_D4}
\begin{array}{ccccccccccccc}
   &       & & Y_3 & & & &
   & & Y^1_3 & & \\
   & & & \cap\shortmid & & & &
   & & \cap\shortmid & &   \\
   & Y_1 & \subseteq
   & X_0 & \supseteq
   & Y_2 & \hspace{3mm}
   & Y^1_1 & \subseteq
   & X^1_0 & \supseteq
   & Y^1_2 \\
   & & & \cup\shortmid & & & &
   & & \cup\shortmid & &   \\
   & \rho      & & Y_4 & & & &
   \Phi^+\rho  & & Y^1_4 & & \\ \vspace{5mm}
\end{array}
\end{equation}

The space $X^1_0$ --- the space of the representation $\Phi^+\rho$
--- is
$$
   X^1_0 = \{ (\eta_1, \eta_2, \eta_3, \eta_4) \mid \eta_i \in Y_{i},
   \quad \sum\eta_i = 0 \},
$$
where $i \in \{1,2,3,4\}$. As in the case of $D^{2,2,2}$
\S\ref{cox_functor}, we set
$$
    R = \bigoplus\limits_\text{$i = 1,2,3,4$}Y_i,
$$
i.e.,
$$
  R = \{ (\eta_1, \eta_2, \eta_3, \eta_4) \mid \eta_i \in Y_i, \quad i = 1,2,3 \}.
$$
Then, $X^1_0 \subseteq R$.

 For the case of $D^4$, the spaces $G_i$ and $G{'}_i$
are introduced as follows:
\begin{equation}
  \begin{split}
       & G_1 = \{(\eta_1, 0, 0, 0) \mid \eta_1 \in Y_1\}, \quad
         G'_1 = \{(0, \eta_2, \eta_3, \eta_4) \mid \eta_i \in Y_i\}, \\
       & G_2 = \{(0, \eta_2, 0, 0) \mid \eta_2 \in Y_2\}, \quad
         G'_2 = \{(\eta_1, 0, \eta_3, \eta_4) \mid \eta_i \in Y_i\}, \\
       & G_3 = \{(0, 0, \eta_3, 0) \mid \eta_3 \in Y_3\},  \quad
         G'_3 = \{(\eta_1, \eta_2, 0, \eta_4) \mid \eta_i \in Y_i\}, \\
       & G_4 = \{(0, 0, 0, \eta_4) \mid \eta_4 \in Y_4\},  \quad
         G'_4 = \{(\eta_1, \eta_2, \eta_3, 0) \mid \eta_i \in Y_i\}. \\
  \end{split}
\end{equation}

For details, see \cite[p.43]{GP74}\footnote{Compare with
definition (\ref{spaces_G_H}) of the spaces $G_i, G'_i, H_i, H'_i$
in the case of $D^{2,2,2}$; see also Table
 \ref{compare_notions}.}.

The associated representations $\nu_0$, $\nu_1$ in $R$ are defined
by Gelfand and Ponomarev \cite[eq.(7.2)]{GP74}\footnote{Compare
with definitions for the case of $D^{2,2,2}$ in
\S\ref{seq_assoc}.}:
\begin{equation}
 \begin{split}
   \label{assoc_repr_D4}
        & \nu^0(e_i) = X_0^{1} + G_i, \indent i = 1,2,3, \\
        & \nu^1(e_i) = X_0^{1}G'_i,  \indent i = 1,2,3.
  \end{split}
\end{equation}

Following \cite{GP74}, we introduce the elementary maps
$\varphi_i$:
 $$
 \varphi_i : X_0^{1}  \longrightarrow X_0,
  (\eta_1, \eta_2, \eta_3, \eta_4) \longmapsto \eta_i.
 $$
 From the definition we have
\begin{equation}  \label{sum_fi_D4}
        \varphi_1  + \varphi_2 + \varphi_3 + \varphi_4  = 0.
\end{equation}

\begin{table}[h]
  \renewcommand{\arraystretch}{1.7}
  \begin{tabular} {|| c | c | c ||}
   \hline \hline
       Notions & $D^{2,2,2}$ & $D^4$  \\
   \hline \hline 
       Generators
       & $\{x_1 \subseteq y_1,
            x_2 \subseteq y_2, x_3 \subseteq y_3\}$
       & $\{e_1, e_2, e_3, e_4\}$ \\
     \hline
       Atomic
       & $a_n^{ij} = x_i + y_j{a}_{n-1}^{jk}$,
       & $a_n^{ij} = e_i + e_j{a}^{kl}$,  \cr
       elements
       & $A_n^{ij} = y_i + x_j{A}_{n-1}^{ki},$
       & $\{i,j,k,l \} = \{1,2,3,4\}$ \cr
       & $\{i,j,k\} = \{1,2,3\}$ &  \\
     \hline
       Representation $\rho$
       & $\rho(x_i) = X_i, \rho(y_i) = Y_i$
       & $\rho(e_i) = Y_i$ \\
     \hline
       Space $X^1_0$
       & $X^1_0 =
        \{ (\eta_1, \eta_2, \eta_3) \mid \eta_i \in Y_{i} \},$
       & $X^1_0 =
          \{ (\eta_1, \eta_2, \eta_3, \eta_4) \mid \eta_i \in Y_{i} \},$ \cr
       & where $\eta_1 + \eta_2 + \eta_3 = 0$
       & where $\eta_1 + \eta_2 + \eta_3 + \eta_4 = 0$ \\
     \hline
       Representation $\Phi^+\rho$
       & $\rho(x_i) = X^1_i, \rho(y_i) = Y^1_i$
       & $\rho(y_i) = Y^1_i$ \\
     \hline
       Spaces $G_i$, $H_i$
       & $G_1 = \{(\eta_1, 0, 0) \mid \eta_1 \in Y_1\}$,
       & $G_1 = \{(\eta_1, 0, 0, 0) \mid \eta_1 \in Y_1\}$ \cr
       & $H_1 = \{(\xi_1, 0, 0) \mid \xi_1 \in X_1\}$,
       & $G_2 = \{(0, \eta_2, 0, 0) \mid \eta_2 \in Y_2\}$ \cr
       & $G_2 = \{(0, \eta_2, 0) \mid \eta_2 \in Y_2\}$,
       & $G_3 = \{(0, 0, \eta_3, 0) \mid \eta_3 \in Y_3\}$ \cr
       & $H_2 = \{(0, \xi_2, 0) \mid \xi_2 \in X_2\}$,
       & $G_4 = \{(0, 0, 0, \eta_4) \mid \eta_4 \in Y_4\}$ \cr
       & $G_3 = \{(0, 0, \eta_3) \mid \eta_3 \in Y_3\}$,
       & \cr
       & $H_3 = \{(0, 0, \xi_3) \mid \xi_3 \in X_3\}$,
       & \\
     \hline
       Joint maps $\psi_i$
       & $\psi_{i}(a) = X_0^{1} + G_{i}(H_{i}^{'} + \nu^{1}(a))$
       & $\psi_{i}(a) = X_0^{1} + G_{i}(G_{i}^{'} + \nu^{1}(a))$ \\
     \hline
       Quasi-
       & $\psi_i(a)\psi_i(b)$ =
       & $\psi_i(a)\psi_i(b)$ = \cr
       multiplicativity
       & $\psi_i((a + e_i)(b + x_j{x}_k))$
       & $\psi_i((a + e_i)(b + e_j{e}_k{e}_l)$ \\
     \hline
       Elementary maps $\varphi_i$
       & $\varphi_i : X_0^{1}  \longrightarrow X_0$,
       & $\varphi_i : X_0^{1}  \longrightarrow X_0$, \cr
       & $(\eta_1, \eta_2, \eta_3) \longmapsto \eta_i$
       & $(\eta_1, \eta_2, \eta_3, \eta_4) \longmapsto \eta_i$ \\
     \hline
       Fundamental properties
       & $\varphi_i\varphi_j\varphi_i + \varphi_i\varphi_k\varphi_i = 0$,
       & $\varphi_i\varphi_k\varphi_j + \varphi_i\varphi_l\varphi_j = 0$ \cr
       of elementary maps $\varphi_i$
       & $\varphi_i^3 = 0$
       & $\varphi_i^2 = 0$ \\
     \hline
       Fundamental properties & &  \cr
       of indices
       & \fbox{$iji = iki$}
       & \fbox{$ikj = ilj$} \cr
       (admissible sequences) && \\
     \hline \hline
  \end{tabular}
  \vspace{2mm}
\caption{\hspace{3mm}Comparison of notions in $D^{2,2,2}$ and
$D^4$}
  \label{compare_notions}
  \vspace{3.7mm} 
\end{table}

\section{Basic relations for the elementary and joint maps in the case of $D^4$}
 \label{sect_basic_rel_D4}

 We define {\it joint} maps $\psi_i\colon{D}_4 \longrightarrow
\mathcal{L}(R)$ as in the case of $D^{2,2,2}$:
\begin{equation}
  \label{psi_D4}
    \psi_{i}(a) = X_0^{1} + G_i(G'_i + \nu^1(a)).
\end{equation}

\begin{proposition}
 \label{basic_eq_D4}
In the case $D^4$ the joint maps $\psi_i$ satisfy the following
basic relations\footnote{Compare with Proposition \ref{basic_eq}
for $D^{2,2,2}$}:
\begin{enumerate}
 \item $\psi_i(e_i) = X_0^1$, \vspace{2mm}
 \item $\psi_i(e_j) = \nu^0(e_i(e_k + e_l))$, \vspace{2mm}
 \item $\psi_i(I) = \nu^0(e_i(e_j + e_k + e_l))$, \vspace{2mm}
 \item $\psi_i(e_k{e}_l) = \psi_j(e_k{e}_l) =
     \nu^0(e_i{e}_j)$. \vspace{2mm}
\end{enumerate}
\end{proposition}

\PerfProof 1) From (\ref{psi_D4}) and (\ref{assoc_repr_D4}) we
have $\psi_i(e_i) = X_0^1 + G_i{G'_i} = X_0^1$. \qedsymbol
\vspace{2mm}

2)We have
   $\psi_i(y_j) = X_0^1 + G_i(G'_i + X_0^1{G'_j})$.
       From $G_i \subseteq G'_j$ for $i \neq j$ and by
       the permutation property (\ref{permutation1}), we get
$$
       \psi_i(y_j) = X_0^1 + G_i(G'_i{G'_j} + X_0^1).
$$
       Since $G'_iG'_j =  G_k + G_l$, where $i,j,k,l$ are distinct indices,
       we have
$$
       \psi_i(y_j) = X_0^1 + G_i(X_0^1 + G_k + G_l) =
       (X_0^1 + G_i)(X_0^1 + G_k + G_l) = \nu^0(e_i(e_k + e_l)).
       \qed \vspace{2mm}
$$

3) Again,
\begin{equation}
 \begin{split}
     & \psi_i(I) = X_0^1 + G_i(G'_i + X_0^1) = \\
     & X_0^1 + G_i(G_j + G_k + G_l + X_0^1) =
      \nu^0(e_i(e_j + e_k + e_l)).
     \qed \vspace{2mm}
 \end{split}
\end{equation}

4) Since $G_i \subseteq G'_k$ for all $i \neq k$, we have
 \begin{equation}
   \label{e_k__e_l}
 \begin{split}
     & \psi_i(e_k{e}_l) = X_0^1 + G_i(G'_i + X_0^1G'_kG'_l) = \\
     & X_0^1 + G_i(G'_i{G'_k}G'_l + X_0^1).
 \end{split}
\end{equation}
Since $G_j(G_i + G_k + G_l) = 0$, then
\begin{equation}
   \label{intersect_G_2}
    G'_i{G'_j} = G_k + G_l.
\end{equation}
Indeed,
 \begin{equation}
   \label{intersect_G_3}
 \begin{split}
     & G'_i{G'_j} = (G_j + G_k + G_l)(G_i + G_k + G_l) = \\
     & G_k + G_l +  G_j(G_i + G_k + G_l) = G_k + G_l.
 \end{split}
\end{equation}
From (\ref{e_k__e_l}) and (\ref{intersect_G_2}), we see that
 \begin{equation*}
 \begin{split}
     & \psi_i(e_k{e}_l) = X_0^1 + G_i(G'_i{G'_k}G'_l + X_0^1) = \\
     & X_0^1 + G_i(X_0^1 + G_j) = \nu_0(e_i{e}_j).
     \qed \vspace{2mm}
 \end{split}
\end{equation*}

The main relation between $\varphi_i$ and $\psi_i$
 (Proposition \ref{phi_and_psi}) holds also for the case of $D^4$.
 Namely, let $a,b,c \subseteq D^{2,2,2}$, then

\begin{equation}
  \label{phi_and_psi_D4}
\begin{split}
    & \text{ (i) If } \psi_i(a) = \nu^0(b), \text{ then }
      \varphi_i\Phi^+\rho(a) = \rho(b). \\
    & \text{(ii) If } \psi_i(a) = \nu^0(b) \text{ and }
      \psi_i(ac) = \psi_i(a)\psi_i(c), \text{ then }\\
    & \qquad
      \varphi_i\Phi^+\rho(ac) = \varphi_i\Phi^+\rho(a)\varphi_i\Phi^+\rho(c).
 \end{split}
\end{equation}

From Proposition \ref{basic_eq_D4} and eq.(\ref{phi_and_psi_D4})
we have

\begin{corollary}
 \label{cor_psi_D4}
 For the elementary map $\varphi_i$ the following basic
 relations hold\footnote{Compare with Corollary \ref{cor_psi} --- a similar proposition
    for $D^{2,2,2}$.}:
\begin{enumerate}
 \item $\varphi_i\Phi^+\rho(e_i) = 0$, \vspace{2mm}
 \item $\varphi_i\Phi^+\rho(e_j) = \rho(e_i(e_k + e_l))$, \vspace{2mm}
 \item $\varphi_i\Phi^+\rho(I) = \rho(e_i(e_j + e_k + e_l))$, \vspace{2mm}
 \item $\varphi_i\Phi^+\rho(e_k{y}_l) = \rho(e_i{e}_j)$. \vspace{2mm}
 \end{enumerate}
\end{corollary}

\section{Additivity and multiplicativity of the joint maps in the case of $D^4$}
 \label{add_multi_D4}

\begin{proposition}
  \label{additivity_D4}
   The map $\psi_i$ is additive and quasimultiplicative
   with respect to the
   lattice operations + and $\cap$,
   namely\footnote{Compare with Proposition \ref{additivity},
     case $D^{2,2,2}$}:
\begin{enumerate}
   \item $\psi_i(a) + \psi_i(b) = \psi_i(a+b)$, \vspace{2mm}
   \item $\psi_i(a)\psi_i(b)$ =
      $\psi_i((a + e_i)(b + x_j{x}_k{x}_l))$, \vspace{2mm}
   \item $\psi_i(a)\psi_i(b)$ =
      $\psi_i(a(b + e_i + x_j{x}_k{x}_l))$. \vspace{2mm}
\end{enumerate}
\end{proposition}
\PerfProof 1) By the modular law (\ref{modular_law})
\begin{equation*}
\begin{split}
 &  \psi_i(a) + \psi_i(b) =
   X_0^1 + G_i(G'_i + \nu^1(a)) +  G_i(G'_i + \nu^1(b)) = \\
 & X_0^1 + G_i
   \left (G'_i + \nu^1(a) + G_i(G'_i + \nu^1(b)) \right ) =
   X_0^1 + G_i \left ((G'_i + \nu^1(b))(G'_i + G_i) + \nu^1(a) \right ).
\end{split}
\end{equation*}
   Since $G_i^{'} + G_i = R$, it follows that
$$
   \psi_i(a) + \psi_i(b) =
       X_0^1 + G_i \left (G'_i + \nu^1(b) + \nu^1(a) \right ) =
       X_0^1 + G_i(G'_i + \nu^1(b + a) ).
     \qed \vspace{2mm}
$$

2) By definition
   (\ref{nu1}) $\nu^1(b) \subseteq X_0^1$, and
   by the permutation property (\ref{permutation2}) we have
$$
   X_0^1 + G_i(G'_i + \nu^1(a)) =   X_0^1 + G'_i(G_i + \nu^1(a)).
$$
   By the modular law (\ref{modular_law}) and by (\ref{permutation2})
\begin{equation*}
\begin{split}
   & \psi_i(a)\psi_i(b) = \\
   & X_0^1 + G_i(G'_i + \nu^1(a))
     \left (X_0^1 + G'_i(G_i + \nu^1(b)) \right ) = \\
   & X_0^1 + G_i \left (X_0^1(G'_i + \nu^1(a)) +
         G'_i(G_i + \nu^1(b)) \right ).
\end{split}
\end{equation*}
   Since
$$
    X_0^1(G'_i + \nu^1(a)) = X_0^1{G'_i} + \nu^1(a)
    \indent \text{ and } \indent X_0^1{G'_i} = \nu^1(e_i),
$$
   we see that
 \begin{equation}
   \label{mul_D4}
     \psi_i(a)\psi_i(b) =
        X_0^1 + G_i \left (\nu^1(e_i) + \nu^1(a) +
          G'_i(G_i + \nu^1(b)) \right ).
 \end{equation}
    By the permutation property (\ref{permutation1}) and by (\ref{mul}) we have
 \begin{equation}
   \label{mul1_D4}
     \psi_i(a)\psi_i(b) =
       X_0^1 + G_i \left (G'_i +
          (\nu^1(e_i) + \nu^1(a))(G_i + \nu^1(b)) \right ).
 \end{equation}
    Since
$$
   G_i = G'_j{G'_k}G'_l \indent  \text{ and } \indent
      \nu^1(e_i) + \nu^1(a) = \nu^1(e_i + a) =
      X_0^1(\nu^1(e_i + a)),
$$
   it follows that
\begin{equation*}
\begin{split}
      & \psi_i(a)\psi_i(b) =
          X_0^1 + G_i \left (G'_i +
          \nu^1(e_i + a)(X_0^1{G'_j}{G'_k}G'_l + \nu^1(b)) \right ) = \\
      & X_0^1 + G_i \left (G'_i +
          \nu^1(e_i + a)(\nu^1(e_j{e}_k{e}_l) + \nu^1(b)) \right ) = \\
      & X_0^1 + G_i \left (G'_i +
          \nu^1(e_i + a)\nu^1(e_j{e}_k{e}_l + b) \right ) =
          \psi_i((a + e_i)(b + e_j{e}_k{e}_l)). \qed \vspace{2mm}
\end{split}
\end{equation*}

3) From (\ref{mul_D4}) and since
  $\nu^1(e_i) = X_0^1{G'_i} \subseteq G'_i$, we have
\begin{equation}
 \label{eq_mul_D4}
     \psi_i(a)\psi_i(b) = X_0^1 + G_i \left (\nu^1(a) +
      G'_i(G_i + \nu^1(b) + \nu^1(e_i)) \right ).
 \end{equation}
 Again, by (\ref{permutation1}) we have
\begin{equation*}
\begin{split}
   & \psi_i(a)\psi_i(b) = X_0^1 + G_i \left (G'_i +
          \nu^1(a)(G_i + \nu^1(b) + \nu^1(e_i)) \right ) = \\
   &  X_0^1 + G_i \left ( G'_i + \nu^1(a)(X_0^1{G'_j}{G'_k}G'_l
             + \nu^1(b) + \nu^1(e_i)) \right ).
\end{split}
\end{equation*}
 Thus,  $\psi_i(a)\psi_i(b) = \psi_i(a(b + x_i + x_j{x}_k{x}_l))$.
       \qedsymbol \vspace{2mm}

We need the following corollary (atomic multiplicativity) from
Proposition \ref{additivity_D4}
\begin{corollary}
  \label{cor_mul_d4}
  1) Suppose one of the following
  inclusions holds\footnote{Compare with Corollary \ref{cor_mul}, case
  $D^{2,2,2}$.}:
\begin{equation*}
\begin{split}
   & {\rm(i)} \indent e_i + e_j{e}_k{e}_l \subseteq a, \\
   & {\rm (ii)} \indent e_i + e_j{e}_k{e}_l \subseteq b,  \\
   & {\rm (iii)} \indent e_i \subseteq a, \indent
              e_j{e}_k{e}_l \subseteq b, \\
   & {\rm (iv)} \indent e_i \subseteq b, \indent
              e_j{e}_k{e}_l \subseteq a. \\
\end{split}
\end{equation*}
 Then the joint map $\psi_i$ operates as a homomorphism on the elements
 $a$ and $b$ with respect to the
 lattice operations $+$ and $\cap$, i.e.,
$$
   \psi_i(a) + \psi_i(b) = \psi_i(a+b), \hspace{3mm}
    \psi_i(a)\psi_i(b) = \psi_i(a)\psi_i(b).
$$

2) The joint map ${\psi_i}$ applied to the following atomic
elements
    is the intersection preserving map, i.e., multiplicative with respect
    to the operation $\cap$:
\begin{equation}
     \label{homo1_D4}
     \psi_i(ba_n^{ij}) = \psi_i(b)\psi_i(a_n^{ij})
             \text{ for every } b \subseteq D^4.
\end{equation}
\end{corollary}

\section{The action of maps $\psi_i$ and $\varphi_i$
            on the atomic elements in $D^4$}
\begin{proposition}
   \label{action_psi_D4}
  The joint maps $\psi_i$ applied to the atomic elements $a^{ij}_n$
 satisfy the
 following relations\footnote{Compare with Proposition \ref{action_psi},
 the modular lattice $D^{2,2,2}$.}
 \begin{enumerate}
   \item $\psi_i(a^{ij}_n) = \nu^0(e_i{a}^{kl}_n)$, \vspace{2mm}
   \item $\psi_j(a^{ij}_n) = \nu^0(e_j(e_k + e_l))$, \vspace{2mm}
   \item $\psi_j(e_i{a}^{kl}_n) = \nu^0(e_j{a}^{kl}_{n+1})$. \vspace{2mm}
  \end{enumerate}
\end{proposition}

\PerfProof 1) Since $\psi_i(e_i) = X^1_0$ (Proposition
\ref{basic_eq_D4}, heading (1)), we have
\begin{equation*}
    \psi_i(a^{ij}_n) = \psi_i(e_i + e_j{a}^{kl}_{n-1}) =
    \psi_i(e_i) + \psi_i(e_j{a}^{kl}_{n-1}) =
   \psi_i(e_j{a}^{kl}_{n-1}).
   \vspace{2mm}
\end{equation*}
We suppose that heading (3) of Proposition \ref{action_psi_D4} for
$n-1$ is true (induction hypothesis), and we get
$$
   \psi_i(a^{ij}_n) = \psi_i(e_j{a}^{kl}_{n-1}) =
   \nu^0(e_j{a}^{kl}_n).
   \qed \vspace{2mm}
$$

2) Here,
\begin{equation*}
   \psi_j(a^{ij}_n) = \psi_j(e_i + e_j{a}^{kl}_{n-1}) =
   \psi_j(e_i) + \psi_j(e_j{a}^{kl}_{n-1}) =
   \psi_j(e_i). \vspace{2mm}
\end{equation*}
Further, by Proposition \ref{basic_eq_D4}, heading (2), we have
$$
   \psi_j(a^{ij}_n) = \psi_j(e_i) = \nu^0(e_j(e_k + e_l)).
   \qed \vspace{2mm}
$$

3) For convenience, without loss of generality, we will show that
\begin{equation}
  \label{e2_a_power_34}
  \psi_1(e_2{a}^{34}_n) = \nu^0(e_1{a}^{34}_{n+1}). \vspace{2mm}
\end{equation}
By permutation property (\ref{permutation1}) and since
  $\nu^1(e_4{a}^{12}_{n-1}) \subseteq X^1_0$,
we have
\begin{equation}
 \label{expand_nu_1}
 \begin{split}
    \nu^1(e_2{a}^{34}_n) = & \nu^1(e_2(e_3 + e_4{a}^{12}_{n-1})) = \\
    & X_0^1{G'_2}(X_0^1{G'_3} + \nu^1(e_4{a}^{12}_{n-1})) =
    X_0^1{G'_2}(G'_3 + \nu^1(e_4{a}^{12}_{n-1})).
    \vspace{2mm}
 \end{split}
\end{equation}
By (\ref{expand_nu_1}) we have
\begin{equation*}
 \label{expand_nu_2}
 \begin{split}
    \psi_1(e_2{a}^{34}_n) =
      & X^1_0 + G_1(G'_1 + \nu^1(e_2{a}^{34}_n)) = \\
    & X^1_0 + G_1(G'_1 +
       X_0^1{G'_2}(G'_3 + \nu^1(e_4{a}^{12}_{n-1})).
    \vspace{2mm}
 \end{split}
\end{equation*}
Since $G_1 \subseteq G'_2$ and $G_1 \subseteq G'_3$, by
permutation property (\ref{permutation1}), we have
\begin{equation*}
 \begin{split}
    \psi_1(e_2{a}^{34}_n) =
     & X^1_0 + G_1(G'_1{G'_2} +
       X_0^1(G'_3 + \nu^1(e_4{a}^{12}_{n-1})) = \\
    & X^1_0 + G_1(G_3 + G_4 +
       X_0^1(G'_3 + \nu^1(e_4{a}^{12}_{n-1})) = \\
    & X^1_0 + G_1((G_3 + G_4)(G'_3 + \nu^1(e_4{a}^{12}_{n-1})) +
       X_0^1).
    \vspace{2mm}
 \end{split}
\end{equation*}
Since $G_4 \subseteq G'_3$, we have
\begin{equation*}
 \begin{split}
    \psi_1(e_2{a}^{34}_n) =
    & X^1_0 + G_1( G_4 + G_3(G'_3 + \nu^1(e_4{a}^{12}_{n-1})) +
       X_0^1) = \\
    & X^1_0 + G_1((X^1_0 + G_4) + (X^1_0 + G_3(G'_3 +
      \nu^1(e_4{a}^{12}_{n-1}))).
    \vspace{2mm}
 \end{split}
\end{equation*}
Pay attention to the fact
$$
  \psi_3(e_4{a}^{12}_{n-1}) =
    X^1_0 + G_3(G'_3 + \nu^1(e_4{a}^{12}_{n-1})),
$$
and therefore
\begin{equation}
 \label{expand_nu_4}
  \begin{split}
     \psi_1(e_2{a}^{34}_n) =
    & X^1_0 + G_1((X^1_0 + G_4) + \psi_3(e_4{a}^{12}_{n-1})) = \\
    & (X^1_0 + G_1)((X^1_0 + G_4) + \psi_3(e_4{a}^{12}_{n-1})) = \\
    &  \nu^0(e_1)(\nu^0(e_4) + \psi_3(e_4{a}^{12}_{n-1})).
  \end{split}
    \vspace{2mm}
\end{equation}

By induction hypothesis for $n-1$ in heading (3) we have
$$
   \psi_3(e_4{a}^{12}_{n-1}) = \nu^0(e_3{a}^{12}_n),
$$
and by (\ref{expand_nu_4}) we have
$$
     \psi_1(e_2{a}^{34}_n) =
      \nu^0(e_1)(\nu^0(e_4) + \nu^0(e_3{a}^{12}_n)) =
      \nu^0(e_1(e_4 + e_3{a}^{12}_n)) =
      \nu^0(e_1{a}^{43}_{n+1}).
      \qed \vspace{2mm}
$$

 From Propositions (\ref{phi_and_psi}) and
(\ref{action_psi_D4}) we have
\begin{corollary}
   \label{act_runner_map_D4}
 The elementary maps $\varphi_i$ applied to the
atomic elements $a^{ij}_n$, where $n \geq 1$, satisfy the
following relations\footnote{Compare with Corollary
\ref{act_runner_map}, case $D^{2,2,2}$}
   \item \vspace{2mm}
   \begin{enumerate}
   \item 
     $\varphi_i\Phi^+\rho(a^{ij}_n) = \rho(e_i{a}^{kl}_n)$, \vspace{2mm}
   \item 
     $\varphi_j\Phi^+\rho(a^{ij}_n) = \rho(e_j(e_k + e_l))$, \vspace{2mm}
   \item 
     $\varphi_j\Phi^+\rho(e_i{a}^{kl}_n) =
                         \rho(e_j{a}^{kl}_{n+1})$. \vspace{2mm}
   \end{enumerate}
\end{corollary}

\section{The fundamental property of the elementary maps}
 \label{sect_fundam_prop_D4}

\begin{proposition}
 \label{motiv_admis_D4}
  For $\{i,j,k,l\} = \{1,2,3,4\}$ the following
  relations hold\footnote{Compare with Proposition \ref{motiv_admis}, case
$D^{2,2,2}$}
\begin{equation}
  \label{fund_varphi_rel_D4}
     \varphi_i\varphi_k\varphi_j +  \varphi_i\varphi_l\varphi_j = 0, \vspace{2mm}
\end{equation}
\begin{equation}
     \varphi_i^2 = 0.
\end{equation}
\end{proposition}
\PerfProof
  For every vector $v \in X_0^1$, by definition of $\varphi_i$
  we have $(\varphi_i + \varphi_j + \varphi_k + \varphi_l)(v) = 0$,
  see eq. (\ref{sum_fi_D4}).
  In other words, $\varphi_i + \varphi_j + \varphi_k + \varphi_l$ = 0. Therefore,
$$
    \varphi_i\varphi_k\varphi_j +  \varphi_i\varphi_l\varphi_j =
    \varphi_i(\varphi_i + \varphi_j)\varphi_j =
    \varphi_i^2\varphi_j + \varphi_i\varphi_j^2.
$$
  So, it suffices to prove that $\varphi_i^2$ = 0.
  For every
   $z \subseteq D^4$, by Corollary \ref{cor_psi_D4}, headings (3) and (1), we have
\begin{equation*}
\begin{split}
  & \varphi_i^2((\Phi^+)^2\rho)(z) \subseteq
      \varphi_i(\varphi_i((\Phi^+)^2\rho(I))) =  \\
  & \varphi_i(\Phi^+\rho(e_i(e_j + e_k + e_l)))
      \subseteq \varphi_i(\Phi^+\rho(e_i)) = 0.
  \qed \vspace{2mm}
\end{split}
\end{equation*}

 \begin{corollary}
   \label{basic_rel_D4}
     The relation
 \begin{equation}
   \label{fund_D4}
      \varphi_i\varphi_k\varphi_j(B)
       = \varphi_i\varphi_l\varphi_j(B)
 \end{equation}
 takes place\footnote{Compare with Corollary \ref{basic_rel}, case
$D^{2,2,2}$}
 for every subspace $B \subseteq X_0^2$,
 where $X_0^2$ is the representation space
 of $(\Phi_i^+)^2\rho$.
 \end{corollary}

Essentially, relations (\ref{fund_varphi_rel_D4}) and
(\ref{fund_D4}) are fundamental and motivate the construction of
the admissible sequences satisfying the following relation:
\begin{equation}
  \label{fund_indices}
   ikj = ilj,
\end{equation}
where indices $i,j,k,l$ are all distinct, see Table
 \ref{compare_notions} and \S\ref{subsect_adm_seq_D4}.

\section{The $\varphi_i-$homomorphic elements in $D^4$}
  \label{subs_homom_D4}
By analogy with the modular lattice $D^{2,2,2}$ (see
\S\ref{psi_homom_L6}),
 we introduce now $\varphi_i-$homomorphic polynomials in $D^4$.

\index{$\varphi_i-$homomorphic polynomial}

   An element $a \subseteq D^4$ is said to be
    {\it $\varphi_i-$homomorphic}, if
    \begin{equation}
       \varphi_i\Phi^+\rho(ap) =
          \varphi_i\Phi^+\rho(a)\varphi_i\Phi^+\rho(p)
            \text{ for all } p \subseteq D^4.
    \end{equation}

    An element $a \subseteq D^4$ is said to be
    {\it $(\varphi_i, e_k)-$homomorphic}, if
    \begin{equation}
        \varphi_i\Phi^+\rho(ap) =
           \varphi_i\Phi^+\rho(e_k{a})\varphi_i\Phi^+\rho(p)
            \text{ for all } p \subseteq e_k.
    \end{equation}

\begin{theorem}
  \label{th_homomorhism_D4}
  1) The polynomials $a^{ij}_n$
  are $\varphi_i-$homomorphic\footnote{Compare with
    Theorem \ref{th_homomorhism}, case $D^{2,2,2}$.}.
    \vspace{2mm}

  2) The polynomials $a^{ij}_n$ are ($\varphi_j,e_k)-$homomorphic
     for distinct indices $\{i,j,k\}$.
            \vspace{2mm}
\end{theorem}
\PerfProof
  1) Follows from multiplicativity (\ref{homo1_D4}),
  Proposition \ref{action_psi_D4} and
  the analog of Proposition \ref{phi_and_psi} for $D^4$.

  2) For convenience, without loss of generality, we will show that
\begin{equation}
  \label{p_a12_D4}
    \varphi_2\Phi^+\rho(a^{12}_n{p}) =
           \varphi_2\Phi^+\rho(e_3{a})\varphi_2\Phi^+\rho(p)
            \text{ for all } p \subseteq e_3.
\end{equation}
  Let $p \subseteq e_3$. First,
\begin{equation*}
      p{a}^{12}_s = p{e}_3{a}^{12}_s = p{e}_3{a}^{21}_s =
      p{a}^{21}_s.
\end{equation*}
  By Corollary \ref{cor_mul_d4} and Proposition \ref{action_psi_D4}, heading (3)
  we have
\begin{equation}
  \label{p_aij_D4_2}
  \begin{split}
     & \psi_2(p{a}^{12}_n) = \psi_2(p{e}_3{a}^{21}_n) =
       \psi_2(p(e_2 + e_1{a}^{34}_{n-1})) = \\
     & \psi_2(p(e_1{a}^{34}_{n-1})) =
       \psi_2(p)\psi_2(e_1{a}^{34}_{n-1}) = \\
     & \psi_2(p)\nu^0(e_2{a}^{34}_n).
  \end{split}
\end{equation}
Since $e_2{a}^{34}_n = e_2(e_1 + e_3 + e_4){a}^{34}_n$ and
\begin{equation*}
  \begin{split}
    & \nu^0(e_2{a}^{34}_n) = \nu^0(e_2(e_1 + e_3 + e_4){a}^{34}_n) = \\
    & \nu^0(e_2(e_1 + e_3 + e_4))\nu^0({a}^{34}_n) =
      \psi_2(I)\nu^0({a}^{34}_n),
  \end{split}
\end{equation*}
by (\ref{p_aij_D4_2}) we have
\begin{equation}
 \label{p_aij_D4_3}
    \psi_2(p{a}^{12}_n) = \psi_2(p)\psi_2(I)\nu^0({a}^{34}_n) =
     \psi_2(p)\nu^0({a}^{34}_n).
\end{equation}
We have $\psi_2(p) \subseteq \psi_2(e_3)$ together with
 $p \subseteq e_3$, and therefore
   $\psi_2(e_3{p}) = \psi_2(e_3)\psi_2(p)$.
 From (\ref{p_aij_D4_3}) we get
\begin{equation*}
  \begin{split}
  &  \psi_2(p{a}^{12}_n) =
     \psi_2(p)\psi_2(e_3)\nu^0({a}^{34}_n) = \\
  &  \psi_2(p)\nu^0(e_2(e_1 + e_4))\nu^0({a}^{34}_n) = \\
  &  \psi_2(p)\nu^0(e_2(e_1 + e_4)({a}^{34}_n)) =
     \psi_2(p)\psi_2(e_3{a}^{21}_n),
  \end{split}
\end{equation*}
i.e.,
\begin{equation}
  \label{before_nabla_D4}
     \psi_2(p{a}^{12}_n) = \psi_2(p)\psi_2(e_3{a}^{21}_n).
     \vspace{2mm}
\end{equation}
 Applying projection $\nabla$ to (\ref{before_nabla_D4}) as in Theorem
\ref{th_homomorhism} we get (\ref{p_a12_D4}) and theorem is
proved.
 \qedsymbol \vspace{2mm}

\section{The theorem on the classes of admissible elements in $D^4$}
 \label{sect_adm_classes_D4}
\begin{theorem}
   \label{th_adm_classes_D4}
   Let $\alpha = i_n{i}_{n-1}\dots{1}$
   be an admissible sequence for $D^4$ and $i \neq i_n$.
   Then ${i}\alpha$ is admissible and, for
    $z_\alpha = e_\alpha$ or $f_{\alpha0}$ from
    Table \ref{table_adm_elem_D4},
    the following relation holds\footnote{Compare with
    Theorem \ref{th_adm_classes}, case $D^{2,2,2}$}:
 \begin{equation}
   \label{adm_classes_D4}
    \varphi_i\Phi^+\rho(z_\alpha) = \rho(z_{i\alpha}).
 \end{equation}
\end{theorem}

For the proof of the theorem on admissible elements in $D^4$, see
\S\ref{proof_adm_D4}. The proof repeatedly uses the basic
properties of the admissible sequences in $D^4$ considered in
\S\ref{basic_adm_D4}, Lemma \ref{homom_polynom_P}.

\begin{table}[h]
 \renewcommand{\arraystretch}{1.3}
  \begin{tabular} {||c|c|c|c||}
  \hline \hline
     & Admissible
     & Admissible & Admissible  \cr
     & sequence $\alpha$ & polynomial $e_\alpha$
     & polynomial $f_{{\alpha}0}$ \\
  \hline         
    $F21$  & $(21)^t(41)^r(31)^s = $
     & $e_2{a}^{31}_{2s}a^{41}_{2r}a^{34}_{2t-1}$
     & $e_\alpha(e_2{a}^{34}_{2t} + a^{41}_{2r+1}a^{31}_{2s-1}) = $ \cr
     & $(21)^t(31)^s(41)^r $ &
     & $e_\alpha({a}^{43}_{2t} + e_1a^{24}_{2r}a^{23}_{2s})$ \\
  \hline         
    $F31$  & $(31)^s(41)^r(21)^t = $
     & $e_3{a}^{21}_{2t}a^{41}_{2r}a^{24}_{2s-1}$
     & $e_\alpha(e_3{a}^{42}_{2s} + a^{41}_{2r+1}a^{21}_{2t-1}) = $ \cr
     & $(31)^s(21)^t(41)^r  $ &
     & $e_\alpha({a}^{42}_{2s} + e_1a^{34}_{2r}a^{32}_{2t})$ \\
  \hline      
    $F41$  & $(41)^r(31)^s(21)^t = $
     & $e_4{a}^{21}_{2t}a^{31}_{2s}a^{32}_{2r-1}$
     & $e_\alpha(e_4{a}^{32}_{2r} + a^{31}_{2s+1}a^{21}_{2t-1}) = $ \cr
     & $(41)^r(21)^t(31)^s$ &
     & $e_\alpha({a}^{32}_{2r} + e_1a^{43}_{2s}a^{42}_{2t})$ \\
  \hline \hline        
    $G11$  & $1(41)^r(31)^s(21)^t$ = & &  \cr
     &   \hspace{3mm} $1(31)^s(41)^r(21)^t$ =
     & $e_1{a}^{24}_{2s}a^{34}_{2t}a^{32}_{2r}$
     & $e_\alpha(e_1{a}^{32}_{2r+1} + a^{24}_{2s+1}a^{34}_{2t-1}) = $ \cr
     &   \hspace{6mm} $1(21)^t(31)^s(41)^r$ &
     & $e_\alpha({a}^{32}_{2r+1} + e_4a^{21}_{2s}a^{31}_{2t})$   \\
  \hline         
    $G21$  & $2(41)^r(31)^s(21)^t =$
     & $e_2{a}^{34}_{2t}a^{31}_{2s+1}a^{14}_{2r-1} =$
     & $e_\alpha(e_2{a}^{14}_{2r} + a^{31}_{2s+2}a^{34}_{2t-1}) =$ \cr
     & $2(31)^{s+1}(41)^{r-1}(21)^t$
     & $e_2{a}^{34}_{2t}a^{31}_{2s-1}a^{14}_{2r+1}$
     & $e_\alpha({a}^{14}_{2r} + e_3a^{21}_{2s+1}a^{24}_{2t})$  \\
  \hline         
     $G31$  & $3(41)^r(21)^t(31)^s =$
     & $e_3{a}^{24}_{2s}a^{21}_{2t+1}a^{14}_{2r-1} =$
     & $e_\alpha(e_3{a}^{14}_{2r} + a^{21}_{2t+2}a^{24}_{2s-1}) =$ \cr
     & $3(21)^{t+1}(41)^{r-1}(31)^s$
     & $e_3{a}^{24}_{2s}a^{21}_{2t-1}a^{14}_{2r+1}$
     & $e_\alpha({a}^{14}_{2r} + e_2a^{31}_{2t+1}a^{34}_{2s})$ \\
  \hline      
    $G41$  & $4(21)^t(31)^s(41)^r =$
     & $e_4{a}^{32}_{2r}a^{31}_{2s+1}a^{12}_{2t-1} =$
     & $e_\alpha(e_4{a}^{12}_{2t} + a^{31}_{2s}a^{32}_{2r+1}) = $ \cr
     & $4(31)^{s+1}(21)^{t-1}(41)^r$
     & $e_4{a}^{32}_{2r}a^{31}_{2s-1}a^{12}_{2t+1}$
     & $e_\alpha({a}^{12}_{2t} + e_3a^{41}_{2s+1}a^{42}_{2r})$ \\
  \hline \hline     
    $H11$  & $(14)^r(31)^s(21)^t$ =
     & & \cr
     & $(14)^r(21)^{t+1}(31)^{s-1}$ =
     & $e_1a^{23}_{2r+1}{a}^{24}_{2s-1}a^{34}_{2t-1} = $
     & $e_\alpha(e_1a^{23}_{2r} + a^{24}_{2s}a^{34}_{2t}) = $ \cr
     & $(13)^s(41)^r(21)^t =$
     & $e_1a^{23}_{2r-1}{a}^{24}_{2s-1}a^{34}_{2t+1} = $
     & $e_\alpha(e_1a^{24}_{2s} + a^{34}_{2t}a^{23}_{2r}) = $ \cr
     & $(13)^s(21)^{t+1}(41)^{r-1}$
     & $e_1a^{23}_{2r-1}{a}^{24}_{2s+1}a^{34}_{2t-1}$
     & $e_\alpha(e_1a^{34}_{2t} + a^{23}_{2r}a^{24}_{2s})$ \cr
     & $(12)^{t+1}(41)^r(31)^{s-1}$ = & & \cr
     & $(12)^{t+1}(31)^{s-1}(41)^r$ & & \\
  \hline      
    $H21$  & $2(14)^r(31)^s(21)^t$
     & $e_2{a}^{14}_{2r-1}a^{31}_{2s-1}a^{34}_{2t+2} = $
     & $e_\alpha(e_2a^{14}_{2r} + a^{31}_{2s}a^{34}_{2t+1}) = $ \cr
     &
     & $e_2{a}^{14}_{2r-1}a^{31}_{2s+1}a^{34}_{2t}$
     & $e_\alpha(e_2a^{34}_{2t+1} + a^{31}_{2s}a^{14}_{2r})$ \\
  \hline      
    $H31$  & $3(14)^r(31)^s(21)^t$
     & $e_3{a}^{14}_{2r-1}a^{21}_{2t-1}a^{24}_{2s + 2} = $
     & $e_\alpha(e_3a^{14}_{2r} + a^{24}_{2s+1}a^{21}_{2t}) = $ \cr
     &
     & $e_3{a}^{14}_{2r-1}a^{21}_{2t+1}a^{24}_{2s}$
     & $e_\alpha(e_3a^{24}_{2s+1} + a^{14}_{2r}a^{21}_{2t})$ \\
  \hline      
    $H41$  & $4(14)^r(31)^s(21)^t$
     &  $e_4{a}^{13}_{2s-1}a^{21}_{2t-1}a^{23}_{2r+2} = $
     &  $e_\alpha(e_4a^{23}_{2r+1} + a^{13}_{2s}a^{21}_{2t}) = $ \cr
     &
     &  $e_4{a}^{13}_{2s-1}a^{21}_{2t+1}a^{23}_{2r}$
     &  $e_\alpha(e_4a^{21}_{2t} + a^{23}_{2r+1}a^{13}_{2s})$ \\
  \hline  \hline
  \end{tabular}
  \vspace{2mm}
  \caption{\hspace{3mm}Admissible polynomials in the modular lattice $D^4$}
  \footnotesize
  \begin{tabular}{l}
    Notes to Table: \cr
  1) For more details about admissible sequences,
     see Proposition \ref{full_adm_seq_D4} and
     Table \ref{table_admissible_ExtD4}. \cr
  2) For relations given in two last columns (definitions of admissible
     polynomials $e_\alpha$ and $f_{{\alpha}0}$), \cr
     see Lemma \ref{homom_polynom_P}. \cr
  3) In each line, each low index should be non-negative. For
     example, for Line $F21$, \cr
     we have:
     $s \geq 0, r \geq 0, t \geq 1$;
     for Line $G21$, we have: $s \geq 0, r \geq 0, t \geq 0$.
  \label{table_adm_elem_D4}
  \end{tabular}
 \normalsize
\end{table}

\subsection{Basic properties of admissible elements in $D^4$}
  \label{basic_adm_D4}
We prove here a number of
 basic properties\footnote{Compare with Table \ref{table_atomic} and
 \S\ref{basic_atomic_L6}, case $D^{2,2,2}$.}
of the atomic elements in $D^4$ used in the proof of the theorem
on admissible elements (Theorem \ref{th_adm_classes_D4}). In
particular, in some cases the lower indices of polynomials
$a^{ij}_s$ entering in the admissible elements $f_{{\alpha}0}$ can
be transformed as in the following
\begin{lemma}
 \label{homom_polynom_P}
 1) Every polynomial $f_{{\alpha}0}$ from Table
 \ref{table_adm_elem_D4} can be represented as an intersection of
$e_\alpha$ and $P$. For every $i \neq i_n$ (see \S\ref{adm_D4}),
 we select $P$ to be some $\varphi_i-$homomorphic
polynomial.

2) The lower indices of polynomials $a^{ij}_s$ entering in the
admissible elements $f_{{\alpha}0}$ can be equalized as follows:
\begin{equation}
 \label{sym_forms_P}
  \begin{split}
      f_{{\alpha}0} =
    & f_{(21)^t(41)^r(31)^s0} = \\
    & e_\alpha(e_2{a}^{34}_{2t} + a^{41}_{2r+1}a^{31}_{2s-1}) =
      e_\alpha({a}^{43}_{2t} + e_1a^{24}_{2r}a^{23}_{2s}).
  \end{split}
\end{equation}
The generic relation\footnote{Throughout this lemma
    we suppose that $\{i,j,k,l\} = \{1,2,3,4\}$.}
 is the following:
\begin{equation}
 \label{sym_forms_P_2}
  \begin{split}
    e_i(e_i{a}^{jl}_{t} + & a^{kj}_{r+1}a^{kl}_{s-1}) = \\
    & e_i({a}^{jl}_{t} + e_i{a}^{kj}_{r+1}a^{kl}_{s-1}) =
      e_i({a}^{jl}_{t} + e_k{a}^{ij}_{r}a^{il}_{s}).
  \end{split}
\end{equation}

3) The substitution
\begin{equation}
 \label{subst_r_s}
  \begin{cases}
     r \mapsto r-2, \\
     s \mapsto s+2
  \end{cases}
\end{equation}
does not change the polynomial $e_i(e_i{a}^{jl}_{t} +
a^{kj}_{r+1}a^{kl}_{s-1})$, namely:
\begin{equation}
 \label{sym_forms_P_3}
    e_i(e_i{a}^{jl}_{t} + a^{kj}_{r+1}a^{kl}_{s-1}) =
    e_i(e_i{a}^{jl}_{t} + a^{kj}_{r-1}a^{kl}_{s+1}).
\end{equation}

4) The substitution (\ref{subst_r_s}) does not change the
polynomial $e_ia^{kj}_{s}a^{kl}_{r}$:
\begin{equation}
 \label{sym_forms_P_4}
    e_ia^{kj}_{r+1}a^{kl}_{s-1} =
    e_ia^{kj}_{r-1}a^{kl}_{s+1}.
\end{equation}

\end{lemma}
\PerfProof 1) We consider here only Line $F21$ from Table
 \ref{table_adm_elem_D4}, all other cases are similarly considered.
Every polynomial $f_{{\alpha}0}$ from Table
 \ref{table_adm_elem_D4} is the intersection $e_\alpha$ and $P$,
where $P$ is the sum contained in the parantheses, and polynomial
$f_{{\alpha}0}$ can be represented in the different equivalent
forms, such that only the polynomial $P$ is changed. For Line
$F21$, we have
\begin{equation}
 \label{diff_forms_P}
  \begin{array}{ccc}
      f_{{\alpha}0} = & f_{(21)^t(41)^r(31)^s0} = & \vspace{2mm} \\
      e_\alpha(e_2{a}^{34}_{2t} + a^{41}_{2r+1}a^{31}_{2s-1}) = &
      e_\alpha(e_2{a}^{43}_{2t} + a^{41}_{2r+1}a^{31}_{2s-1}) = &
    \qquad \qquad \text{ (i)} \vspace{2mm} \\
      e_\alpha({a}^{34}_{2t} + e_2a^{41}_{2r+1}a^{31}_{2s-1}) = &
      e_\alpha({a}^{43}_{2t} + e_2a^{41}_{2r+1}a^{31}_{2s-1}) = &
    \qquad \qquad \text{ (ii)} \vspace{2mm} \\
      e_\alpha({a}^{34}_{2t} + e_2a^{14}_{2r+1}a^{13}_{2s-1}) = &
      e_\alpha({a}^{43}_{2t} + e_2a^{14}_{2r+1}a^{13}_{2s-1}) = &
    \qquad \qquad \text{ (iii)} \vspace{2mm} \\
      e_\alpha(e_2{a}^{34}_{2t} + a^{14}_{2r+1}a^{13}_{2s-1}) = &
      e_\alpha(e_2{a}^{43}_{2t} + a^{14}_{2r+1}a^{13}_{2s-1}) = &
    \qquad \qquad \text{ (iv)} \vspace{2mm} \\
      e_\alpha(e_2{a}^{34}_{2t}a^{13}_{2s-1} + a^{14}_{2r+1}) = &
      e_\alpha(e_2{a}^{43}_{2t}a^{13}_{2s-1} + a^{14}_{2r+1}) = &
  \end{array}
\end{equation}
Relation (i) in (\ref{diff_forms_P}) follows from the definition
of Table \ref{table_adm_elem_D4} and (\ref{main_atomic_D4}).
Relations (ii) in (\ref{diff_forms_P}) follow from the inclusion
$e_\alpha \subseteq e_2$ and permutation property
(\ref{permutation1}). Relations (iii) and (iv) follow from
(\ref{main_atomic_D4}). Further, by (\ref{incl_atomic_D4})
$$
  e_\alpha \subseteq a^{31}_{2s} \subseteq a^{31}_{2s-1}
$$
and by permutation property (\ref{permutation1})
 we have the relation (v).

   The polynomial $P$ in (iv) is $\varphi_1-$homomorphic since
$$
   e_1 \subseteq a^{14}_{2r+1}a^{13}_{2s-1} \text{ and }
   e_2e_3e_4 \subseteq e_2a^{34}_{2t}.
$$

   The polynomial $P$ in (iii) is $\varphi_3-$homomorphic since
$$
   e_3 \subseteq a^{34}_{2t}
    \text{ and }
   e_1e_2e_4 \subseteq e_1e_2 \subseteq e_2a^{14}_{2r+1}a^{13}_{2s-1}.
$$

   The polynomial $P$ in (iii) is $\varphi_4-$homomorphic since
$$
   e_4 \subseteq a^{43}_{2t}
    \text{ and }
   e_1e_2e_3 \subseteq e_1e_2 \subseteq e_2a^{14}_{2r+1}a^{13}_{2s-1}.
$$

Essentially, heading 1) allows to select an appropriate form of
the polynomial $P$ for every given $\varphi_i$.
 \qedsymbol \vspace{2mm}

 2) Let us prove (\ref{sym_forms_P}). By
(\ref{equalize_D4}) we have
\begin{equation*}
  \begin{split}
      & e_\alpha(e_2{a}^{34}_{2t} + a^{41}_{2r+1}a^{31}_{2s-1}) =
        e_\alpha({a}^{43}_{2t} + e_2a^{41}_{2r+1}a^{31}_{2s-1}) = \\
      & e_\alpha({a}^{43}_{2t} + (e_4 + e_2a^{41}_{2r+1}a^{31}_{2s-1})) =
        e_\alpha({a}^{34}_{2t} + e_4 + e_1a^{24}_{2r}a^{23}_{2s}) = \\
      & e_\alpha({a}^{34}_{2t} + e_1a^{24}_{2r}a^{23}_{2s}).
  \end{split}
\end{equation*}

Now consider (\ref{sym_forms_P_2}).
 Since $e_j \subseteq a^{jl}_{t}$, by (\ref{equalize_D4})
we have
\begin{equation*}
  \begin{split}
      & e_i({a}^{jl}_{t} + e_i{a}^{kj}_{r+1}a^{kl}_{s-1}) =
        e_i({a}^{jl}_{t} + (e_j  + e_i{a}^{kj}_{r+1}a^{kl}_{s-1})) = \\
      & e_i({a}^{jl}_{t} + e_k{a}^{ij}_{r+1}a^{il}_{s-1})).
 \qed \vspace{2mm}
  \end{split}
\end{equation*}

3) By (\ref{sym_forms_P_2}) we have
\begin{equation*}
    e_i(e_i{a}^{jl}_{t} + a^{kj}_{r+1}a^{kl}_{s-1}) =
      e_i({a}^{jl}_{t} + e_k{a}^{ij}_{r}a^{il}_{s}),
\end{equation*}
On the other hand
\begin{equation*}
  \begin{split}
   & e_i(e_i{a}^{jl}_{t} + a^{kj}_{r-1}a^{kl}_{s+1}) = \\
   & e_i(e_i{a}^{lj}_{t} + a^{kl}_{s+1}a^{kj}_{r-1}) =
     e_i({a}^{lj}_{t} + e_k{a}^{il}_{s}a^{ij}_{r}),
  \end{split}
\end{equation*}
i.e.,
$$
    e_i(e_i{a}^{jl}_{t} + a^{kj}_{r+1}a^{kl}_{s-1}) =
    e_i(e_i{a}^{jl}_{t} + a^{kj}_{r-1}a^{kl}_{s+1}).
$$
and (\ref{sym_forms_P_3}) is true.

4) By (\ref{equalize_D4}) we have
\begin{equation*}
  \begin{split}
   & e_ia^{kj}_{r+1}a^{kl}_{s-1} = \\
   & e_i(e_k + e_j{a}^{il}_r)a^{kl}_{s-1} =
     e_i(e_k + e_j{a}^{il}_r{a}^{kl}_{s-1}) = \\
   & e_i(e_k + e_l{a}^{ij}_{r-1}{a}^{kj}_s) =
     e_i(e_k + e_l{a}^{ij}_{r-1}){a}^{kj}_{s} =
     e_i{a}^{kl}_r{a}^{kj}_s,
  \end{split}
\end{equation*}
i.e.,
\begin{equation}
  \label{subst_r_s_1}
     e_i{a}^{kj}_{r+1}{a}^{kl}_{s-1} =
     e_i{a}^{kl}_r{a}^{kj}_s.
\end{equation}
By repeating (\ref{subst_r_s_1}), we see that
\begin{equation*}
     e_i{a}^{kj}_{r+1}{a}^{kl}_{s-1} =
     e_i{a}^{kj}_{r-1}{a}^{kl}_{s+1}.
\end{equation*}
 The lemma is proved.
 \qedsymbol \vspace{2mm}

\subsection{Coincidence with the Gelfand-Ponomarev polynomials}
  \label{coinc_GP_D4}
Let $e_\alpha, f_{{\alpha}0}$ be admissible elements constructed
in this work and  $\tilde{e}_\alpha, \tilde{f}_{{\alpha}0}$ be the
admissible elements constructed by Gelfand and Ponomarev
\cite{GP74}. Since Theorem \ref{th_adm_classes_D4} takes place for
both $e_\alpha, f_{{\alpha}0}$ and for $\tilde{e}_\alpha,
\tilde{f}_{{\alpha}0}$ \cite[Th.7.2, Th.7.3]{GP74} we have
\begin{proposition}
 \label{coincidence_GP}
  The elements $e_\alpha$, (resp. $f_{{\alpha}0}$) and
  $\tilde{e}_\alpha$ (resp. $\tilde{f}_{{\alpha}0}$)
  coincide$\mod\theta$.
\end{proposition}

We will prove that this coincidence
 (without restriction$\mod\theta$) takes place for the admissible sequences of the
small length. Recall definitions of $\tilde{e}_\alpha$ and
$\tilde{f}_{{\alpha}0}$ from \cite{GP74}. \vspace{2mm}

\underline{The definition of $\tilde{e}_\alpha$,
\cite[p.6]{GP74}}.

\begin{equation}
  \begin{split}
   & \tilde{e}_{i_n{i}_{n-1}\dots{i}_2{i}_1} =
     \tilde{e}_{i_n}\sum\limits_{\beta \in \Gamma_{e}(\alpha)}\tilde{e}_\beta,
  \end{split}
\end{equation}
where
\begin{equation}
  \begin{split}
   \Gamma_{e}(\alpha) = \{ \beta = (k_{n-1},\dots,k_2,k_1) \mid~
     & k_{n-1} \notin \{i_n,{i}_{n-1}\},\dots, k_1 \notin
     \{i_2,{i}_1\},\text{ and } \\
    &  k_1 \neq k_2, \dots, k_{n-2} \neq k_{n-1}  \}.
  \end{split}
\end{equation}

\underline{The definition of $\tilde{f}_{{\alpha}0}$,
\cite[p.53]{GP74}}.

\begin{equation}
  \begin{split}
   & \tilde{f}_{i_n{i}_{n-1}\dots{i}_2{i}_1{0}} =
     \tilde{e}_{i_n}\sum\limits_{\beta \in \Gamma_{f}(\alpha)}\tilde{e}_\beta,
  \end{split}
\end{equation}
where
\begin{equation}
  \begin{split}
   \Gamma_{f}(\alpha) = \{ \beta = (k_{n},\dots,k_2,k_1) \mid~
     & k_n \notin \{i_n,{i}_{n-1}\},\dots, k_2 \notin
     \{i_2,{i}_1\}, k_1 \notin \{i_1\} \text{ and } \\
    &  k_1 \neq k_2, \dots, k_{n-2} \neq k_{n-1}  \}.
  \end{split}
\end{equation}

\begin{proposition}[The elements $\tilde{e}_\alpha$]
  \label{coincidence_E}
 Consider elements $\tilde{e}_\alpha$ for $\alpha = 21, 121, 321,
2341$ (see \S\ref{examples_D4}). The relation
$$
   e_\alpha = \tilde{e}_\alpha
$$
takes place without restriction$\mod\theta$.
\end{proposition}
\PerfProof \underline{For $n=2$}: $\alpha = 21$. We have
\begin{equation}
  \tilde{e}_{21} = e_2\sum\limits_{j \neq 1,2}{e_j} =
    e_2(e_3 + e_4).
\end{equation}
According to \S\ref{examples_D4} we see that $e_{21} =
\tilde{e}_{21}$.
 \qedsymbol \vspace{2mm}

 \underline{For $n=3$}: 1) $\alpha = 121$,
\begin{equation}
  \begin{split}
  & \Gamma_{e}(\alpha) = \{ (k_2 k_1) \mid k_2 \in \{3,4 \}, k_1 \in
  \{3, 4\}, k_1 \neq k_2 \}, \\
  & \tilde{e}_{121} = e_1\sum\limits_{\beta \in \Gamma_{e}(\alpha)}e_\beta =
     e_1(e_{34} + e_{43}) = e_1(e_3(e_1 + e_2) + e_4(e_1 + e_2)) = e_1a^{34}_2.
  \end{split}
\end{equation}
By \S\ref{examples_D4} we have $e_{121} = \tilde{e}_{121}$.
\qedsymbol \vspace{2mm}

 2) $\alpha = 321 = 341$. We have
\begin{equation*}
  \begin{split}
   \Gamma_{e}(\alpha) = & \{ (k_2 k_1) \mid k_2 \in \{3,2 \}, k_1 \in
  \{2, 1\}, k_1 \neq k_2 \} = \{ 14, 13, 43\}. \\
  \tilde{e}_{321} = & e_1\sum\limits_{\beta \in \Gamma_{e}(\alpha)}e_\beta =
     e_3(e_{14} + e_{13} + e_{43}) = \\
     & e_3(e_1(e_2 + e_4) + e_1(e_2 + e_3) + e_4(e_1 + e_2)) = \\
  &  e_3((e_1 + e_2)(e_2 + e_4)(e_1 + e_4) + e_1(e_2 + e_3)) = \\
  &  e_3((e_1 + e_2)(e_1 + e_4)(e_2 + e_4 + e_1(e_2 + e_3)) = \\
  &  e_3((e_1 + e_2)(e_1 + e_4)(e_4 + (e_1 + e_2)(e_2 + e_3)) = \\
  &  e_3((e_1 + e_2)(e_1 + e_4)(e_4(e_2 + e_3) + (e_1 + e_2)).
  \end{split}
\end{equation*}
Since $e_1 + e_2 \subseteq e_4(e_2 + e_3) + e_1 + e_2$, we have
\begin{equation}
   \tilde{e}_{321} = e_3(e_{14} + e_{13} + e_{43}) =
     e_3(e_1 + e_2)(e_1 + e_4).
\end{equation}
Since $\tilde{e}_{321}$ is symmetric with respect to transposition
$2 \leftrightarrow 4$ we have
\begin{equation}
   \tilde{e}_{321} = \tilde{e}_{341} =
   e_3(e_{14} + e_{13} + e_{43}) =
   e_3(e_{12} + e_{13} + e_{23}) =
     e_3(e_1 + e_2)(e_1 + e_4).
\end{equation}
By \S\ref{examples_D4} we have $e_{321} = \tilde{e}_{321}$.
\qedsymbol \vspace{2mm}

 \underline{For $n=4$}:  $\alpha = 2341 = 2321 = 2141$. We have
\begin{equation}
  \begin{split}
   \Gamma_{e}(\alpha) = & \{ (k_3 k_2 k_1) \mid k_3 \in \{1,4 \}, k_2 \in
  \{1, 2\}, k_1 \in \{ 2,3 \},\quad k_1 \neq k_2 , k_2 \neq k_3 \} = \\
  & \{ (123), (412), (413) = (423) \},
  \end{split}
\end{equation}
and
\begin{equation}
  \begin{split}
   & \tilde{e}_{2341} = e_2\sum\limits_{\beta \in \Gamma_{e}(\alpha)}e_\beta =
     e_2(e_{123} + e_{412} + e_{413}) = \\
   & e_2(e_1(e_2 + e_3)(e_4 + e_3) + e_4(e_1 + e_2)(e_3 + e_2) +
         e_4(e_1 + e_3)(e_2 + e_3)) = \\
   & e_2(e_4 + e_3)(e_1(e_2 + e_3) + e_4(e_1 + e_2)(e_3 + e_2) +
         e_4(e_1 + e_3)(e_3 + e_2)) =  \\
   & e_2(e_4 + e_3)(e_1(e_2 + e_3) + e_4 +
         e_4(e_1 + e_3)(e_1 + e_2)(e_3 + e_2)) =  \\
   & e_2(e_4 + e_3)(e_1(e_2 + e_3) + e_4).
  \end{split}
\end{equation}
By \S\ref{examples_D4} we have $e_{2341} = \tilde{e}_{2341}$.
 The proposition is proved.
\qedsymbol \vspace{2mm}

\begin{proposition}[The elements $\tilde{f}_{{\alpha}0}$]
  \label{coincidence_F}
Consider elements $\tilde{f}_{{\alpha}0}$ for $\alpha = 21, 121,
321$ (see \S\ref{examples_D4}). The following relation
$$
   f_{{\alpha}0} = \tilde{f}_{{\alpha}0}
$$
takes place without restriction$\mod\theta$.
\end{proposition}
\PerfProof
 \underline{For $n=2$}: $\alpha = 21$. We have
\begin{equation}
 \begin{split}
  & \Gamma_{f}(\alpha) = \{ (k_2 k_1) \mid k_2 \in \{3,4 \}, k_1 \in
  \{2, 3, 4\}, k_1 \neq k_2 \}, \\
  & \tilde{f}_{210} = e_2\sum\limits_{\beta \in \Gamma_{f}(\alpha)}{e_\beta} =
    e_2(e_{32} + e_{34} + e_{42} + e_{43}) = \\
  & e_2(e_3(e_1 + e_4) + e_3(e_2 + e_1) +
    e_4(e_1 + e_3) + e_4(e_1 + e_2)) = \\
  & e_2((e_3 + e_4)(e_1 + e_4)(e_1 + e_3) + e_3(e_1 + e_2) + e_4(e_1 + e_2)) = \\
  & e_2(e_3 + e_4)((e_1 + e_4)(e_1 + e_3) + e_3(e_1 + e_2) + e_4(e_1 + e_2)) = \\
  & e_2(e_3 + e_4)(e_1 + e_4(e_1 + e_3) + e_3(e_1 + e_2) + e_4(e_1 + e_2)) = \\
  & e_2(e_3 + e_4)(e_4(e_1 + e_3) + (e_1 + e_3)(e_1 + e_2) + e_4(e_1 + e_2)) = \\
  & e_2(e_3 + e_4)(e_4(e_1 + e_3)(e_1 + e_2) + (e_1 + e_3)(e_1 + e_2) + e_4) = \\
  & e_2(e_3 + e_4)(e_4 + (e_1 + e_3)(e_1 + e_2)) = \\
  & e_2(e_3 + e_4)(e_4 + e_1 + e_3(e_1 + e_2)).
 \end{split}
\end{equation}
By \S\ref{examples_D4} we have
 $f_{210} = \tilde{f}_{210}$.
 \qedsymbol \vspace{2mm}

\underline{For $n=3$}: 1) $\alpha = 121$.
 For this case, we have
\begin{equation}
   \Gamma_{f}(\alpha) = \{ (k_3 k_2 k_1) \mid k_3 \in \{3,4 \}, k_2 \in
  \{3, 4\}, k_1 \in \{2,3,4\}, k_1 \neq k_2, k_2 \neq k_3 \}, \\
\end{equation}
and
\begin{equation}
 \begin{split}
  & \tilde{f}_{1210} = e_1\sum\limits_{\beta \in \Gamma_{f}(\alpha)}{e_\beta} =
    e_1(e_{342} + e_{343} + e_{432} + e_{434}) = \\
    e_1[& e_3(e_1 + e_2)(e_4 + e_2) + e_3(e_1(e_4 + e_3) + e_2(e_4 + e_3))+ \\
       & e_4(e_1 + e_2)(e_3 + e_2) + e_4(e_1(e_4 + e_3) + e_2(e_4 + e_3))]= \\
    e_1[& e_3(e_1 + e_2)(e_4 + e_2 + e_3(e_1 + e_2(e_4 + e_3))) + \\
        & e_4(e_1 + e_2)(e_3 + e_2 + e_4(e_1 + e_2(e_4 + e_3)))] = \\
    e_1[& e_3(e_1 + e_2)(e_4 + e_2 + e_1(e_3 + e_2(e_4 + e_3))) + \\
        & e_4(e_1 + e_2)(e_3 + e_2 + e_1(e_4 + e_2(e_4 + e_3)))] = \\
    e_1[& (e_4 + e_2 + e_1(e_3 + e_2)(e_4 + e_3))]  \\
       [& (e_3 + e_2 + e_1(e_4 + e_2)(e_4 + e_3))]  \\
       [& e_3(e_1 + e_2) + e_4(e_1 + e_2)] = \\
    e_1[& (e_4 + e_2)(e_4 + e_3) + e_1(e_3 + e_2))]  \\
       [& (e_3 + e_2 + e_1(e_4 + e_2)(e_4 + e_3))]  \\
       [& e_3(e_1 + e_2) + e_4(e_1 + e_2)] = \\
    e_1[& (e_4 + e_2)(e_4 + e_3) + e_1(e_3 + e_2))]  \\
       [& e_1(e_3 + e_2) + (e_4 + e_2)(e_4 + e_3))]  \\
       [& e_3(e_1 + e_2) + e_4(e_1 + e_2)].
 \end{split}
\end{equation}
Since the first two intersection polynomials in the last
expression of $\tilde{f}_{1210}$ coincide with
$$
   (e_4 + e_2)(e_4 + e_3) + e_1(e_3 + e_2)),
$$
we have
\begin{equation}
 \begin{split}
   \tilde{f}_{1210} = & \\
  & e_1((e_4 + e_2)(e_4 + e_3) + e_1(e_3 + e_2)))(e_3 + e_4(e_1 + e_2)).
 \end{split}
\end{equation}
By \S\ref{examples_D4} we have
 $f_{1210} = \tilde{f}_{1210}$.
 \qedsymbol \vspace{2mm}

2) $\alpha = 321$. Here we have
\begin{equation}
 \begin{split}
  & \Gamma_{f}(\alpha) = \{ (k_3 k_2 k_1) \mid k_3 \in \{1,4 \}, k_2 \in
    \{3, 4\}, k_1 \in \{2,3,4\}, k_1 \neq k_2, k_2 \neq k_3 \} =  \\
  & \{ (132) = (142), (134), (143, (432), (434) \}
 \end{split}
\end{equation}
and
\begin{equation}
 \label{f_3210_transf_1}
 \begin{split}
   \tilde{f}_{3210} = & e_1\sum\limits_{\beta \in \Gamma_{f}(\alpha)}{e_\beta} =
    e_3(e_{132} + e_{134} + e_{143} + e_{432} + e_{434}) = \\
    e_3[& e_1(e_4 + e_2)(e_3 + e_2) + e_1(e_4 + e_2)(e_3 + e_4) + \\
        & e_1(e_3 + e_2)(e_3 + e_4) + e_4(e_3 + e_2)(e_1 + e_2) + \\
        & e_4(e_1 + e_2(e_3 + e_4))].
 \end{split}
\end{equation}
Since
\begin{equation}
 \begin{split}
  & e_1(e_4 + e_2)(e_3 + e_4) + e_4(e_1 + e_2(e_3 + e_4)) = \\
  & (e_3 + e_4)(e_4 + e_2)(e_1 + e_4(e_1 + e_2(e_3 + e_4)) = \\
  & (e_3 + e_4)(e_4 + e_2)(e_1 + e_4)(e_1 + e_2(e_3 + e_4)),
 \end{split}
\end{equation}
by (\ref{f_3210_transf_1}) we have
\begin{equation}
 \label{f_3210_transf_2}
 \begin{split}
   \tilde{f}_{3210} = & \\
    e_3[& e_1(e_4 + e_2)(e_3 + e_2) + \\
        & (e_3 + e_4)(e_4 + e_2)(e_1 + e_4)(e_1 + e_2(e_3 + e_4)) + \\
        &  e_1(e_3 + e_2)(e_3 + e_4) +
           e_4(e_3 + e_2)(e_1 + e_2)] = \\
    e_3(&e_1 + e_2)(e_1 + e_4)[e_1(e_4 + e_2)(e_3 + e_2) + \\
        & (e_3 + e_4)(e_4 + e_2)(e_1 + e_2(e_3 + e_4)) + \\
        & e_1(e_3 + e_2)(e_3 + e_4) +
          e_4(e_3 + e_2)] = \\
    e_3(&e_1 + e_2)(e_1 + e_4)[e_1(e_4 + e_2)(e_3 + e_2)(e_3 + e_4) + \\
        & (e_3 + e_4)(e_4 + e_2)(e_1 + e_2(e_3 + e_4)) + \\
        & e_1(e_3 + e_2) +
          e_4(e_3 + e_2)].
 \end{split}
\end{equation}
Since
$$
  e_1(e_4 + e_2)(e_3 + e_2)(e_3 + e_4) \subseteq e_1(e_3 + e_2),
$$
by (\ref{f_3210_transf_2}) we have
\begin{equation}
 \label{f_3210_transf_3}
 \begin{split}
   \tilde{f}_{3210} = & \\
    e_3(&e_1 + e_2)(e_1 + e_4)
        [(e_3 + e_4)(e_4 + e_2)(e_1 + e_2(e_3 + e_4)) + \\
        & e_1(e_3 + e_2) + e_4(e_3 + e_2)] = \\
    e_3(&e_1 + e_2)(e_1 + e_4)
        [e_1(e_3 + e_4)(e_4 + e_2) + e_2(e_3 + e_4)) + \\
        & e_1(e_3 + e_2) + e_4(e_3 + e_2)] = \\
    e_3(&e_1 + e_2)(e_1 + e_4)
        [e_1(e_3 + e_4)(e_4 + e_2)(e_3 + e_2) + e_2(e_3 + e_4)) + \\
        & e_1 + e_4(e_3 + e_2)].
 \end{split}
\end{equation}
Since
$$
   e_1(e_3 + e_4)(e_4 + e_2)(e_3 + e_2) \subseteq e_1,
$$
by (\ref{f_3210_transf_3}) we have
\begin{equation}
 \label{f_3210_transf_4}
 \begin{split}
   \tilde{f}_{3210} = & \\
    e_3(&e_1 + e_2)(e_1 + e_4)
        (e_2(e_3 + e_4)) + e_1 + e_4(e_3 + e_2)) = \\
    e_3(&e_1 + e_2)(e_1 + e_4)
        (e_1 + (e_2 + e_4)(e_3 + e_2)(e_3 + e_4)).
 \end{split}
\end{equation}

By \S\ref{examples_D4} we have
 $f_{3210} = \tilde{f}_{3210}$.
 The proposition is proved.
\qedsymbol \vspace{2mm}

\begin{conjecture}
 \label{conj_4}{\rm
  For every admissible sequence $\alpha$,
  the elements $e_\alpha$ (resp. $f_{{\alpha}0}$) and
  $\tilde{e}_\alpha$ (resp. $\tilde{f}_{{\alpha}0}$) coincide
  without restriction$\mod\theta$
  (see Proposition \ref{coincidence_GP})}.
\end{conjecture}

In Propositions \ref{coincidence_E} and \ref{coincidence_F}, this
conjecture was proven for small values of lengths of the
admissible sequence $\alpha$.

\begin{appendix}
\chapter{\sc\bf Modular and linear lattices}
 \label{on_lattices}

\epigraph{The notion of modular lattice misled logicians and
mathematicians for decades. \dots most examples of modular
lattices occurring in algebra and combinatorics enjoy the stronger
property of being linear lattices. Unlike modular lattices, it is
not known whether linear lattices can not be defined by identities
alone \dots This lack of an abstract definition is perhaps the
reason why in the past the theory of linear lattices was subsumed
into the theory of modular lattices.}{Mainetti~M., Yan~C.~H.,
 \cite[p.12]{MY99}.}

\section{Distributive lattices and Boolean algebras}
  \label{lattice_1}
  \index{lattice}
  \index{sum}
  \index{intersection}
  \index{inclusion}
  \index{idempotent operation}
  \index{absorbtion law}

  A {\it lattice} is a set $L$ with two commutative and
  associative operations: {\it a sum} and {\it an intersection}.
  If $a, b \in L$, then we denote the intersection
  by $ab$ and the sum by $a + b$. More frequently
  the sum of $a$ and $b$  is denoted by
  $a\cup{b}$ or $a\vee{b}$ and
  the intersection is denoted by $a\cap{b}$ or
  $a\wedge{b}$.
  Both operations are {\it idempotent} and satisfy the
  {\it absorbtion law}:
\begin{equation*}
\begin{array}{lll}
    & a{a} = a,  & \text{ idempotency of intersection}, \\
    & a + a = a, & \text{ idempotency of sum}, \\
    & a(a + b) = a, & \text{ absorbtion law}, \\
    & a + a{b} = a, & \text{ absorbtion law}. \\
\end{array}
\end{equation*}
The sum and intersection of $n$ elements $a_1$,\dots,$a_n$ is
denoted, respectively, by
\begin{equation*}
   \sum\limits_{i=1}^n{a}_i  \text{ and }
   \bigcap\limits_{i=1}^n{a}_i
\end{equation*}
On the lattice $L$, an order relation $\subseteq$ ({\it an
inclusion}) is defined:
$$
    a \subseteq b \Longleftrightarrow a{b} = a
    \text{ or } a + b = b.
$$

 \index{distributive lattice}
 \index{lattice!- distributive}
 \index{distributive laws}

A lattice $L$ is said to be {\it distributive} if the following
two {\it distributive laws} hold
\begin{equation}
 \label{distributive_laws}
 \begin{array}{lll}
    & a(b + c) = ab + ac
       & \text{ for every } a, b, c \subseteq L, \\
    & a + bc = (a + b)(a + c)
       & \text{ for every } a, b, c \subseteq L.\\
 \end{array}
\end{equation}

 \index{Boolean algebra}
 \index{atom in Boolean algebra}
 \index{coatom in Boolean algebra}
 \index{lattice!- Boolean algebra}

A distributive lattice $L$ is said to be a {\it Boolean algebra}
if $L$ contains {\it a minimal} element $O$ and {\it a maximal}
element $I$ and, for each element $a \in L$, there exists a unique
{\it complement} ${a}'$ such that
\begin{equation}
\begin{split}
   & a{a}' = O,  \\
   & a + a' = I,   \\
   & (a')'  = a,  \\
   & (ab)' = a' + b', \\
   & (a + b)' = a'{b}'. \\
\end{split}
\end{equation}
An {\it atom} $a$ of a Boolean algebra is a minimal nonzero
element, i.e., if $x \subseteq a$ then, either $x = O$ or
 $x = a$. An element $y$ is {\it atomless} if there is no atom
 $x$ such that $x \subseteq y$. A Boolean algebra is {\it atomic} if there are no
atomless elements.

\begin{fact}
 Some infinite Boolean algebras do not contain
any atoms. All finite Boolean algebras are atomic,
\cite[p.135]{BS81}, \cite{Hal63}.
\end{fact}

A {\it coatom} $a$ of a Boolean algebra is a maximal nonzero
element, i.e., if $x \supseteq a$, then either $x = I$ or $x = a$.

\begin{fact}[Stone, \cite{Stn36}]
Every finite Boolean algebra is isomorphic to the Boolean algebra
of all subsets of a finite set; in particular, the number of
elements of every finite Boolean algebra is a power of two,
\cite[Ch.2]{Gr98}.
\end{fact}

\section{Modular lattices}
  \label{modular_1}

 \index{modular lattice}
 \index{lattice!- modular}
 \index{modular law(=Dedekind's law)}

  \label{modularity}
 A lattice is said to be {\it modular} if,
 for every $b, a, c \in L$,
\begin{equation}
  \label{modular_law}
       a \subseteq b \Longrightarrow b(a + c) = a + b{c}.
\end{equation}

The relation (\ref{modularity}) is said to be {\it the modular
law}. The {\it modular law} is sometimes referred to as {\it
Dedekind's law}. Every distributive lattice is modular but not
conversely.

A set of subspaces of the given vector space $R$ (and, of course,
submodules of the given module $M$) form a modular lattice, namely
if $A, B, C \subseteq R$ and $A \subseteq B$, then
$$
   B(A + C) = A + B{C}.
$$
Indeed, the inclusion $A + BC \subseteq  B(A + C)$ holds, because
 $A + BC \subseteq A + C$ and $A + BC \subseteq B$. Conversely, let
 $x \in B(A + C)$, i.e., $x \in B$ and $x = x_a + x_c$,
where vectors $x_a$ and $x_c$ are components of the vector
 $x$, such that $x_a \in A$ and $x_c \in C$.
Then $x_c = x - x_a \in B$, i.e., $x_c \in BC$ and
 $x_a + x_c \in A + BC$.

The main properties of the modular lattices were derived by
R.~Dedekind \cite{De1897}, \cite{Bir48}. We need the two features
of the modular lattice, which are called {\it permutation
properties of the modular lattice} or just {\it permutation
properties}.
 \index{permutation properties of modular lattice}

\begin{proposition} [Permutation properties]
   In any modular lattice $L$,
   for every four elements $A,B,C,D \subseteq L$,
   the following two permutation properties hold
\begin{equation}
   \label{permutation1}
      A(B+CD) = A(CB+D) = A(CB+CD)
         \text{ for every } A \subseteq C.
\end{equation}
\begin{equation}
  \label{permutation2}
   A+D(B+C) = A+B(D+C) = A+(D+C)(B+C)
      \text{ for every } C \subseteq A.
\end{equation}
\end{proposition}
\PerfProof By the modular law (\ref{modularity}) we have
\begin{equation*}
\begin{split}
   & A(B+CD) = AC(B+CD) = A(CB+CD) \text{ and } \\
   & A(CB+D) = AC(CB+D) = A(CB+CD),
\end{split}
\end{equation*}
which proves (\ref{permutation1}). Similarly,
\begin{equation*}
\begin{split}
   & A+D(B+C) = A+C+D(B+C) = A+(D+C)(B+C) \text{ and } \\
   & A+B(D+C) = A+C+B(D+C) = A+(B+C)(D+C),
\end{split}
\end{equation*}
which proves (\ref{permutation2}). \qedsymbol

\section{Hasse diagrams}
 \index{Hasse diagram}
 \index{upward drawing}

\begin{figure}[h]
\includegraphics{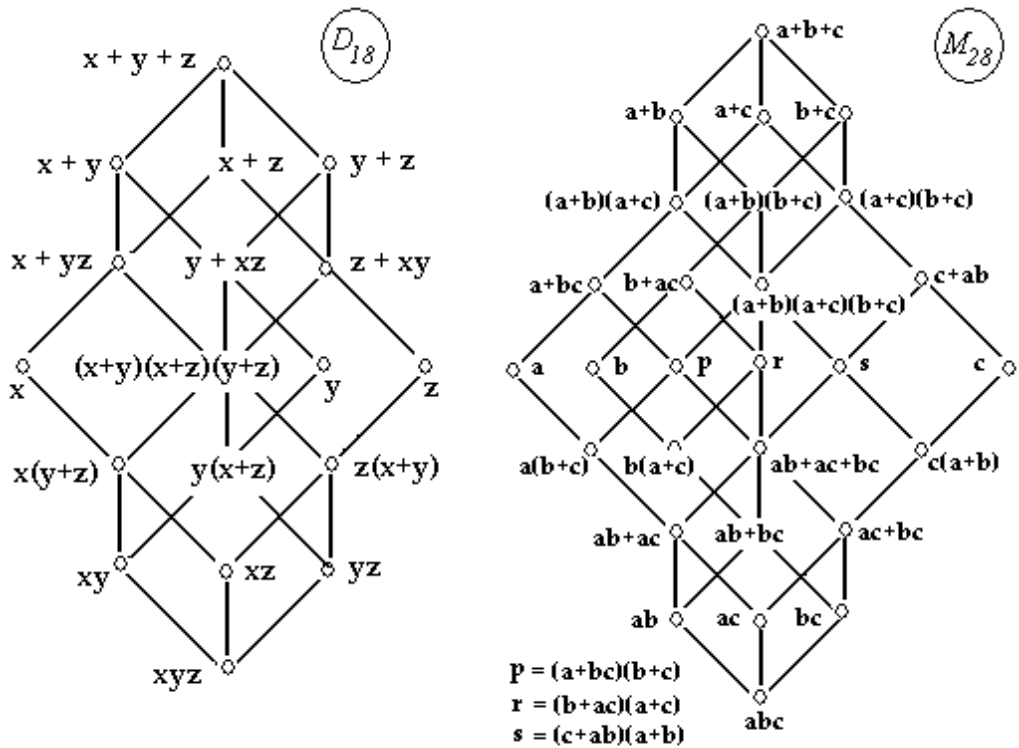}
\caption[\hspace{3mm}Hasse diagrams of $D_{18}$ and $M_{28}$]
{\hspace{3mm}Hasse diagrams of two lattices: the free distributive
3-generated lattice $D_{18}$ and the free modular 3-generated
lattice $M_{28}$}
\label{D18_M28}
\end{figure}

Hasse diagrams (also called {\it upward drawings}) are constructed
to clear up the lattice by eliminating obvious edges. The vertices
are arranged vertically, so the directional arrows are implied
(always upwards) and omitted. Any edge implied by transitivity is
not shown, i.e., if there is some element $x$ such that $a
\subseteq x \subseteq b$, then the edge $[a, b]$ is not shown on
the Hasse diagram.
 On Fig. \ref{D18_M28}, two Hasse diagrams are depicted: free
distributive lattice $D_{18}$ generated by three generators and
free modular lattice $M_{28}$ generated by three generators.
Lattices $D_{18}$ and  $M_{28}$ are finite, $D_{18}$ contains $18$
elements and $M_{28}$ contains $28$ elements.

 \index{modular lattice!- $M_{28}$}
 \index{distributive lattice!- $D_{18}$}

\begin{figure}[h]
\includegraphics{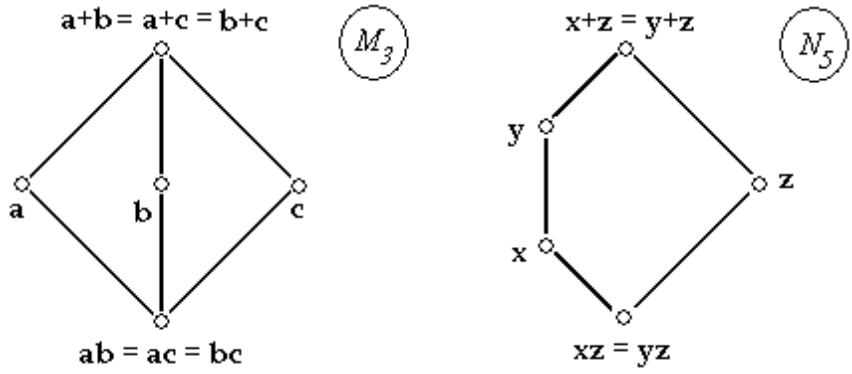}
\caption{\hspace{3mm}Two Hasse diagrams: the diamond $M_3$ and
 the pentagon $N_5$}
\label{M3_N5}
\end{figure}

\begin{figure}[th]
\includegraphics{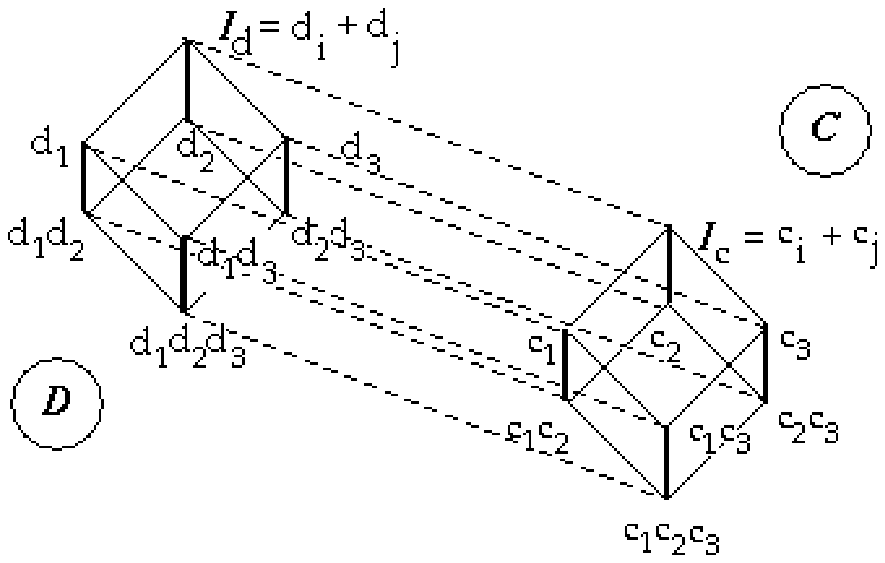}
\caption{\hspace{3mm}The $4-$dimensional cube $D$ $\bigcup$ $C$}
\label{16elem_C_D}
\end{figure}

 \index{diamond $M_3$}
 \index{lattice!- diamond $M_3$}
 \index{pentagon $N_5$}
 \index{lattice!- pentagon $N_5$}
For other examples of Hasse diagrams, see

 a) The $64$-element distributive lattice $H^+(n)$, Fig. \ref{cubic64},

 b) Perfect sublattices in $D^4$ and $D^{2,2,2}$, Fig. \ref{cube_comparison},

 c) The $16$-element Boolean algebra $U_n \cup{V}_{n+1}$, Fig. \ref{boolean16},

 d) The Boolean algebras $U_n$ and $V_{n+1}$, Fig. \ref{8elem_U_V},

 e) The diamond $M_3$ and the pentagon $N_5$, Fig. \ref{M3_N5}.

 f) $4-$dimensional cube $D \cup C$, Fig. \ref{16elem_C_D}.

\section{A characterization of distributive and modular lattices}
  \label{sect_char_distr_modul}
Every distributive lattice $L$ is also modular. Indeed, if
 $a \subseteq b$, then
$$
  b(a + c) =  ba + bc = a + bc,
$$
i.e., (\ref{modular_law}) holds. The converse is false. Consider
the diamond lattice $M_3$, see Fig. \ref{M3_N5}. It easily follows
from the modular law (\ref{modular_law}) that the lattice $M_3$ is
modular. However, $M_3$ is not distributive, because
\begin{equation}
 \label{diamond}
 c(a + b) = c, \text{ but } ca + cb = ab,
 \text{ i.e., } c(a + b) \neq ca + cb,
\end{equation}
see Fig. \ref{M3_N5}. Further, the pentagon $N_5$ is not a modular
lattice, because
\begin{equation}
  \label{pentagon}
  x \subseteq y, \indent  y(x + z) = y, \indent
  x + yz = x, \text{ i.e., }  y(x + z) \neq x + yz.
\end{equation}

 The following characterization of modular lattices is well-known
 \cite{Bir48}.
\begin{proposition}
 \index{pentagon $N_5$}
 \label{ch_modular}
   A lattice $L$ is modular if and only if
   it does not contain pentagon $N_5$ as a sublattice.
\end{proposition}
\PerfProof According to (\ref{pentagon}) $L$ does not contain
$N_5$, since $N_5$ is not modular. \\
\begin{minipage}{6cm}
  \psfig{file=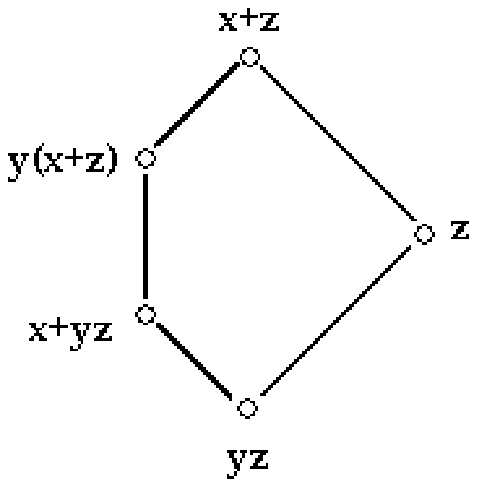, width=2in}
\end{minipage}
\begin{minipage}{10cm}
   We will show that, conversely, if $L$ is not modular,
   then it must contain a pentagon. Let us suppose that there are
   three elements $x, y, z$ such that
\begin{equation}
\label{contra_N5}
\begin{split}
  & x \subseteq y, \\
  & x + yz \neq  y(x + z), \\
  & x + yz \subseteq  y(x + z).
\end{split}
\end{equation}
The last inclusion in (\ref{contra_N5}) is true in every lattice.
\end{minipage}

Consider the sublattice generated by 5 the following elements:
$$
    yz \subseteq x + yz  \subseteq y(x + z) \subseteq x + z
      \text{ and } z.
$$
We will show that this sublattice is isomorphic to $N_5$. For
sums, we have
$$
 x + z \subseteq (x + yz) + z  \subseteq y(x + z) + z
   \subseteq x + z,
$$
and therefore
$$
    y(x + z) + z = x + z, \indent (x + yz) + z = x + z.
$$
For intersections, we have
$$
  zy \subseteq z(x + yz) \subseteq zy
$$
since $x \subseteq y$. Therefore $z(x + yz) = zy$.
 Besides, $z(y(x + z)) = zy$. Finally, let us show that
$$
   yz \neq x + yz.
$$
Suppose that $yz = x + yz$; then $x \subseteq y$ and
$$
   yz = x + yz \subseteq y(x + z) \subseteq y(yz + z) = yz,
$$
i.e., $x + yz = y(x + z)$ which contradicts the hypothesis
(\ref{contra_N5}). \qedsymbol \vspace{2mm}

\begin{proposition}[\cite{Bir48}]
 \label{ch_distributive}
 \index{diamond $M_3$}
   A lattice $L$ is distributive if and only if it does not
contain sublattices $M_3$ and $N_5$.
\end{proposition}
\PerfProof  According to (\ref{pentagon}), (\ref{diamond}) any
distributive lattice $L$ does not contain $M_3$ and $N_5$, since
$M_3$ and $N_5$ are not distributive. Suppose, conversely, that
$L$ does not contain $N_5$. Then by Proposition \ref{ch_modular}
the lattice $L$ is modular. Let $L$ be not distributive.

\begin{minipage}{8cm}
  \psfig{file=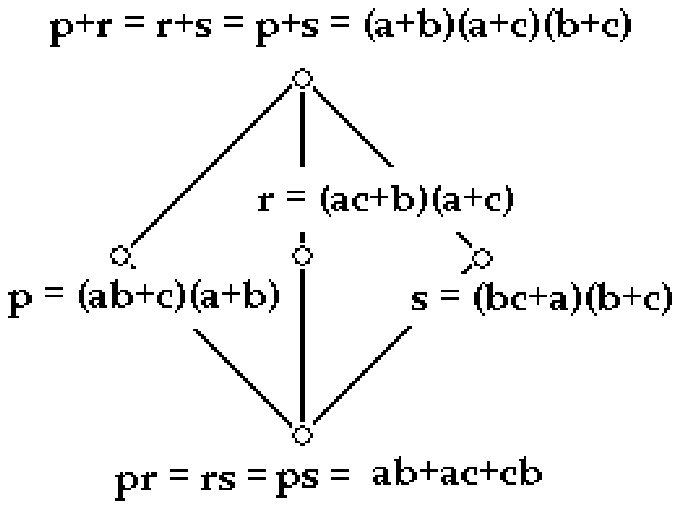, width=3in}
\end{minipage}
\begin{minipage}{7.5cm}
   We will show that, if $L$ is not distributive,
   then it must contain a diamond. Suppose that there are
   three elements $a, b, c$ that sublattice
   $\{a, b, c \}$ is not distributive.
   We consider the sublattice generated by the following elements $p, r, s$:
\begin{equation}
\label{contra_M3}
\begin{split}
  & p = (ab + c)(a + b), \\
  & r = (ac + b)(a + c), \\
  & s = (bc + a)(b + c). \\
\end{split}
\end{equation}
\end{minipage}

Since the modular law holds in $L$, we have
\begin{equation*}
\begin{split}
    p + r = &
     (ab + c)(a + b) + (ac + b)(a + c) = (a + b)(a + c)
           (ab + c + ac + b) = \\
             & (a + b)(a + c)(c + b). \\
\end{split}
\end{equation*}
Similarly,
\begin{equation}
   p + s = r + s = p + r = (a + b)(a + c)(b + c).
\end{equation}
In just the same way,
\begin{equation*}
\begin{split}
    pr = &
     (ab + c)(a + b)(ac + b)(a + c) = (ab + c)(ac + b) = \\
   & ac + b(ab + c) = ac + ab + bc, \\
\end{split}
\end{equation*}
and
\begin{equation}
   ps = pr  = rs =
   ac + ab + bc.
\end{equation}
We will show now that the elements $p, r, s$ are distinct.
Suppose, for example, that $p = r$. Then
$$
   p + r = p = r = pr
$$
and
\begin{equation}
  \label{top_bottom}
   (a + b)(a + c)(b + c) =  ac + ab + bc
\end{equation}
Consider now the intersection of both sides of the
(\ref{top_bottom}) with $a$:
$$
    a(b + c) = a(ac + ab + bc) = ac + ab + abc = ac + ab.
$$
Similarly,
$$
   b(a + c) = ba + bc, \indent
   c(a + b) = ca + cb,
$$
which contradicts the hypothesis that $\{a, b, c\}$ is not
distributive. \qedsymbol \vspace{2mm}

 The following proposition shows
the equivalence of the two distributivity laws.
\begin{proposition}
\label{distrib_equiv} Two following distributivity laws
 (see (\ref{distributive_laws})) are equivalent for an arbitrary lattice
$L$:
\begin{equation*}
\begin{array}{lll}
    \text{ (i)  } &  a(b + c) = ab + ac & \text{ for every } a, b, c \in L, \\
    \text{ (ii) } & a + bc = (a + b)(a + c) & \text{ for every } a, b, c \in L
\end{array}
\end{equation*}
\end{proposition}
\PerfProof Let us prove the implication (i) $\Longrightarrow$ (ii):
$$
   (a + b)(a + c) = (a + b)a + (a + b)c =
    a + (a + b)c = a + ac + bc = a + bc.
$$
Conversely, (ii) $\Longrightarrow$ (i):
$$
    ab + ac = (ab + c)(ab + a) =
    (ab + c)a = (a + c)(b + c)a = a(b + c).
    \qed \vspace{2mm}
$$

 The following criterion of distributivity of the
modular lattice is proved by B.~Jonnson:
\begin{proposition}{\em (\cite[Th.5]{Jo55})}
 \label{jonsson}
  Let $A$ be a modular lattice, $p$ a positive integer, and
  $X_1, X_2,\dots,X_p$ nonempty
  linearly ordered subsets of $A$.  For the sublattice
  of $A$ generated by the set
$$
  X_1\cup{X}_1\cup\dots\cup{X}_p
$$
to be distributive it is necessary and sufficient that, for any
$$
  x_1 \in X_1,x_2 \in X_2,\dots,x_p \in X_p,
$$
the sublattice of $A$ generated by the set $\{x_1,x_2,...,x_p\}$
be distributive. \qedsymbol \vspace{2mm}
\end{proposition}

We use this result for $p = 3$ in the proof of the distributivity
of the sublattice of perfect elements $H^+(n)$ in Proposition
\ref{distrib_Hn}.

\section{When the union of two 3D cubes is a 4D cube?}
 \index{$3D$ cube}
 \index{$4-$dimensional cube($=4D$ cube)}

\begin{proposition}
  \label{DC_lemma}
    Let $D$ and $C$ be two $8$-element Boolean algebras with
coatoms $d_1, d_2, d_3 \in D$, coatoms $c_1, c_2, c_3 \in C$,
maximal elements $I_d \in D$, $I_c \in C$ and minimal elements
$O_d \in D$, $O_c \in C$. If the following inclusions hold:
\begin{equation} \label {c_less_d}
      c_i \subseteq d_i, \indent i = 1,2,3,
\end{equation}
\begin{equation} \label {di_less_ci_dj}
      d_i \subseteq c_i + d_j,  \indent
            i,j \in \{1,2,3\},
\end{equation}
\begin{equation} \label {dicj_less_ci}
      d_i{c}_j \subseteq c_i, \indent i \neq j,
            i,j \in \{1,2,3\},
\end{equation}
then the union $D \cup C$ is a $16$-element Boolean algebra with
coatoms $d_1, d_2, d_3$ and $I_c$.
\end{proposition}
   \PerfProof Let us prove that the \underline{sum does not lead out of $D$}.
  \begin{equation}
     c \in C, \indent d \in D \Longrightarrow c + d \in D.
  \end{equation}
\underline{Case $c = c_i$, 1) -- 4)}

1) By (\ref{di_less_ci_dj}) $d_i + d_j \subseteq
   c_i + d_j = c_i + d_i{d}_j + d_j =
   d_i(c_i + d_j) + d_j \subseteq d_i + d_j$, i.e.,
\begin{equation} \label {di_dj}
    c_i + d_j = d_i + d_j = I_d.
\end{equation}

2) From (\ref{di_dj}) we have
\begin{equation} \label {ci_didj}
    c_i + d_i{d}_j = d_i(c_i + d_j) = d_i.
\end{equation}

3) Since $c_j{c}_k \subseteq d_j{d}_k$, we get
\begin{equation*}
  \begin{split}
   & c_i + d_j{d}_k = c_i{c}_j + c_i + c_j{c}_k + d_j{d}_k = \\
   & c_j(c_i + c_k) + c_i + d_j{d}_k =
      c_j + c_i + d_j{d}_k = c_i + (c_j + d_j{d}_k).
  \end{split}
\end{equation*}
 According to (\ref{ci_didj}) we have
\begin{equation}
 \label{ci_djdk}
    c_i + d_j{d}_k = c_i + d_j = I_d.
\end{equation}

 \index{$4-$dimensional cube}

4) From  (\ref{ci_djdk}) for $O_d = d_i{d}_j{d}_k$ we have
\begin{equation} \label {ci_didjdk}
    c_i + d_i{d}_j{d}_k = d_i(c_i + d_j{d}_k) = d_i.
\end{equation}

\underline{Case $c = c_i{c}_j$, 5) -- 8)}

5) Since $c_j \subseteq d_j$, we get
\begin{equation}
  \label {cicj_dj}
    c_i{c}_j + d_j = d_j.
\end{equation}

6) From (\ref{cicj_dj}) we obtain
\begin{equation}
  \label {cicj_dk}
    c_i{c}_j + d_k =
    c_i{c}_j + c_i{c}_k + d_k =
    c_i + d_k = I_d.
\end{equation}

7) Similarly,
\begin{equation} \label {cicj_didk}
    c_i{c}_j + d_i{d}_k =
    d_i(c_i{c}_j + d_k) = d_i.
\end{equation}

8) For $O_d = d_i{d}_j{d}_k$, we have
\begin{equation} \label {cicj_didjdk}
    c_i{c}_j + d_i{d}_j{d}_k = d_i{d}_j(c_i{c}_j + d_k) =
    d_i{d}_j
\end{equation}

\underline{Case $c = c_i{c}_j{c}_k = O_c$}.

 In this case we have
 $c \subseteq d$ \text{ for every } $d \in D$ and
\begin{equation}
   \label {cicjck_d}
    c_i{c}_j{c}_k + d = d
    \text{ for every } d \in D.
\end{equation}

\underline{Case $c = I_c$}.

In this case $c$ is the sum of two coatoms $c_i + c_j$, here we
have also
\begin{equation}
  \label {ic_d}
    I_c + d \in D \text{ for every } d \in D.
\end{equation}

 We will prove that \underline{the intersection does not lead out
 of $C$}.
\begin{equation}
    c \in C, \indent d \in D
   \Longrightarrow dc \in C.
\end{equation}

\underline{Cases $c = c_i, c = I_c, c = c_i{c}_j$, 9) -- 11)} .

9)  From (\ref{dicj_less_ci}) we have
 $d_i{c}_j \subseteq c_i{c}_j$. From (\ref{c_less_d}) we have
 $c_i{c}_j \subseteq d_i{c}_j$, so
\begin{equation} \label {dicj_cicj}
   d_i{c}_j = c_i{c}_j.
\end{equation}

10) From (\ref{dicj_cicj}) we deduce that
\begin{equation}
  \label{diIc_ci}
   d_i{I}_c = d_i(c_i + c_j) =
   c_i + d_i{c}_j = c_i + c_i{c}_j = c_i.
\end{equation}

11) Similarly,
\begin{equation}
  \label {dicjck_cicjck}
   d_i{c}_j{c}_k  = (d_i{c}_j)(d_i{c}_k) =
   c_i{c}_j{c}_k = O_c.
\end{equation}

12) Finally,
\begin{equation}
  \label {didjck_cicjck}
   d_i{d}_j{c}_k  =  (d_i{c}_k)(d_j{c}_k) =
   c_i{c}_j{c}_k = O_c.
\end{equation}
For all other $d = d_i{d}_j$ or for $d = d_i{d}_j{d}_k$ we also
have $d{c} \in C$ \text{ for every } $c \in C$. \vspace{2mm}

\underline{Complementarity in $D \cup C$}.

 From (\ref{ci_djdk}) and (\ref{didjck_cicjck}) we obtain
\begin{equation}
 \label {complement_1}
   {c}'_i = d_j{d}_k.
\end{equation}
From (\ref{cicj_dk}) and (\ref{dicjck_cicjck}) we have
\begin{equation}
  \label {complement_2}
   {d}'_i = c_j{c}_k.
 \qed \vspace{2mm}
\end{equation}

\section{Representations of the modular lattices}
 \label{repres_lat}

In this section we mostly follow definitions of Gelfand and
 Ponomarev \cite{GP74}.

 \index{indecomposable representation $\rho$}
 \index{representation $\rho$ of a modular lattice}
 \index{linear representation}
 \index{representation!- indecomposable}
 \index{representation!- of modular lattice}

Let $L$ be a modular lattice, $X$ a finite dimensional vector
space, $\mathcal{L}(X)$ the modular lattice of linear subspaces in
$L$. A morphism  $\rho: L \longrightarrow \mathcal{L}(X)$ is
called a {\it linear representation} of $L$ in the space $X$.
Thus, the representation $\rho: L \longrightarrow \mathcal{L}(X)$
maps every element $a \in L$ to the subspace
 $\rho(a) \subseteq X$ such that, for every $a, b \in L$, we have:
$$
   \rho(ab) = \rho(a)\rho(b),  \indent
   \rho(a + b) = \rho(a) + \rho(b).
$$

Let $\rho_1: \longrightarrow \mathcal{L}(U)$ and $\rho_2:
\longrightarrow \mathcal{L}(V)$ be representations of $L$ in the
spaces $U$ and $V$. We put
$$
   \rho(a) = \rho_1(a) \oplus \rho_2(a)
$$
for every $a \in L$, where  the subspace
 $\rho(a) \subseteq U \oplus V$ is the set of all pairs
$$
 \{ (u, v) \mid u \in \rho_1(a) \text{ and } v \in \rho_2(a)\}.
$$
This correspondence gives the representation $\rho$ in the space
$U \oplus V$. The representation $\rho$ is called the {\it direct
sum of representations} $\rho_1 \oplus \rho_2$.

 \index{representation!- decomposable}
The representation $\rho$ is said to be {\it decomposable} if it
isomorphic to a direct sum $\rho_1 \oplus \rho_2$ of non-zero
representations $\rho_1$ and $\rho_2$. It is easy to prove that
the representation $\rho$ in the space $X$ is decomposable if and
only if there are subspaces $V$ and $U$ such that
 $X = U \oplus V$ and
$$
    \rho(a) = U\rho(a) + V\rho(a)
    \indent \text{ for every}\footnote{Recall that $U\rho(a)$
     denotes the intersection of
     subspaces: $U\rho(a) = U \cap \rho(a)$.}
     \indent a \in L.
$$

 \index{representation!- indecomposable}

The representation $\rho$ is said to be {\it indecomposable} if
there exists an element $a \in L$ such that
$$
   \rho(a) \neq U\rho(a) + V\rho(a)
$$
for every decomposition $X$ into the direct sum $X = U \oplus V$.

 \index{perfect element}
Following \cite{GP74} we introduce now the notions of {\it perfect
elements} and {\it linear equivalence relation}.
 The element $a \in L$ is said to be {\it perfect} if
either $\rho(a) = 0$ or $\rho(a) = X$ for each indecomposable
representation $\rho$ of the modular lattice $L$, where $X$ is the
representation space of $\rho$.

Obviously, the sum and intersection of perfect elements is also a
 perfect element. Thus, perfect elements form a sublattice in the
lattice $L$.

  \index{modulo linear equivalence$\mod\theta$}
  \index{linear equivalence relation $\theta$}
  \index{$\mod\theta$, modulo linear equivalence}
  \index{$\theta$, linear equivalence relation}
 Two elements $a, b \in L$ are called {\it linearly
equivalent} if $\rho(a) = \rho(b)$ for all indecomposable
representations $\rho: L \longrightarrow \mathcal{L}(X)$ and we
write
\begin{equation}
  \label{lin_equiv}
   a \equiv b\mod\theta.
\end{equation}
The relation $\theta$ is called the {\it linear equivalence
relation}.

We say that a property of $L$ is true {\it modulo linear
equivalence} if the property is true for each indecomposable
representation $\rho$ of $L$.

\index{subspace admissible with respect to representation}
\index{admissible subspace with respect to representation}

Let $\rho$ be the representation of the modular lattice $L$ in the
vector space $X$. The subspace $U \subseteq X$ is said to be
 {\it admissible with respect to } $\rho$,
if one of following two relations hold
\begin{equation}
\label{distrib_eq_2}
\begin{split}
    & \text{ (i) }  U(\rho(a) + \rho(b)) = U\rho(a) + U\rho(b), \\
    & \text{ (ii) } U + \rho(a)\rho(b) = (U + \rho(a))(U + \rho(b))
\end{split}
\end{equation}
for every $a, b \in L$. The equivalence of (i) and (ii) follows
from Proposition \ref{distrib_equiv}.

 \index{representation!- in the subspace}
 \index{representation!- in the quotient space}
\begin{proposition}
  \label{admis_space}
  {\em (\cite[Prop.2.1]{GP74})}
  Let $\rho$ be a representation of the lattice $L$ in a
  vector space $X$, let $U$ be the subspace of $X$, and
  $\nabla : X \longrightarrow X/U$ the canonical map.
  Then the following relations are equivalent:
  \begin{enumerate}
   \item The subspace $U$ is admissible with respect
        to $\rho$. \vspace{2mm}
   \item The map $x \longmapsto U\rho(x)$
        defines a representation in the subspace $U$. \vspace{2mm}
   \item The mapping $x \longmapsto \nabla\rho(x)$
        defines a representation in the quotient space $X/U$. \vspace{2mm}
  \end{enumerate}
\end{proposition}
 \PerfProof Let us prove the implication $(1) \Longrightarrow (2)$.
 If $U$ is admissible with respect to $\rho$, then
 $$
     U\rho(a + b) = U(\rho(a) + \rho(b)) = U\rho(a) + U\rho(b).
 $$
Besides,
$$
     U\rho(ab) = U\rho(a)\rho(b) = (U\rho(a))(U\rho(b)).
$$
Thus, we have a representation in the subspace $U$. The
representation in the subspace $U$ is called a
 {\it subrepresentation} of $\rho$ and is denoted by $\rho_U$.
 \qedsymbol \vspace{2mm}

 Conversely, $(2) \Longrightarrow (1)$:
\begin{equation*}
\begin{split}
   &  \rho_U(a) + \rho_U(b) = \rho_U(a + b), \\
   &  U\rho(a) + U\rho(b) = U\rho(a + b) = U(\rho(a) + \rho(b)).
\end{split}
\end{equation*}
So, we have (\ref{distrib_eq_2},(i)). \qedsymbol \vspace{2mm}

 Now, consider $(1) \Longrightarrow (3)$. Set
$$
   \nu(a) = \rho(a) + U,
$$
Then (\ref{distrib_eq_2},(ii)) is equivalent to the relation
\begin{equation}
 \label{repr_nu}
     \nu(a)\nu(b) = \nu(ab).
\end{equation}
Since, $\ker\nabla = U$, we have $\nabla\rho(a) =
 \nabla\nu(a)$. Further,
$$
  \nabla\nu(a + b) = \nabla\rho(a + b) =  \nabla(\rho(a) + \rho(b)) =
  \nabla\rho(a) + \nabla\rho(b) =  \nabla\nu(a) + \nabla\nu(b),
$$
and (3) from Proposition \ref{admis_space} is equivalent to the
following relation:
\begin{equation}
 \label{repr_nabla_nu}
     (\nabla\nu(a))(\nabla\nu(b)) = \nabla\nu(ab).
\end{equation}
Thus, we just need to prove that (\ref{repr_nabla_nu}) follows
from (\ref{repr_nu}). The inclusion
$$
  (\nabla\nu(a))(\nabla\nu(b)) \supseteq \nabla\nu(ab)
$$
is obvious, since $\nu(ab) \subseteq \nu(a), \nu(b)$. So, it is
sufficient to prove that
\begin{equation}
 \label{nabla_nu_inclusion}
     (\nabla\nu(a))(\nabla\nu(b)) \subseteq \nabla\nu(ab).
\end{equation}
Let $z \in \nabla\nu(a)$ and $z \in \nabla\nu(b)$. Then there
exist vectors $v_a \in \nu(a)$ and $v_a \in \nu(b)$ such that
$$
   \nabla(v_a) = \nabla(v_b) = z.
$$
Then, $w = v_a - v_b  \in \ker\nabla = U$, so
$$
   v_b = v_a - w \in \nu(a) + U = \nu(a),
$$
i.e., $v_b  \in \nu(a)$, so
$$
  v_b  \in \nu(a)\nu(b) = \nu(ab)
$$
and $z \in \nabla\nu(ab)$, and hence (\ref{nabla_nu_inclusion})
and (\ref{repr_nabla_nu}) are proved. \qedsymbol \vspace{2mm}

 Finally, consider the implication
  $(3) \Longrightarrow (1)$. We need to prove that
(\ref{repr_nu}) follows from (\ref{repr_nabla_nu}). Let
 $v \in \nu(a)\nu(b)$, i.e.,
$$
    v \in \nu(a) = \rho(a) + U, \text{ and }
    v \in \nu(b) = \rho(b) + U.
$$
By (\ref{repr_nabla_nu}) there exists $w \in \nu(ab)$ such that
$\nabla(v) = \nabla(w)$. Therefore, $v - w \in U$ and
 $v \in \nu(ab) + U = \nu(ab)$. Thus,
$\nu(a)\nu(b) \subseteq \nu(ab)$. The inverse inclusion is obvious
and (\ref{repr_nu}) holds. \qedsymbol

\section[Representations of graphs and lattices]
        {A relation between representations of graphs
         and representation of lattices}
\label{sect_proj_repr}

In a sense, representations of lattices are equivalent to the
representations of graphs, \cite{BGP73}, \cite{Gab72}. We
demonstrate this fact on the example of the graph $\tilde{E}_6$.
\begin{proposition}
 \label{repr_graphs}
Let $\rho$ be the indecomposable representation of the diagram
$\tilde{E}_6$, (\ref{arbitr_repr_rho}).
\begin{equation}
\label{arbitr_repr_rho}
\begin{array}{rccccccccccccl}
   & g_1  & & f_1 & &  f_2 & &  g_2 & &  \\
     X_1 & \longrightarrow
   & Y_1 & \longrightarrow
   & X_0 & \longleftarrow
   & Y_2 & \longleftarrow
   & X_2 &  \\
   & & & & \uparrow  & f_3 \\
   & &  \rho    & & Y_3  \\
   & & & & \uparrow & g_3 \\
   & & & & X_3 &
\end{array}
\end{equation}
If $X_0 \neq 0$, then maps $f_i$ and $g_i$ are monomorphisms.
\end{proposition}
\PerfProof 1) Let $g_1$ be not a monomorphism and
 $G =\ker{g}_1$. Consider the subspace $M$ in $X$ complementary
to $G$ in the direct sum $X_1 = G \oplus M$. Then $\rho$ is the
direct sum of the following two subrepresentations:
\begin{equation}
\label{decompose_g_1}
\begin{array}{rccccccccccccl}
     M & Y_1 & X_0  & Y_2 & X_2 &    \hspace{1cm} &
     G & 0   & 0    & 0   & 0   &  \hspace{1cm} \\
             &      & Y_3  &         &         &
                  \hspace{1cm}\oplus\hspace{1cm} &
             &      & 0    &         &         &  \\
             &      & X_3  &         &         &  \hspace{1cm} &
             &      & 0    & \\
\end{array}
\end{equation}
Since $\rho$ is an indecomposable representation, we see that
$\rho$ coincides with either the first or the second summand in
(\ref{decompose_g_1}). Since $X_0 \neq 0$, we deduce that $\rho$
coincides with the first summand in (\ref{decompose_g_1}) and
 $G = 0$ in the second summand. Thus, $g_1$ is a monomorphism.

2)  Let us prove that $f_1$ is a monomorphism.
\begin{figure}[h]
\includegraphics{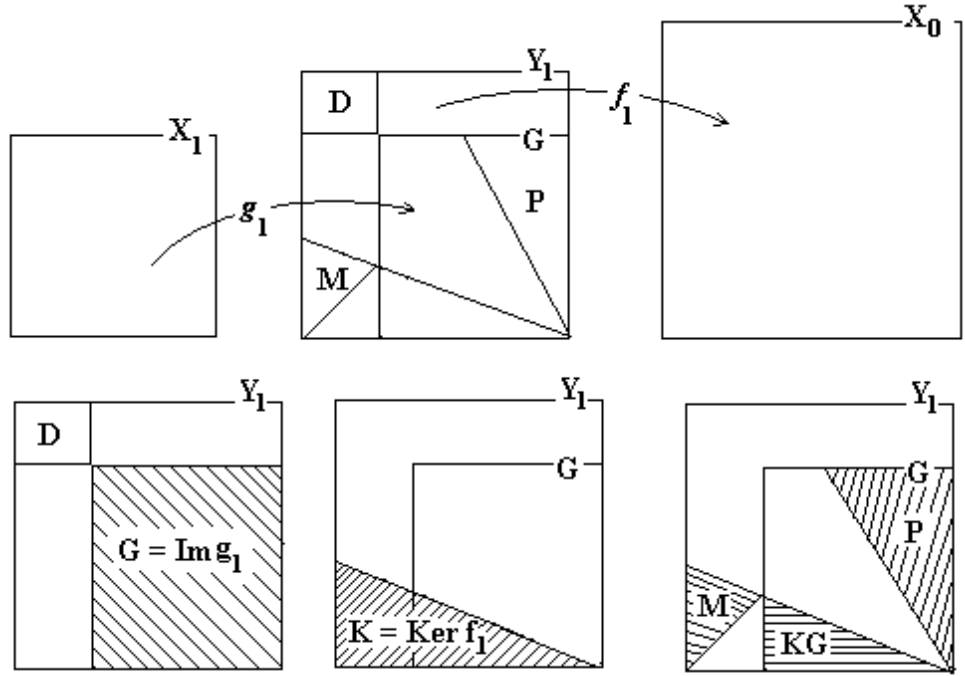}
\caption{\hspace{3mm}Maps $f_1$ and $g_1$.
  $K = \ker{f}_1$, $G = \Im{g}_1$}
\label{monomorphisms}
\end{figure}
 Let $K = \ker{f}_1$ and $G = \Im{g}_1$.
 Consider a complementary subspace $M$ to $KG$ in $K$
 and a complementary subspace $P$ to $KG$ in $G$:
\begin{equation}
 \label{complem_KG}
    K = KG \oplus M, \indent G = KG \oplus P.
\end{equation}
Then, we have
\begin{equation}
 \label{complem_KG_1}
    MG = MKG = 0, \indent PK = PKG = 0.
\end{equation}
Since $PK = P(KG + M) = 0$, we see that the subspaces $P, M$ and
 $KG$ form the direct sum
\begin{equation}
 \label{complem_KG_2}
     M \oplus  P \oplus KG = K + G
\end{equation}
Further, there is a complementary subspace $D$ to $K + G$ in
$Y_1$:
\begin{equation}
\begin{split}
 \label{complem_KG_3}
    & Y_1 = (K + G) \oplus D, \text{ i.e., } \\
    & Y_1 = M \oplus  P \oplus KG  \oplus D  \text{ or } \\
    & Y_1 = M \oplus  G  \oplus D.
\end{split}
\end{equation}


3) Since we already proved that $g_1$ is monomorphism and
$g^{-1}_1(G \oplus D) = X_1$, we see that $\rho$ in the direct sum
of following two subrepresentations:
\begin{equation}
\label{decompose_g_3}
\begin{array}{rccccccccccccl}
     0   & M & 0  & 0 & 0    & \hspace{1cm} &
     X_1 & G \oplus D & X_0  & Y_2 & X_2 &  \hspace{1cm} \\
         &            & 0    &     &     &
                  \hspace{1cm}\oplus\hspace{1cm} &
         &            & Y_3  &     &     &  \\
         &            & 0    &     &     &  \hspace{1cm} &
         &            & X_3  & \\
\end{array}
\end{equation}
Since the representation $\rho$ is indecomposable with
 $X_0 \neq 0$, it follows that $M = 0$. From (\ref{complem_KG}) and
(\ref{complem_KG_3}) we see that $K = KG$ and $K \subseteq G$,
i.e., $\ker{f}_1 \subseteq \Im{g}_1$, and therefore
\begin{equation}
\label{complem_KG_4}
\begin{split}
        & G = K \oplus P, \\
        & Y_1 = K \oplus P \oplus D.
\end{split}
\end{equation}
Since
$$
   X_1 =  g^{-1}(K \oplus P) =
                g^{-1}(K) \oplus g^{-1}(P),
$$
we see that the representation $\rho$ is the direct sum of the
following two subrepresentations:
\begin{equation}
\label{decompose_g_4}
\begin{array}{rccccccccccccl}
     g^{-1}(K) & K          & 0   & 0 & 0 &  \hspace{1cm} &
     g^{-1}(P) & P \oplus D & X_0 & 0 & 0 &  \hspace{1cm} \\
               &            & 0  &         &         &
                               \hspace{1cm}\oplus\hspace{1cm} &
               &            & Y_3    &         &         &  \\
               &            & 0      &         &         &  \hspace{1cm} &
               &            & X_3    & \\
\end{array}
\end{equation}
As above, the first summand is zero, $K = \ker{f}_1 = 0$,
 $G = P = \Im{g}_1$ and $f_1$ is monomorphism. \qedsymbol

 An analog of this proposition for an arbitrary star
graph with orientation in which all arrows directed to the branch
vertex space $X_0$ is also true, as one can show.
\begin{corollary}
 \label{corol_repr_graphs}
  The category of representations of the graph $\widetilde{E}_6$ with
  orientation (\ref{arbitr_repr_rho}) coincides
  with the category of representations of the lattice $D^{2,2,2}$ in non-zero
  spaces:
$$
     \rho : L \longrightarrow \mathcal{L}(X_0), \indent
     X_0 \neq 0.
$$
\end{corollary}

 \index{Coxeter functors $\Phi^+$, $\Phi^-$}
 \index{Gelfand-Ponomarev!- Coxeter functor $\Phi^+$}

In particular, the Coxeter functors $\Phi^+$, $\Phi^-$ which
appeared in the representation theory of graphs (\cite{BGP73},
\cite{GP74}) can be applied to the representations of lattices.
The reflection functors from \cite{BGP73} can not, however, be
used because these functors change the orientation of the graph.

 \index{representation!- projective}
 \index{representation!- injective}
 \index{injective representation}
 \index{projective representation}
The representation $\rho$ for which $\Phi^+\rho$ = 0 (resp.
$\Phi^-\rho$ = 0) is called the {\it projective} (resp. {\it
injective}). For every indecomposable representation $\rho$, the
new indecomposable representation $\Phi^+\rho$ (resp.
$\Phi^-\rho$) can be constructed except for the case where $\rho$
is {\it projective} (resp. {\it injective}).

For example, we can construct a new indecomposable representation
$\Phi^+\rho$ of $\widetilde{E}_6$ except for seven indecomposable
representations $\rho_{x_0}$ and $\rho_{y_i}$, $\rho_{x_i}$, $i =
1,2,3$; for them $\Phi^+\rho$ = 0, see Table \ref{proj_repr_E6}.

By \cite[Prop. 8,9]{GP79} {\it the projective indecomposable
representations of any graph are naturally enumerated by the
vertices of the graph and can be recovered from the orientation of
the graph} .

If $\rho$ is indecomposable and not {\it projective}, i.e.,
$\Phi^+\rho \neq 0$, then $\rho = \Phi^-\Phi^+\rho$.

If $\rho$ is indecomposable and not {\it injective}, i.e.,
$\Phi^-\rho \neq 0$, then $\rho = \Phi^+\Phi^-\rho$.

If $\rho$ is indecomposable and $(\Phi^+)^k\rho \neq 0$, then
$\rho = (\Phi^-)^k(\Phi^+)^k\rho$.

If $\rho$ is indecomposable and $(\Phi^-)^k\rho \neq 0$, then
$\rho = (\Phi^+)^k(\Phi^-)^k\rho$.

 \index{representation!- preprojective}
 \index{preprojective representation}

The representation $\rho$ is called {\it preprojective} if, for
some {\it projective} representation $\tilde\rho$,
$$
    (\Phi^+)^k\rho = (\Phi^+)^k(\Phi^-)^k\tilde\rho = \tilde\rho,
    \indent
    (\Phi^+)^{k+1}\rho =  \Phi^+\tilde\rho = 0.
$$

 \index{representation!- preinjective}
 \index{preinjective representation}

The representation $\rho$ is called {\it preinjective} if, for
some {\it injective} representation $\tilde\rho$,
$$
    (\Phi^-)^k\rho = (\Phi^-)^k(\Phi^+)^k\tilde\rho = \tilde\rho,
    \indent
    (\Phi^-)^{k+1}\rho =  \Phi^-\tilde\rho = 0.
$$

 \index{representation!- regular}
 \index{regular representation}

The representation $\rho$ is called {\it regular}
if $(\Phi^+)^k\rho \neq 0$ and $(\Phi^+)^k\rho \neq 0$ for every
$k \in \mathbb{Z}$.

\begin{table}[h]
  \label{proj_repr_E6}
  \renewcommand{\arraystretch}{1.35}
  \begin{tabular} {|| c || c | c | c ||}
   \hline \hline
       &  $\rho$ & $\Phi^-\rho$ & $(\Phi^-)^2\rho$  \\
   \hline \hline 
       $\rho_{x_0}$
       &
       $\begin{array}{rccccl}
     {\bf 0} & {\bf 0} & {\bf 1}  & {\bf 0} & {\bf 0} &  \\
             &         & {\bf 0}  &         &         &  \\
             &         & {\bf 0}  &         &         &
        \end{array}$
       &
       $\begin{array}{rccccl}
     {\bf 0} & {\bf 1} & {\bf 2}  & {\bf 1} & {\bf 0} &  \\
             &         & {\bf 1}  &         &         &  \\
             &         & {\bf 0}  &         &         &
        \end{array}$
       &
       $\begin{array}{rccccl}
     {\bf 1} & {\bf 2} & {\bf 4}  & {\bf 2} & {\bf 1} &   \\
             &         & {\bf 2}  &         &         &   \\
             &         & {\bf 1}  &
        \end{array}$ \\
   \hline \hline 
       $\rho_{y_1}$
       &
       $\begin{array}{rccccl}
     {\bf 0} & {\bf 1} & {\bf 1}  & {\bf 0} & {\bf 0} &  \\
             &         & {\bf 0}  &         &         &  \\
             &         & {\bf 0}  &         &         &
        \end{array}$
       &
       $\begin{array}{rccccl}
     {\bf 1} & {\bf 1} & {\bf 2}  & {\bf 1} & {\bf 0} &  \\
             &         & {\bf 1}  &         &         &  \\
             &         & {\bf 0}  &         &         &
        \end{array}$
       &
       $\begin{array}{rccccl}
     {\bf 0} & {\bf 1} & {\bf 3}  & {\bf 2} & {\bf 1} &   \\
             &         & {\bf 2}  &         &         &   \\
             &         & {\bf 1}  &
        \end{array}$ \\
   \hline  
       $\rho_{y_2}$
       &
       $\begin{array}{rccccl}
     {\bf 0} & {\bf 0} & {\bf 1}  & {\bf 1} & {\bf 0} &  \\
             &         & {\bf 0}  &         &         &  \\
             &         & {\bf 0}  &         &         &
        \end{array}$
       &
       $\begin{array}{rccccl}
     {\bf 0} & {\bf 1} & {\bf 2}  & {\bf 1} & {\bf 1} &  \\
             &         & {\bf 1}  &         &         &  \\
             &         & {\bf 0}  &         &         &
        \end{array}$
       &
       $\begin{array}{rccccl}
     {\bf 1} & {\bf 1} & {\bf 3}  & {\bf 2} & {\bf 0} &   \\
             &         & {\bf 2}  &         &         &   \\
             &         & {\bf 1}  &
        \end{array}$ \\
   \hline  
       $\rho_{y_3}$
       &
       $\begin{array}{rccccl}
     {\bf 0} & {\bf 0} & {\bf 1}  & {\bf 0} & {\bf 0} &  \\
             &         & {\bf 1}  &         &         &  \\
             &         & {\bf 0}  &         &         &
        \end{array}$
       &
       $\begin{array}{rccccl}
     {\bf 0} & {\bf 1} & {\bf 2}  & {\bf 1} & {\bf 0} &  \\
             &         & {\bf 1}  &         &         &  \\
             &         & {\bf 1}  &         &         &
        \end{array}$
       &
       $\begin{array}{rccccl}
     {\bf 1} & {\bf 1} & {\bf 3}  & {\bf 2} & {\bf 1} &   \\
             &         & {\bf 2}  &         &         &   \\
             &         & {\bf 0}  &
        \end{array}$ \\
    \hline \hline 
       $\rho_{x_1}$
       &
       $\begin{array}{rccccl}
     {\bf 1} & {\bf 1} & {\bf 1}  & {\bf 0} & {\bf 0} &  \\
             &         & {\bf 0}  &         &         &  \\
             &         & {\bf 0}  &         &         &
        \end{array}$
       &
       $\begin{array}{rccccl}
     {\bf 0} & {\bf 0} & {\bf 1}  & {\bf 1} & {\bf 0} &  \\
             &         & {\bf 1}  &         &         &  \\
             &         & {\bf 0}  &         &         &
        \end{array}$
       &
       $\begin{array}{rccccl}
     {\bf 0} & {\bf 1} & {\bf 2}  & {\bf 1} & {\bf 1} &   \\
             &         & {\bf 1}  &         &         &   \\
             &         & {\bf 1}  &
        \end{array}$ \\
    \hline 
       $\rho_{x_2}$
       &
       $\begin{array}{rccccl}
     {\bf 0} & {\bf 0} & {\bf 1}  & {\bf 1} & {\bf 1} &  \\
             &         & {\bf 0}  &         &         &  \\
             &         & {\bf 0}  &         &         &
        \end{array}$
       &
       $\begin{array}{rccccl}
     {\bf 0} & {\bf 1} & {\bf 1}  & {\bf 0} & {\bf 0} &  \\
             &         & {\bf 1}  &         &         &  \\
             &         & {\bf 0}  &         &         &
        \end{array}$
       &
       $\begin{array}{rccccl}
     {\bf 1} & {\bf 1} & {\bf 2}  & {\bf 1} & {\bf 0} &   \\
             &         & {\bf 1}  &         &         &   \\
             &         & {\bf 1}  &
        \end{array}$ \\
    \hline 
       $\rho_{x_3}$
       &
       $\begin{array}{rccccl}
     {\bf 0} & {\bf 0} & {\bf 1}  & {\bf 0} & {\bf 0} &  \\
             &         & {\bf 1}  &         &         &  \\
             &         & {\bf 1}  &         &         &
        \end{array}$
       &
       $\begin{array}{rccccl}
     {\bf 0} & {\bf 1} & {\bf 1}  & {\bf 1} & {\bf 0} &  \\
             &         & {\bf 0}  &         &         &  \\
             &         & {\bf 0}  &         &         &
        \end{array}$
       &
       $\begin{array}{rccccl}
     {\bf 1} & {\bf 1} & {\bf 2}  & {\bf 1} & {\bf 1} &   \\
             &         & {\bf 1}  &         &         &   \\
             &         & {\bf 0}  &
        \end{array}$ \\
     \hline \hline
  \end{tabular}
  \vspace{2mm}
\caption[\hspace{3mm}Preprojective representations]
{\hspace{3mm}Preprojective representations. Seven representations
in the first column are projective}
  \label{RepOrder12}
\end{table}

\section{The Coxeter functor $\Phi^+$ for $D^{2,2,2}$}
 \label{correctness}

Let $\rho$ be an indecomposable representation of
$\widetilde{E}_6$ in a space $X_0$ and $\rho^1 = \Phi^+\rho$ an
indecomposable representation of $\widetilde{E}_6$ in the space
$X^1_0$. According to \cite{BGP73}, the Coxeter functor $\Phi^+$
is constructed as the sequence of reflection functors $F^+_z$,
where $z$ runs over the vertices of $\widetilde{E}_6$:
\begin{equation*}
     \Phi^+\rho =
      F^+_{x_1}F^+_{x_2}F^+_{x_3}
      F^+_{y_1}F^+_{y_2}F^+_{y_3}
      F^+_{x_0}\rho\hspace{0.5mm}.
\end{equation*}
\index{reflection functor $F^+_{x_i}$}

\begin{figure}[h]
\includegraphics{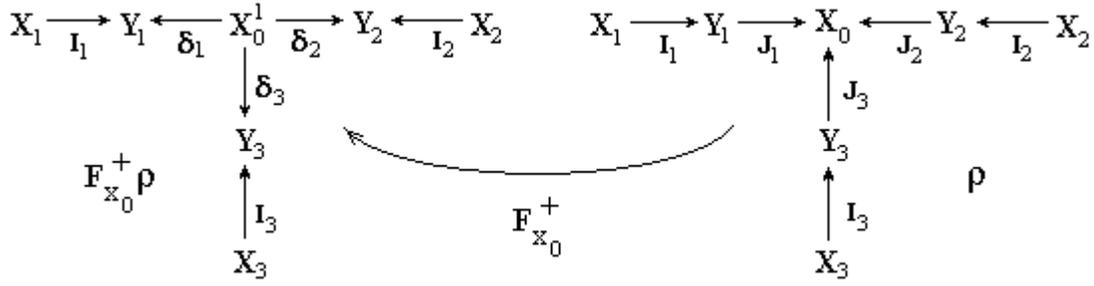}
\caption{\hspace{3mm}The reflection functor $F^+_{x_0}$}
\label{reflect1}
\end{figure}

The reflection functor $F^+_{x_0}$ changes only the space $X_0$ to
$X^1_0$ and maps $Y_i \longrightarrow X_0$ to maps
 $\delta_i: X^1_0 \longrightarrow Y_i$ for each $i=1,2,3$. Let
$$
   \nabla:\{(\eta_1, \eta_2, \eta_3) \mid \eta_i \in Y_i\}
   \longrightarrow \sum{J_i({\eta_i})},
$$
see Fig. \ref{reflect1}. Then by \cite{BGP73}, we have $X^1_0$ =
$\ker\nabla$ from the exact sequence
\begin{equation}   \label{es1}
0 \longrightarrow \ker\nabla
  \longrightarrow \oplus Y_i
  \stackrel{\nabla}{\longrightarrow}
   X_0 \longrightarrow 0
\end{equation}
and
\begin{equation}
 \delta_i : \{(\eta_1, \eta_2, \eta_3) \mid \eta_i \in Y_i,
   \hspace{1mm} \Sigma{J_i(\eta_i)} = 0 \}
   \longrightarrow \eta_i.
\end{equation}

Now, we have
    $X^1_0$ $\stackrel{\delta_i}\longrightarrow$ $Y_i$
         $\stackrel{I_i}\longleftarrow$ $X_i$.
Applying the reflection functor $F^+_{y_i}$, we only change the
space $Y_i$ to the space $Y^1_i$, the corresponding maps
$\delta_i$ to $\mu_i$ and $I_i$ to $P_i$ for each $i =1,2,3$. So
\begin{equation*}
 \begin{split}
  \nabla_i : \{((\eta_1, \eta_2, \eta_3), \xi_i) \mid
    \eta_i \in Y_i, & \Sigma{J_i(\eta_i)} = 0, \xi_i \in X_i \}
       \longrightarrow \\
       & \longrightarrow \delta_i(\eta_1, \eta_2, \eta_3) + I_i(\xi) =
          \eta_i + I_i(\xi_i).
 \end{split}
\end{equation*}
Then $Y^1_i = \ker\nabla_i$ from the following exact sequences,
where $i = 1,2,3$.
\begin{equation}  \label{es2}
0 \longrightarrow \ker\nabla_i
  \longrightarrow X^1_0\oplus X_i
  \stackrel{\nabla_i}\longrightarrow
   Y_i \longrightarrow 0.
\end{equation}
We have,
\begin{equation}
 \label{es2_1}
 \begin{split}
 Y^1_i = & \ker\nabla_i = \\
 & \{((\eta_1, \eta_2, \eta_3), \xi_i) \mid
   \eta_i \in Y_i, \quad
   \Sigma{J_k(\eta_k)} = 0, \quad
   \xi_i \in X_i , \quad
    \eta_i + I_i(\xi_i) = 0\}
 \end{split}
\end{equation}
and $X^1_0 \stackrel{\mu_i}\longleftarrow Y^1_i
         \stackrel{P_i}\longrightarrow X_i$.
Here $\mu_i : Y^1_i \longrightarrow X^1_0$,
\begin{equation} \label{es2_2}
\begin{split}
& \mu_1((\eta_1, \eta_2, \eta_3), \xi_1) =
      ((\eta_1, \eta_2, \eta_3) \mid
       \eta_k \in Y_k, \quad
       \eta_1 + I_1(\xi_1) = 0, \quad
       \Sigma{J_k(\eta_k)} = 0)), \\
& \mu_2((\eta_1, \eta_2, \eta_3), \xi_2) =
      ((\eta_1, \eta_2, \eta_3) \mid
      \eta_k \in Y_k, \quad
       \eta_2 + I_1(\xi_2) = 0, \quad
       \Sigma{J_k(\eta_k)} = 0)), \\
& \mu_3((\eta_1, \eta_2, \eta_3), \xi_3) =
      ((\eta_1, \eta_2, \eta_3) \mid
      \eta_k \in Y_k, \quad
       \eta_3 + I_1(\xi_3) = 0, \quad
       \Sigma{J_k(\eta_k)} = 0)).
\end{split}
\end{equation}
\\
\begin{figure}[h]
\includegraphics{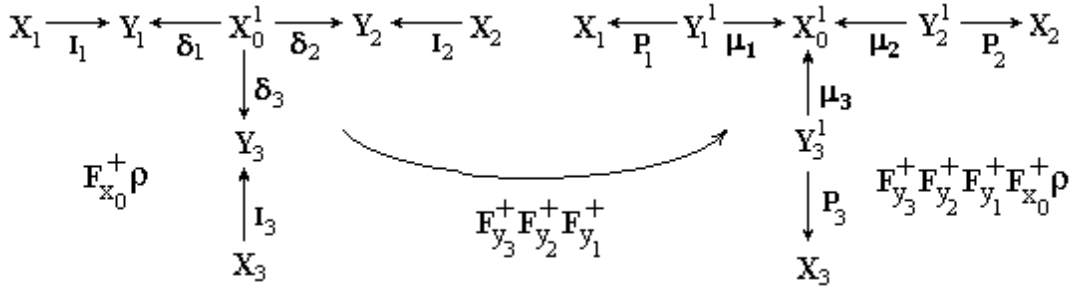}
\caption{\hspace{3mm}The reflection functor $F^+_{y_3}F^+_{y_2}F^+_{y_1}F^+_{x_0}$}
\label{reflect2}
\end{figure}
\begin{remark}
  {\rm In any representation of the lattice $D^{2,2,2}$, the maps
  $I_k$ and $J_k$ are monomorphisms and we have  $\eta_i = -\xi_i$ in (\ref{es2_1}),
  i.e., $\eta_i \in X_i$. In other words,
\begin{equation}
 \label{es2_2_mu}
 \begin{split}
 & \mu_1((\eta_1, \eta_2, \eta_3), \xi_1) =
      ((-\xi_1, \eta_2, \eta_3) \mid
        \eta_k \in Y_k, \quad \xi_1 \in X_1, \quad
       \Sigma{\eta_k} = 0)), \\
 & \mu_2((\eta_1, \eta_2, \eta_3), \xi_2) =
      ((\eta_1, -\xi_2, \eta_3) \mid
        \eta_k \in Y_k, \quad \xi_2 \in X_2, \quad
       \Sigma{\eta_k} = 0)), \\
 & \mu_3((\eta_1, \eta_2, \eta_3), \xi_3) =
      ((\eta_1, \eta_2, -\xi_3) \mid
        \eta_k \in Y_k, \quad \xi_3 \in X_3, \quad
       \Sigma{\eta_k} = 0)).
 \end{split}
\end{equation}
  and \Im$\mu_i$ = $G{'}_i{X}^1_0$.
     This proves that (\ref{def_X}) from \S\ref{cox_functor} is well-defined.}
\end{remark}
Further,
\begin{equation}  \label{es3}
0 \longrightarrow \ker{P}_i
  \stackrel{Q_i}\longrightarrow Y^1_i
  \stackrel{P_i}\longrightarrow
   X_i \longrightarrow 0  \\
\end{equation}
and $X^1_i = \ker{P}_i$, where
\begin{equation} \label{es3_1}
\begin{split}
& P_1((\eta_1, \eta_2, \eta_3), \xi_1) =
      (\xi_1 \mid \eta_k \in Y_k, \quad
       \eta_1 + I_1(\xi_1) = 0, \quad \Sigma{J_k(\eta_k)} = 0)), \\
& P_2((\eta_1, \eta_2, \eta_3), \xi_2) =
      (\xi_2 \mid \eta_k \in Y_k, \quad
       \eta_2 + I_1(\xi_2) = 0, \quad \Sigma{J_k(\eta_k)} = 0)), \\
& P_3((\eta_1, \eta_2, \eta_3), \xi_3) =
      (\xi_3 \mid \eta_k \in Y_k, \quad
       \eta_3 + I_1(\xi_3) = 0, \quad \Sigma{J_k(\eta_k)} = 0)). \\
\end{split}
\end{equation}

\begin{remark}
  {\rm For any representation of the lattice $D^{2,2,2}$,
  the maps $J_k$ are monomorphisms, therefore
\begin{equation} \label{es3_2}
\begin{split}
& \ker{P}_1 =
  \{((0, \eta_2, \eta_3), 0) \mid \eta_k \in Y_k, \quad \Sigma\eta_k = 0)\}, \\
& \ker{P}_2 =
  \{((\eta_1, 0, \eta_3), 0) \mid \eta_k \in Y_k, \quad \Sigma\eta_k = 0)\}, \\
& \ker{P}_3 =
 \{((\eta_1, \eta_2, 0), 0) \mid \eta_k \in Y_k, \quad \Sigma\eta_k = 0)\}.
 \vspace{1mm} \\
\end{split}
\end{equation}
\begin{figure}[h]
\includegraphics{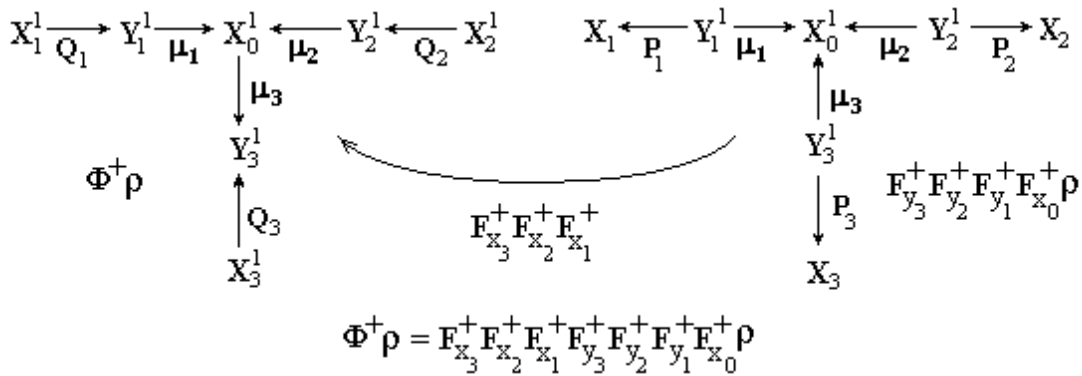}
\caption{\hspace{3mm}The Coxeter functor $\Phi^+$}
\label{reflect3}
\end{figure}
Since
\begin{equation} \label{es3_3}
\begin{split}
& \mu_1Q_1(\ker{P}_1) =
   \{(0, \eta_2, \eta_3) \mid \eta_k \in Y_k, \quad \Sigma\eta_k = 0)\}, \\
& \mu_2Q_2(\ker{P}_2) =
   \{(\eta_1, 0, \eta_3) \mid \eta_k \in Y_k, \quad \Sigma\eta_k = 0)\}, \\
& \mu_3Q_3(\ker{P}_3) =
   \{(\eta_1, \eta_2, 0) \mid \eta_k \in Y_k, \quad \Sigma\eta_k = 0)\} \\
\end{split}
\end{equation}
and $\Im(\mu_i{Q}_i) = H{'}_i{X}^1_0$, it follows that
(\ref{def_Y}) from \S\ref{cox_functor} is well-defined.}
\end{remark}

\section{Linear lattices}
 \label{sect_lin_lat}
 \index{partition of set}
 \index{linear lattice}
 \index{partition lattice}
 \index{lattice!- of partitions}
 \index{blocks of partition}
 In this section, we mostly follow the definitions of the {\it linear lattices}
from works of M.~Haiman \cite{Hai85} and M.~Mainetti, C.~H.~Yan
 \cite{MY99}.

A set $\pi$ of non-empty pairwise disjoint subsets such that their
union coincides with $A$ is said to be a {\it partition} of the
set $A$. The elements of $\pi$ are called {\it blocks}
 of the partition $\pi$. If elements $a, b \in A$ lie in the same
 block $\pi$, we write
$$
    a \equiv b \mod\pi, \text{ or } a \equiv b (\pi), \text{ or } a\pi{b}.
$$
\begin{proposition}
 \label{part_eqiuv_rel}
  There exists a bijection between partitions $\pi$ of the
  set $A$ and equivalence relations $R_{\pi}$ on the set $A$.
\end{proposition}
 \PerfProof
Given a partition $\pi$, we define the equivalence relation
$R_\pi$ on the set $A$ as follows:
\begin{equation}
 \begin{array}{ccccc}
   & <a, b> \in R_{\pi} & \Longleftrightarrow & a\pi{b}, & \text{ or }\\
   & aR_{\pi}b & \Longleftrightarrow & a\pi{b}. &\\
 \end{array}
\end{equation}
Conversely, if $R$ is an equivalence relation, then the partition
$\pi$ is defined to be:
\begin{equation}
   \pi = \{ [a]R  \mid  a \in A \}.
\qed \vspace{2mm}
\end{equation}

\begin{proposition}
 \label{lattice_of_part}
The set of partitions ${\rm Part}(A)$ of set $A$ constitutes a
lattice called the partition lattice.
\end{proposition}

\PerfProof
 Let $\pi_1$, $\pi_2$ be partitions of $A$:
\begin{equation*}
 \begin{split}
  & \pi_1 = \bigcup\limits_{s \in S}\pi_1^s,
     \text{ where } \pi_1^{s_1} \cap \pi_1^{s_2} =
      \emptyset \text{ for all } s_1,s_2 \in S, \\
  & \pi_2 = \bigcup\limits_{t \in T}\pi_1^t,
     \text{ where } \pi_1^{t_1} \cap \pi_1^{t_2} =
      \emptyset \text{ for all } t_1,t_2 \in T. \\
 \end{split}
\end{equation*}
For two partitions $\pi_1$ and $\pi_2$, their intersection,
inclusion and union are defined as follows. Blocks of intersection
are pairwise intersections of blocks of $\pi_1$ and $\pi_2$.
\begin{equation}
 \label{part_inters}
   \pi_1 \cap \pi_2 = \bigcup\limits_{s \in S, t \in T}(\pi_1^s \cap \pi_2^t),
     \text{ for all } s \in S, t \in T. \\
\end{equation}
Every block of $\pi_1 \subseteq \pi_2$ is included in some block
of $\pi_2$.
\begin{equation*}
   \pi_1 \subseteq \pi_2 =
     \{ x \equiv y(\pi_1) \Longrightarrow x \equiv y(\pi_2) \}. \\
\end{equation*}
The union $\pi_1 \cup \pi_2$ is the intersection of all partitions
$\pi$ containing both $\pi_1$ and $\pi_2$.
\begin{equation*}
   \pi_1 \cup \pi_2 = \bigcap\limits_{\pi_1, \pi_2 \subseteq \pi}\pi.
  \qed \vspace{2mm}
\end{equation*}
For more details concerning the {\it partition lattice}, see
\cite[Ch.4]{Gr98}.

\begin{corollary}
  \label{lattice_eq_rel}
   The set of equivalence relations $\{ R_{\pi} \}$ constitutes
   the lattice of equivalence relations ${\rm Rel}(A)$
   isomorphic to the partition lattice
   \begin{equation}
       {\rm Rel}(A) \simeq {\rm Part}(A).
   \end{equation}
\end{corollary}
\PerfProof It follows from Propositions \ref{part_eqiuv_rel} and
\ref{lattice_of_part} \qedsymbol

 \index{lattice!- of commuting equivalence relations}
 \index{lattice!- linear}
 \index{commuting equivalence relations}

Consider now {\it commuting equivalence relations}. Composition
$R\circ{T}$ of relations $R$  and $T$ is defined as follows:
\begin{equation}
  x(R\circ{T})y \text{ if and only if there exists } z \in A \text{ such that }
            xRz \text{ and } zTy .
\end{equation}

Two equivalence relations $R$ and $T$ are said to
 {\it commute}\footnote{It is easy to give an example of non-commuting
equivalence relations: if $R$ is ``father of'' and $T$ is ``mother
of'', then $R$ and $T$ do not commute, since $RT$ is ``maternal
grandfather of'' and $TR$ is ``paternal grandmother of'',
\cite{Bz2001}.} if
$$
 R\circ{T} = T\circ{R}.
$$

\begin{definition}[\cite{MY99}, \cite{MY2000}]
 \label{def_linear_lat}
{\rm A {\it linear lattice} is a sublattice of the lattice of
equivalence relations ($\simeq$ partition lattice) on a set, with
the property that any two equivalence relations in the sublattice
commute, in the sense of composition of relations.}
\end{definition}

The classic example of the linear lattice is the lattice of
subspaces of a vector space $V$. Any subspace $S \subseteq V$
determines an equivalence relation $R_S$ on our vector space,
given by letting
\begin{equation}
  \label{rel_subspace}
  x(R_S)y \text{ if and only if } x-y \in S.
\end{equation}
\begin{proposition}
 \label{subsp_lin_lat}
  1) All equivalence relations produced, as in (\ref{rel_subspace}),
  by subspaces of the vector space $V$, commute.

  2) Commuting equivalence relations given by (\ref{rel_subspace})
  constitute the linear lattice.
\end{proposition}
\PerfProof 1) We have
\begin{equation}
 \label{rel_x_z}
   x(R_{S_1}\circ{R_{S_2}})z \text{ if and only if }
    x-y \in S_1 \text{ and } y-z \in S_2 \text{ for some } y.
\end{equation}
By (\ref{rel_x_z}), we have
\begin{equation}
 \label{rel_x_z_2}
   x-z \in S_1 + S_2.
\end{equation}
Conversely, suppose $x - z = s_1 + s_2$, where $s_1 \in S_1$ and
 $s_2 \in S_2$. Let $y = z - s_1$. Then
\begin{equation}
 \label{rel_x_z_3}
 \begin{split}
  & x - y = x - z + s_1 = s_2 \in S_2, \quad\text{ and } \\
  & y - z = z - s_1 - z  = -s_1 \in S_1.
 \end{split}
\end{equation}
Thus, by (\ref{rel_x_z_2}), (\ref{rel_x_z_2}) and
(\ref{rel_x_z_3}), we see that
\begin{equation*}
   x(R_{S_1}\circ{R_{S_2}})z \text{ if and only if }
    x-z \in S_1 + S_2.
\end{equation*}
Similarly,
\begin{equation*}
   x(R_{S_2}\circ{R_{S_1}})z \text{ if and only if }
    x-z \in S_2 + S_1,
\end{equation*}
and therefore
\begin{equation*}
   x(R_{S_1}\circ{R_{S_2}})z \text{ if and only if }
   x(R_{S_2}\circ{R_{S_1}})z,
\end{equation*}
i.e.,
\begin{equation}
  \label{rel_sum}
  R_{S_1}\circ{R_{S_2}} = R_{S_2}\circ{R_{S_1}} = R_{S_1 + S_2}.
  \qed \vspace{2mm}
\end{equation}

2) By definition (\ref{def_linear_lat}) it suffices to prove that
the commuting equivalence relations defined by
(\ref{rel_subspace}) constitute a lattice. The equivalence
relation corresponding to the sum of subspaces $S_1 + S_2$ is
defined by (\ref{rel_sum}).

By (\ref{def_linear_lat}) blocks of the partition $\pi_S$,
corresponding to the equivalence relation $R_S$ are cosets $V/S$.
In other words, the partition $\pi_S$ consists of the blocks
$\pi^\alpha_S$ as follows
$$
  \pi_S = \{ \pi^\alpha_S \mid x - y \in \pi^\alpha_S
    \Longleftrightarrow x - y \in S \},
$$
By (\ref{part_inters})
the partition corresponding to the intersection of two relations
$S_1$ and $S_2$ consists of cosets $V/(S_1\cap{S_2})$, and
$$
   \pi_{S_1} \cap \pi_{S_2} = \pi_{S_1\cap{S_2}},
$$
i.e.,
$$
   R_{S_1} \cap R_{S_2} = R_{S_1\cap{S_2}}.
   \qed \vspace{2mm}
$$
\begin{remark}
 \label{lots_lin_lat}
{\rm
There are lots of linear lattices. Proposition \ref{subsp_lin_lat}
can be generalized to groups. Given any group $G$, each normal
subgroup $H$ determines an equivalence relation on $G$, which we
call $R_H$, such that
$$
   x(R_H)y \quad \text{ if and only if }\quad  xy^{-1} \in H.
$$
So the lattice of normal subgroups of $G$ becomes a linear
lattice. Similarly, lattices of ideals of a ring are also linear
lattices, see Fig. \ref{hierar_lattices}.}
\end{remark}

 \index{Arguesian law}
 \index{Arguesian lattice}
 \index{lattice!- Arguesian}
 \index{Desargues's theorem}

\begin{definition}[\cite{Jo53}, \cite{Sch45}]
 {\rm
  Let $a_0, a_1, a_2, b_0, b_1, b_2$ be elements of the modular
  lattice $L$, and let
 \begin{equation}
  \label{arg_terms}
  \begin{split}
   & z_{ij} = (a_i + a_j)(b_i + b_j) \text{ for } 0 \leq i < j < 3, \\
   & z = z_{01}(z_{02} + z_{12}).
  \end{split}
 \end{equation}
 The identity
 \begin{equation}
  \label{arg_law}
    (a_0 + b_0)(a_1 + b_1)(a_2 + b_2)
     \subseteq a_0(z + a_1) + b_0(z + b_1)
 \end{equation}
 is called the Arguesian law and
 lattices satisfying the Arguesian law are called {\it Arguesian
 lattices}. For another version of the Arguesian law,
 see \cite[p.3]{MY2000} and \cite{Hai85}.
 Identity (\ref{arg_law}) is a lattice-theoretic form of the
 Desargues theorem of projective geometry, see
 \cite{Sch45}, \cite[p.2]{MY2000}, \cite[Ch.4]{Gr98}.
 }
\end{definition}

\begin{remark}{\rm
 \label{arguesian_lin_lat}
  Every Arguesian lattice is modular, see \cite[Th.1.9]{Jo54}.
  There exists modular lattices which are not Arguesian, see
  \cite[Ch.4.\S5]{Gr98}.
  Linear lattices are Arguesian, but the converse is
  not true: in \cite{Hai91} M.~Haiman  constructed a family of
  Arguesian lattices which cannot be represented as linear
  lattices, see \cite{KY2003}, see Fig. \ref{hierar_lattices}.
  }
\end{remark}

\section{Grassmann-Cayley algebras and modular lattices}
 \label{GC_algebra}

\subsection{A Peano space and the bracket}

 \index{Grassmann-Cayley algebra}
 \index{join $\vee$}
 \index{meet $\wedge$}
 \index{progressive and regressive products}
 \index{Ausdehnungslehre}
 \index{geometric calculus}
 \index{Grassmann's geometric calculus}
 \index{double algebra}

A Grassmann-Cayley algebra (or {\it double algebra}) is
essentially the exterior algebra of a vector space, equipped with
{\bf two} operations, {\it join}\footnote{The classical Grassmann
algebra is the exterior algebra of a vector space, equipped with
only one operation, the wedge product $\wedge$. Barnabei, Brini
and Rota \cite[p. 127]{BBR85} write ``The product in the exterior
algebra is called the {\it join} and is denoted by the symbol
$\vee$. We note that the usage is at variance with ordinary usage,
where exterior multiplication is denoted by the wedge $\wedge$.
This departure from common usage will be amply justified in the
sequel.'' For motivation of the usage the symbol $\vee$, see
\S\ref{motiv_vee}. }
 $\vee$ and {\it meet} $\wedge$, which are the algebraic
analogs of the {\it union}\footnote{This lattice-theoretic
operation union $\cup$ does not coincide with the set-theoretic
operation union $\sqcup$. For vector spaces, the operation union
$\cup$ is the linear span of the union $\sqcup$ as follows:
\begin{equation*}
    A \cup B = \overline{A \sqcup B}.
\end{equation*}
 The equivalent notation frequently used instead the operation union $\cup$
 is the operation {\it sum} $+$.}
$\cup$ and {\it intersection} $\cap$ of subspaces of a vector
space:
\begin{equation}
 \label{Rota_corresp}
  \begin{array}{cccc}
    & \text{ the join } \vee \quad & \longleftrightarrow \quad
      & \text{ the union } \cup (\text{ the sum } + )\\
    & \text{ the meet } \wedge \quad & \longleftrightarrow \quad
      & \text{ the intersection } \cap
  \end{array}
\end{equation}
These operations were introduced by Grassmann \cite{Gra11}, under
the name of {\it progressive} and {\it regressive} products. In
works of G.-C.Rota and his coworkers, the Grassmann work has been
recognized and extended, see \cite{DRS76}, \cite{BBR85}.

 Barnabei, Brini and Rota \cite[p. 121]{BBR85} write ``To the best
of our knowledge of published work, the first mathematicians to
understand, albeit imperfectly, the program {\it Ausdehnungslehre}
\footnote{Grassmann's program {\it Ausdehnungslehre} refers
 to $1844$, see \cite{Gra11}. The English translation of {\it
Ausdehnungslehre} is {\it geometric calculus}. Grassmann's idea
was to develop a calculus for the {\it join} and {\it meet}. In
\cite{BBR85} authors consider Grassmann's geometric calculus and
its relevance in invariant theory.} were Clifford and Schr\"oder,
first, and later A.~N.~Whitehead, \'Elie Cartan, and best
 of all Giuseppe Peano. ... It was Schr\"oder, in an appendix to his
 {\it Algebra and Logik}, who first stressed the analogy between
 the algebra of {\it progressive} and {\it regressive} products,
 and the algebra of sets with {\it union} and {\it intersection}.''

 By analogy to a {\it Hilbert space} endowed
 with a symmetric bilinear form $(x, y)$, called the inner product, and
 to a {\it symplectic space} endowed with an antisymmetric
 bilinear form $[x,y]$, a {\it Peano space} is a vector space $V$
 of dimension $n$, endowed with an antisymmetric
 non-degenerate multilinear form
 $[x_1, \dots, x_n]$ defined for $x_i \in V$, called the {\it
 bracket} \cite{BBR85}. The group of all linear transformations preserving
 the bracket is $SL(V)$, the special linear group.
 The algebra of $SL(V)$-invariants is generated by
 brackets (the first fundamental theorem of
 invariant theory, see \cite[Th. II.6.A]{Weyl46}) or \cite[Th. 2.1]{Dol2003}).

 Let $K$ be the ground field. A {\it bracket} (of step $n$) over the space $V$ is a
 non-degenerate alternating $n$-linear form, i.e.,  a function
 $$
    x_1,\dots, x_n \longmapsto [x_1,\dots,x_n] \in K,
 $$
 where vectors $x_1,\dots, x_n$ lie in the vector space $V$,
 with the following properties:
\begin{equation}
 \begin{split}
   & (1) \quad [x_1,\dots,x_n] = 0 \quad\text{ if at least two among the $x$'s coincide;} \\
   & (2) \quad  \text{ for every }  x,y \in V \quad\text{ and } \alpha, \beta \in K  \\
   & \qquad  [x_1,\dots,x_{i-1}, \alpha{x} + \beta{y}, x_{i+1},\dots,x_n] =  \\
   & \qquad \qquad \qquad 
      \alpha[x_1,\dots,x_{i-1}, x, x_{i+1},\dots,x_n] +
      \beta[x_1,\dots,x_{i-1}, y, x_{i+1},\dots,x_n]; \\
   & (3) \quad \text{there exists a basis }  \{ b_1,\dots,b_n \} \text{ of } V
      \text{ such that } [b_1,\dots,b_n] \neq 0.
 \end{split}
\end{equation}

 \index{exterior algebra}

Let $F(V)$ be the free associative algebra with unity over $K$
generated by the elements of a Peano space $V$. The {\it exterior
algebra $E(V)$ of the Peano space $V$} is obtained as the quotient
of $F(V)$ modulo the ideal generated by $v^2$, for all $v \in V$,
see, for example, \cite[Def. 2.1]{MY99} or \cite[Prop.
3.6]{BBR85}. For further properties of the exterior algebra see,
e.g., \cite{Gre67}. For exterior Grassmann algebra in the context
of supersymmetries, see, e.g., \cite{Ber87} or \cite{L05}.

Let
$$
     \Phi : F(V) \rightarrow E(V),
$$
denote the canonical projection. If
 $x_1x_2\dots{x}_k$ is a word in $F(V)$, we denote its image by
 \index{extensor}
 \index{join of two extensors}
 \index{decomposable antisymmetric tensor}
 \index{decomposable $k$-vector}
 \index{antisymmetric tensor of step $k$}
\begin{equation}
   \Phi(x_1x_2\dots{x}_k) = x_1\vee{x}_2\vee\dots\vee{x}_k.
\end{equation}
The element $x_1\vee{x}_2\vee\dots\vee{x}_k$ is called an {\it
extensor of step $k$}, or {\it decomposable antisymmetric tensor}
or {\it decomposable $k$-vector}. The join of two extensors is an
extensor. A linear combination of extensors of step $k$ is called
an {\it antisymmetric tensor of step $k$}, \cite{BBR85}.

\begin{proposition}{\em (Properties of the join, \cite[Prop. 3.1-3.4]{BBR85})}
  \label{join_prop}
   1) For every permutation $\sigma$ of $\{ 1,2,\dots,k \}$ we have
\begin{equation}
      x_{\sigma(1)}\vee{x}_{\sigma(2)}\dots
      \vee{x}_{\sigma(k)} =
      {\rm sign}(\sigma)x_1\vee{x}_2\dots\vee{x}_k.
\end{equation}

   2) Let $E_k(V)$ be the subspace of $E(V)$ generated by all
      extensors of step $k$. Let $A \in E_j(V), B \in E_k(V)$.
      Then
\begin{equation}
      B \vee A = (-1)^{jk}A \vee B.
\end{equation}

   3) $\dim E_k(V)=\displaystyle\binom{n}{k}$ for $0 \leq k \leq n$, and
   hence $\dim E(V)=2^n$.

   4) Let $B$ be a dimension $k > 0$ subspace of $V$, let
      $\{x_1, x_2, \dots, x_k \}$ and
      $\{y_1, y_2, \dots, y_k \}$ are two bases of $B$. Then
\begin{equation}
      x_1\vee{x}_2\dots\vee{x}_k = Cy_1\vee{y}_2\dots\vee{y}_k
\end{equation}
for some non-zero scalar $C$.
\end{proposition}
  \index{associated extensor}
    By Proposition \ref{join_prop}, heading 4) there exists one-to-one correspondence
    (up to non-zero scalar multiple) between non-trivial subspaces
    of $V$ and non-zero extensor uniquely representing this
    subspace. Such a representing extensor is said to be {\it associated} to
    the subspace.

\subsection{Motivation of the usage the symbol $\vee$}
  \label{motiv_vee}

  Here we formulate the main property of the operation join
  $\vee$ connecting it with the
  lattice-theoretic operation sum $+$ (or, equivalently,
  with the union $\cup$), see \cite[Prop. 3.5]{BBR85}.
\begin{proposition}
  \label{join_and_union}
  Let $A$, $B$ be two subspaces of $V$ associated to extensors
  $F$ and $G$, respectively.
  Then

  (i) $F \vee G = 0$ if and only if $A \cap B \neq {0}$;

  (ii) if $A \cap B = {0}$, then the extensor $F \vee G$ is assosiated
  to the subspace generated by $A + B$.
\end{proposition}

  Reformulate Proposition \ref{join_and_union} as follows. Let
\begin{equation}
  \label{two_extensors}
 \begin{split}
  & F = f_1\vee{f}_2\vee\dots\vee{f}_j, \\
  & G = g_1\vee{g}_2\vee\dots\vee{g}_k,
 \end{split}
\end{equation}
and $\overline{F}$ (resp. $\overline{G}$)
denotes the linear span of extensor $F$ (resp. $G$). \\

  If vectors
$$
    f_1,{f}_2,\dots,{f}_j,g_1,{g}_2,\dots,{g}_k
$$
are distinct and linearly independent, i.e.
 $\overline{F} \cap \overline{G} = 0$,  then
\begin{equation}
  \label{F_vee_G}
   \overline{F\vee{G}} = \overline{F} + \overline{G}.
\end{equation}
Thus, the join operation $\vee$ on extensors represents the
geometric operation of the lattice sum of subspaces of $V$ at
least for independent vectors. This is one of motivations to use
the symbol $\vee$.

Another motivation is a discovery of a class of identities that
hold in the exterior algebras, which turn into the classical
theorems of projective geometry on the incidence of subspaces of
projective spaces, and therefore to identities in modular
lattices, see (\ref{Rota_corresp}).

  For more details concerning these identities,
  see \cite{Haw94}, \cite{Haw96}, \cite{MY2000}, \cite{MY99}.

 \index{Desargues's theorem}
 \index{Rota's rendering}

Theorems such as Desargues, Pappus, Bricard, and Fonten\'e of
classical projective geometry are revealed to be expressible as
simple and elegant identities holding among joins and meets of
extensors in the framework of the Grassmann-Cayley algebra,
\cite{BBR85}, \cite{MY99}, \cite{MY2000}.
 We suggest to call the transformation
 of modular lattice polynomials to Grassmann-Cayley algebra
 polynomials and {\it vice versa} by the substitution
 (\ref{Rota_corresp}) the {\it Rota rendering}. The Rota rendering establishes
 a connection between the theory of Grassmann-Cayley algebras
 and the modular (and even linear) lattices.

\subsection{The meet}

A similar operation on extensors should be defined for the lattice
operation intersection of two subspaces. This is the {\it meet}
operation. Let extensors $F$ and $G$ are given by
(\ref{two_extensors}) with $j + k \geq n$, where $n$ is dimension
of a {\it Peano space} $V$. Then

\begin{equation}
  \label{def_meet}
   F\wedge{G} = \sum\limits_{\sigma}\text{sign}(\sigma)
     [a_{\sigma(1)},\dots,a_{\sigma(d-k)},b_1,\dots,b_k]
     a_{\sigma(d-k+1)}\dots{a}_{\sigma(j)}.
\end{equation}

The meet is associative and anti-commutative in the following
 sense (\cite[Prop. 4.4]{BBR85}):
\begin{equation}
  \label{anti_comm}
   F\wedge{G} = (-1)^{(d-k)(d-j)}G\wedge{F}.
\end{equation}

Both the join and meet are defined on extensors only, but the
definitions are extended to arbitrary elements of $E(V)$ by
distributivity.

\begin{proposition}{\em (\cite[Prop. 4.3]{BBR85}, \cite[Prop.
2.5]{MY99})}
1) We have
$$
   F\wedge{G} = 0 \text{ if and only if } F + G \text{ does not span } V.
$$

2) If $F + G$ spans $V$, then\footnote{Compare with
(\ref{F_vee_G}).}
  \index{integral in the Grassmann-Cayley algebra}
\begin{equation}
   \overline{F\wedge{G}} = \overline{F} \cap \overline{G}.
\end{equation}

\end{proposition}
Thus the meet operation corresponds to the intersection $\cap$.
  \index{unimodular basis}
A basis $\{e_1,\dots,e_n\}$ in the space $V$ is said to be
 {\it unimodular} if $[e_1,\dots,e_n] = 1$. The extensor
$$
    E = e_1 \vee e_2 \vee \dots \vee e_n
$$
in $E(V)$ is called the {\it integral}. Indeed, it corresponds to
the Berezin integral on supermanifolds, cf. \cite{Ber87},
\cite{L05}. The integral is well defined and does nor depend on
the choice of a unimodular basis.

\subsection{Note on Howe-Huang projective invariants of quadruples}

 \index{projective invariants of quadruples}
 \index{Howe-Huang invariants}

Projective invariants of quadruples were obtained by Howe~R.,
Huang~R. in \cite{HH96}. They gave an explicit set of generators
for the ring of invariants of a set of four subspaces in a
projective space.
\begin{remark}{\rm For other results in this
direction, see \cite{Gro98}, \cite{Rin80}. In  \cite{Gro98},
Grosshans gave a complete description of the relationship between
invariants and quadruples. Invariants of the tame quivers were
obtained by Ringel in \cite{Rin80}.}
\end{remark}
 \index{Cayley factorization}
  Among projective invariants of quadruples
obtained by Howe-Huang, the invariants of type III (\cite[Th.
14]{HH96}) can be transformed to the polynomials $f_{GC}(E_1, E_2,
E_3, E_4)$ in four subspaces $E_1, E_2, E_3, E_4$ under the join
$\vee$ and meet $\wedge$ in the framework of the Grassmann-Cayley
algebra. The transformation of Howe-Huang projective invariants to
polynomials $f_{GC}(E_1, E_2, E_3, E_4)$ is performed by
 means of the {\it Cayley factorization} in the sense of
 \cite{Wh91}.

 Further, applying the Rota rendering (\ref{Rota_corresp}) to
$f_{GC}(E_1, E_2, E_3, E_4)$, we get polynomials $f_{ML}(E_1, E_2,
E_3, E_4)$ in the modular lattice $D^4$. It may be shown that {\it
polynomials $f_{ML}(E_1, E_2, E_3, E_4)$ are admissible elements
in $D^4$}.

\chapter{\sc\bf A proof of the theorem on admissible classes}
 \label{sect_proof_adm}

We will briefly write
\begin{equation}
\begin{split}
 & \OperatorFi{a}{i}{b}
 \quad \text{ instead of } \quad
 \varphi_i\Phi^+\rho(a)  = \rho(b), \\
 & \OperatorFiFj{a}{i}{b}{j}{c}
\quad \text{ instead of } \quad \varphi_i\Phi^+\rho(a)  = \rho(b),
\quad
\varphi_j\Phi^+\rho(b) = \rho(c),  \\
 & \ldots  \\
 & \OperatorFi{a}{i}{b} \ldots  \OperatorFi{c}{i}{d}
                            \quad \text{ instead of } \quad
 \varphi_i\Phi^+\rho(a) = \rho(b), \ldots,
 \varphi_j\Phi^+\rho(c) = \rho(d),
\end{split}
\end{equation}
see (\ref{rel_elem_map_0}), (\ref{rel_elem_map_1}).

\section{A proof of the theorem on admissible elements in lattice $D^{2,2,2}$}

 We will now prove the main property of admissible elements.

{\bf Theorem \ref{th_adm_classes}}
  {\it Let $\alpha = i_n{i}_{n-1}\dots{1}$
   be an admissible sequence and $i \neq i_n$. Then
   $i\alpha$ is admissible and,
   for $z_\alpha = f_\alpha, e_\alpha, g_{\alpha0}$ from the
   Table \ref{table_adm_elems}, the following relation holds:}
\begin{equation} \label{adm_classes_rep}
    \varphi_i\Phi^+\rho(z_\alpha) = \rho(z_{i\alpha}).
\end{equation}

{\it Proof of Theorem \ref{th_adm_classes}} For the cases of
 Lines 1--6, it suffices to prove only
\begin{equation} \label{line12}
\OperatorFi{z_\gamma}{3}{z_{3\gamma}} \hspace{3mm} \text{ and }
\hspace{3mm} \OperatorFi{z_{3\gamma}}{1}{z_{(13)\gamma}},
\end{equation}
because the relations
$$
   \OperatorFiFj{z_{(13)\gamma}}{2}{z_{(213)\gamma}}{3}{z_{3(213)\gamma}}
   \stackrel{\varphi_1}{\longmapsto}
   \OperatorFi{z_{13(213)\gamma}}{2}{z_{(213)(213)\gamma}}
$$
follow from (\ref{line12}) as the permutation $1 \longrightarrow 3
\longrightarrow 2 \longrightarrow 1$ is applied to Lines 1--3 of
Table \ref{table_adm_elems} from \S\ref{sect_adm_classes}.

\subsection{Case $\OperatorFi{z_\gamma}{3}{z_{3\gamma}}$}
   \label{case_line_1_2}
~\\

 \underline{\text{(F) Line 1}}. We will show that
\begin{equation}
  \label{f_line1}
\begin{split}
  & \OperatorFi{f_\gamma}{3}{f_{3\gamma}}, \\
  & f_\gamma =
   \OperatorFi{y_1y_2a^{13}_qA^{32}_{k-1}}{3}
              {y_3(x_1 + x_2)A^{23}_qa^{31}_{k-1}} =
   f_{3\gamma}.
\end{split}
\end{equation}
 Applying Corollary \ref{cor_psi} to $y_1y_2$ and
 Corollary \ref{act_runner_map} to $a^{13}_q$ and $A^{32}_{k-1}$
 we have
 \begin{equation} \label{f_line1_comp}
 \begin{split}
        & \OperatorFi{y_1y_2} {3} {y_3(y_1 + y_2)}, \\
        & \OperatorFi{a^{13}_q} {3} {y_3A^{23}_q}, \\
        & \OperatorFi{A^{32}_{k-1}} {3} {y_3(y_1 + y_2)a^{31}_{k-1}}.
 \end{split}
 \end{equation}
 Since $a^{13}_q$ and $A^{32}_{k-1}$ are
 the $\varphi_3-$homomorphic polynomials
 (Theorem  \ref{th_homomorhism}), we see that
 ``intersection of relations''
 (\ref{f_line1_comp}) gives us relation (\ref{f_line1}).
 If $q=1$, then
 $$
    \OperatorFi{a^{13}_q} {3} {y_3(y_1+y_2)A^{23}_q}
    \subseteq y_3(y_1+y_2)
 $$
 and we have the same ``intersection of relations''.
\qedsymbol

 \underline{\text{(E) Line 1}}. We want to show that
\begin{equation} \label{f_line2}
\begin{split}
  &  \OperatorFi{e_\gamma}{3}{e_{3\gamma}},  \\
  &  e_\gamma = \OperatorFi{y_2A^{32}_{k-1}a^{21}_qA^{13}_ka^{13}_q}{3}
                 {y_3a^{31}_{k-1}A^{12}_{q+1}a^{12}_kA^{23}_q} =
    e_{3\gamma}.
\end{split}
\end{equation}
By Corollary \ref{act_runner_map} we have
\begin{equation} \label{f_line2_comp}
 \begin{split}
       & \OperatorFi {y_2a^{21}_q}{3}{y_3A^{12}_{q+1}}, \\
       & \OperatorFi {A^{32}_{k-1}}{3}{y_3(y_1 + y_2)a^{31}_{k-1}}, \\
       & \OperatorFi {A^{13}_k}{3}{y_3(y_1 + x_2)a^{12}_k}, \\
       & \OperatorFi {a^{13}_q}{3}{y_3A^{23}_q}.    \\
 \end{split}
\end{equation}
  Since  $a^{13}_q, A^{32}_{k-1}$ are $\varphi_3-$homomorphic
  and $A^{13}_k$ is $(\varphi_3, y_2)-$homomorphic
  (Theorem \ref{th_homomorhism}), and
  $y_3{A}^{12}_{q+1}\subseteq y_3(y_1 + x_2)$,
  we see that (\ref{f_line2}) comes out from (\ref{f_line2_comp}).
  As in the previous case, we have the same result for $q=1$.
\qedsymbol

\underline{\text{(G) Line 1}}. It will be shown that
\begin{equation} \label{f_line3}
\begin{split}
 & \OperatorFi{g_{\gamma0}}{3}{g_{3\gamma0}}, \\
 & g_{\gamma0} = \OperatorFi{e_\gamma(x_1+a^{32}_qA^{32}_k)}{3}
                 {e_{3\gamma}(y_2y_3+A^{12}_qa^{31}_k)} =
   g_{3\gamma0}.
\end{split}
\end{equation}
By (\ref{f_line2}) we have $e_\gamma
\stackrel{\varphi_3}{\longmapsto}
 e_{3\gamma}$. By basic relations (\ref{basic_eq}) and atomic
multiplicativity (\ref{cor_mul}) we have
$$
   \psi_3(x_1+a^{32}_qA^{32}_k) = \psi_3(x_1) +
   \psi_3(a^{32}_q)\varphi_3(A^{32}_k).
$$
According to Proposition \ref{action_psi} (the joint map actions)
we obtain
$$
  \varphi_3(x_1+a^{32}_qA^{32}_k) =
  \nu^0(y_2y_3) + \nu^0(y_3A^{12}_q(y_3(y_1 + y_2)a^{31}_k)) =
  \nu^0(y_2y_3 + y_3A^{12}_qa^{31}_k).
$$
Therefore by Proposition \ref{phi_and_psi} (the connection) we see
that
$x_1+a^{32}_qA^{32}_k$ is $\varphi_3-$homomorphic and   \\
$$
  \OperatorFi{x_1+a^{32}_qA^{32}_k}{3}{y_2y_3 + y_3A^{12}_qa^{31}_k} =
    y_3(y_2y_3 + A^{12}_qa^{31}_k).
$$
By (\ref{f_line2}) $e_{3\gamma} \subseteq  y_3$, hence we get
(\ref{f_line3}). \qedsymbol

\subsection{Case \OperatorFi{z_{3\gamma}}{1}{z_{(13)\gamma}}}
   \label{case_line_2_3}
~\\

\underline{\text{(F) Line 2}}. We will show that
\begin{equation} \label{f_line4}
\begin{split}
 & \OperatorFi{f_{3\gamma}}{1}{f_{(13)\gamma}}, \\
 & f_{3\gamma} =
   \OperatorFi{y_3(x_1 + x_2)A^{23}_qa^{31}_{k-1}}{1}
              {y_3y_1a^{32}_{q+1}A^{21}_{k-1}}.
               =
   f_{(13)\gamma}
\end{split}
\end{equation}
 By Corollary \ref{act_runner_map}
 \begin{equation} \label{f_line4_comp}
 \begin{split}
        & \OperatorFi{x_1 + x_2}{1}{y_1y_3} , \\
        & \OperatorFi{y_3A^{23}_q} {1} {y_1a^{32}_{q+1}} , \\
        & \OperatorFi{a^{31}_{k-1}} {1} {y_1A^{21}_{k-1}}.
 \end{split}
 \end{equation}
 As above, $a^{31}_{k-1}$ and $x_1 + x_2$ are
 $\varphi_1-$homomorphic
 (Theorem  \ref{th_homomorhism})
  then ``intersection of relations''
 (\ref{f_line4_comp}) gives us relation (\ref{f_line4}). For
 $k-1=1$ we have
 \OperatorFi{a^{31}_{k-1}} {1} {y_1(y_2+y_3)A^{21}_{k-1}}. Since
 $a^{32}_{q+1} \subseteq y_2+y_3$, we have the same result.
\qedsymbol

\underline{\text{(E) Line 2}}. It will be shown that
\begin{equation} \label{f_line5}
\begin{split}
  & \OperatorFi{e_{3\gamma}}{1}{e_{(13)\gamma}}, \\
  &  e_{3\gamma} = \OperatorFi{y_3a^{31}_{k-1}A^{12}_{q+1}a^{12}_kA^{23}_q}{1}
                 {y_1A^{21}_{k-1}a^{13}_{q+1}A^{32}_ka^{32}_{q+1}} =
    e_{(13)\gamma}.
\end{split}
\end{equation}
By Corollary \ref{act_runner_map} we have
\begin{equation} \label{f_line5_comp}
 \begin{split}
       & \OperatorFi {y_3A^{23}_q}{1}{y_1a^{32}_{q+1}}, \\
       & \OperatorFi {a^{31}_{k-1}}{1}{y_1A^{21}_{k-1}}, \\
       & \OperatorFi {a^{12}_k}{1}{y_1A^{32}_k} \subseteq y_1(y_2 + y_3), \\
       & \OperatorFi {A^{12}_{q+1}}{1}{y_1(y_2 +
               y_3)a^{13}_{q+1}}.
 \end{split}
\end{equation}
  The atomic elements $a^{13}_{q-1}, a^{12}_k$ and $A^{12}_{q+1}$ are
  $\varphi_1-$homomorphic (Theorem \ref{th_homomorhism}). Then
  (\ref{f_line5}) follows from (\ref{f_line5_comp}). As in the
  previous case for $k-1=1$ we have the same result.
\qedsymbol

\underline{\text{(G) Line 2}}. It will be shown that
\begin{equation} \label{f_line6}
\begin{split}
  & \OperatorFi{g_{3\gamma0}}{1}{g_{(13)\gamma0}}, \\
 &  g_{3\gamma} = \OperatorFi{e_{3\gamma}(y_2y_3 + A^{12}_qa^{31}_k)}{1}
                 {e_{(13)\gamma}(x_3 + a^{21}_{q+1}A^{21}_k)} =
    g_{(13)\gamma}.
\end{split}
\end{equation}
Relation (\ref{f_line6}) follows from
\begin{equation} \label{f_line6_comp}
 \begin{split}
       & \OperatorFi {y_2y_3}{1}{y_1(x_2 + x_3)},   \\
       & \OperatorFi {A^{12}_q}{1}{y_1(y_2 + y_3)a^{13}_q}, \\
       & \OperatorFi {a^{31}_k}{1}{y_1A^{21}_k}
 \end{split}
\end{equation}
since
\begin{equation*}
\begin{split}
  & y_1(x_2 + x_3) + y_1(y_2 + y_3)a^{13}_qA^{21}_k =
    y_1(y_2 + y_3)(x_2 + x_3 + y_1a^{13}_qA^{21}_k) = \\
  & y_1(y_2 + y_3)(A^{21}_k(x_2 + y_1a^{13}_q) + x_3) =
    y_1(y_2 + y_3)(A^{21}_ka^{21}_{q+1} + x_3)
\end{split}
\end{equation*}
   and
   $e_{(13)\gamma} \subseteq y_1(y_2 + y_3)$.
   Again, as in two previous cases if $k=1$ for $a^{31}_k$, then
 $$
   \OperatorFi {a^{31}_k}{1}{y_1(y_2+y_3)A^{21}_k}
 $$
  with the same result.
\qedsymbol

We will omit the case \OperatorFi
{a^{ij}_1}{j}{y_j(y_i+y_k)A^{kj}_1} (Corollary
\ref{act_runner_map}, (2), $k=1$) for all further cases of the
theorem.
 As above, for the cases of Lines 9 -- 14 it suffices to prove only
\begin{equation} \label{line9_10}
\OperatorFi{z_\beta}{3}{z_{3\beta}} \hspace{3mm} and \hspace{3mm}
\OperatorFi{z_{3\beta}}{i}{z_{(13)\beta}},
\end{equation}
because the following relation
$$
   \OperatorFiFj{z_{(13)\beta}}{2}{z_{(213)\beta}}{3}{z_{3(213)\beta}}
   \stackrel{\varphi_1}{\longmapsto}
   \OperatorFi{z_{13(213)\beta}}{2}{z_{(213)(213)\beta}}
$$
follows from (\ref{line9_10}) as the permutation $1
\longrightarrow 3 \longrightarrow 2 \longrightarrow 1$ is applied
to Lines 9 -- 11 of Table \ref{table_adm_elems}. We will omit
references to the $\varphi_i-$homomorphic Theorem (Theorem
\ref{th_homomorhism}) and ``intersection of relations'' in
\S\S\hspace{1mm}\ref{case_line_1_2}-\ref{case_line_2_3}.

\subsection{Case \OperatorFi{z_\beta}{3}{z_{3\beta}}}
~\\

\underline{\text{(F) Line 9}}.
 It should be shown that
\begin{equation} \label{f_line9_1}
\begin{split}
  &
  \OperatorFi{f_\beta}{3}{f_{3\beta}}, \\
  & f_\beta =
   \OperatorFi{y_2y_3a^{21}_qA^{13}_k}{3}
              {x_3A^{12}_{q+1}a^{12}_k} =
   f_{3\beta}.
\end{split}
\end{equation}
The following relations are true (Corollary \ref{act_runner_map}):
 \begin{equation} \label{f_line9_1_comp}
 \begin{split}
        & \OperatorFi{y_3} {3} {x_3(y_1 + y_2)}, \\
        & \OperatorFi{y_2a^{21}_q} {3} {y_3A^{12}_{q+1}} \subseteq x_2 + y_1, \\
        & \OperatorFi{y_2A^{13}_k} {3} {y_3(x_2 + y_1)a^{12}_k}.
 \end{split}
 \end{equation}
 Then \ref{f_line9_1} follows from \ref{f_line9_1_comp}
\qedsymbol \vspace{2mm}

\underline{\text{(E) Line 9}}.  We want to prove that
\begin{equation} \label{f_line9_2}
\begin{split}
  & \OperatorFi{e_\beta}{3}{e_{3\beta}}, \\
 &  e_\beta =
   \OperatorFi{y_2A^{32}_ka^{21}_qA^{13}_ka^{13}_{q+1}}{3}
              {y_3a^{31}_kA^{12}_{q+1}a^{12}_k}A^{23}_{q+1} =
   e_{3\beta}.
\end{split}
\end{equation}
 It follows from the following relations (Corollary
 \ref{act_runner_map}):
 \begin{equation} \label{f_line9_2_comp}
 \begin{split}
        & \OperatorFi{A^{32}_k} {3} {y_3(y_1 + y_2)a^{31}_k}, \\
        & \OperatorFi{y_2a^{21}_q} {3} {y_3A^{12}_{q+1}} \subseteq x_2 + y_1, \\
        & \OperatorFi{y_2A^{13}_k} {3} {y_3(x_2 + y_1)a^{12}_k}, \\
        & \OperatorFi{a^{13}_{q+1}} {3} {y_3A^{23}_{q+1}}.
\qed \vspace{2mm}
 \end{split}
 \end{equation}

\underline{\text{(G) Line 9}}. We need to show
\begin{equation} \label{f_line9_3}
\begin{split}
  & \OperatorFi{g_\beta0}{3}{f_{3\beta0}}, \\
 & g_{\beta0} =
 \OperatorFi{e_\beta(x_2 + a^{13}_qA^{13}_{k+1})}{3}
              {e_{3\beta}(y_1y_3 + A^{23}_qa^{12}_{k+1})} =
   g_{3\beta0}.
\end{split}
\end{equation}
We have
$$
   e_\beta(x_2 + a^{13}_qA^{13}_{k+1}) = e_\beta(x_2 + y_2a^{13}_qA^{13}_{k+1})
$$
and
$$
  e_{3\beta}(y_1y_3 + A^{23}_qa^{12}_{k+1}) =
     e_{3\beta}(y_1y_3 + y_3(x_2 + y_1)A^{23}_qa^{12}_{k+1}),
$$
and the following relations (Corollary \ref{act_runner_map}) are
true:
\begin{equation}
  \label{f_line9_3_comp}
 \begin{split}
        & \OperatorFi{x_2} {3} {y_1y_3}, \\
        & \OperatorFi{a^{13}_q} {3} {y_3A^{23}_q}, \\
        & \OperatorFi{y_2A^{13}_{k+1}} {3} {y_3(x_2 + y_1)a^{12}_{k+1}}
               \subseteq A^{12}_1 \subseteq e_{3\beta}.
 \end{split}
\end{equation}
 Then relation (\ref{f_line9_3}) follows from the (\ref{f_line9_3_comp}).
\qedsymbol

\subsection{Case \OperatorFi{z_{3\beta}}{1}{z_{(13)\beta}}}
~\\

\underline{\text{(F) Line 10}}. It will be shown that
\begin{equation} \label{f_line10_1}
\begin{split}
  & \OperatorFi{f_{3\beta}}{1}{f_{(13)\beta}}, \\
 &  f_{3\beta} =
   \OperatorFi {x_3A^{12}_{q+1}a^{12}_k}{1}
               {y_1y_2a^{13}_{q+1}A^{32}_k}  =
   f_{(13)\beta}.
\end{split}
\end{equation}
The following relations are true (Corollary \ref{act_runner_map}):
 \begin{equation} \label{f_line10_1_comp}
 \begin{split}
        & \OperatorFi{x_3} {1} {y_1y_2}, \\
        & \OperatorFi{A^{12}_{q+1}} {1} {y_1(y_2 + y_3)A^{13}_{q+1}}, \\
        & \OperatorFi{a^{12}_k} {1} {A^{32}_k} \subseteq y_2 + y_3.
 \end{split}
 \end{equation}
 Then (\ref{f_line10_1}) follows from (\ref{f_line10_1_comp}).
\qedsymbol

\underline{\text{(E) Line 10}}. It will be shown that
\begin{equation} \label{f_line10_2}
\begin{split}
  & \OperatorFi{e_{3\beta}}{1}{e_{(13)\beta}}, \\
 &  e_{3\beta} =
   \OperatorFi {y_3a^{31}_kA^{12}_{q+1}a^{12}_kA^{23}_{q+1}} {1}
               {y_1A^{21}_ka^{13}_{q+1}A^{32}_ka^{32}_{q+2}}  =
   e_{(13)\beta}.
\end{split}
\end{equation}
The following relations are true (Corollary \ref{act_runner_map})
 \begin{equation} \label{f_line10_2_comp}
 \begin{split}
        & \OperatorFi{a^{31}_k} {1} {y_1A^{21}_k}, \\
        & \OperatorFi{A^{12}_{q+1}} {1} {y_1(y_2 + y_3)A^{13}_{q+1}}, \\
        & \OperatorFi{a^{12}_k} {1} {A^{32}_k} \subseteq y_2 + y_3, \\
        & \OperatorFi{y_3A^{23}_{q+1}}{1}{y_1a^{32}_{q+2}}.
 \end{split}
 \end{equation}
 Then (\ref{f_line10_2}) follows from (\ref{f_line10_2_comp}).
\qedsymbol

\underline{\text{(G) Line 10}}. It will be shown that
\begin{equation} \label{f_line10_3}
\begin{split}
  & \OperatorFi{g_{3\beta0}}{1}{g_{(13)\beta0}}, \\
 &  g_{3\beta0} =
   \OperatorFi {e_{3\beta}(y_1y_3 + A^{23}_qa^{12}_{k+1})} {1}
               {e_{(13)\beta}(x_1 + a^{32}_{q+1}A^{32}_{k+1})}  =
   g_{(13)\beta0}.
\end{split}
\end{equation}
  We have
$$
  e_{3\beta}(y_1y_3 + A^{23}_qa^{12}_{k+1}) =
     e_{3\beta}(y_1y_3 + y_3A^{23}_qa^{12}_{k+1})
$$
and
$$
    e_{(13)\beta}(x_1 + a^{32}_{q+1}A^{32}_{k+1}) =
    e_{(13)\beta}(x_1(x_3 + y_2) + y_1a^{32}_{q+1}A^{32}_{k+1})),
$$
     so the following relations (Corollary \ref{act_runner_map})
     are true:
 \begin{equation} \label{f_line10_3_comp}
 \begin{split}
        & \OperatorFi{y_1y_3} {1} {x_1(y_2 + y_3)y_1(x_3 + y_2)}
                 = x_1(x_3 + y_2), \\
        & \OperatorFi{y_3A^{23}_q}{1}{y_1a^{32}_{q+1}}
                  \subseteq x_3 + y_2, \\
        & \OperatorFi{a^{12}_{k+1}} {1} {y_1A^{32}_{k+1}}.
 \end{split}
 \end{equation}
 Then (\ref{f_line10_3}) follows from (\ref{f_line10_3_comp}).
\qedsymbol

Now, we will prove
\begin{equation}
 \begin{split}
   &  \OperatorFiFj{\underline{\text{Line 1}}(p = 0)}
       {1}{\underline{\text{Line 7}}}{2}{\underline{\text{Line 9}}(p = 0)}
      \stackrel{\varphi_1}{\longmapsto}
      \OperatorFi{\underline{\text{Line 8}}}{2}{\underline{\text{Line 1}}(p=0)}, \\
   &  \OperatorFiFj{z_{(21)^{2k}}}{1}{z_{1(21)^{2k}}}{2}{z_{(21)^{2k+1}}}
      \stackrel{\varphi_1}{\longmapsto}
      \OperatorFi{z_{1(21)^{2k+1}}}{2}{z_{(21)^{2k+2}}} \hspace{2mm}
      {\rm for} \hspace{2mm} z_\alpha = f_\alpha, e_\alpha, g_{\alpha0}.
 \end{split}
\end{equation}

\subsection{Case \OperatorFi {z_{(21)^{2k}}} {1} {z_{1(21)^{2k}}}}
~\\

\underline{\text{(F)} \OperatorFi{\text{ Line
1}(p=0)}{1}{\text{Line 7}}}. We will prove that
 \begin{equation}  \label{C1_f}
 \begin{split}
   & \OperatorFi {f_{(21)^{2k}}} {1} {f_{1(21)^{2k}}}, \\
   &  {f_{(21)^{2k}}} =
     \OperatorFi {y_1y_2a^{13}_kA^{32}_{k-1}} {1} {x_1A^{23}_ka^{23}_k} =
     {f_{1(21)^{2k}}}.
\end{split}
\end{equation}
The following relations are true (Corollary \ref{act_runner_map}):
 \begin{equation} \label{C1_f_comp}
 \begin{split}
        & \OperatorFi{y_1} {1} {x_1(y_2+y_3)}, \\
        & \OperatorFi{a^{13}_k} {1} {y_1A^{23}_k} \subseteq y_2+y_3, \\
        & \OperatorFi{y_2A^{32}_k}{1}{y_1A^{23}_k}.
 \end{split}
 \end{equation}
 Then (\ref{C1_f}) follows from (\ref{C1_f_comp}).
\qedsymbol

\underline{\text{(E)} \OperatorFi{\text{ Line 1
}(p=0)}{1}{\text{Line 7}}}. We will prove that
 \begin{equation}  \label{C2_f}
 \begin{split}
 & \OperatorFi {e_{(21)^{2k}}} {1} {e_{1(21)^{2k}}}, \\
 & {e_{(21)^{2k}}} =
     \OperatorFi {y_2A^{32}_{k-1}a^{21}_kA^{13}_ka^{13}_k} {1}
         {y_1A^{31}_ka^{12}_kA^{23}_ka^{23}_k} =
     {e_{1(21)^{2k}}}.
 \end{split}
 \end{equation}
The following relations are true (Corollary \ref{act_runner_map}):
 \begin{equation} \label{C2_f_comp}
 \begin{split}
        & \OperatorFi{y_2A^{32}_{k-1}} {1} {y_1a^{23}_k}, \\
        & \OperatorFi{a^{21}_k} {1} {y_1A^{31}_k},  \\
        & \OperatorFi{A^{13}_k} {1} {y_1(y_2+y_3)a^{12}_k},  \\
        & \OperatorFi{a^{13}_k}{1}{y_1A^{23}_k} \subseteq y_2+y_3.
 \end{split}
 \end{equation}
 Then (\ref{C2_f}) follows from (\ref{C2_f_comp}).
\qedsymbol

\underline{\text{(G)} \OperatorFi{\text{ Line 1
}(p=0)}{1}{\text{Line 7}}}. We need to show
 \begin{equation}  \label{C3_g}
 \begin{split}
 & \OperatorFi {g_{(21)^{2k}0}} {1} {g_{1(21)^{2k}0}}, \\
 & {g_{(21)^{2k}0}} =
   \OperatorFi {e_{(21)^{2k}}(x_1 + a^{32}_kA^{32}_k)} {1}
         {e_{1(21)^{2k}}(y_2a^{21}_k+y_3a^{31}_k)} =
     {g_{1(21)^{2k}0}}.
 \end{split}
 \end{equation}
Since $A^{32}_k = y_3 + x_2A^{13}_{k-1}$ and $x_2 \subseteq
a^{32}_k$, we see that $a^{32}_kA^{32}_k = x_2A^{13}_{k-1} +
y_3a^{32}_k$  and
$$
  x_1 + a^{32}_kA^{32}_k = (x_1+y_3a^{32}_k) + x_2A^{13}_{k-1} =
  a^{13}_{k+1} + x_2A^{13}_{k-1}.
$$
The transformation is well-defined because $k > 0$ (see Table
 \ref{table_adm_elems}). Further, since $e_{(21)^{2k}} \subseteq
A^{13}_{k-1}$ (by Table \ref{table_adm_elems}), we have
\begin{equation} \label{C3_transform}
e_{(21)^{2k}}(x_1 + a^{32}_kA^{32}_k) = e_{(21)^{2k}}(a^{13}_{k+1}
+ x_2A^{13}_{k-1}) = e_{(21)^{2k}}(x_2 +
a^{13}_{k+1}A^{13}_{k-1}).
\end{equation}
Apply now $\varphi_1$.
The polynomial $x_2 +
a^{13}_{k+1}A^{13}_{k-1}$ is $\varphi_1-$homomorphic because $x_1
\subseteq a^{13}_{k+1}A^{13}_{k-1}$. We get
\begin{equation}
 \label{C3_g_1}
 \begin{split}
    &  \OperatorFi {e_{(21)^{2k}}} {1} {e_{1(21)^{2k}}}, \\
    &  \OperatorFi {x_2} {1} {y_1y_3},                      \\
    &  \OperatorFi {a^{13}_{k+1}} {1} y_1{A^{23}_{k+1}} \subseteq y_2 + y_3, \\
    &  \OperatorFi {A^{13}_{k-1}} {1} {y_1(y_2+y_3)a^{12}_{k-1}}.
 \end{split}
\end{equation}
Since $e_{1(21)^{2k+1}} \subseteq y_1$, we get
\begin{equation} \label{C3_g_2}
     \OperatorFi {e_{(21)^{2k}}(x_1 + a^{32}_kA^{32}_k)} {1}
        {e_{1(21)^{2k+1}}(y_1y_3 + y_1A^{23}_{k+1}a^{12}_{k-1})} =
        {e_{1(21)^{2k+1}}(y_1y_3 + A^{23}_{k+1}a^{12}_{k-1})}.
\end{equation}
Since $e_{1(21)^{2k+1}} \subseteq a^{12}_{k-1}$, we have
\begin{equation*}
\begin{split}
& e_{1(21)^{2k+1}}(y_1y_3 + A^{23}_{k+1}a^{12}_{k-1}) =
  e_{1(21)^{2k+1}}(A^{23}_{k+1} + y_1y_3a^{12}_{k-1}) =  \\
& e_{1(21)^{2k+1}}(y_2 + x_3A^{12}_k + y_1y_3a^{12}_{k-1}) =
  e_{1(21)^{2k+1}}(y_2 + y_3A^{12}_k(x_3 + y_1a^{12}_{k-1})) =  \\
& e_{1(21)^{2k+1}}(y_2 + y_3A^{12}_ka^{31}_k).
\end{split}
\end{equation*}
  Again, by Table \ref{table_adm_elems}
$e_{1(21)^{2k+1}} \subseteq y_1 \subseteq A^{12}_k$, therefore the
last polynomial is equal to
$$
  e_{1(21)^{2k+1}}(y_2A^{12}_k + y_3a^{31}_k),
$$
 i.e.,
\begin{equation} \label{C3_g_3}
        e_{1(21)^{2k+1}}(y_1y_3 + A^{23}_{k+1}a^{12}_{k-1}) =
        e_{1(21)^{2k+1}}(y_2A^{12}_k + y_3a^{31}_k).
\end{equation}
(\ref{C3_g}) follows from (\ref{C3_g_3}) and from the next lemma.
\begin{lemma}
  The following relation takes place
  \begin{equation} \label{C3_g_4}
        e_{1(21)^{2k+1}}(y_2A^{12}_k + y_3a^{31}_k) =
        e_{1(21)^{2k+1}}(y_2a^{21}_k+y_3a^{31}_k).
  \end{equation}
\end{lemma}
\PerfProof Indeed,
\begin{equation*}
\begin{split}
 & e_{1(21)^{2k+1}}(y_2A^{12}_k + y_3a^{31}_k) =
   e_{1(21)^{2k+1}}(y_2(y_1 + x_2A^{31}_{k-1}) + y_3a^{31}_k) =  \\
 & e_{1(21)^{2k+1}}(y_2y_1 + x_2A^{31}_{k-1} + y_3a^{31}_k).
\end{split}
\end{equation*}
Since $e_{1(21)^{2k+1}} \subseteq A^{31}_{k-1}$, we have
\begin{equation*}
\begin{split}
 & e_{1(21)^{2k+1}}(y_2A^{12}_k + y_3a^{31}_k) =
   e_{1(21)^{2k+1}}(x_2 + A^{31}_{k-1}(y_1y_2 + y_3a^{31}_k)) = \\
 & e_{1(21)^{2k+1}}(x_2 + y_1y_2A^{31}_{k-1} + y_3a^{31}_k).
\end{split}
\end{equation*}
Further, by basic relations (Table \ref{table_atomic}, 4.1) we
have
\begin{equation*}
\begin{split}
 & e_{1(21)^{2k+1}}(y_2A^{12}_k + y_3a^{31}_k) =
   e_{1(21)^{2k+1}}(x_2 + y_1y_2a^{13}_{k-1} + y_3a^{31}_k) = \\
 & e_{1(21)^{2k+1}}(y_2(x_2 + y_2a^{13}_{k-1}) + y_3a^{31}_k) =
   e_{1(21)^{2k+1}}(y_2a^{21}_k + y_3a^{31}_k).
\end{split}
\end{equation*}
Thus, the lemma is proven, (\ref{C3_g_4}) and so (\ref{C3_g}) are
true. \qedsymbol

Action $\varphi_2$ transforms Line 7 to Line 9 with $p = 0$. For
$z_\alpha = f_\alpha, e_\alpha, g_{\alpha0}$, we are going to
prove

\subsection{Case \OperatorFi {z_{1(21)^{2k}}} {2}
{z_{(21)^{2k+1}}}}
~\\

 \underline{\text{(F)} \OperatorFi{\text{
Line 7 }}{2}{\text{Line 9 }(p=0)}}. We want to prove that
\begin{equation}  \label{D1_f}
\begin{split}
 & \OperatorFi{f_{1(21)^{2k}}}{2}{f_{(21)^{2k+1}}}, \\
 & {f_{1(21)^{2k}}} =
     \OperatorFi {x_1A^{23}_ka^{23}_k} {2} {y_2y_3a^{21}_kA^{13}_k} =
     {f_{(21)^{2k+1}}}.
\end{split}
\end{equation}
The following relations are true (Corollary \ref{act_runner_map}):
\begin{equation} \label{D1_f_comp}
\begin{split}
        & \OperatorFi{x_1} {2} {y_2y_3}, \\
        & \OperatorFi{A^{23}_k} {2} {y_2(y_1 + y_3)a^{21}_k}, \\
        & \OperatorFi{a^{23}_k}{2}{y_2A^{13}_k}  \subseteq y_1 + y_3.
\end{split}
\end{equation}
 Then (\ref{D1_f}) follows from (\ref{D1_f_comp}).
\qedsymbol

\underline{\text{(E)} \OperatorFi{\text{ Line 7 }}{2}{\text{Line 9
}(p=0)}}. We want to prove that
 \begin{equation}  \label{D2_e}
 \begin{split}
   & \OperatorFi{e_{1(21)^{2k}}}{2}{e_{(21)^{2k+1}}}, \\
   & {e_{1(21)^{2k}}} =
     \OperatorFi {y_1A^{31}_ka^{12}_kA^{23}_ka^{23}_k} {2}
                 {y_2A^{32}_ka^{21}_kA^{13}_k}a^{13}_{k+1} =
     {e_{(21)^{2k+1}}}.
 \end{split}
 \end{equation}
The following relations are true (Corollary \ref{act_runner_map}):
 \begin{equation} \label{D2_e_comp}
 \begin{split}
        & \OperatorFi{y_1A^{31}_k} {2} {y_1a^{13}_{k+1}}, \\
        & \OperatorFi{a^{12}_k} {2} {y_2A^{32}_k}, \\
        & \OperatorFi{A^{23}_k} {2} {y_2(y_1 + y_3)a^{21}_k}, \\
        & \OperatorFi{a^{23}_k}{2}{y_2A^{13}_k}  \subseteq y_1 + y_3.
 \end{split}
 \end{equation}
Thus, (\ref{D2_e}) follows from (\ref{D2_e_comp}). \qedsymbol

\underline{\text{(G)} \OperatorFi{\text{ Line 7 }}{2}{\text{Line 9
}(p=0)}}. It will be shown that
\begin{equation}  \label{D3_g}
 \begin{split}
   & \OperatorFi{g_{1(21)^{2k}0}}{2}{g_{(21)^{2k+1}0}}, \\
   & {g_{1(21)^{2k}0}} =
     \OperatorFi {e_{1(21)^{2k+1}}(y_2a^{21}_k+y_3a^{31}_k)}
         {2} {e_{(21)^{2k+1}}(x_2 + a^{13}_kA^{13}_k)}  =
     {g_{(21)^{2k+1}0}}.
 \end{split}
 \end{equation}
 The following relations are true:
 \begin{equation}  \label{D3_g_1}
 \begin{split}
     & \OperatorFi {y_2a^{21}} {2} {x_2(y_1 + y_3)a^{13}y_2kA^{31}_k} =
                x_2A^{31}_k, \\
     & \OperatorFi {y_3a^{31}} {2} {y_2A^{13}_k}.
 \end{split}
 \end{equation}
 Therefore,
 \begin{equation}  \label{D3_g_2}
      \OperatorFi{y_2a^{21}_k+y_3a^{31}_k} {2}
         {y_2(x_2A^{31}_k + A^{13}_k)} =
           y_2(x_2a^{13}_k + A^{13}_k) = x_2a^{13}_k + y_2A^{13}_k.
 \end{equation}
 Since
 \begin{equation}  \label{D3_g_3}
      \OperatorFi{e_{1(21)^{2k+1}}} {2} {e_{(21)^{2k+1}}}
           \subseteq a^{13}_{k+1} \subseteq a^{13}_k,
 \end{equation}
 we see by (\ref{D3_g_2}) that
 \begin{equation}  \label{D3_g_4}
      \OperatorFi{e_{1(21)^{2k+1}}(y_2a^{21}_k+y_3a^{31}_k)} {2}
          {e_{(21)^{2k+1}}(x_2a^{13}_k + y_2A^{13}_k)} =
           e_{(21)^{2k+1}}(x_2 + y_2a^{13}_kA^{13}_k).
\qed \vspace{2mm}
 \end{equation}

Action $\varphi_1$ transforms Line 9 with $p = 0$ to Line 8 with
$p = 0$. For $z_\alpha = f_\alpha, e_\alpha, g_{\alpha0}$, we are
going to prove

\subsection{Case \OperatorFi {z_{(21)^{2k+1}}}{1}{z_{1(21)^{2k+1}}}}
~\\

 \underline{\text{(F)} \OperatorFi{\text{ Line 9
}(p=0)}{1}{\text{Line 8}}}. It will be shown that
\begin{equation}
\label{line9_8_f}
\begin{split}
 & \OperatorFi{f_{(21)^{2k+1}}}{1}{f_{1(21)^{2k+1}}}, \\
 & {f_{(21)^{2k+1}}} =
     \OperatorFi {y_2y_3a^{21}_kA^{13}_k} {1} {y_1(x_2 + x_3)A^{31}_ka^{12}_k} =
     {f_{1(21)^{2k+1}}}.
\end{split}
\end{equation}
The following relations are true (Corollary \ref{act_runner_map}):
\begin{equation} \label{line9_8_f_comp}
\begin{split}
        & \OperatorFi{y_2y_3} {1} {y_1(x_2 + x_3)} \subseteq y_2 + y_3, \\
        & \OperatorFi{a^{21}_k} {1} {y_1A^{31}_k}, \\
        & \OperatorFi{A^{13}_k}{1}{y_1(y_2+y_3)a^{12}_k}.
\end{split}
\end{equation}
 Then (\ref{line9_8_f}) follows from (\ref{line9_8_f_comp}).
\qedsymbol

\underline{\text{(E)} \OperatorFi{\text{ Line 9
}(p=0)}{1}{\text{Line 8}}}. Here, it should be proven that
 \begin{equation}  \label{line9_8_e}
 \begin{split}
 & \OperatorFi{e_{(21)^{2k+1}}}{1}{e_{1(21)^{2k+1}}}, \\
 & {e_{(21)^{2k+1}}} =
     \OperatorFi {y_2A^{32}_ka^{21}_kA^{13}_ka^{13}_{k+1}} {1}
                 {y_1A^{31}_ka^{12}_kA^{23}_{k+1}}a^{23}_{k+1} =
     {e_{1(21)^{2k+1}}}.
 \end{split}
 \end{equation}
The following relations are true (Corollary \ref{act_runner_map}):
 \begin{equation} \label{line9_8_e_comp}
 \begin{split}
        & \OperatorFi{y_2A^{32}_k} {1} {y_1a^{23}_{k+1}} \subseteq y_2 + y_3, \\
        & \OperatorFi{a^{21}_k} {1} {y_1A^{31}_k}, \\
        & \OperatorFi{A^{13}_k} {1} {y_1(y_2 + y_3)a^{12}_k}, \\
        & \OperatorFi{a^{13}_{k+1}}{1}{y_1A^{23}_{k+1}}.
 \end{split}
 \end{equation}
 Therefore, (\ref{line9_8_e}) follows from (\ref{line9_8_e_comp}).
\qedsymbol

\underline{\text{(G)} \OperatorFi{\text{ Line 9
}(p=0)}{1}{\text{Line 8}}}. We will show that
\begin{equation}  \label{line9_8_g}
 \begin{split}
   & \OperatorFi{g_{(21)^{2k+1}0}}{1}{g_{1(21)^{2k+1}0}}, \\
   & {g_{(21)^{2k+1}0}} =
     \OperatorFi {e_{(21)^{2k+1}}(x_2 + a^{13}_kA^{13}_{k+1})}
         {1} {e_{1(21)^{2k+1}}(y_1a^{12}_{k+1} + y_1a^{13}_{k+1})}  =
     {g_{1(21)^{2k+1}0}}.
 \end{split}
 \end{equation}
 First, from the relations
 \begin{equation}  \label{line9_8_g1}
 \begin{split}
    & \OperatorFi{e_{(21)^{2k+1}}}{1}{e_{1(21)^{2k+1}}} \subseteq y_1, \\
    & \OperatorFi{x_2}{1}{y_1y_3}, \\
    & \OperatorFi{a^{13}_k}{1}{y_1A^{23}_k} \subseteq y_2+y_3, \\
    & \OperatorFi{A^{13}_{k+1}}{1}{y_1(y_2+y_3)a^{12}_{k+1}}
 \end{split}
 \end{equation}
 we have
 \begin{equation}
  \label{line9_8_g2}
  \begin{split}
   & \OperatorFi {e_{(21)^{2k+1}}(x_2 + a^{13}_kA^{13}_{k+1})}
         {1} {e_{1(21)^{2k+1}}(y_1y_3 + y_1A^{23}_ka^{12}_{k+1})} = \\
   &  e_{1(21)^{2k+1}}(y_1y_3 + A^{23}_ka^{12}_{k+1}).
  \end{split}
 \end{equation}
 Transform the last polynomial from (\ref{line9_8_g2}):
 \begin{equation*}
 \begin{split}
   &  y_1y_3 + A^{23}_ka^{12}_{k+1} = y_1y_3 + A^{23}_k(x_1 + y_2a^{23}_k) = \\
   &  y_1y_3 + x_1A^{23}_k + y_2a^{23}_k =
      y_1(y_3 + x_1A^{23}_k) +  y_2a^{23}_k =
      y_1A^{31}_{k+1} + y_2a^{23}_k.
 \end{split}
 \end{equation*}
 Thus,
 \begin{equation}  \label{line9_8_g3}
 \OperatorFi {e_{(21)^{2k+1}}(x_2 + a^{13}_kA^{13}_{k+1})}
         {1} {e_{1(21)^{2k+1}}(y_1A^{31}_{k+1} + y_2a^{23}_k)}.
 \end{equation}
 Then (\ref{line9_8_g}) follows from the next lemma.
 \begin{lemma}
   We have
 \begin{equation}
    e_{1(21)^{2k+1}}(y_1A^{31}_{k+1} + y_2a^{23}_k) =
    e_{1(21)^{2k+1}}(y_1a^{12}_{k+1} + y_1a^{13}_{k+1}).
 \end{equation}
 \end{lemma}
 \PerfProof Indeed,
 \begin{equation*}
 \begin{split}
  &  e_{1(21)^{2k+1}}(y_1A^{31}_{k+1} + y_2a^{23}_k) =
     e_{1(21)^{2k+1}}(y_1(y_3 + x_1A^{23}_k) + y_2a^{23}_k) = \\
  &  e_{1(21)^{2k+1}}(y_1y_3 + x_1A^{23}_k + y_2a^{23}_k).
 \end{split}
 \end{equation*}
 Since $e_{1(21)^{2k+1}} \subseteq A^{23}_{k+1} \subseteq A^{23}_k$,
 it follows that
 \begin{equation*}
 \begin{split}
  & e_{1(21)^{2k+1}}(y_1A^{31}_{k+1} + y_2a^{23}_k) =
    e_{1(21)^{2k+1}}(x_1 + A^{23}_k(y_1y_3 + y_2a^{23}_k)) =  \\
  & e_{1(21)^{2k+1}}(x_1 + y_1y_3A^{23}_k + y_2a^{23}_k)).
 \end{split}
 \end{equation*}

 By basic relations of atomic elements (Table \ref{table_atomic}, 4.1)
 we have
 \begin{equation*}
 \begin{split}
  & e_{1(21)^{2k+1}}(x_1 + A^{23}_k(y_1y_3 + y_2a^{23}_k)) =
    e_{1(21)^{2k+1}}(x_1 + y_1y_3a^{32}_k + y_2a^{23}_k)) =  \\
  & e_{1(21)^{2k+1}}(y_1(x_1 + y_3a^{32}_k) + y_2a^{23}_k)) =
    e_{1(21)^{2k+1}}(y_1a^{13}_{k+1} + y_2a^{23}_k)).
 \end{split}
 \end{equation*}

 Now,
 \begin{equation*}
 \begin{split}
  & e_{1(21)^{2k+1}}(y_1a^{13}_{k+1} + y_2a^{23}_k)) =
    e_{1(21)^{2k+1}}(y_1a^{13}_{k+1} + x_1 + y_2a^{23}_k)) = \\
  & e_{1(21)^{2k+1}}(y_1a^{13}_{k+1} + a^{12}_{k+1}) =
    e_{1(21)^{2k+1}}(y_1a^{13}_{k+1} + y_1a^{12}_{k+1}).
 \end{split}
 \end{equation*}
 The lemma is proved.
\qedsymbol

Action $\varphi_2$ transforms Line 8 to the Line 1 with $p = 0$.
For $z_\alpha = f_\alpha, e_\alpha, g_{\alpha0}$ \\
we are going to prove

\subsection{Case \OperatorFi{z_{1(21)^{2k+1}}}{2}{z_{(21)^{2k+2}}}}
~\\

\underline{\text{(F)} \OperatorFi{\text{ Line 8 }}{2}{\text{Line 1
}(p=0)}}. We want to prove that
\begin{equation}
 \label{line8_1_f}
\begin{split}
 & \OperatorFi{f_{1(21)^{2k+1}}}{2}{f_{(21)^{2k+2}}}, \\
 & {f_{1(21)^{2k+1}}} =
     \OperatorFi {y_1(x_2 + x_3)A^{31}_ka^{12}_k} {1} {y_1y_2a^{13}_{k+1}A^{32}_k} =
     {f_{(21)^{2k+2}}}.
\end{split}
\end{equation}
The following relations are true (Corollary \ref{act_runner_map}):
\begin{equation} \label{line8_1_f_comp}
\begin{split}
        & \OperatorFi{x_2+x_3} {2} {y_1y_2},  \\
        & \OperatorFi{y_1A^{31}_k} {2} {y_2a^{13}_{k+1}}, \\
        & \OperatorFi{a^{12}_k}{2}{y_1A^{32}_k}.
\end{split}
\end{equation}
 Then (\ref{line8_1_f}) follows from (\ref{line8_1_f_comp}).
\qedsymbol

\underline{\text{(E)} \OperatorFi{\text{ Line 8 }}{2}{\text{Line 1
}(p=0)}}. We need to show
\begin{equation}  \label{line8_1_e}
 \begin{split}
   & \OperatorFi{e_{1(21)^{2k+1}}}{2}{e_{(21)^{2k+2}}}, \\
   & {e_{1(21)^{2k+1}}} =
     \OperatorFi {y_1A^{31}_ka^{12}_kA^{23}_{k+1}a^{23}_{k+1}} {2}
                 {y_2A^{32}_ka^{21}_{k+1}A^{13}_{k+1}}a^{13}_{k+1} =
     {e_{(21)^{2k+2}}}.
 \end{split}
 \end{equation}
The following relations are true (Corollary \ref{act_runner_map}):
 \begin{equation} \label{line8_1_e_comp}
 \begin{split}
        & \OperatorFi{y_1A^{31}_k} {2} {y_2a^{13}_{k+1}} \subseteq y_1 + y_3, \\
        & \OperatorFi{a^{12}_k} {2} {y_2A^{32}_k}, \\
        & \OperatorFi{A^{23}_{k+1}} {2} {y_2(y_1 + y_3)a^{21}_{k+1}}, \\
        & \OperatorFi{a^{23}_{k+1}}{2}{y_1A^{13}_{k+1}}.
 \end{split}
 \end{equation}
 Thus, (\ref{line8_1_e}) follows from (\ref{line8_1_e_comp}).
\qedsymbol

\underline{\text{(G)} \OperatorFi{\text{ Line 8 }}{2}{\text{Line 1
}(p=0)}}. It will be shown that
 \begin{equation}  \label{line8_1_g}
 \begin{split}
   & \OperatorFi{g_{1(21)^{2k+1}0}}{2}{g_{(21)^{2k+2}0}}, \\
   & {g_{1(21)^{2k+1}0}} =
     \OperatorFi {e_{1(21)^{2k+1}}(y_1a^{12}_{k+1} + y_1a^{13}_{k+1}))}
         {2} {e_{(21)^{2k+2}}(x_1 + a^{32}_{k+1}A^{32}_{k+2}})  =
     {g_{(21)^{2k+2}0}}.
 \end{split}
 \end{equation}
 First, since $e_{1(21)^{2k+1}} \subseteq y_1$, we have
 $$
   e_{1(21)^{2k+1}}(y_1a^{12}_{k+1} + y_1a^{13}_{k+1})) =
   e_{1(21)^{2k+1}}(x_1 + y_2a^{23}_k + y_1y_3a^{32}_k).
 $$
 Since  $x_2 \subseteq y_2a^{23}_k$ and the polynomial
 $x_1 + y_2a^{23}_k + y_1y_3a^{32}_k$
 meets the requirements of atomic multiplicativity
 (Corollary \ref{cor_mul}), we have
\begin{equation*}  \label{line8_1_g1}
    \OperatorFi
       {e_{1(21)^{2k+1}}(x_1 + y_2a^{23}_k + y_1y_3a^{32}_k)}{2}
{e_{(21)^{2k+2}}(y_2y_3 + x_2(y_1 + y_3)A^{13}_k +
y_2(x_1+x_3)A^{12}_k)}.
\end{equation*}
We transform the last polynomial. Since $e_{(21)^{2k+2}} \subseteq
y_2$, it follows that
\begin{equation*}
\begin{split}
 & e_{(21)^{2k+2}}(y_2y_3 + x_2(y_1 + y_3)A^{13}_k + y_2(x_1+x_3)A^{12}_k) = \\
 & e_{(21)^{2k+2}}(y_2y_3 + x_2A^{13}_k + (x_1+x_3)A^{12}_k) =  \\
 & e_{(21)^{2k+2}}(y_2y_3 + x_2A^{13}_k + x_1+x_3A^{12}_k) =
   e_{(21)^{2k+2}}(y_2y_3 + x_2A^{13}_k + x_1+x_3a^{12}_k).
\end{split}
\end{equation*}
Since $e_{(21)^{2k+2}} \subseteq$  $a^{12}_k$, we see by
permutation property (\ref{permutation1}) that
\begin{equation*}
\begin{split}
 & e_{(21)^{2k+2}}(y_2y_3 + x_2A^{13}_k + x_1+x_3a^{12}_k) =
   e_{(21)^{2k+2}}(x_1 + x_3 + a^{12}_k(y_2y_3 + x_2A^{13}_k)) = \\
 & e_{(21)^{2k+2}}(x_1 + x_3 + x_2A^{13}_k + y_2y_3a^{12}_k )) =
   e_{(21)^{2k+2}}(x_1 + y_3(x_3 + y_2a^{12}_k) + x_2A^{13}_k) = \\
 & e_{(21)^{2k+2}}(x_1 + y_3a^{32}_{k+1} + x_2A^{13}_k) =
   e_{(21)^{2k+2}}(x_1 + a^{32}_{k+1}(y_3 + x_2A^{13}_k) =
   e_{(21)^{2k+2}}(x_1 + a^{32}_{k+1}A^{32}_{k+1})),
\end{split}
\end{equation*}
which was to be proved. \qedsymbol

The proof of Theorem \ref{th_adm_classes} is finished. \qedsymbol

\section{A proof of the theorem on admissible elements in lattice $D^4$}
 \label{proof_adm_D4}

Every element $e_\alpha$ from Table \ref{table_adm_elem_D4} is the
intersection of three polynomials. In all cases, two of these
polynomials are $\varphi_i-$homomorphic under action $\psi_i$. The
third one is not $\varphi_i-$homomorphic, but satisfies the
following relation:
$$
   \OperatorFi{e_j{e}^{kl}_n}{i}{e_i{e}^{kl}_{n+1}},
$$
see Corollary \ref{act_runner_map_D4}, heading(3).

In this proof, we repeatedly use the basic properties of
admissible elements from Lemma \ref{homom_polynom_P},
\S\ref{basic_adm_D4}.

\subsection{Line $F21$, $\OperatorFi{e_{\alpha}}{3}{e_{3\alpha}}$,
  $\OperatorFi{f_{{\alpha}0}}{3}{f_{{3\alpha}0}}$}

 Here,
$$
  \alpha = (21)^t(41)^r(31)^s, e_\alpha = e_{(21)^t(41)^r(31)^s},
$$
and
\begin{equation}
 \label{Line_F21_act_3}
   e_\alpha =
   \OperatorFi{e_2a^{31}_{2s}a^{41}_{2r}a^{34}_{2t-1}}{3}
              {e_3a^{24}_{2s}a^{21}_{2t-1}a^{41}_{2r+1}} =
   e_{3\alpha},
\end{equation}
where $3\alpha$ is obtained by eq.(\ref{act_psi_4_F21}). By
Corollary \ref{act_runner_map_D4} we have
\begin{equation*}
 \begin{split}
  & \OperatorFi{a^{31}_{2s}}{3}{e_3{a}^{24}_{2s}}, \\
  & \OperatorFi{a^{34}_{2t-1}}{3}{e_3{a}^{21}_{2t-1}}, \\
  & \OperatorFi{e_2a^{41}_{2r}}{3}{e_3a^{41}_{2r+1}}. \\
 \end{split}
\end{equation*}
Since $a^{31}_{2s}$ and $a^{34}_{2t-1}$ are
$\varphi_3-$homomorphic (Corollary \ref{cor_mul_d4}), we get
(\ref{Line_F21_act_3}).
 The polynomial $e_{3\alpha}$ corresponds to Line $G31$
 under the substitution
\begin{equation}
  \label{subst_F21_act_3}
  \begin{split}
    & r \mapsto r+1, \\
    & t \mapsto t-1,
  \end{split}
\end{equation}
see relation (\ref{act_psi_4_F21}).
 For the case $f_{{\alpha}0}$, we have
$$
  e_\alpha(e_2{a}^{34}_{2t} + a^{14}_{2r+1}a^{13}_{2s-1}) =
  e_\alpha({a}^{34}_{2t} + e_2a^{14}_{2r+1}a^{13}_{2s-1}),
$$
and we select
$$
  P = {a}^{34}_{2t} + e_2a^{14}_{2r+1}a^{13}_{2s-1},
$$
see Lemma \ref{homom_polynom_P}.
 By Corollary \ref{act_runner_map_D4} we have
\begin{equation*}
 \begin{split}
  & \OperatorFi{a^{34}_{2t}}{3}{e_3{a}^{21}_{2t}}, \\
  & \OperatorFi{a^{13}_{2s-1}}{3}{e_3a^{42}_{2s-1}}. \\
  & \OperatorFi{e_2a^{14}_{2r+1}}{3}{e_3{a}^{14}_{2r+2}}, \\
 \end{split}
\end{equation*}
and we get
\begin{equation}
 \label{Line_F21_act_3_f}
   P =
   \OperatorFi{{a}^{34}_{2t} + e_2a^{14}_{2r+1}a^{13}_{2s-1}}{3}
              {e_3a^{21}_{2t} + e_3a^{14}_{2r+2}a^{42}_{2s-1}} =
   \tilde{P}.
\end{equation}
Since $a^{42}_{2s} \subseteq a^{42}_{2s-1}$, we have
\begin{equation*}
  \begin{split}
  & e_{3\alpha}\tilde{P} =
    e_3a^{24}_{2s}a^{21}_{2t-1}a^{41}_{2r+1}
    (e_3a^{21}_{2t} + e_3a^{14}_{2r+2}a^{42}_{2s-1}) = \vspace{2mm} \\
  & e_3a^{24}_{2s}a^{21}_{2t-1}a^{41}_{2r+1}
    (a^{14}_{2r+2} + e_3a^{21}_{2t}a^{42}_{2s-1}),
  \end{split}
\end{equation*}
 which corresponds to Line $G31$ with
 substitution (\ref{subst_F21_act_3}), see relation
 (\ref{act_psi_4_F21}).
\qedsymbol  \vspace{2mm}

\subsection{Line $F21$, $\OperatorFi{e_{\alpha}}{4}{e_{4\alpha}}$,
 $\OperatorFi{f_{{\alpha}0}}{4}{f_{{4\alpha}0}}$}
 \label{subs_line_21}
We have $\alpha = (21)^t(41)^r(31)^s$,
 $e_\alpha = e_{(21)^t(41)^r(31)^s}$,
\begin{equation}
 \label{Line_F21_act_4}
   e_\alpha =
   \OperatorFi{e_2a^{31}_{2s}a^{41}_{2r}a^{34}_{2t-1}}{4}
              {e_4a^{31}_{2s+1}a^{32}_{2r}a^{12}_{2t-1}} =
   e_{4\alpha},
\end{equation}
because by Corollary \ref{act_runner_map_D4} we have
\begin{equation*}
 \begin{split}
  & \OperatorFi{a^{41}_{2r}}{4}{e_4{a}^{32}_{2r}}, \\
  & \OperatorFi{a^{43}_{2t-1}}{4}{e_4{a}^{12}_{2t-1}}, \\
  & \OperatorFi{e_4a^{31}_{2s}}{4}{e_4a^{31}_{2s+1}}. \\
 \end{split}
\end{equation*}
 The polynomial $e_{4\alpha}$ corresponds to Line $G41$.
 For the case $f_{{\alpha}0}$,
 we select
$$
  P = {a}^{43}_{2t} + e_2a^{14}_{2r+1}a^{13}_{2s-1}.
$$
 We have
\begin{equation}
 \label{Line_F21_act_4_f}
  P =
   \OperatorFi{{a}^{43}_{2t} + e_2a^{14}_{2r+1}a^{13}_{2s-1}}{4}
              {e_4a^{21}_{2t} + e_4a^{23}_{2r+1}a^{13}_{2s}} =
   \tilde{P},
\end{equation}
because by Corollary \ref{act_runner_map_D4} we have
\begin{equation*}
 \begin{split}
  & \OperatorFi{a^{41}_{2r+1}}{4}{e_4{a}^{32}_{2r+1}}, \\
  & \OperatorFi{a^{43}_{2t}}{4}{e_4{a}^{12}_{2t}}, \\
  & \OperatorFi{e_4a^{31}_{2s-1}}{4}{e_4a^{31}_{2s}}. \\
 \end{split}
\end{equation*}
 The polynomial $\tilde{P}$ corresponds to Line $G41$.
\qedsymbol \vspace{2mm}

\subsection{Line $F21$, $\OperatorFi{e_{\alpha}}{1}{e_{1\alpha}}$,
 $\OperatorFi{f_{{\alpha}0}}{1}{f_{{1\alpha}0}}$}
The admissible sequence $\alpha$ and polynomial $e_\alpha$ are
given as in \S\ref{subs_line_21}. Here,
\begin{equation}
 \label{Line_F21_act_1}
   e_\alpha =
   \OperatorFi{e_2a^{31}_{2s}a^{41}_{2r}a^{34}_{2t-1}}{1}
              {e_1a^{42}_{2s}a^{23}_{2r}a^{34}_{2t}} =
   e_{1\alpha},
\end{equation}
because by Corollary \ref{act_runner_map_D4} we have
\begin{equation*}
 \begin{split}
  & \OperatorFi{a^{31}_{2s}}{1}{e_1{a}^{42}_{2s}}, \\
  & \OperatorFi{a^{41}_{2r}}{1}{e_1{a}^{23}_{2r}}, \\
  & \OperatorFi{e_2a^{34}_{2t-1}}{1}{e_1a^{34}_{2t}}. \\
 \end{split}
\end{equation*}
 The polynomial $e_{1\alpha}$ corresponds to Line $G11$.
 For the case $f_{{\alpha}0}$,
 we select
$$
  P = e_2{a}^{43}_{2t} + a^{14}_{2r+1}a^{13}_{2s-1},
$$
 and we have
\begin{equation}
 \label{Line_F21_act_1_f}
  P =
   \OperatorFi{e_2{a}^{43}_{2t} + a^{14}_{2r+1}a^{13}_{2s-1}}{1}
              {e_1a^{43}_{2t+1} + e_1a^{23}_{2r+1}a^{24}_{2s-1}} =
   \tilde{P},
\end{equation}
because by Corollary \ref{act_runner_map_D4} we have
\begin{equation*}
 \begin{split}
  & \OperatorFi{a^{41}_{2r+1}}{1}{e_1{a}^{32}_{2r+1}}, \\
  & \OperatorFi{a^{31}_{2s-1}}{1}{e_1a^{24}_{2s-1}}. \\
  & \OperatorFi{e_2a^{43}_{2t}}{1}{e_1{a}^{43}_{2t+1}}, \\
 \end{split}
\end{equation*}
Since $a^{24}_{2s} \subseteq a^{24}_{2s-1}$ (see Line $G11$), by
(\ref{sym_forms_P_2}) we have
\begin{equation*}
 \begin{split}
   & e_{1\alpha}\tilde{P} =
     e_1a^{42}_{2s}a^{23}_{2r}a^{34}_{2t}(e_1a^{43}_{2t+1} + a^{23}_{2r+1}a^{24}_{2s-1}) = \\
   & e_1a^{42}_{2s}a^{23}_{2r}a^{34}_{2t}(e_1a^{43}_{2t+1}a^{24}_{2s-1} + a^{23}_{2r+1}) =
     e_{1\alpha}(a^{43}_{2t-1}a^{24}_{2s+1} + e_1a^{23}_{2r+1}),
\end{split}
\end{equation*}
 which corresponds to Line $G11$.
\qedsymbol \vspace{2mm}

\subsection{Line $F31$, $\OperatorFi{e_{\alpha}}{1}{e_{1\alpha}}$,
   $\OperatorFi{f_{{\alpha}0}}{1}{f_{{1\alpha}0}}$}

First, we have
\begin{equation}
 \label{Line_F31_act_1_e}
  e_\alpha =
    \OperatorFi{e_3{a}^{21}_{2t}a^{41}_{2r}a^{24}_{2s-1}}{1}
               {e_1{a}^{43}_{2t}a^{23}_{2r}a^{24}_{2s}} =
          e_{1\alpha}.
\end{equation}
 For the case $f_{{\alpha}0}$, select
$$
  P = e_3{a}^{42}_{2s} + a^{14}_{2r+1}a^{12}_{2t-1},
$$
see Lemma \ref{homom_polynom_P}. We have
\begin{equation}
 \label{Line_F31_act_1_f}
   P =
   \OperatorFi{e_3{a}^{42}_{2s} + {a}^{14}_{2r+1}a^{12}_{2t-1}}{1}
              {e_1a^{42}_{2s+1} + e_1{a}^{23}_{2r+1}a^{34}_{2t-1}}.
\end{equation}
Since ${a}^{34}_{2t} \subseteq {a}^{34}_{2t-1}$, we get
\begin{equation*}
 \begin{split}
  & e_{1\alpha}(e_1a^{42}_{2s+1} + e_1a^{23}_{2r+1}a^{34}_{2t-1}) =
   e_1{a}^{43}_{2t}a^{23}_{2r}a^{24}_{2s}
    (e_1a^{42}_{2s+1} + e_1a^{23}_{2r+1}a^{34}_{2t-1}) = \vspace{2mm} \\
  & e_1{a}^{43}_{2t}a^{23}_{2r}a^{24}_{2s}
    (e_1a^{42}_{2s+1}a^{34}_{2t-1} + e_1a^{23}_{2r+1}) =
    e_{1\alpha}(e_1a^{23}_{2r+1} + a^{42}_{2s+1}a^{34}_{2t-1}) =
   \tilde{P},
 \end{split}
\end{equation*}
we get Line $G11$.
 \qedsymbol \vspace{2mm}

 Cases $F31$ (actions
$\psi_2, \psi_4$), and $F41$ (actions $\psi_1, \psi_2, \psi_3$)
are similarly proved.
 \qedsymbol \vspace{2mm}

\subsection{Line $G11$, $\OperatorFi{e_{\alpha}}{2}{e_{2\alpha}}$,
  $\OperatorFi{f_{{\alpha}0}}{2}{f_{{2\alpha}0}}$}
 \label{subs_line_g11}

We have
$$
   \alpha = 1(21)^t(31)^r(41)^s, \quad
   e_\alpha = e_{1(21)^t(31)^r(41)^s}.
$$
Here $2\alpha = (21)^{t+1}(31)^r(41)^s$, and by Corollary
\ref{act_runner_map_D4} we get
\begin{equation*}
 \begin{split}
  & \OperatorFi{a^{24}_{2s}}{2}{e_2{a}^{13}_{2s}}, \\
  & \OperatorFi{a^{32}_{2r}}{2}{e_2{a}^{14}_{2r}}, \\
  & \OperatorFi{e_1a^{34}_{2t}}{2}{e_2a^{34}_{2t+1}}, \\
 \end{split}
\end{equation*}
hence,
\begin{equation}
 \label{Line_G11_act_2}
   e_\alpha =
   \OperatorFi{e_1a^{24}_{2s}a^{34}_{2t}a^{32}_{2r}}{2}
              {e_2a^{13}_{2s}a^{34}_{2t+1}a^{14}_{2r}} =
   e_{2\alpha}.
\vspace{2mm}
\end{equation}
 The polynomial $e_{2\alpha}$ corresponds to Line $G21$
 under the substitution $t \mapsto t+1$.
 For the case $f_{{\alpha}0}$,
\begin{equation*}
 \begin{split}
   & e_\alpha(e_1{a}^{32}_{2r+1} + a^{24}_{2s+1}a^{34}_{2t-1}) = \\
   & e_\alpha(e_1{a}^{23}_{2r+1} + a^{24}_{2s+1}a^{34}_{2t-1}) =
     e_\alpha({a}^{23}_{2r+1} + e_1a^{24}_{2s+1}a^{34}_{2t-1}),
 \end{split}
\end{equation*}
and we select
$$
  P = {a}^{23}_{2r+1} + e_1a^{24}_{2s+1}a^{34}_{2t-1}.
$$
 We have
\begin{equation}
 \label{Line_G11_act_1_f}
   P =
   \OperatorFi{e_1{a}^{23}_{2r+1} + e_1a^{24}_{2s+1}a^{34}_{2t-1}}{2}
              {e_2a^{14}_{2r+1} + e_2{a}^{13}_{2s+1}a^{34}_{2t}} =
              \tilde{P}.
\end{equation}
Since
$$
  e_{2\alpha} \subseteq  e_2{a}^{13}_{2s} \subseteq
  e_2{a}^{13}_{2s-1},
$$
 by the equalizing property (\ref{sym_forms_P_3}) we have
\begin{equation*}
 \begin{split}
   & e_{2\alpha}\tilde{P} =
     e_{2\alpha}(e_2a^{14}_{2r+1} + e_2{a}^{13}_{2s+1}a^{34}_{2t}) = \\
   & e_{2\alpha}(e_2a^{14}_{2r+1} + e_2{a}^{13}_{2s-1}a^{34}_{2t+2})=
     e_{2\alpha}(a^{34}_{2t+2} + e_2{a}^{13}_{2s-1}a^{14}_{2r+1}),
 \end{split}
\end{equation*}
i.e., we get Line $F21$ under substitution the $t \mapsto t+1$.
\qedsymbol \vspace{2mm}

 Cases $G11$ (actions $\psi_3, \psi_4$) are similarly proved.
 \qedsymbol \vspace{2mm}

\subsection{Line $G21$, $\OperatorFi{e_{\alpha}}{3}{e_{2\alpha}}$,
 $\OperatorFi{f_{{\alpha}0}}{3}{f_{{3\alpha}0}}$}
 \label{subs_line_G21}

We have
$$
  \alpha = 2(41)^r(31)^s(21)^t, \quad
  e_\alpha = e_{2(41)^r(31)^s(21)^t}.
$$
By eq.(\ref{act_psi_3_F31}) $3\alpha = (31)^{s+1}(41)^r(21)^t$,
and by Corollary \ref{act_runner_map_D4} we have
\begin{equation*}
 \begin{split}
  & \OperatorFi{a^{34}_{2t}}{3}{e_3{a}^{12}_{2t}}, \\
  & \OperatorFi{a^{31}_{2s+1}}{3}{e_3{a}^{24}_{2s+1}}, \\
  & \OperatorFi{e_2a^{14}_{2r-1}}{3}{e_3a^{14}_{2r}}, \\
 \end{split}
\end{equation*}
i.e.,
\begin{equation}
 \label{Line_G21_act_3}
   e_\alpha =
   \OperatorFi{e_2a^{34}_{2t}a^{31}_{2s+1}a^{14}_{2r-1}}{3}
              {e_3a^{12}_{2t}a^{24}_{2s+1}a^{14}_{2r}} =
   e_{3\alpha}.
\end{equation}
 The polynomial $e_{3\alpha}$ corresponds to Line $F31$ for
 substitution $s \mapsto s+1$.
 For the case $f_{{\alpha}0}$,
 we select
$$
   P = e_2a^{14}_{2r} + {a}^{31}_{2s+2}a^{34}_{2t-1},
$$
and we have
\begin{equation}
 \label{Line_G21_act_3_f}
   P =
   \OperatorFi{e_2a^{14}_{2r} + {a}^{31}_{2s+2}a^{34}_{2t-1}}{3}
              {e_3a^{14}_{2r+1} + {a}^{24}_{2s+2}a^{12}_{2t-1}} =
   \tilde{P}.
\end{equation}
Since
 $e_{3\alpha} \subseteq a^{12}_{2t} \subseteq a^{12}_{2t-1}$,
we get
\begin{equation}
  \label{Line_G21_act_3_e}
 \begin{split}
  & e_{3\alpha}\tilde{P} =
    e_{3\alpha}(e_3a^{14}_{2r+1} + {a}^{24}_{2s+2}a^{12}_{2t-1}) =\\
  & e_{3\alpha}(a^{14}_{2r+1}a^{12}_{2t-1} + e_3{a}^{24}_{2s+2})
 \end{split}
\end{equation}
which corresponds to Line $F31$. By (\ref{Line_G21_act_3_f}) and
(\ref{Line_G21_act_3_e}) $f_{3\alpha}$ and $e_{3\alpha}$
correspond to Line $F31$ under substitution the $s \mapsto s+1$.
\qedsymbol \vspace{2mm}

\subsection{Line $G21$, $\OperatorFi{e_{\alpha}}{4}{e_{4\alpha}}$,
  $\OperatorFi{f_{{\alpha}0}}{4}{f_{{4\alpha}0}}$}
The admissible sequence $\alpha$ and the polynomial
 $e_\alpha$ are given as above in \S\ref{subs_line_G21}.
 \vspace{2mm}

By eq.(\ref{act_psi_4_F31}) $4\alpha = (41)^r(31)^{s+1}(21)^t$,
and by Corollary \ref{act_runner_map_D4} we have
\begin{equation*}
 \begin{split}
  & \OperatorFi{a^{34}_{2t}}{4}{e_4{a}^{12}_{2t}}, \\
  & \OperatorFi{a^{14}_{2r-1}}{4}{e_4{a}^{23}_{2r-1}}, \\
  & \OperatorFi{e_2a^{31}_{2s+1}}{4}{e_4a^{31}_{2s+2}}, \\
 \end{split}
\end{equation*}
and therefore
\begin{equation}
 \label{Line_G21_act_4}
   e_\alpha =
   \OperatorFi{e_2a^{34}_{2t}a^{31}_{2s+1}a^{14}_{2r-1}}{4}
              {e_4a^{12}_{2t}a^{23}_{2r-1}a^{31}_{2s+2}} =
   e_{4\alpha}.
\end{equation}
The polynomial $e_{4\alpha}$ corresponds to Line $F41$ under the
 substitution $s \mapsto s + 1$ . For the case $f_{{\alpha}0}$,
 we have
\begin{equation*}
 \begin{split}
 & e_\alpha(a^{14}_{2r} + e_2{a}^{31}_{2s+2}a^{34}_{2t-1}) =
   e_\alpha(e_2a^{14}_{2r} + e_2{a}^{31}_{2s+2}a^{34}_{2t-1}) = \\
 & e_\alpha(e_2a^{41}_{2r} + e_2{a}^{31}_{2s+2}a^{34}_{2t-1}) =
   e_\alpha(a^{41}_{2r} + e_2{a}^{31}_{2s+2}a^{34}_{2t-1}),
 \end{split}
\end{equation*}
and we select
$$
  P = a^{41}_{2r} + e_2{a}^{31}_{2s+2}a^{34}_{2t-1}.
$$
 Here, we get
\begin{equation}
 \label{Line_G21_act_4_f}
   P =
   \OperatorFi{a^{41}_{2r} + e_2{a}^{31}_{2s+2}a^{34}_{2t-1}}{4}
              {e_4a^{23}_{2r} + e_4{a}^{31}_{2s+3}a^{12}_{2t-1}} =
              \tilde{P}.
\end{equation}
Since
$$
  e_{4\alpha}\tilde{P} =
  e_{4\alpha}(e_4a^{32}_{2r} + {a}^{31}_{2s+3}a^{12}_{2t-1}),
$$
 we see that $e_{4\alpha}\tilde{P}$ corresponds to Line
$F41$ under the substitution $s \mapsto s + 1$.
 \qedsymbol \vspace{2mm}

\subsection{Line $G21$, \OperatorFi{e_{\alpha}}{1}{e_{1\alpha}},
    $\OperatorFi{f_{{\alpha}0}}{1}{f_{{1\alpha}0}}$}

 The admissible sequence $\alpha$ and the polynomial
 $e_\alpha$ are given as above in \S\ref{subs_line_G21}.
 \vspace{2mm}

By eq.(\ref{act_psi_1_F31}) $1\alpha = (14)^r(31)^{s+1}(21)^t$,
and by Corollary \ref{act_runner_map_D4} we have
\begin{equation*}
 \begin{split}
  & \OperatorFi{a^{31}_{2s+1}}{1}{e_1{a}^{24}_{2s+1}}, \\
  & \OperatorFi{a^{14}_{2r-1}}{1}{e_1{a}^{23}_{2r-1}}, \\
  & \OperatorFi{e_2a^{34}_{2t}}{1}{e_1a^{34}_{2t+1}}, \\
 \end{split}
\end{equation*}
and therefore
\begin{equation}
 \label{Line_G21_act_1_e}
   e_\alpha =
   \OperatorFi{e_2a^{31}_{2s+1}a^{14}_{2r-1}a^{34}_{2t}}{1}
              {e_1a^{24}_{2s+1}a^{23}_{2r-1}a^{34}_{2t+1}} =
   e_{1\alpha}.
\end{equation}
The polynomial $e_{1\alpha}$ corresponds to Line $H11$ under
 the substitution $t \mapsto t+1$.
 For the case $f_{{\alpha}0}$, we have
\begin{equation*}
 \begin{split}
  & e_\alpha{P} =
    e_2a^{31}_{2s+1}a^{14}_{2r-1}a^{34}_{2t}
    (e_2a^{14}_{2r} + {a}^{31}_{2s+2}a^{34}_{2t-1}) = \\
  & e_2a^{31}_{2s+1}a^{14}_{2r-1}a^{34}_{2t}
    (a^{14}_{2r} + e_2{a}^{31}_{2s+2}a^{34}_{2t-1}),
 \end{split}
\end{equation*}
and we select
$$
  P = a^{14}_{2r} + e_2{a}^{31}_{2s+2}a^{34}_{2t-1}.
$$
Thus,
\begin{equation}
 \label{Line_G21_act_1_f}
   P =
   \OperatorFi{a^{14}_{2r} + e_2{a}^{31}_{2s+2}a^{34}_{2t-1}}{1}
              {e_1a^{23}_{2r} + e_1{a}^{24}_{2s+2}a^{34}_{2t}} =
              \tilde{P}.
\end{equation}
By (\ref{sym_forms_P_4})
$$
  e_1{a}^{24}_{2s+2}a^{34}_{2t} = e_1{a}^{24}_{2s+2}a^{34}_{2t},
$$
and by (\ref{Line_G21_act_1_f}), we get
\begin{equation}
 \label{Line_G21_act_1_f_cont}
   \tilde{P} = e_1a^{23}_{2r} + e_1{a}^{24}_{2s}a^{34}_{2t+2}.
\end{equation}
 By (\ref{Line_G21_act_1_f_cont}) and (\ref{Line_G21_act_1_e})
$f_{1\alpha}$ and $e_{1\alpha}$ correspond to Line $H11$ under
 the substitution $t \mapsto t+1$.
 \qedsymbol \vspace{2mm}

 Cases $G31$ (actions $\psi_1, \psi_2, \psi_4$), and $G41$ (actions
$\psi_1, \psi_2, \psi_3$) are similarly proved.
 \qedsymbol \vspace{2mm}

\subsection{Line $H11$, $\OperatorFi{e_{\alpha}}{2}{e_{2\alpha}}$,
  $\OperatorFi{f_{{\alpha}0}}{2}{f_{{2\alpha}0}}$}
 Here,
$$
  \alpha = (14)^r(31)^s(21)^t,
  e_\alpha = e_{(23)^{2r-1}(24)^{2s-1}(34)^{2t+1}}.
  \vspace{2mm}
$$
Since $2\alpha = 2(14)^r(31)^s(21)^t$, and by Corollary
\ref{act_runner_map_D4} we have
\begin{equation*}
 \begin{split}
  & \OperatorFi{a^{23}_{2r-1}}{2}{e_2{a}^{14}_{2r-1}}, \\
  & \OperatorFi{a^{24}_{2s-1}}{2}{e_2{a}^{31}_{2s-1}}, \\
  & \OperatorFi{e_1a^{34}_{2t+1}}{2}{e_2a^{34}_{2t+2}}, \\
 \end{split}
\end{equation*}
i.e.,
\begin{equation}
 \label{Line_H11_act_2}
   e_\alpha =
   \OperatorFi{e_1a^{23}_{2r-1}a^{24}_{2s-1}a^{34}_{2t+1}}{2}
              {e_2a^{14}_{2r-1}a^{31}_{2s-1}a^{34}_{2t+2}} =
   e_{2\alpha}.
\end{equation}
The polynomial $e_{2\alpha}$ corresponds to Line $H21$.
 For the case $f_{{\alpha}0}$,
 select
$$
  P = e_1a^{34}_{2t} + {a}^{24}_{2s}{a}^{23}_{2r},
$$
 and we have
\begin{equation}
 \label{Line_H11_act_2_f}
   P =
   \OperatorFi{e_1{a}^{34}_{2t} + {a}^{24}_{2s}{a}^{23}_{2r}}{2}
              {e_2{a}^{34}_{2t+1} + e_2{a}^{13}_{2s}a^{14}_{2r}} =
    \tilde{P}.
\end{equation}
By (\ref{Line_H11_act_2_f}) and (\ref{Line_H11_act_2})
$f_{2\alpha}$ and $e_{2\alpha}$ correspond to Line $H21$.
 \qedsymbol \vspace{2mm}

 Cases $H11$ (actions $\psi_3, \psi_4$) are similarly proved.
 \qedsymbol \vspace{2mm}

\subsection{Line $H21$, $\OperatorFi{e_{\alpha}}{1}{e_{1\alpha}}$,
   $\OperatorFi{f_{{\alpha}0}}{1}{f_{{1\alpha}0}}$}
 \label{subs_line_H21}
 We have
$$
  \alpha = 2(14)^r(31)^s(21)^t,  e_\alpha = e_{2(14)^r(31)^s(21)^t}.
  \vspace{2mm}
$$

Since $1\alpha = 12(14)^r(31)^s(21)^t = (14)^r(31)^s(21)^{t+1}$
(see heading (10) of  Proposition \ref{relations_1_14}), and by
Corollary \ref{act_runner_map_D4} we have
\begin{equation*}
 \begin{split}
  & \OperatorFi{a^{14}_{2r-1}}{1}{e_1{a}^{23}_{2r-1}}, \\
  & \OperatorFi{a^{13}_{2s-1}}{1}{e_1{a}^{24}_{2s-1}}, \\
  & \OperatorFi{e_2a^{34}_{2t+2}}{1}{e_1a^{34}_{2t+3}}, \\
 \end{split}
\end{equation*}
i.e.,
\begin{equation}
 \label{Line_H21_act_1}
   e_\alpha =
   \OperatorFi{e_2a^{14}_{2r-1}a^{13}_{2s-1}a^{34}_{2t+2}}{1}
              {e_1a^{23}_{2r-1}a^{24}_{2s-1}a^{34}_{2t+3}} =
   e_{1\alpha}.
\end{equation}
The polynomial $e_{2\alpha}$ corresponds to Line $H11$ under
 the substitution $t \mapsto t+1$.
 For the case $f_{{\alpha}0}$, we select
$$
  P = {a}^{14}_{2r} + e_2{a}^{31}_{2s}a^{34}_{2t+1},
$$
 and we have
\begin{equation}
 \label{Line_H21_act_1_f}
   P =
   \OperatorFi{{a}^{14}_{2r} + e_2{a}^{31}_{2s}a^{34}_{2t+1}}{1}
              {e_1{a}^{23}_{2r} + e_1{a}^{24}_{2s}a^{34}_{2t+2}} =
     \tilde{P}.
\end{equation}
By (\ref{Line_H21_act_1_f}) and (\ref{Line_H21_act_1})
$f_{1\alpha}$ and $e_{1\alpha}$ correspond to Line $H11$ under
 the substitution $t \mapsto t+1$.
 \qedsymbol \vspace{2mm}

\subsection{Line $H21$, $\OperatorFi{e_{\alpha}}{3}{e_{1\alpha}}$,
  $\OperatorFi{f_{{\alpha}0}}{3}{f_{{3\alpha}0}}$}
The admissible sequence $\alpha$ and the polynomial
 $e_\alpha$ are given as above in \S\ref{subs_line_H21}.
 \vspace{2mm}

Since $3\alpha = 32(14)^r(31)^s(21)^t = (31)^s(21)^{t+1}(41)^r$
(see heading (12) of  Proposition \ref{relations_1_14}), and by
Corollary \ref{act_runner_map_D4} we have
\begin{equation*}
 \begin{split}
  & \OperatorFi{a^{13}_{2s-1}}{3}{e_3{a}^{24}_{2s-1}}, \\
  & \OperatorFi{a^{34}_{2t+2}}{3}{e_3a^{12}_{2t+2}}, \\
  & \OperatorFi{e_2a^{14}_{2r-1}}{3}{e_3{a}^{14}_{2r}}, \\
 \end{split}
\end{equation*}
i.e.,
\begin{equation}
 \label{Line_H21_act_3}
   e_\alpha =
   \OperatorFi{e_2a^{14}_{2r-1}a^{13}_{2s-1}a^{34}_{2t+2}}{3}
              {e_3a^{14}_{2r}a^{24}_{2s-1}a^{12}_{2t+2}} =
   e_{3\alpha}.
\end{equation}
The polynomial $e_{3\alpha}$ corresponds to Line $F31$ under
 the substitution $t \mapsto t+1$.
 For the case $f_{{\alpha}0}$, we select
$$
  P = e_2{a}^{14}_{2r} + {a}^{31}_{2s}a^{34}_{2t+1},
$$
 and we have
\begin{equation}
 \label{Line_H21_act_3_f}
   P =
   \OperatorFi{e_2{a}^{14}_{2r} + {a}^{31}_{2s}a^{34}_{2t+1}}{3}
              {e_3{a}^{14}_{2r+1} + e_3{a}^{24}_{2s}a^{12}_{2t+1}}
              = \tilde{P}.
\end{equation}
By (\ref{Line_H21_act_3})
$$
  e_{3\alpha}\tilde{P} =
  e_{3\alpha}(e_3{a}^{14}_{2r+1} + {a}^{24}_{2s}a^{12}_{2t+1})
    = e_{3\alpha}(e_3{a}^{24}_{2s} + {a}^{14}_{2r+1}a^{12}_{2t+1}),
$$
and therefore, by (\ref{Line_H21_act_3_f}) and
(\ref{Line_H21_act_3}), we see that $f_{1\alpha}$ and
$e_{1\alpha}$ correspond to Line $F31$ under the
 substitution $t \mapsto t+1$.
 \qedsymbol \vspace{2mm}

\subsection{Line $H21$, $\OperatorFi{e_{\alpha}}{4}{e_{4\alpha}}$,
    $\OperatorFi{f_{{\alpha}0}}{4}{f_{{4\alpha}0}}$}
The admissible sequence $\alpha$ and the polynomial
 $e_\alpha$ are given as above in \S\ref{subs_line_H21}.
 \vspace{2mm}

Since $4\alpha = 42(14)^r(31)^s(21)^t = (41)^s(21)^{t+1}(31)^r$
(see heading (13) of  Proposition \ref{relations_1_14}), and by
Corollary \ref{act_runner_map_D4} we have
\begin{equation*}
 \begin{split}
  & \OperatorFi{a^{34}_{2t+2}}{4}{e_4a^{12}_{2t+2}}, \\
  & \OperatorFi{a^{14}_{2r-1}}{4}{e_4{a}^{23}_{2r-1}}, \\
  & \OperatorFi{e_2a^{13}_{2s-1}}{4}{e_4{a}^{13}_{2s}}, \\
 \end{split}
\end{equation*}
i.e.,
\begin{equation}
 \label{Line_H21_act_4}
   e_\alpha =
   \OperatorFi{e_2a^{14}_{2r-1}a^{13}_{2s-1}a^{34}_{2t+2}}{4}
              {e_4a^{23}_{2r-1}a^{13}_{2s}a^{12}_{2t+2}} =
   e_{4\alpha}.
\end{equation}
The polynomial $e_{4\alpha}$ corresponds to Line $F41$ under
 the substitution $t \mapsto t+1$.
 For the case $f_{{\alpha}0}$, we select
$$
  P = {a}^{41}_{2r} + e_2{a}^{31}_{2s}a^{34}_{2t+1},
$$
 and we have
\begin{equation}
 \label{Line_H21_act_4_f}
   P =
   \OperatorFi{{a}^{41}_{2r} + e_2{a}^{31}_{2s}a^{34}_{2t+1}}{4}
              {e_4{a}^{23}_{2r} + e_4{a}^{31}_{2s+1}a^{12}_{2t+1}}.
\end{equation}
By (\ref{Line_H21_act_4_f}) and (\ref{Line_H21_act_4}) we see that
$f_{4\alpha}$ and $e_{4\alpha}$ correspond to Line $F41$ under
 the substitution $t \mapsto t+1$.
 \qedsymbol \vspace{2mm}

Cases $H31$, (actions $\psi_1, \psi_3, \psi_4$) and $H41$,
(actions $\psi_1, \psi_2, \psi_4$) are similarly proved.
 \qedsymbol \vspace{2mm}

  Theorem \ref{th_adm_classes_D4} on the admissible elements in
  $D^4$ is proved for all polynomials $e_\alpha$,
  see Table \ref{table_adm_elem_D4}.
 \qedsymbol \vspace{2mm}

\chapter{\sc\bf The inclusion theorem and properties of $H^{+}(n)$}
 \label{sect_inclusion_theorem}

\section{Proof of the inclusion theorem}
\subsection{More properties of the atomic elements}

The atomic elements $a_n^{ij}$ and $A_n^{ji}$ can frequently
substitute each other. This has already been shown in
 Lines 3.1 -- 5.2 of Table \ref{table_atomic}.
\begin{proposition}
  \label{more_eq}
 The following relations hold for $n-2 \leq m$: \vspace{2mm}
\begin{enumerate}
  \item
     $y_i{a}_m^{ik}a_n^{kj} = y_i{a}_m^{ik}A_n^{jk}$, \vspace{2mm}
  \item
     $y_i{A}_m^{ji}a_n^{jk} = y_i{A}_m^{ji}A_n^{kj}$. \vspace{2mm}
\end{enumerate}
\end{proposition}
\PerfProof 1) Since
$$
  y_i{a}_m^{ik}a_n^{kj} =
     y_i{a}_m^{ik}(x_k + y_j(x_j + y_i{a}_{n-2}^{ik}))
$$
and $y_i{a}_{n-2}^{ik} \supseteq y_i{a}_m^{ik}$ for $n-2 \leq m$,
it follows that
$$
   y_i{a}_m^{ik}a_n^{kj} =
      y_i{a}_m^{ik}(y_j +
         x_k(x_j + y_i{a}_{n-2}^{ik})).
$$
Since $a_{n-2}^{ik} \supseteq x_k$, we have
$$
   y_i{a}_m^{ik}{a}_n^{kj} =
      y_i{a}_m^{ik}(y_j + x_k(y_i + x_j{a}_{n-2}^{ik})).
$$
Further, by Table \ref{table_atomic}, Line 3.1, we have
  $x_j{a}_{n-2}^{ik} = x_j{A}_{n-2}^{ki}$ and
$$
  y_i{a}_m^{ik}a_n^{kj} =
     y_i{a}_m^{ik}(y_j + x_k(y_i + x_j{A}_{n-2}^{ki})) =
           y_i{a}_m^{ik}A_n^{jk}.
 \qed \vspace{2mm}
$$

 2)  Since
$$
  y_i{A}_m^{ji}a_n^{jk} =
     y_i{A}_m^{ji}(x_j + y_k(x_k + y_i{a}_{n-2}^{ij})) =
        y_i{A}_m^{ji}(x_j +
          x_k + y_k{y}_i{a}_{n-2}^{ij})
$$
and by Table \ref{table_atomic}, Line 4.1, we have
  $y_k{y}_i{a}_{n-2}^{ij} =
         y_k{y}_i{A}_{n-2}^{ji}$, and
$$
  y_i{A}_m^{ji}a_n^{jk} =
       y_i{A}_m^{ji}(x_j + x_k + y_k{y}_i{A}_{n-2}^{ji}).
$$
Since $A_{n-2}^{ji} \supseteq A_m^{ji}$, we see that
\begin{equation*}
 \begin{split}
     & y_i{A}_m^{ji}a_n^{jk} =
         y_i{A}_m^{ji}(y_k{y}_i + A_{n-2}^{ji}(x_j + x_k)) = \\
     & y_i{A}_m^{ji}(y_k{y}_i + x_j + x_k{A}_{n-2}^{ji}) =
       y_i{A}_m^{ji}(x_j + y_k(y_i + x_k{A}_{n-2}^{ji})).
  \end{split}
\end{equation*}
     Further,
$$
     y_i + x_k{A}_{n-2}^{ji} \supseteq y_i{A}_m^{ji},
$$
and we have
$$     y_i{A}_m^{ji}{a}_n^{jk} =
       y_i{A}_m^{ji}(y_k + x_j(y_i + x_k{A}_{n-2}^{ji})) =
       y_i{A}_m^{ji}A_n^{kj}.
\qed \vspace{5mm}
$$
\subsection{The $a$-form and $A$-form of admissible elements}

\begin{proposition}
 \label{subst_a_A}
    For all admissible elements $f_\alpha$ and $e_\alpha$ of
    Table \ref{table_adm_elems},
 \begin{enumerate}
 \item
 Every symbol ``$A$'' can be substituted by the symbol ``$a$''. \vspace{2mm}
 \item
 Every symbol ``$a$'' can be substituted by the symbol ``$A$''. \vspace{2mm}
 \end{enumerate}
\end{proposition}
 \index{admissible elements! - $a$-form, $A$-form}

   The substitutions $A \longleftrightarrow a$
   are collected in Table \ref{table_two_forms}.
   The original form of  admissible elements
   is given by Table \ref{table_adm_elems}.
   Two other forms of admissible elements are given by
   Proposition \ref{subst_a_A}.
   The substitution (1) (resp. (2)) from Proposition \ref{subst_a_A}
   is said to be an {\it $a$-form} (resp. {\it $A$-form}).
   (Table \ref{table_two_forms}).
    \vspace{2mm}

\small
\begin{table}
 \renewcommand{\arraystretch}{1.35}
 \begin{tabular} {|c||c|c|c|c|c|c|}
  \hline \hline
    $N$ & $\alpha$
            & $f_\alpha$,
            & $e_\alpha$,
            & $g_{\alpha0}$    \cr
            &
            & $a$-form = $A$-form
            & $a$-form = $A$-form
            &     \\
  \hline  \hline
    1 & $(213)^{2p}(21)^{2k}$
      & $y_1y_2a^{13}_qa^{23}_{k-1}$ =
      & $y_2a^{23}_{k-1}a^{21}_qa^{31}_ka^{13}_q$ =
      & $e_\alpha(x_1+a^{32}_qA^{32}_k)$ \cr
      &
      & $y_1y_2A^{31}_qA^{32}_{k-1}$
      & $y_2A^{32}_{k-1}A^{12}_qA^{13}_kA^{31}_q$
      & \\
  \hline
    2 & $3(213)^{2p}(21)^{2k}$
      & $y_3(x_1+x_2)a^{32}_qa^{31}_{k-1}$ =
      & $y_3a^{31}_{k-1}a^{21}_{q+1}a^{12}_ka^{32}_q$ =
      & $e_\alpha(y_2y_3+A^{12}_qa^{31}_k)$ \cr
      &
      & $y_3(x_1+x_2)A^{23}_qA^{13}_{k-1}$
      & $y_3A^{13}_{k-1}A^{12}_{q+1}A^{21}_kA^{23}_q$
      & \\
  \hline
    3 & $13(213)^{2p}(21)^{2k}$
      & $y_3y_1a^{32}_{q+1}a^{12}_{k-1}$ =
      & $y_1a^{12}_{k-1}a^{13}_{q+1}a^{23}_ka^{32}_{q+1}$ =
      & $e_\alpha(x_3+a^{21}_{q+1}A^{21}_k)$    \cr
      &
      & $y_3y_1A^{23}_{q+1}A^{21}_{k-1}$
      & $y_1A^{21}_{k-1}A^{31}_{q+1}A^{32}_kA^{23}_{q+1}$
      & \\
  \hline
    4 & $(213)^{2p+1}(21)^{2k}$
      & $y_2(x_3+x_1)a^{21}_{q+1}a^{23}_{k-1}$ =
      & $y_2a^{23}_{k-1}a^{13}_{q+2}a^{31}_ka^{21}_{q+1}$ =
      & $e_\alpha(y_1y_2+A^{31}_{q+1}a^{23}_k)$ \cr
      &
      & $y_2(x_3+x_1)A^{12}_{q+1}A^{32}_{k-1}$
      & $y_2A^{32}_{k-1}A^{31}_{q+2}A^{13}_kA^{12}_{q+1}$
      & \\
  \hline
    5 & $3(213)^{2p+1}(21)^{2k}$
      & $y_2y_3a^{21}_{q+2}a^{31}_{k-1}$ =
      & $y_3a^{31}_{k-1}a^{32}_{q+2}a^{12}_ka^{21}_{q+2}$ =
      & $e_\alpha(x_2+a^{13}_{q+2}A^{13}_k)$    \cr
      &
      & $y_2y_3A^{12}_{q+2}A^{13}_{k-1}$
      & $y_3A^{13}_{k-1}A^{23}_{q+2}A^{21}_kA^{12}_{q+2}$
      & \\
  \hline
    6 & $13(213)^{2p+1}(21)^{2k}$
      & $y_1(x_2+x_3)a^{13}_{q+2}a^{12}_{k-1}$ =
      & $y_1a^{12}_{k-1}a^{32}_{q+3}a^{23}_ka^{13}_{q+2}$ =
      & $e_\alpha(y_1y_3+A^{23}_{q+2}a^{12}_k)$ \cr
      &
      & $y_1(x_2+x_3)A^{31}_{q+2}A^{21}_{k-1}$
      & $y_1A^{21}_{k-1}A^{23}_{q+3}A^{32}_kA^{31}_{q+2}$
      & \\
  \hline \hline
    7 & $1(21)^{2k}$
      & $x_1a^{32}_ka^{23}_k$ =
      & $y_1a^{13}_ka^{12}_ka^{32}_ka^{23}_k$ =
      & $e_\alpha(y_2a^{21}_k+y_3a^{31}_k)$ \cr
      &
      & $x_1A^{23}_kA^{32}_k$
      & $y_1A^{31}_kA^{21}_kA^{23}_kA^{32}_k$
      & \\
  \hline
    8 & $1(21)^{2k+1}$
      & $y_1(x_2+x_3)a^{13}_ka^{12}_k$ =
      & $y_1a^{13}_ka^{12}_ka^{32}_{k+1}a^{23}_{k+1}$ =
      & $e_\alpha(y_1a^{12}_{k+1}+y_1a^{13}_{k+1})$ \cr
      &
      & $y_1(x_2+x_3)A^{31}_kA^{21}_k$
      & $y_1A^{31}_kA^{21}_kA^{23}_{k+1}A^{32}_{k+1}$
      & \\
  \hline  \hline
    9 & $(213)^{2p}(21)^{2k+1}$
      & $y_2y_3a^{21}_qa^{31}_k$ =
      & $y_2a^{23}_ka^{21}_qa^{31}_ka^{13}_{q+1}$ =
      & $e_\alpha(x_2+a^{13}_qA^{13}_{k+1})$  \cr
      &
      & $y_2y_3A^{12}_qA^{13}_k$
      & $y_2A^{32}_kA^{12}_qA^{13}_kA^{31}_{q+1}$
      & \\
  \hline
   10 & $3(213)^{2p}(21)^{2k+1}$
      & $x_3a^{21}_{q+1}a^{12}_k$ =
      & $y_3a^{31}_ka^{21}_{q+1}a^{12}_ka^{32}_{q+1}$ =
      & $e_\alpha(y_1y_3+A^{23}_qa^{12}_{k+1})$ \cr
      &
      & $x_3A^{12}_{q+1}A^{21}_k$
      & $y_3A^{13}_kA^{12}_{q+1}A^{21}_kA^{23}_{q+1}$
      & \\
  \hline
   11 & $13(213)^{2p}(21)^{2k+1}$
      & $y_1y_2a^{13}_{q+1}a^{23}_k$ =
      & $y_1a^{12}_ka^{13}_{q+1}a^{23}_ka^{32}_{q+2}$ =
      & $e_\alpha(x_1+a^{32}_{q+1}A^{32}_{k+1})$    \cr
      &
      & $y_1y_2A^{31}_{q+1}A^{32}_k$
      & $y_1A^{21}_kA^{31}_{q+1}A^{32}_kA^{23}_{q+2}$
      & \\
  \hline
   12 & $(213)^{2p+1}(21)^{2k+1}$
      & $x_2a^{13}_{q+2}a^{31}_k$ =
      & $y_2a^{23}_ka^{13}_{q+2}a^{31}_ka^{21}_{q+2}$ =
      & $e_\alpha(y_2y_3+A^{12}_{q+1}a^{31}_{k+1})$ \cr
      &
      & $x_2A^{31}_{q+2}A^{13}_k$
      & $y_2A^{32}_kA^{31}_{q+2}A^{13}_kA^{12}_{q+2}$
      & \\
  \hline
   13 & $3(213)^{2p+1}(21)^{2k+1}$
      & $y_3y_1a^{32}_{q+2}a^{12}_k$ =
      & $y_3a^{31}_ka^{32}_{q+2}a^{12}_ka^{21}_{q+3}$ =
      & $e_\alpha(x_3+a^{21}_{q+2}A^{21}_{k+1})$    \cr
      &
      & $y_3y_1A^{23}_{q+2}A^{21}_k$
      & $y_3A^{13}_kA^{23}_{q+2}A^{21}_kA^{12}_{q+3}$
      & \\
  \hline
    14 & $13(213)^{2p+1}(21)^{2k+1}$
      & $x_1a^{32}_{q+3}a^{23}_k$ =
      & $y_1a^{12}_ka^{32}_{q+3}a^{23}_ka^{13}_{q+3}$ =
      & $e_\alpha(y_1y_2+A^{31}_{q+2}a^{23}_{k+1})$ \cr
      &
      & $x_1A^{23}_{q+3}A^{32}_k$
      & $y_1A^{21}_kA^{23}_{q+3}A^{32}_kA^{31}_{q+3}$
      &  \\
  \hline \hline
  \end{tabular}
  \vspace{2mm}
  \caption{\hspace{3mm}The admissible elements, $a$-form and $A$-form}
  \label{table_two_forms}
  \begin{center}
    In Lines 1 -- 6: $k > 0$ and  $p \geq 0$.
    In Lines 7 -- 8: $k > 0$. \\
    In Lines 9 -- 14: $k \geq 0, p \geq 0$.
    In all Lines: $q=k+3p$.
  \end{center}
\end{table}
\normalsize

{\it Proof of Proposition \ref{subst_a_A}.}
  For elements $f_\alpha$, both statements (1) and (2) of
  Proposition \ref{subst_a_A} follow from
  Table \ref{table_atomic}, headings (3.1), (4.1) and (5.2).
  For convenience, we give these properties here:
\begin{equation*}
\begin{array}{lll}
     & x_k{A}_n^{ij} = x_k{a}_n^{ji}
         & \text{ (Table \ref{table_atomic}, Line 3.1)},  \\
     & y_i{y}_j{A}_n^{ki} =
                 y_i{y}_j{a}_n^{ik}
         & \text{ (Table \ref{table_atomic}, Line 4.1)},   \\
     & y_i(x_j + x_k){A}_n^{ij} =
            y_i(x_j + x_k)a_n^{ji}
         & \text{ (Table \ref{table_atomic}, Line 5.2)}.
\end{array}
\end{equation*}
  Consider elements $e_\alpha$. We will use
  Proposition \ref{subst_a_A} and property 5.1 from Table \ref{table_atomic}.
  We have
  \begin{equation*}
     a_n^{ij}a_m^{jk} = A_n^{ji}a_m^{jk}
        \text{ for } n \leq m+1
        \qquad \text{ (Table \ref{table_atomic}, Line 5.1)}.
  \end{equation*}

  \underline{Line 1, Line 9 of Table \ref{table_adm_elems}}.
  Take Line 1:
$$
  e_\alpha = y_2A^{32}_{k-1}a^{21}_qA^{13}_ka^{13}_q.
$$
  According to Line 5.1 of Table \ref{table_atomic},  we have
$$
  a^{21}_qa^{13}_q = A^{12}_qa^{13}_q  \text{ and }
  e_\alpha = y_2A^{32}_{k-1}A^{12}_qA^{13}_ka^{13}_q.
$$
  Further, by Proposition \ref{more_eq}, (2) we have
$$
  y_2A^{12}_ka^{13}_q = y_2A^{12}_kA^{31}_q
  \text{ and }
  e_\alpha =  y_2A^{32}_{k-1}A^{12}_qA^{13}_kA^{31}_q.
$$
  This is the $A$-form of $e_\alpha$.

  On the other hand, by Proposition \ref{more_eq}, (2) we have
$$
  y_2A^{32}_{k-1}A^{13}_k = y_2A^{32}_{k-1}a^{31}_k \text{ and }
  e_\alpha = y_2A^{32}_{k-1}a^{21}_qa^{31}_ka^{13}_q.
$$
  Again, by Line 5.1 of Table \ref{table_atomic}, we have
$$
  A^{32}_{k-1}a^{31}_k = a^{23}_{k-1}a^{31}_k \text{ and }
  e_\alpha = y_2a^{23}_{k-1}a^{21}_qa^{31}_ka^{13}_q.
$$
  It is the $a$-form of $e_\alpha$.
  The same for Line 9 of Table \ref{table_adm_elems}.
\qedsymbol \vspace{2mm}

  \underline{Line 2, Line 10 of Table \ref{table_adm_elems}}.
  For Line 2, we have
$$
   e_\alpha = y_3a^{31}_{k-1}A^{12}_{q+1}a^{12}_kA^{23}_q.
$$
   By Table \ref{table_atomic}, Line 5.1 we have
$$
   a^{31}_{k-1}a^{12}_k = A^{13}_{k-1}a^{12}_k \text{ and }
   e_\alpha = y_3A^{13}_{k-1}A^{12}_{q+1}a^{12}_kA^{23}_q.
$$
   By proposition \ref{more_eq}, (2) we obtain
$$
   y_3A^{13}_{k-1}a^{12}_k = y_3A^{13}_{k-1}A^{21}_k \text{ and }
   e_\alpha = y_3A^{13}_{k-1}A^{12}_{q+1}A^{21}_kA^{23}_q \quad
   \text{($A$-form)}.
$$

   On the other hand,
   by proposition \ref{more_eq}, (2) we have
$$
   y_3A^{23}_qA^{12}_{q+1} = y_3A^{23}_qa^{21}_{q+1} \text{ and }
   e_\alpha = y_3a^{31}_{k-1}a^{21}_{q+1}a^{12}_kA^{23}_q.
$$
   By Table \ref{table_atomic}, Line 5.1 we get
$$
   A^{23}_qa^{21}_{q+1} = a^{32}_qa^{21}_{q+1} \text{ and }
   e_\alpha = y_3a^{31}_{k-1}a^{21}_{q+1}a^{12}_ka^{32}_q, \indent
   \text{($a$-form)}.
$$
   The same for Line 10 of Table \ref{table_adm_elems}.
\qedsymbol \vspace{2mm}

 \underline{Lines 3--6, Lines 11--14 of Table \ref{table_adm_elems}}.
   The proof for Line 3 follows from
   Line 1 under substitution
   $1 \longrightarrow 3 \longrightarrow 2 \longrightarrow 1$.
   Under the same substitution the case of Line 11 follows
   from the case of Line 10. Similarly, by the same substitution
   we have
   \begin{equation*}
   \begin{split}
     & \text{Line 1} \Longrightarrow \text{Line 3}
     \Longrightarrow \text{Line 5},  \\
     & \text{Line 9}
     \Longrightarrow \text{Line 11}
     \Longrightarrow \text{Line 13}, \\
     & \text{Line 2}
     \Longrightarrow \text{Line 4}
     \Longrightarrow \text{Line 6},     \\
     & \text{Line 10}
     \Longrightarrow \text{Line 12}
     \Longrightarrow \text{Line 14}.
\qed \vspace{2mm}
   \end{split}
   \end{equation*}

 \underline{Lines 7-8, Table \ref{table_adm_elems}}.
  Let us take Line 7, where
  $e_\alpha = y_1A^{31}_ka^{12}_kA^{23}_ka^{23}_k$.
   By Proposition \ref{more_eq}, (2) we have
$$
   y_1A^{23}_ka^{23}_k = y_1A^{23}_kA^{32}_k \text{ and }
   e_\alpha = y_1A^{31}_ka^{12}_kA^{23}_kA^{32}_k.
$$
   Again,
$$
   y_1A^{32}_ka^{12}_k = y_1A^{32}_kA^{21}_k, \text{ so }
   e_\alpha = y_1A^{31}_kA^{21}_kA^{23}_kA^{32}_k \quad
   \text{($A$-form)}.
$$

   On the other hand,
   by Proposition \ref{more_eq}, (2) we have
$$
   y_1A^{31}_kA^{23}_k = y_1A^{31}_ka^{32}_k \text{ and }
   e_\alpha = y_1A^{31}_ka^{12}_ka^{32}_ka^{23}_k.
$$
   By Table \ref{table_atomic}, Line 5.1 we have
$$
   A^{31}_ka^{32}_k = a^{13}_ka^{32}_k, \text{ so }
   e_\alpha = y_1a^{13}_ka^{12}_ka^{32}_ka^{23}_k \quad
   \text{($a$-form)}.
$$
   The same for Line 8 of Table \ref{table_adm_elems}.
\qedsymbol

\footnotesize
\begin{table}[ht]
  \renewcommand{\arraystretch}{1.35}
  \begin{tabular} {|c||c|c|c|c|}
  \hline \hline
    $N$  & $\alpha$1
             & $\alpha$ (index of  $g_{\alpha0}$) =
             & $e_{\alpha1}$
             & $g_{\alpha0}$ =
                  \cr
             & index of $e_{\alpha1}$
             & $T^{123}[\tilde{\alpha}]$
                or $T^{12}[\tilde{\alpha}]$
             & $A$-form
             & $T^{123}[g_{\tilde{\alpha}0}]$
                or $T^{12}[g_{\tilde{\alpha}0}]$ \\
  \hline  \hline  
    $1$  & $(213)^{2p}(21)^{2k}$
     & $21(321)^{2p-1}(32)^{2k}$ =
     & $y_2A_{k-1}^{32}A_q^{12}A_k^{13}A_q^{31}$
     & $y_2A_{k-1}^{32}A_q^{31}A_k^{13}A_{q-1}^{12}$ \cr
     &
     & $T^{123}[13(213)^{2p-1}(21)^{2k}]$
     &
     & $\bigcap(y_2y_1 + A_{q-1}^{31}a_k^{23})$ \\
  \hline         
    $1{'}$ &    $(21)^{2k}$
     &  $2(12)^{2k-1}$ =
     & $y_2A_{k-1}^{32}A_k^{12}A_k^{13}A_k^{31}$
     & $y_2A_{k-1}^{32}A_{k-1}^{12}A_k^{13}A_k^{31}$    \cr
     &
     & $T^{12}[1(21)^{2k-1}]$
     &
     & $\bigcap(y_2a_k^{21}+y_2a_k^{23})$ \\
  \hline         
   $2$  & $3(213)^{2p}(21)^{2k}$
     & $(321)^{2p}(32)^{2k}$ =
     & $y_3A_{k-1}^{13}A_{q+1}^{12}A_k^{21}A_q^{23}$
     & $y_3A_{k-1}^{13}A_q^{23}A_k^{21}A_q^{12}$      \cr
     &
     & $T^{123}[(213)^{2p}(21)^{2k}]$
     &
     & $\bigcap(x_2 + a_q^{13}A_k^{13})$ \\
  \hline        
   $3$  & $13(213)^{2p}(21)^{2k}$
     & $1(321)^{2p}(32)^{2k}$ =
     & $y_1A_{k-1}^{21}A_{q+1}^{31}A_k^{32}A_{q+1}^{23}$
     & $y_1A_{k-1}^{21}A_{q+1}^{23}A_k^{32}A_q^{31}$ \cr
     &
     & $T^{123}[3(213)^{2p}(21)^{2k}]$
     &
     & $\bigcap(y_3y_1 + A_q^{23}a_k^{12})$ \\
  \hline         
   $4$  & $(213)^{2p+1}(21)^{2k}$
     & $21(321)^{2p}(32)^{2k}$ =
     & $y_2A_{k-1}^{32}A_{q+2}^{31}A_k^{13}A_{q+1}^{12}$
     & $y_2A_{k-1}^{32}A_{q+1}^{12}A_k^{13}A_{q+1}^{31}$ \cr
     &
     & $T^{123}[13(213)^{2p}(21)^{2k}]$
     &
     & $\bigcap(x_1 + a_{q+1}^{32}A_k^{32})$ \\
  \hline         
   $5$  & $3(213)^{2p+1}(21)^{2k}$
     & $(321)^{2p+1}(32)^{2k}$ =
     & $y_3A_{k-1}^{13}A_{q+2}^{23}A_k^{21}A_{q+2}^{12}$
     & $y_3A_{k-1}^{13}A_{q+2}^{12}A_k^{21}A_{q+1}^{23}$ \cr
     &
     & $T^{123}[(213)^{2p+1}(21)^{2k}]$
     &
     & $\bigcap(y_2y_3 + A_{q+1}^{12}a_k^{31})$ \\
  \hline         
   $6$  & $13(213)^{2p+1}(21)^{2k}$
     & $1(321)^{2p+1}(32)^{2k}$ =
     & $y_1A_{k-1}^{21}A_{q+3}^{23}A_k^{32}A_{q+2}^{31}$
     & $y_1A_{k-1}^{21}A_{q+2}^{31}A_k^{32}A_{q+2}^{23}$ \cr
     &
     & $T^{123}[3(213)^{2p+1}(21)^{2k}]$
     &
     & $\bigcap(x_3 + a_{q+2}^{21}A_k^{21})$ \\
  \hline \hline   
   $7$  & $1(21)^{2k}$
     & $(12)^{2k}$ =
     & $y_1A_k^{31}A_k^{21}A_k^{23}A_k^{32}$
     & $y_1A_{k-1}^{31}A_k^{21}A_k^{23}A_k^{32}$ \cr
     &
     & $T^{12}[(21^{2k}]$
     &
     & $\bigcap(x_2 + a_k^{31}A_k^{31})$  \\
  \hline          
   $8$  & $1(21)^{2k+1}$
     & $(12)^{2k+1}$ =
     & $y_1A_k^{31}A_k^{21}A_{k+1}^{23}A_{k+1}^{32}$
     & $y_1A_k^{31}A_k^{21}A_k^{23}A_{k+1}^{32}$    \cr
     &
     & $T^{12}[(21)^{2k}]$
     &
     & $\bigcap(x_1 + a_k^{23}A_{k+1}^{23})$ \\
  \hline  \hline  
   $9$  & $(213)^{2p}(21)^{2k+1}$
     & $21(321)^{2p-1}(32)^{2k+1}$ =
     & $y_2A_k^{32}A_q^{12}A_k^{13}A_{q+1}^{31}$
     & $y_2A_k^{32}A_q^{31}A_k^{13}A_q^{12}$ \cr
     &
     & $T^{123}[13(213)^{2p-1}(21)^{2k+1}]$
     &
     & $\bigcap(y_2y_3 + A_{q-1}^{12}a_{k+1}^{31})$ \\
  \hline          
   $9{'}$  & $(21)^{2k+1}$
     & $2(12)^{2k}$ =
     &  $y_2A_k^{32}A_k^{12}A_k^{13}A_{k+1}^{31}$
     &  $y_2A_k^{32}A_k^{12}A_k^{13}A_k^{31}$ \cr
     &
     & $T^{12}[1(21)^{2k}]$
     &
     & $\bigcap(y_2a_k^{21} + y_2a_k^{23})$ \\
  \hline         
   $10$  & $3(213)^{2p}(21)^{2k+1}$
     & $(321)^{2p}(32)^{2k+1}$ =
     & $y_3A_k^{13}A_{q+1}^{12}A_k^{21}A_{q+1}^{23}$
     & $y_3A_k^{13}A_q^{23}A_k^{21}A_{q+1}^{12}$    \cr
     &
     & $T^{123}[(213)^{2p}(21)^{2k+1}]$
     &
     & $\bigcap(x_3 + a_q^{21}A_{k+1}^{21})$ \\
  \hline        
   $11$  & $13(213)^{2p}(21)^{2k+1}$
     & $1(321)^{2p}(32)^{2k+1}$ =
     & $y_1A_k^{21}A_{q+1}^{31}A_k^{32}A_{q+2}^{23}$
     & $y_1A_k^{21}A_{q+1}^{23}A_k^{32}A_{q+1}^{31}$ \cr
     &
     & $T^{123}[3(213)^{2p}(21)^{2k+1}]$
     &
     & $\bigcap(y_2y_1 + A_q^{31}a_{k+1}^{23})$ \\
  \hline          
   $12$  & $(213)^{2p+1}(21)^{2k+1}$
     & $21(321)^{2p}(32)^{2k+1}$ =
     & $y_2A_k^{32}A_{q+2}^{31}A_k^{13}A_{q+2}^{12}$
     & $y_2A_k^{32}A_{q+1}^{12}A_k^{13}A_{q+2}^{31}$ \cr
     &
     & $T^{123}[(213)^{2p+1}(21)^{2k+1}]$
     &
     & $\bigcap(x_2 + a_{q+1}^{13}A_{k+1}^{13})$ \\
  \hline         
   $13$  & $3(213)^{2p+1}(21)^{2k+1}$
     & $(321)^{2p+1}(32)^{2k+1}$ =
     & $y_3A_k^{13}A_{q+2}^{23}A_k^{21}A_{q+3}^{12}$
     & $y_3A_k^{13}A_{q+2}^{12}A_k^{21}A_{q+2}^{23}$ \cr
     &
     & $T^{123}[(213)^{2p+1}(21)^{2k+1}]$
     &
     & $\bigcap(y_3y_1 + A_{q+1}^{23}a_{k+1}^{12})$ \\
  \hline         
   $14$  & $13(213)^{2p+1}(21)^{2k+1}$
     & $1(321)^{2p+1}(32)^{2k+1}$ =
     & $y_1A_k^{21}A_{q+3}^{23}A_k^{32}A_{q+3}^{31}$
     & $y_1A_k^{21}A_{q+2}^{31}A_k^{32}A_{q+3}^{23}$ \cr
     &
     & $T^{123}[3(213)^{2p+1}(21)^{2k+1}]$
     &
     & $\bigcap(x_1 + a_{q+2}^{32}A_{k+1}^{32})$ \\
  \hline \hline
  \end{tabular}
  \vspace{2mm}
  \caption{
   \hspace{3mm}$\alpha1$, $\alpha$, $e_{\alpha1}$ and $g_{\alpha0}$}
  \label{table_alpha}
  \begin{center}
    For Line 1 and Line 9: $p > 0$.
    For Lines 1--6: $k > 0, p \geq 0$. \\
    For Lines 7--8: $k > 0$.
    For Lines 9--14: $k \geq 0, p \geq 0$.
    For all lines: $q = k + 3p$.\\
  \end{center}
 \end{table}
 \normalsize

 \begin{remark} [to Table \ref{table_alpha}]
 {\rm
 Transformation $T^{123}$ substitutes
$1 \rightarrow 2 \rightarrow 3 \rightarrow 1$ and $T^{12}$
substitutes $1 \rightarrow 2 \rightarrow 1$. These transformations
are needed for getting corresponding admissible elements
$g_{\alpha0}$ from the Table \ref{table_adm_elems} or their
$a$-forms and $A$-forms from Table \ref{table_two_forms}. For
example, consider Line 2. In this case $\alpha =
(321)^{2p}(32)^{2k}$. This sequence can be obtained by applying
$T^{123}$ to the sequence $\tilde{\alpha} = (213)^{2p}(21)^{2k}$
from Table \ref{table_adm_elems}, i.e.,
$$
  \alpha =
   (321)^{2p}(32)^{2k} = T^{123}[\tilde{\alpha}] =
   T^{123}[(213)^{2p}(21)^{2k}].
$$
}
\end{remark}

\subsection{Admissible sequences ${\alpha}1$ and $\alpha$. The
 inclusion theorem.}
  \label{proof_incl_th}
   Our intention is to prove the following theorem. \\

{\bf Theorem \ref{inclusion}}
  {\it For every admissible sequence
  $\alpha$ or ${\alpha}i$, where $i = 1,2,3$ from
  Table \ref{table_admissible}, the following inclusion holds:}
  \begin{equation*}
      e_{{\alpha}i} \subseteq g_{\alpha0}, \indent i=1,2,3.
  \end{equation*}

   Without loss of generality it suffices
   to prove Theorem \ref{inclusion} for $i = 1$.
   Thus, we need description of admissible
   sequences ${\alpha}1$ and $\alpha$ in the form of sequences
   of Table \ref{table_adm_elems}.
\begin{proposition}  \label{apl_alp1}
  For the corresponding indices ${\alpha}1$ and $\alpha$ from
  Table \ref{table_alpha} (columns $2$ and $3$), we have
  \begin{equation*}
      \varphi_{\alpha}(1) = \alpha1.
  \end{equation*}
\end{proposition}

\PerfProof \underline{Line 1}. Let $p > 0$. We have
\begin{equation*}
\begin{split}
  & \varphi_{21(321)^{2p-1}(32)^{2k}}(1) =
    21(321)^{2p-1}(32)^{2k}1 =   \\
  & 21(321)(321)^{2p-2}(32)^{2k}1 =
    213((21)(321)^{2p-2}(32)^{2k}1).
\end{split}
\end{equation*}
By induction, from Line 4, Table \ref{table_alpha} we get
$$
   \varphi_{21(321)^{2p-1}(32)^{2k}}(1) =
   213[(213)^{2p-1}(21)^{2k}] = (213)^{2p}(21)^{2k} =
   \alpha1.
$$
If $p=0$, then
   $2(12)^{2k}1 = (21)^{2k+1}$.
   The same for Line 9.
\qedsymbol

 \underline{Line 4}. Here,
$$
   21(321)^{2p}(32)^{2k}1 =
   21(321)(321)^{2p-1}(32)^{2k}1 =
   213((21)(321)^{2p-1}(32)^{2k}1).
$$
   By induction, from heading 1) we get
$$
    21(321)^{2p}(32)^{2k}1 =
    (213)(213)^{2p}(21)^{2k} =
    (213)^{2p+1}(21)^{2k} = \alpha1.
$$
   The same for Line 12.
\qedsymbol

   Other lines of Table \ref{table_alpha} are considered
   by analogy with these two cases.
\qedsymbol

 \index{admissible elements! - $S$-form (symmetric) of
$g_{{\alpha}0}$}

\subsection{An $S$-form of $g_{{\alpha}0}$}

\begin{definition}
{\rm The form of the elements $g_{\alpha0}$ from the right column
of Table \ref{table_symmetric} is said to be an {\it $S$-form}
(symmetric).}
\end{definition}

\begin{proposition}
Both forms of the elements $g_{\alpha0}$ from Table
 \ref{table_symmetric} coincide.
\end{proposition}
\large
\begin{table}[h]
 \renewcommand{\arraystretch}{1.35}
  \begin{tabular} {|c||c|c|c|c|c|}
  \hline \hline
    $N$  & $g_{\alpha0}$
             & $g_{\alpha0}$ \cr
             & from Table \ref{table_alpha}
             & $S$-form \\
  \hline  \hline  
    $1$
     & $y_2A_{k-1}^{32}A_q^{31}A_k^{13}A_{q-1}^{12}
        (y_2y_1 + A_{q-1}^{31}a_k^{23})$
     & $y_2A_{k-1}^{32}A_q^{31}A_k^{13}A_{q-1}^{12}
        (a^{21}_q + a^{12}_{k+1})$
     \\
  \hline         
    $1{'}$
     & $y_2A_{k-1}^{32}A_{k-1}^{12}A_k^{13}A_k^{31}
       (y_2a_k^{21}+y_2a_k^{23})$
     & $y_2A_{k-1}^{32}A_{k-1}^{12}A_k^{13}A_k^{31}
        (a^{21}_k + a^{12}_{k+1})$
     \\
  \hline         
     $2$
     & $y_3A_{k-1}^{13}A_q^{23}A_k^{21}A_q^{12}
        (x_2 + a_q^{13}A_k^{13})$
     & $y_3A_{k-1}^{13}A_q^{23}A_k^{21}A_q^{12}
       (A^{12}_{q+1} + A^{13}_k)$
     \\
  \hline        
    $3$
     & $y_1A_{k-1}^{21}A_{q+1}^{23}A_k^{32}A_q^{31}
        (y_3y_1 + A_q^{23}a_k^{12})$
     &$y_1A_{k-1}^{21}A_{q+1}^{23}A_k^{32}A_q^{31}
       (a^{13}_{q+1} + a^{31}_{k+1})$
     \\
  \hline         
    $4$
     & $y_2A_{k-1}^{32}A_{q+1}^{12}A_k^{13}A_{q+1}^{31}
        (x_1 + a_{q+1}^{32}A_k^{32})$
     & $y_2A_{k-1}^{32}A_{q+1}^{12}A_k^{13}A_{q+1}^{31}
       (A^{31}_{q+2} + A^{32}_k)$
     \\
  \hline         
    $5$
     & $y_3A_{k-1}^{13}A_{q+2}^{12}A_k^{21}A_{q+1}^{23}
        (y_2y_3 + A_{q+1}^{12}a_k^{31})$
     & $y_3A_{k-1}^{13}A_{q+2}^{12}A_k^{21}A_{q+1}^{23}
       (a^{32}_{q+2} + a^{23}_{k+1})$
     \\
  \hline         
     $6$
     & $y_1A_{k-1}^{21}A_{q+2}^{31}A_k^{32}A_{q+2}^{23}
        (x_3 + a_{q+2}^{21}A_k^{21})$
     & $y_1A_{k-1}^{21}A_{q+2}^{31}A_k^{32}A_{q+2}^{23}
       (A^{23}_{q+3} + A^{21}_k)$
     \\
  \hline \hline   
    $7$
     & $y_1A_{k-1}^{31}A_k^{21}A_k^{23}A_k^{32}
        (x_2 + a_k^{31}A_k^{31})$
     & $y_1A_{k-1}^{31}A_k^{21}A_k^{23}A_k^{32}
      (A^{32}_{k+1} + A^{31}_k)$
     \\
  \hline          
    $8$
     & $y_1A_k^{31}A_k^{21}A_k^{23}A_{k+1}^{32}
        (x_1 + a_k^{23}A_{k+1}^{23})$
     & $y_1A_k^{31}A_k^{21}A_k^{23}A_{k+1}^{32}
      (A^{21}_{k+1} + A^{23}_{k+1})$
     \\
  \hline  \hline  
    $9$
     & $y_2A_k^{32}A_q^{31}A_k^{13}A_q^{12}
        (y_2y_3 + A_{q-1}^{12}a_{k+1}^{31})$
     & $y_2A_k^{32}A_q^{31}A_k^{13}A_q^{12}
      (a^{13}_{q+1} + a^{31}_{k+1})$
     \\
  \hline          
    $9{'}$
     & $y_2A_k^{32}A_k^{12}A_k^{13}A_k^{31}
        (y_1a_k^{12} + y_3a_k^{32})$
     & $y_2A_k^{32}A_k^{12}A_k^{13}A_k^{31}
      (a^{31}_{k+1} + a^{13}_{k+1})$
     \\
  \hline         
   $10$
     & $y_3A_k^{13}A_q^{23}A_k^{21}A_{q+1}^{12}
        (x_3 + a_q^{21}A_{k+1}^{21})$
     & $y_3A_k^{13}A_q^{23}A_k^{21}A_{q+1}^{12}
      (A^{23}_{q+1} + A^{21}_{k+1})$
     \\
  \hline        
    $11$
     & $y_1A_k^{21}A_{q+1}^{23}A_k^{32}A_{q+1}^{31}
          (y_2y_1 + A_q^{31}a_{k+1}^{23})$
     & $y_1A_k^{21}A_{q+1}^{23}A_k^{32}A_{q+1}^{31}
      (a^{32}_{q+2} + a^{23}_{k+1})$
     \\
  \hline          
   $12$
     & $y_2A_k^{32}A_{q+1}^{12}A_k^{13}A_{q+2}^{31}
          (x_2 + a_{q+1}^{13}A_{k+1}^{13})$
     & $y_2A_k^{32}A_{q+1}^{12}A_k^{13}A_{q+2}^{31}
      (A^{12}_{q+2} + A^{13}_{k+1})$
     \\
  \hline         
   $13$
     & $y_3A_k^{13}A_{q+2}^{12}A_k^{21}A_{q+2}^{23}
          (y_3y_1 + A_{q+1}^{23}a_{k+1}^{12})$
     & $y_3A_k^{13}A_{q+2}^{12}A_k^{21}A_{q+2}^{23}
      (a^{21}_{q+3} + a^{12}_{k+1})$
     \\
  \hline         
   $14$
     & $y_1A_k^{21}A_{q+2}^{31}A_k^{32}A_{q+3}^{23}
          (x_1 + a_{q+2}^{32}A_{k+1}^{32})$
     & $y_1A_k^{21}A_{q+2}^{31}A_k^{32}A_{q+3}^{23}
      (A^{31}_{q+3} + A^{32}_{k+1})$
     \\
  \hline \hline
  \end{tabular}
  \vspace{2mm}
  \caption{
   \hspace{3mm}The $S$-forms of $g_{\alpha0}$}
  \label{table_symmetric}
\end{table}

\normalsize

\PerfProof   1) Since $e_{\alpha} \subseteq A_q^{31} \subseteq
A_{q-1}^{31}$, we see that
$$
  e_{\alpha}(y_2y_1 + A_{q-1}^{31}a_k^{23}) =
      e_{\alpha}(y_2y_1A_{q-1}^{31} + a_k^{23}).
$$

Since $y_2y_1A_{q-1}^{31} = y_2y_1a_{q-1}^{13}$ and
      $e_{\alpha} \subseteq y_2$,
it follows that
$$
      e_{\alpha}(y_2y_1 + A_{q-1}^{31}a_k^{23}) =
      e_{\alpha}(y_2y_1a_{q-1}^{13} + a_k^{23}) =
      e_{\alpha}(y_1a_{q-1}^{13} + y_2a_k^{23}).
$$
Finally, since $x_2 \subseteq y_2a_k^{23}$ and
      $x_1  \subseteq y_1a_{q-1}^{13}$, it follows that
$$
      e_{\alpha}(y_2y_1 + A_{q-1}^{31}a_k^{23}) =
      e_{\alpha}((x_2 + y_1a_{q-1}^{13}) + (x_1 + y_2a_k^{23})) =
      e_{\alpha}(a_q^{21} + a_{k+1}^{12}).
      \qed \vspace{2mm}
$$

 1$'$) Since $e_{\alpha} \subseteq y_2$, we have
$$
      e_{\alpha}(y_2a^{21}_k + y_2a_k^{23}) =
      e_{\alpha}(a^{21}_k + y_2a_k^{23}).
$$
      Since $x_1 \subseteq a_k^{21}$, we have
$$
      e_{\alpha}(y_2a^{21}_k + y_2a_k^{23}) =
      e_{\alpha}(a^{21}_k + (x_1 + y_2a_k^{23})) =
      e_{\alpha}(a^{21}_k + a_{k+1}^{12}).
      \qed \vspace{2mm}
$$

 2) Since
$$
  a^{13}_q = x_1 + y_3a^{32}_{q-1} \supseteq y_3a^{32}_{q-1}
     \supseteq y_3a^{32}_q \supseteq e_{\alpha}
$$
(Table \ref{table_two_forms}, Line 2, $a$-form), it follows that
$$
      e_{\alpha}(x_2 + a^{13}_qA^{13}_k) =
      e_{\alpha}(x_2a^{13}_q + A^{13}_k) =
      e_{\alpha}(x_2A^{31}_q + A^{13}_k).
$$
      Further, since $y_1 \subseteq  A^{13}_k$, it follows that
$$
      e_{\alpha}(x_2 + a^{13}_qA^{13}_k) =
      e_{\alpha}(y_1 + x_2A^{31}_q + A^{13}_k) =
      e_{\alpha}(A^{12}_{q+1} + A^{13}_k).
      \qed \vspace{2mm}
$$

 3) Since $e_{\alpha} \subseteq A_{q+1}^{23} \subseteq  A_q^{23}$,
we see that
$$
      e_{\alpha}(y_3y_1 + A_q^{23}a_k^{12}) =
      e_{\alpha}(y_3y_1A_q^{23} + a_k^{12}) =
      e_{\alpha}(y_3y_1a_q^{32} + a_k^{12}).
$$
Since $e_{\alpha} \subseteq y_1$, it follows that
\begin{equation*}
 \begin{split}
     & e_{\alpha}(y_3y_1a_q^{32} + a_k^{12}) =
       e_{\alpha}(y_3a_q^{32} + y_1a_k^{12}) = \\
     & e_{\alpha}(x_1 + y_3a_q^{32} + x_3 + y_1a_k^{12}) =
       e_{\alpha}(a_{q+1}^{13} + a_{k+1}^{31}).
       \qed \vspace{2mm}
 \end{split}
\end{equation*}

 4) Since $a^{32}_{q+1} = x_3 + y_2a^{21}_q \supseteq
y_2a^{21}_{q+1}
     \supseteq e_{\alpha}$ (Table \ref{table_two_forms}, Line 4, $a$-form),
      we see that
\begin{equation*}
\begin{split}
   &  e_{\alpha}(x_1 + a^{32}_{q+1}A^{32}_k) =
      e_{\alpha}(x_1a^{32}_{q+1} + A^{32}_k) =
      e_{\alpha}(x_1A^{23}_{q+1} + A^{32}_k) =  \\
   &  e_{\alpha}(y_3 + x_1A^{23}_{q+1} + A^{32}_k) =
      e_{\alpha}(A^{31}_{q+2} + A^{32}_k).
      \qed \vspace{2mm}
\end{split}
\end{equation*}

 5) Since $e_{\alpha} \subseteq y_3A^{12}_{q+2}
      \subseteq y_3A^{12}_{q+1}$ and $y_3 \supseteq  e_{\alpha}$,
      we see that
\begin{equation*}
\begin{split}
    & e_{\alpha}(y_2y_3 + A^{12}_{q+1}a^{31}_k) =
      e_{\alpha}(y_2y_3A^{12}_{q+1} + a^{31}_k) =
      e_{\alpha}(y_2a^{21}_{q+1} + y_3a^{31}_k) = \\
    & e_{\alpha}((x_3 + y_2a^{21}_{q+1}) + (x_2 +y_3a^{31}_k)) =
      e_{\alpha}(a^{32}_{q+2} + a^{23}_{k+1}).
      \qed \vspace{2mm}
\end{split}
\end{equation*}

 6) Since $a^{21}_{q+2} = x_2 + y_1a^{13}_{q+1}
\supseteq y_1a^{13}_{q+1}
    \supseteq e_{\alpha}$ (Table \ref{table_two_forms}, Line 6, $a$-form),
      it follows that
\begin{equation*}
\begin{split}
    &  e_{\alpha}(x_3 + a^{21}_{q+2}A^{21}_k) =
       e_{\alpha}(x_3a^{21}_{q+2} + A^{21}_k) =
       e_{\alpha}(x_3A^{12}_{q+2} + A^{21}_k) = \\
    &  e_{\alpha}(y_2 + x_3A^{12}_{q+2} + A^{21}_k) =
       e_{\alpha}(A^{23}_{q+3} + A^{21}_k).
      \qed \vspace{2mm}
\end{split}
\end{equation*}

 7) By (Table \ref{table_two_forms}, Line 7, $a$-form)
      $a^{31}_k = x_3 + y_1a^{12}_{k-1} \supseteq y_1a^{12}_{k-1}
      \supseteq y_1a^{12}_k \supseteq e_{\alpha}$,
      and therefore
\begin{equation*}
\begin{split}
     & e_{\alpha}(x_2 + a^{31}_kA^{31}_k) =
       e_{\alpha}(x_2a^{31}_k + A^{31}_k) =
       e_{\alpha}(x_2A^{13}_k + A^{31}_k) = \\
     & e_{\alpha}(y_3 + x_2A^{13}_k + A^{31}_k) =
       e_{\alpha}(A^{32}_{k+1} + A^{31}_k).
      \qed \vspace{2mm}
\end{split}
\end{equation*}

 8) By (Table \ref{table_two_forms}, Line 8, $a$-form)
      $a^{23}_k  \supseteq a^{23}_{k+1} \supseteq e_{\alpha}$
      and we have
\begin{equation*}
\begin{split}
     & e_{\alpha}(x_1 + a^{23}_kA^{23}_{k+1}) =
       e_{\alpha}(x_1a^{23}_k + A^{23}_{k+1}) = \\
     & e_{\alpha}(y_2 + x_1A^{32}_k + A^{23}_{k+1}) =
       e_{\alpha}(A^{21}_{k+1} + A^{23}_{k+1}).
       \qed \vspace{2mm}
\end{split}
\end{equation*}

 9) Since $A^{12}_{q-1} \supseteq A^{12}_q \supseteq
e_{\alpha}$,
   we see that
\begin{equation*}
\begin{split}
    &  e_{\alpha}(y_2y_3 + A^{12}_{q-1}a^{31}_{k+1}) =
       e_{\alpha}(y_2y_3A^{12}_{q-1} + a^{31}_{k+1}) =
       e_{\alpha}(y_2y_3a^{21}_{q-1} + a^{31}_{k+1}) = \\
    &  e_{\alpha}(x_3 + y_2y_3a^{21}_{q-1} + a^{31}_{k+1}) =
       e_{\alpha}(y_3(x_3 + y_2a^{21}_{q-1}) + a^{31}_{k+1}) =
       e_{\alpha}(y_3a^{32}_q + a^{31}_{k+1}) =  \\
    &  e_{\alpha}(x_1 + y_3a^{32}_q + a^{31}_{k+1}) =
       e_{\alpha}(a^{13}_{q+1} + a^{31}_{k+1}).
      \qed \vspace{2mm}
\end{split}
\end{equation*}

 9$'$) Here,
$$
      e_{\alpha}(y_1a^{12}_k + y_3a^{32}_k) =
      e_{\alpha}((x_3 + y_1a^{12}_k) + (x_1 + y_3a^{32}_k)) =
      e_{\alpha}(a^{31}_{k+1} + a^{13}_{k+1}).
      \qed \vspace{2mm}
$$

 10) By (Table \ref{table_two_forms}, Line 10, $a$-form)
      $a^{21}_q  \supseteq a^{21}_{q+1} \supseteq e_{\alpha}$,
      thus,
\begin{equation*}
\begin{split}
    &  e_{\alpha}(x_3 + a^{21}_qA^{21}_{k+1}) =
       e_{\alpha}(x_3a^{21}_q + A^{21}_{k+1}) =
       e_{\alpha}(x_3A^{12}_q + A^{21}_{k+1}) =  \\
    &  e_{\alpha}(y_2 + x_3A^{12}_q + A^{21}_{k+1}) =
       e_{\alpha}(A^{23}_{q+1} + A^{21}_{k+1}).
      \qed \vspace{2mm}
\end{split}
\end{equation*}

 11) Since $A^{31}_q \supseteq A^{31}_{q+1} \supseteq
e_{\alpha}$,
      we get
\begin{equation*}
\begin{split}
    &  e_{\alpha}(y_2y_1 + A^{31}_qa^{23}_{k+1}) =
       e_{\alpha}(y_2y_1A^{31}_q + a^{23}_{k+1}) =
       e_{\alpha}(y_2y_1a^{13}_q + a^{23}_{k+1}) = \\
    &  e_{\alpha}(x_2 + y_2y_1a^{13}_q + a^{23}_{k+1}) =
       e_{\alpha}(y_2(x_2 + y_1a^{13}_q) + a^{23}_{k+1}) =
       e_{\alpha}(y_2a^{21}_{q+1} + a^{23}_{k+1}) = \\
    &  e_{\alpha}(x_3 + y_2a^{21}_{q+1} + a^{23}_{k+1}) =
       e_{\alpha}(a^{32}_{q+2} + a^{23}_{k+1}).
      \qed \vspace{2mm}
\end{split}
\end{equation*}

 12) By (Table \ref{table_two_forms}, Line 12, $a$-form)
      $a^{13}_{q+1}  \supseteq a^{13}_{q+2} \supseteq e_{\alpha}$
      and we get
\begin{equation*}
\begin{split}
    &  e_{\alpha}(x_2 + a^{13}_{q+1}A^{13}_{k+1}) =
       e_{\alpha}(x_2a^{13}_{q+1} + A^{13}_{k+1}) = \\
    &  e_{\alpha}(y_1 + x_2A^{31}_{q+1} + A^{32}_{k+1}) =
       e_{\alpha}(A^{12}_{q+2} + A^{32}_{k+1}).
      \qed \vspace{2mm}
\end{split}
\end{equation*}

 13) Here,
\begin{equation*}
\begin{split}
   & e_{\alpha}(y_3y_1 + A^{23}_{q+1}a^{12}_{k+1}) =
     e_{\alpha}(y_3y_1A^{23}_{q+1} + a^{12}_{k+1}) =
     e_{\alpha}(y_3y_1a^{32}_{q+1} + a^{12}_{k+1}) =  \\
   & e_{\alpha}(x_1 + y_3y_1a^{32}_{q+1} + a^{12}_{k+1}) =
     e_{\alpha}(y_1(x_1 + y_3a^{32}_{q+1}) + a^{12}_{k+1}) =
     e_{\alpha}(y_1a^{13}_{q+2} + a^{12}_{k+1}) = \\
   & e_{\alpha}((x_2 + y_1a^{13}_{q+2}) + a^{12}_{k+1}) =
     e_{\alpha}(a^{21}_{q+3} + a^{12}_{k+1}).
      \qed \vspace{2mm}
\end{split}
\end{equation*}

 14) By (Table \ref{table_two_forms}, Line 14, $a$-form)
      $a^{32}_{q+2}  \supseteq a^{32}_{q+3} \supseteq e_{\alpha}$
      and we have
\begin{equation*}
\begin{split}
    & e_{\alpha}(x_1 + a^{32}_{q+2}A^{32}_{k+1}) =
      e_{\alpha}(x_1a^{32}_{q+2} + A^{32}_{k+1}) =  \\
    & e_{\alpha}(y_3 + x_1a^{32}_{q+2} + A^{32}_{k+1}) =
      e_{\alpha}(y_3 + x_1A^{23}_{q+2} + A^{32}_{k+1}) =
      e_{\alpha}(A^{31}_{q+3} + A^{32}_{k+1}).
      \qed \vspace{2mm}
\end{split}
\end{equation*}

 After all previous preparations the proof of
Theorem \ref{inclusion} is rather trivial. Instead of inclusion
\begin{equation*}
      e_{{\alpha}i} \subseteq g_{\alpha0},
         \text{ where } i=1,2,3,
\end{equation*}
we will prove a sharper statement
\begin{proposition} \label{strong_eq}
For the lattice elements $e_{{\alpha}i}$ and $g_{\alpha0}$, we
have
\begin{equation*}
      e_{{\alpha}i} = g_{\alpha0}Z, \indent i=1,2,3,
\end{equation*}
where the element $Z$ for $i = 1$ is given in Table
 \ref{tab_instead_inclusion}.
\end{proposition}
\begin{table}
 \renewcommand{\arraystretch}{1.35}
  \begin{tabular} {|c||c|c|c|}
  \hline \hline
    $N$  & $\alpha$1
             & $\alpha$
             & $e_{{\alpha}1} = g_{\alpha0}Z$ \cr
             & index of $e_{\alpha1}$
             & index of $g_{\alpha0}$
             &   \\
  \hline  \hline  
    $1$  & $(213)^{2p}(21)^{2k}$
     & $21(321)^{2p-1}(32)^{2k}$
     & $e_{{\alpha}1} = g_{\alpha0}A^{12}_q$ \\
  \hline         
   $2$  & $3(213)^{2p}(21)^{2k}$
     & $(321)^{2p}(32)^{2k}$
     & $e_{{\alpha}1} = g_{\alpha0}A^{12}_{q+1}$ \\
  \hline        
   $3$  & $13(213)^{2p}(21)^{2k}$
     & $1(321)^{2p}(32)^{2k}$
     & $e_{{\alpha}1} = g_{\alpha0}A^{31}_{q+1}$   \\
  \hline         
   $4$  & $(213)^{2p+1}(21)^{2k}$
     & $21(321)^{2p}(32)^{2k}$
     & $e_{{\alpha}1} = g_{\alpha0}A^{31}_{q+2}$   \\
  \hline         
   $5$  & $3(213)^{2p+1}(21)^{2k}$
     & $(321)^{2p+1}(32)^{2k}$
     & $e_{{\alpha}1} = g_{\alpha0}A^{23}_{q+2}$  \\
  \hline         
   $6$  & $13(213)^{2p+1}(21)^{2k}$
     & $1(321)^{2p+1}(32)^{2k}$
     & $e_{{\alpha}1} = g_{\alpha0}A^{23}_{q+3}$ \\
  \hline \hline   
   $7$  & $1(21)^{2k}$
     & $(12)^{2k}$
     & $e_{{\alpha}1} = g_{\alpha0}A^{31}_k$  \\
  \hline          
   $8$  & $1(21)^{2k+1}$
     & $(12)^{2k+1}$
     & $e_{{\alpha}1} = g_{\alpha0}A^{23}_{k+1}$ \\
  \hline  \hline  
   $9$  & $(213)^{2p}(21)^{2k+1}$
     & $21(321)^{2p-1}(32)^{2k+1}$
     & $e_{{\alpha}1} = g_{\alpha0}A^{31}_{q+1}$  \\
  \hline         
   $10$  & $3(213)^{2p}(21)^{2k+1}$
     & $(321)^{2p}(32)^{2k+1}$
     & $e_{{\alpha}1} = g_{\alpha0}A^{23}_{q+1}$ \\
  \hline        
   $11$  & $13(213)^{2p}(21)^{2k+1}$
     & $1(321)^{2p}(32)^{2k+1}$
     & $e_{{\alpha}1} = g_{\alpha0}A^{23}_{q+2}$ \\
  \hline          
   $12$  & $(213)^{2p+1}(21)^{2k+1}$
     & $21(321)^{2p}(32)^{2k+1}$
     & $e_{{\alpha}1} = g_{\alpha0}A^{12}_{q+2}$ \\
  \hline         
   $13$  & $3(213)^{2p+1}(21)^{2k+1}$
     & $(321)^{2p+1}(32)^{2k+1}$
     & $e_{{\alpha}1} = g_{\alpha0}A^{12}_{q+3}$ \\
  \hline         
   $14$  & $13(213)^{2p+1}(21)^{2k+1}$
     & $1(321)^{2p+1}(32)^{2k+1}$
     & $e_{{\alpha}1} = g_{\alpha0}A^{31}_{q+3}$  \\
  \hline \hline
  \end{tabular}
  \vspace{2mm}
  \caption{
   \hspace{3mm}The relation $e_{{\alpha}1} = g_{\alpha0}Z$}
  \label{tab_instead_inclusion}
  \end{table}

 \PerfProofW of Proposition \ref{strong_eq} follows from
 the $S$-form of elements $g_{\alpha0}$ from Table \ref{table_symmetric}
 and the $a$-form and $A$-form of the elements $e_{{\alpha}1}$
 from Table \ref{table_two_forms}.
 The polynomial $g_{\alpha0}$ from Table \ref{table_symmetric}
 is the intersection of $M$ and $P$:
 \begin{equation}  \label{MP}
           g_{\alpha0} = M \cap P,
 \end{equation}
 where $M$ is the expression located on the left of the brackets in
 the $S$-form of Table \ref{table_two_forms} and
 $P$ is the sum contained in the parantheses of this $S$-form.
 Then, for every line
 from Table \ref{table_two_forms}, we have
 \begin{equation}  \label{MZ}
           MZ = e_{{\alpha}1} \text{ and }
             e_{{\alpha}1} \subseteq P.
 \end{equation}
 From (\ref{MP}) and (\ref{MZ}) it follows that
 \begin{equation*}
         g_{\alpha0}Z =  MPZ = e_{{\alpha}1}P =
             e_{{\alpha}1}.
 \end{equation*}

 For example, consider Lines 1 and 2 of Table \ref{tab_instead_inclusion}.

 \underline{Line 1}. We have
\begin{equation*}
\begin{split}
   & g_{\alpha0} = y_2A_{k-1}^{32}A_q^{31}A_k^{13}A_{q-1}^{12}
        (a^{21}_q + a^{12}_{k+1}), \\
   & M = y_2A_{k-1}^{32}A_q^{31}A_k^{13}A_{q-1}^{12}, \quad
     P = a^{21}_q + a^{12}_{k+1}, \quad
     Z = A_q^{12}.
\end{split}
\end{equation*}
   Since $A_q^{12} \subseteq A_{q-1}^{12}$, it follows that
$$
   MZ = y_2A_{k-1}^{32}A_q^{31}A_k^{13}A_q^{12} = e_{{\alpha}1}
$$
   and by {\it $a$-form} from Table \ref{table_two_forms} we obtain
   $e_{{\alpha}1} \subseteq a^{21}_q \subseteq P$.
\qedsymbol

 \underline{Line 2}. We have
\begin{equation*}
\begin{split}
  & g_{\alpha0} = y_3A_{k-1}^{13}A_q^{23}A_k^{21}A_q^{12}
       (A^{12}_{q+1} + A^{13}_k), \\
  & M = y_3A_{k-1}^{13}A_q^{23}A_k^{21}A_q^{12},  \quad
    P = A^{12}_{q+1} + A^{13}_k,  \quad
    Z = A_{q+1}^{12}.
\end{split}
\end{equation*}
   Since $A_{q+1}^{12} \subseteq A_q^{12}$, we see that
$$
   MZ = y_2A_{k-1}^{32}A_q^{31}A_k^{13}A_{q+1}^{12} = e_{{\alpha}1}.
$$
   In addition,
   $e_{{\alpha}1} \subseteq A^{12}_{q+1} \subseteq P$.
\qedsymbol

Proposition \ref{strong_eq} together with Theorem \ref{inclusion}
(Inclusion Theorem) are proven.

\section{Properties of $H^+(n)$}
   \label{properties}

\subsection{The perfectness of $a_i(1)$}
 \label{ai_perfect}
Without loss of generality we consider $i=1$; we have
\begin{equation*}
\begin{split}
  \rho(a_1(1)) = &
   \rho(x_2(1) + x_3(1)) + \rho(y_2(2) + y_3(2)) = \\
 & \rho(x_2 + x_3) +
    \sum\varphi_i\Phi^+\rho(y_2(1) + y_3(1)) = \\
 & X_2 + X_3 +
    \sum\varphi_i\Phi^+\rho(y_2 + y_3).
\end{split}
\end{equation*}

 \underline{(a) Case $\Phi^+\rho(y_2 + y_3) \neq 0$}. In this case,
$\rho(a_1(1)) \neq 0$. Since $y_2 + y_3$ is perfect
 (Proposition \ref{first_perfect}), we see that
 $\Phi^+\rho(y_2 + y_3) = X_0^1$ and by Corollary \ref{cor_psi} we
obtain
\begin{equation} \label{rhoa1}
\begin{split}
 \rho(a_1(1)) = & X_2 + X_3 + \sum\varphi_i{X}_0^1=
  X_2 + X_3 + \sum\varphi_i\Phi^+\rho(I) = \\
& X_2 + X_3 + (Y_1 + Y_2)(Y_1 + Y_3)(Y_2 + Y_3).
\end{split}
\end{equation}
If $X_2 + X_3 = 0$, then
$$
  \rho(a_1(1)) = (Y_1 + Y_2)(Y_1 + Y_3)(Y_2 + Y_3)
$$
and $a_1(1)$ is perfect because $y_i + y_j$ is perfect.
 If $X_2 + X_3 \neq 0$ then $Y_2 + Y_3 = X_0$,
and therefore
\begin{equation}
 \label{rhoa2}
 \begin{split}
 \rho(a_1(1)) =
   & X_2 + X_3 + (Y_1 + Y_2)(Y_1 + Y_3) = \\
   & (Y_1 + Y_2 + X_3)(Y_1 + Y_3 + X_2) = \\
   & \rho((y_1 + y_2 + x_3)(y_1 + y_3 + x_2)).
 \end{split}
\end{equation}
Thus, $a_1(1)$ is perfect by Corollary \ref{second_perfect}.

  \underline{(b) Case $\Phi^+\rho(y_2 + y_3) = 0$}. Here, we have
 $Y_2^1 + Y_3^1$=0. Then $\Phi^+\rho$ is of the from
$$
 \begin{array}{cccccc}
   X^1_1 & Y^1_1 & X^1_0  & 0 & 0 &   \\
              &  & 0      &   &   &   \\
              &  & 0      & \\
 \end{array}
$$

For any indecomposable representation $\Phi^+\rho$, it can be true
only for
\begin{equation}
  \label{only_this_rho}
   \Phi^+\rho = \left \{
  \begin{array}{l}
     \rho_{x_0} \text{ or } \\
     \rho_{y_1} \text{ or } \\
     \rho_{x_1},
  \end{array}
    \right .
     \hspace{2mm} \text{ i.e., } \quad
   \rho = \left \{
   \begin{array}{l}
      \Phi^{-1}\rho_{x_0} \text{ or }  \\
      \Phi^{-1}\rho_{y_1} \text{ or }  \\
      \Phi^{-1}\rho_{x_1},
   \end{array}
    \right .
\end{equation}
see Table \ref{RepOrder12}. For these 3 representations,
 $\rho(x_2 + x_3)$ = 0, and therefore $\rho(a_1(1))$=0. \qedsymbol

\subsection{The perfectness of $b_i(1)$}
 \label{bi_perfect}

To prove the {\it perfectness} of
 $b_1(1) \supseteq a_1(1)$, it suffices to consider only the
representation $\rho $, for which $\rho(a_1(1)) = 0$, i.e., the
case \S\ref{ai_perfect} (b). In this case, $\rho$ satisfies eq.
(\ref{only_this_rho}). Since $b_i(1) = a_i(1) + x_1(2)$, we get
$$
   \rho(b_1(1)) = \rho(a_1(1) + y_3{y}_2) = Y_3{Y}_2.
$$

 \underline{Case $\rho = \Phi^{-1}\rho_{x_0}$} (similarly, for
$\Phi^{-1}\rho_{y_1}$). In this case, $\rho$ is the direct sum
$$
\vspace{3mm}
\begin{array}{rccccccccccccl}
     0 & Y_1 & Y_1  & 0   & 0 &  \hspace{1cm} &
     0 & 0   & Y_2  & Y_2 & 0 &  \hspace{1cm} \\
             &      & 0   &   &         &
                     \hspace{1cm}\oplus\hspace{1cm} &
             &      & Y_2 &   &         &  \\
             &      & 0   &   &         &  \hspace{1cm} &
             &      & 0   & \\
\end{array}
\vspace{3mm}
 $$
where $\dim{Y}_1 = \dim{Y}_2 = \dim{Y}_3 = 1$, which contradicts
the indecomposability of $\rho$.

 \underline{Case $\rho = \Phi^{-1}\rho_{x_1}$}. Then
 $\rho(b_1(1)) = Y_3{Y}_2  = \mathbb{R}^1$, i.e.,
 $\rho(b_1(1)) = X_0$ and $b_1(1)$ is perfect. \qedsymbol
 \vspace{2mm}

\subsection{The perfectness of $c_i(1)$}
 \label{ci_perfect}
Since $c_1(1) \supseteq b_1(1)$, we see as above, that it suffices
to consider only the case
 $\rho(b_1(1)) = \rho(a_1(1)) = 0$, i.e., again, $\rho$ satisfies
eq. (\ref{only_this_rho}). Since $c_i(1) = a_i(1) + y_1(2)$, it
follows that
$$
 \rho(c_1(1)) = \rho(a_i(1) + y_2(x_1 + y_3) +
  y_3(x_1 + y_2)) =  Y_2(X_1 + Y_3) + Y_3(X_1 + Y_2).
$$
  For $\rho = \Phi^{-1}\rho_{x_0}$ and
      $\rho = \Phi^{-1}\rho_{x_1}$, we have $X_1 = 0$. Then
$$
  Y_2(X_1 + Y_3) + Y_3(X_1 + Y_2) = Y_2{Y}_3
$$
  and we are back to the argumentation of
  \S\ref{bi_perfect}. Let us consider $\rho = \Phi^{-1}\rho_{y_1}$.

  Here, $X_1 + Y_2 = X_0$, otherwise
  $\dim(X_1 + Y_2) = 1$ and $X_1 = Y_2$ and the
  representation $\Phi^{-1}\rho_{y_1}$ is decomposable. Since
  $X_1 + Y_2 = X_0$, we get
  $Y_3(X_1 + Y_2) = Y_3$; by analogy we have
  $Y_2(X_1 + Y_3) = Y_2$, and hence
  $\rho(c_1(1)) = Y_3 + Y_2$. Since $y_2 + y_3$
  is perfect, it follows that $c_1(1)$ is also a perfect element.

  The perfectness of elements
  $a_i(1), b_i(1), c_i(1)$ for $i = 2,3$ is similarly proved. \qedsymbol

\subsection{The perfectness of $a_i(n), b_i(n), c_i(n)$}
  \label{perfectness_abc_n}

Now, we will prove the perfectness of the elements $a_i(1),
b_i(1), c_i(1)$.

{\bf Proposition \ref{perfect_abc}}
  {\it 1) For every element $v_i(n) = a_i(n), b_i(n), c_i(n)$,
  we have
 \begin{equation}
   \label{vn_n_back}
    \rho(v_i(n)) = \sum_\text{$k=1,2,3$}\varphi_k\Phi^+\rho(v_i(n-1)).
 \end{equation}

  2) The elements $a_i(n), b_i(n), c_i(n)$ are perfect
   for every $n \geq 1$ and $i=1,2,3$.}

 \PerfProof  1) follows from the definition of
 the cumulative polynomials from \S\ref{def_cumul},
 Proposition \ref{cumul_polyn} and the
 definition of the perfect elements $v_i(n)$ from \S\ref{repres_lat}.
\qedsymbol

 2) For $n=1$, it was only just proved in
 \S\S\hspace{1mm}\ref{ai_perfect}, \ref{bi_perfect}, \ref{ci_perfect}.
 We will prove, for example, the
 perfectness of $a_1(n)$. From eq. (\ref{vn_n_back}) we have
 \begin{equation}  \label{an_n}
    \rho(a_1(n)) = \sum_\text{$i=1,2,3$}\varphi_i\Phi^+\rho(a_1(n-1)).
 \end{equation}
By (\ref{an_n}), Proposition \ref{new_perfect} and relation
(\ref{f12}) we see that the perfectness of $a_1(n)$ follows from
the perfectness of $a_1(n-1)$. \qedsymbol

\subsection{The distributivity of $H^+(n)$}
  \label{sect_distrib_Hn}

Our proof of the following proposition is founded on the
well-known result of B. Jonsson \cite{Jo55} (see Proposition
\ref{jonsson}):

{\it A modular lattice generated by the chains $s_1, s_2, s_3$
  is distributive if and only if every sublattice
  $[v_1, v_2, v_3]$, where $v_i$ is
  some element from $s_i$, is distributive}.\vspace{2mm}

{\bf Proposition \ref{distrib_Hn}}
 {\it $H^+(n)$ is distributive sublattices for every $n \geq 0$}.

\PerfProof
  By Jonsson's criterion it suffices to prove that
\begin{equation}  \label{distributive_1}
   v_i{v}_k + v_j{v}_k = (v_i + v_j)v_k
\end{equation}
for distinct $i,j,k$. In our case, every two generators from
different chains compose the same sum
\begin{equation}
  \label{distributive_2}
    v_i + v_j = I_n,  \text{ where } i \neq j.
\end{equation}
For $n \geq 1$, we have
\begin{equation}  \label{distributive_sum}
  I_n = \sum_{i=1,2,3}x_i(n) + \sum_{i=1,2,3}y_i(n+1).
\end{equation}
For $n = 0$, we have
\begin{equation}
  I_0 = \sum_{i=1,2,3}y_i.
\end{equation}
By (\ref{distributive_1}) and (\ref{distributive_2}) it suffices
to prove
\begin{equation}  \label{distributive_3}
   v_i{v}_k + v_j{v}_k = v_k
\end{equation}
for distinct $i,j,k$.  We will omit index $n$ in the polynomials
$a_i(n), b_i(n), c_i(n)$.

Suppose that
\begin{equation} \label{distributive_4}
     a_i{a}_j + a_k = I_n.
\end{equation}
  Then, $b_i{a}_j + a_k = I_n$ and
  $c_i{a}_j + a_k = I_n$.
  Further, by the modular law (\ref{modular_law}) we have
\begin{equation}
\begin{split}
  & a_i = a_i(a_i{a}_j + a_k) = a_i{a}_j + a_i{a}_k, \\
  & b_i = b_i(a_i{a}_j + a_k) = b_i{a}_j + b_i{a}_k, \\
  & c_i = c_i(a_i{a}_j + a_k) = c_i{a}_j + c_i{a}_k.
\end{split}
\end{equation}
Now, if $a_i \subseteq v_i$, where
     $v_i = b_i$ or $c_i$,  then
$$
  a_i = a_i{a}_j + a_i{a}_k
    \subseteq a_i{v}_j + a_i{v}_k
        \subseteq a_i(v_j + a_i{v}_k) \subseteq a_i,
$$
     i.e.,  $a_i = a_i{v}_j + a_i{v}_k$.
     The same for $b_i, c_i$:
\begin{equation}
 \label{distributive_5}
 \begin{split}
   & a_i = a_i{v}_j + a_i{v}_k, \\
   & b_i = b_i{v}_j + b_i{v}_k, \\
   & c_i = a_i{v}_j + c_i{v}_k.
 \end{split}
\end{equation}
So, (\ref{distributive_5}) and (\ref{distributive_3}) follow from
(\ref{distributive_4}). Therefore, it suffices to prove
(\ref{distributive_4}). For $n = 0$ we have $a_i = y_j + y_k$ and
\begin{equation*}
 I_n \supseteq a_i{a}_j + a_k =
   (y_j + y_k)(y_i + y_k) + y_i + y_j \supseteq
       y_k + y_i + y_j = I_n.
\end{equation*}
   Set
\begin{equation} \label{designation_t}
     t_i(n) = x_i(n) + y_i(n+1).
\end{equation}
Obviously, $a_i(n) = t_j(n) + t_k(n)$. By
(\ref{distributive_sum}), for $n \geq 1$, we have (parameter $n$
is dropped)
\begin{equation*}
  I_n \supseteq a_i{a}_j + a_k =
  (t_j + t_k)(t_i + t_k) + t_i + t_j
  \supseteq t_k + t_i + t_j = I_n.
\end{equation*}

The distributivity of $H^+(n)$ is proven. \qedsymbol

\subsection{The cardinality of $H^+(n)$}
 \label{proofs_27_64}

 Now we will prove that the distributive lattice $H^+(0)$ is a
$27$-element sublattice in $D^{2,2,2}$ and the distributive
lattice
 $H^+(n) (n \geq 1)$ is the $64$-element sublattice in $D^{2,2,2}$.

 \index{characteristic function $\chi_\rho^v$}

 We introduce {\it characteristic function $\chi_\rho^v$} to be
\begin{equation}
  \label{char_fun}
   \chi_\rho^v =
     \begin{cases}
          1  \text{ for } \rho(v)=X_0, \\
          0  \text{ for } \rho(v)= 0.
        \end{cases}
\end{equation}

{\bf Proposition \ref{27_0_distinct_elem}}
  {\it The lattice $H^+(0)$ contains $27$ distinct elements}.

\PerfProof It suffices to enumerate all $27$ elements in the
lattice $H^+(0)$ and find representation which separates each pair
of elements. Table \ref{table_char_fun0} contains values of the
characteristic function (\ref{char_fun}) on $6$ representations
$\rho_{y_i}, \rho_{x_i}$, where $i=1,2,3$ (see Table
\ref{RepOrder12}), and all $27$ elements of $H^+(0)$. All
signatures in Table \ref{table_char_fun0} are distinct.
 \qedsymbol \vspace{2mm}

\begin{table}[h]
 \renewcommand{\arraystretch}{1.35}
 \begin{tabular} {||c|c||c|c|c||c|c|c||c||}
  \hline \hline
    $N$ & Polynomial
            & $\rho_{y_1}$ & $\rho_{y_2}$ & $\rho_{y_3}$
            & $\rho_{x_1}$ & $\rho_{x_2}$ & $\rho_{x_3}$
            & (Zeros, Units)
     \\
  \hline  \hline
    1 & $a_1(0) = y_2 + y_3 $
      & 0 & 1 & 1
      & 0 & 1 & 1
      & (2,4)   \\
  \hline
    2 & $a_2(0) = y_1 + y_3 $
      & 1 & 0 & 1
      & 1 & 0 & 1
      & (2,4)   \\
  \hline
    3 & $a_3(0) = y_2 + y_1 $
      & 1 & 1 & 0
      & 1 & 1 & 0
      & (2,4)   \\
  \hline
    4 & $b_1(0) = x_1 + y_2 + y_3 $
      & 0 & 1 & 1
      & 1 & 1 & 1
      & (1,5)   \\
  \hline
    5 & $b_2(0) = x_2 + y_1 + y_3 $
      & 1 & 0 & 1
      & 1 & 1 & 1
      & (1,5)   \\
  \hline
    6 & $b_3(0) = x_3 + y_2 + y_1 $
      & 1 & 1 & 0
      & 1 & 1 & 1
      & (1,5)   \\
  \hline
    7 & $a_1(0)a_2(0)$
      & 0 & 0 & 1
      & 0 & 0 & 1
      & (4,2)   \\
  \hline
    8 & $a_1(0)a_3(0)$
      & 0 & 1 & 0
      & 0 & 1 & 0
      & (4,2)   \\
  \hline
    9 & $a_2(0)a_3(0)$
      & 1 & 0 & 0
      & 1 & 0 & 0
      & (4,2)   \\
  \hline
   10 & $b_1(0)b_2(0)$
      & 0 & 0 & 1
      & 1 & 1 & 1
      & (2,4)   \\
  \hline
   11 & $b_1(0)b_3(0)$
      & 0 & 1 & 0
      & 1 & 1 & 1
      & (2,4)   \\
  \hline
   12 & $b_2(0)b_3(0)$
      & 1 & 0 & 0
      & 1 & 1 & 1
      & (2,4)   \\
  \hline
   13 & $a_1(0)b_2(0)$
      & 0 & 0 & 1
      & 0 & 1 & 1
      & (3,3)   \\
  \hline
   14 & $a_1(0)b_3(0)$
      & 0 & 1 & 0
      & 0 & 1 & 1
      & (3,3)   \\
  \hline
   15 & $a_2(0)b_1(0)$
      & 0 & 0 & 1
      & 1 & 0 & 1
      & (3,3)   \\
  \hline
   16 & $a_2(0)b_3(0)$
      & 1 & 0 & 0
      & 1 & 0 & 1
      & (3,3)   \\
  \hline
   17 & $a_3(0)b_1(0)$
      & 0 & 1 & 0
      & 1 & 1 & 0
      & (3,3)   \\
  \hline
   18 & $a_3(0)b_2(0)$
      & 1 & 0 & 0
      & 1 & 1 & 0
      & (3,3)   \\
  \hline
   19 & $a_1(0)a_2(0)b_3(0)$
      & 0 & 0 & 0
      & 0 & 0 & 1
      & (5,1)   \\
  \hline
   20 & $a_1(0)b_2(0)a_3(0)$
      & 0 & 0 & 0
      & 0 & 1 & 0
      & (5,1)   \\
  \hline
   21 & $b_1(0)a_2(0)a_3(0)$
      & 0 & 0 & 0
      & 1 & 0 & 0
      & (5,1)   \\
  \hline
   22 & $a_1(0)b_2(0)b_3(0)$
      & 0 & 0 & 0
      & 0 & 1 & 1
      & (4,2)   \\
  \hline
   23 & $b_1(0)a_2(0)b_3(0)$
      & 0 & 0 & 0
      & 1 & 0 & 1
      & (4,2)   \\
  \hline
   24 & $b_1(0)b_2(0)a_3(0)$
      & 0 & 0 & 0
      & 1 & 1 & 0
      & (4,2)   \\
  \hline
   25 & $a_1(0)a_2(0)a_3(0)$
      & 0 & 0 & 0
      & 0 & 0 & 0
      & (6,0)   \\
  \hline
   26 & $b_1(0)b_2(0)b_3(0)$
      & 0 & 0 & 0
      & 1 & 1 & 1
      & (3,3)   \\
  \hline
   27 & $I_0 = y_1 + y_2 + y_3$
      & 1 & 1 & 1
      & 1 & 1 & 1
      & (6,0)   \\
  \hline \hline
  \end{tabular}
  \vspace{2mm}
  \caption{\hspace{3mm}The characteristic function
   on the $H^+(0)$}
  \label{table_char_fun0}
 \small
 \underline{Notes}: (a) Lines 7 -- 27 are obtained
 by multiplying the appropriate Lines 1 -- 6, since
\begin{equation} \label{chi_products}
 \chi_\rho^{v_1v_2} = \chi_\rho^{v_1}\chi_\rho^{v_1}.
\end{equation}
 (b) We call a sequence of ``zeros'' and ``units''
  in the set $\rho_{y_1},\rho_{y_2},\dots\rho_{x_3}$
  the {\it signature} of the corresponding polynomial.
 \normalsize
  \vspace{7.5mm} 
 \end{table}


{\bf Proposition \ref{64_1_distinct_elem}}
  {\it The lattice $H^+(1)$ contains $64$ distinct elements}.

As above, for the sublattice $H^+(0)$, it suffices to enumerate
all $64$ elements in the lattice $H^+(1)$ and find representation
which separates each pair of elements. Tables
 \ref{table_char_fun1_gen}, \ref{table_char_fun1} contain values of
the characteristic function (\ref{char_fun}) on $9$
representations $\rho_{x_i}, \Phi^{-1}\rho_{y_i},
\Phi^{-1}\rho_{x_i}$, where $i=1,2,3,$ (see Table
 \ref{RepOrder12}) and all $64$ elements of $H^+(1)$.

By (\ref{def_gen_abc}) the generators $a_i, b_i, c_i, i=1,2,3$
from the sublattice $H^+(1)$ are:
\begin{equation} \label{generators_abc}
\begin{split}
   & a_1(1) = x_2 + x_3 +
       y_3(x_2 + y_1) + y_1(x_2 + y_3) + y_1(x_3 + y_2) + y_2(x_3 + y_1), \\
   & a_2(1) = x_1 + x_3 +
       y_3(x_1 + y_2) + y_2(x_1 + y_3) + y_1(x_3 + y_2) + y_2(x_3 + y_1),  \\
   & a_3(1) = x_2 + x_1 +
       y_3(x_2 + y_1) + y_1(x_2 + y_3) + y_3(x_1 + y_2) + y_2(x_1 + y_3),
\end{split}
\end{equation}
\begin{equation}  \label{generators_abc_1}
\begin{split}
   & b_1(1) = a_1(1) + y_2{y}_3,  \hspace{7mm}
     c_1(1) = a_1(1) + y_2(x_1 + y_3) + y_3(x_1 + y_2), \\
   & b_2(1) = a_2(1) + y_1{y}_3,  \hspace{7mm}
     c_2(1) = a_2(1) + y_1(x_2 + y_3) + y_3(x_2 + y_1), \\
   & b_3(1) = a_3(1) + y_1{y}_2,  \hspace{7mm}
     c_3(1) = a_3(1) + y_1(x_3 + y_2) + y_2(x_3 + y_1),
\end{split}
\end{equation}
see Table \ref{table_char_fun1_gen}. By (\ref{chi_products}) Lines
10 -- 64 of Table \ref{table_char_fun1} are obtained by
multiplying of appropriate Lines 1 -- 9 of Table
 \ref{table_char_fun1_gen}.

 No two signatures in
 Tables \ref{table_char_fun1_gen} coincide,
 which was to be proved.\qedsymbol

\small
\begin{table}[h]
 \renewcommand{\arraystretch}{1.35}
 \begin{tabular} {||c|c||c|c|c||c|c|c||c|c|c||c||}
  \hline \hline
      $N$ & Poly- & & & & & & & & & &  (Zeros, \cr
          & nomial
          & $\rho_{x_1}$ & $\rho_{x_2}$ & $\rho_{x_3}$
          & $\Phi^{-1}\rho_{y_1}$ & $\Phi^{-1}\rho_{y_2}$ & $\Phi^{-1}\rho_{y_3}$
          & $\Phi^{-1}\rho_{x_1}$ & $\Phi^{-1}\rho_{x_2}$ & $\Phi^{-1}\rho_{x_3}$
          & Units)
     \\
  \hline  \hline
    1 & $a_1(1) $
      & 0 & 1 & 1
      & 0 & 1 & 1
      & 0 & 1 & 1
      & (3,6)   \\
  \hline
    2 & $a_2(1) $
      & 1 & 0 & 1
      & 1 & 0 & 1
      & 1 & 0 & 1
      & (3,6)   \\
  \hline
    3 & $a_3(1) $
      & 1 & 1 & 0
      & 1 & 1 & 0
      & 1 & 1 & 0
      & (3,6)   \\
  \hline
    4 & $b_1(1) $
      & 0 & 1 & 1
      & 0 & 1 & 1
      & 1 & 1 & 1
      & (2,7)   \\
  \hline
    5 & $b_2(1) $
      & 1 & 0 & 1
      & 1 & 0 & 1
      & 1 & 1 & 1
      & (2,7)   \\
  \hline
    6 & $b_3(1) $
      & 1 & 1 & 0
      & 1 & 1 & 0
      & 1 & 1 & 1
      & (2,7)   \\
  \hline
    7 & $c_1(1) $
      & 0 & 1 & 1
      & 1 & 1 & 1
      & 1 & 1 & 1
      & (1,8)   \\
  \hline
    8 & $c_2(1) $
      & 1 & 0 & 1
      & 1 & 1 & 1
      & 1 & 1 & 1
      & (1,8)   \\
  \hline
    9 & $c_3(1) $
      & 1 & 1 & 0
      & 1 & 1 & 1
      & 1 & 1 & 1
      & (1,8)   \\
  \hline \hline
  \end{tabular}
  \vspace{2mm}
  \caption{\hspace{3mm}The characteristic function on the $H^+(1)$.
  The generators $a_i(1), b_i(1), c_i(1)$}
  \label{table_char_fun1_gen}
 \end{table}
\normalsize \tiny
\begin{table}[h]
 \begin{tabular} {||c|c||c|c|c||c|c|c||c|c|c||c||}
  \hline \hline
     $N$  & Poly- & & & & & & & & & &  (Zeros, \cr
          & nomial
          & $\rho_{x_1}$ & $\rho_{x_2}$ & $\rho_{x_3}$
          & $\Phi^{-1}\rho_{y_1}$ & $\Phi^{-1}\rho_{y_2}$ & $\Phi^{-1}\rho_{y_3}$
          & $\Phi^{-1}\rho_{x_1}$ & $\Phi^{-1}\rho_{x_2}$ & $\Phi^{-1}\rho_{x_3}$
          & Units)
     \\
  \hline  \hline
   10 & $a_1(1)a_2(1) $
      & 0 & 0 & 1
      & 0 & 0 & 1
      & 0 & 0 & 1
      & (6,3)   \\
  \hline
   11 & $a_1(1)a_3(1) $
      & 0 & 1 & 0
      & 0 & 1 & 0
      & 0 & 1 & 0
      & (6,3)   \\
  \hline
   12 & $a_2(1)a_3(1) $
      & 1 & 0 & 0
      & 1 & 0 & 0
      & 1 & 0 & 0
      & (6,3)   \\
  \hline
   13 & $b_1(1)b_2(1) $
      & 0 & 0 & 1
      & 0 & 0 & 1
      & 1 & 1 & 1
      & (4,5)   \\
  \hline
   14 & $b_1(1)b_3(1) $
      & 0 & 1 & 0
      & 0 & 1 & 0
      & 1 & 1 & 1
      & (4,5)   \\
  \hline
   15 & $b_2(1)b_3(1) $
      & 1 & 0 & 0
      & 1 & 0 & 0
      & 1 & 1 & 1
      & (4,5)   \\
  \hline
   16 & $c_1(1)c_2(1) $
      & 0 & 0 & 1
      & 1 & 1 & 1
      & 1 & 1 & 1
      & (2,7)   \\
  \hline
   17 & $c_1(1)c_3(1) $
      & 0 & 1 & 0
      & 1 & 1 & 1
      & 1 & 1 & 1
      & (2,7)   \\
  \hline
   18 & $c_2(1)c_3(1) $
      & 1 & 0 & 0
      & 1 & 1 & 1
      & 1 & 1 & 1
      & (2,7)   \\
  \hline
   19 & $a_1(1)b_2(1) $
      & 0 & 0 & 1
      & 0 & 0 & 1
      & 0 & 1 & 1
      & (5,4)   \\
  \hline
   20 & $a_1(1)b_3(1) $
      & 0 & 1 & 0
      & 0 & 1 & 0
      & 0 & 1 & 1
      & (5,4)   \\
  \hline
   21 & $a_2(1)b_1(1) $
      & 0 & 0 & 1
      & 0 & 0 & 1
      & 1 & 0 & 1
      & (5,4)   \\
  \hline
   22 & $a_2(1)b_3(1) $
      & 1 & 0 & 0
      & 1 & 0 & 0
      & 1 & 0 & 1
      & (5,4)   \\
  \hline
   23 & $a_3(1)b_1(1) $
      & 0 & 1 & 0
      & 0 & 1 & 0
      & 1 & 1 & 0
      & (5,4)   \\
  \hline
   24 & $a_3(1)b_2(1) $
      & 1 & 0 & 0
      & 1 & 0 & 0
      & 1 & 1 & 0
      & (5,4)   \\
  \hline
   25 & $a_1(1)c_2(1) $
      & 0 & 0 & 1
      & 0 & 1 & 1
      & 0 & 1 & 1
      & (4,5)   \\
  \hline
   26 & $a_1(1)c_3(1) $
      & 0 & 1 & 0
      & 0 & 1 & 1
      & 0 & 1 & 1
      & (4,5)   \\
  \hline
   27 & $a_2(1)c_1(1) $
      & 0 & 0 & 1
      & 1 & 0 & 1
      & 1 & 0 & 1
      & (4,5)   \\
  \hline
   28 & $a_2(1)c_3(1) $
      & 1 & 0 & 0
      & 1 & 0 & 1
      & 1 & 0 & 1
      & (4,5)   \\
  \hline
   29 & $a_3(1)c_1(1) $
      & 0 & 1 & 0
      & 1 & 1 & 0
      & 1 & 1 & 0
      & (4,5)   \\
  \hline
   30 & $a_3(1)c_2(1) $
      & 1 & 0 & 0
      & 1 & 1 & 0
      & 1 & 1 & 0
      & (4,5)   \\
  \hline
   31 & $b_1(1)c_2(1) $
      & 0 & 0 & 1
      & 0 & 1 & 1
      & 1 & 1 & 1
      & (3,6)   \\
  \hline
   32 & $b_1(1)c_3(1) $
      & 0 & 1 & 0
      & 0 & 1 & 1
      & 1 & 1 & 1
      & (3,6)   \\
  \hline
   33 & $b_2(1)c_1(1) $
      & 0 & 0 & 1
      & 1 & 0 & 1
      & 1 & 1 & 1
      & (3,6)   \\
  \hline
   34 & $b_2(1)c_3(1) $
      & 1 & 0 & 0
      & 1 & 0 & 1
      & 1 & 1 & 1
      & (3,6)   \\
  \hline
   35 & $b_3(1)c_1(1) $
      & 0 & 1 & 0
      & 1 & 1 & 0
      & 1 & 1 & 1
      & (3,6)   \\
  \hline
   36 & $b_3(1)c_2(1) $
      & 1 & 0 & 0
      & 1 & 1 & 0
      & 1 & 1 & 1
      & (3,6)   \\
  \hline  \hline
   37 & $a_1(1)a_2(1)b_3(1) $
      & 0 & 0 & 0
      & 0 & 0 & 0
      & 0 & 0 & 1
      & (8,1)   \\
  \hline
   38 & $a_1(1)a_3(1)b_2(1) $
      & 0 & 0 & 0
      & 0 & 0 & 0
      & 0 & 1 & 0
      & (8,1)   \\
  \hline
   39 & $a_2(1)a_3(1)b_1(1) $
      & 0 & 0 & 0
      & 0 & 0 & 0
      & 1 & 0 & 0
      & (8,1)   \\
  \hline
   40 & $a_1(1)a_2(1)c_3(1) $
      & 0 & 0 & 0
      & 0 & 0 & 1
      & 0 & 0 & 1
      & (7,2)   \\
  \hline
   41 & $a_1(1)a_3(1)c_2(1) $
      & 0 & 0 & 0
      & 0 & 1 & 0
      & 0 & 1 & 0
      & (7,2)   \\
  \hline
   42 & $a_2(1)a_3(1)c_1(1) $
      & 0 & 0 & 0
      & 1 & 0 & 0
      & 1 & 0 & 0
      & (7,2)   \\
  \hline
   43 & $b_1(1)b_2(1)a_3(1) $
      & 0 & 0 & 0
      & 0 & 0 & 0
      & 1 & 1 & 0
      & (7,2)   \\
  \hline
   44 & $b_1(1)b_3(1)a_2(1) $
      & 0 & 0 & 0
      & 0 & 0 & 0
      & 1 & 0 & 1
      & (7,2)   \\
  \hline
   45 & $b_2(1)b_3(1)a_1(1) $
      & 0 & 0 & 0
      & 0 & 0 & 0
      & 0 & 1 & 1
      & (7,2)   \\
  \hline
   46 & $b_1(1)b_2(1)c_3(1) $
      & 0 & 0 & 0
      & 0 & 0 & 1
      & 1 & 1 & 1
      & (5,4)   \\
  \hline
   47 & $b_1(1)b_3(1)c_2(1) $
      & 0 & 0 & 0
      & 0 & 1 & 0
      & 1 & 1 & 1
      & (5,4)   \\
  \hline
   48 & $b_2(1)b_3(1)c_1(1) $
      & 0 & 0 & 0
      & 1 & 0 & 0
      & 1 & 1 & 1
      & (5,4)   \\
  \hline
   49 & $c_1(1)c_2(1)a_3(1) $
      & 0 & 0 & 0
      & 1 & 1 & 0
      & 1 & 1 & 0
      & (5,4)   \\
  \hline
   50 & $c_1(1)c_3(1)a_2(1) $
      & 0 & 0 & 0
      & 1 & 0 & 1
      & 1 & 0 & 1
      & (5,4)   \\
  \hline
   51 & $c_2(1)c_3(1)a_1(1) $
      & 0 & 0 & 0
      & 0 & 1 & 1
      & 0 & 1 & 1
      & (5,4)   \\
  \hline
   52 & $c_1(1)c_2(1)b_3(1) $
      & 0 & 0 & 0
      & 1 & 1 & 0
      & 1 & 1 & 1
      & (4,5)   \\
  \hline
   53 & $c_1(1)c_3(1)b_2(1) $
      & 0 & 0 & 0
      & 1 & 0 & 1
      & 1 & 1 & 1
      & (4,5)   \\
  \hline
   54 & $c_2(1)c_3(1)b_1(1) $
      & 0 & 0 & 0
      & 0 & 1 & 1
      & 1 & 1 & 1
      & (4,5)   \\
  \hline
   55 & $a_1(1)b_2(1)c_3(1) $
      & 0 & 0 & 1
      & 0 & 0 & 1
      & 0 & 1 & 1
      & (7,2)   \\
  \hline
   56 & $a_1(1)b_3(1)c_2(1) $
      & 0 & 0 & 0
      & 0 & 1 & 0
      & 0 & 1 & 1
      & (7,2)   \\
  \hline
   57 & $a_2(1)b_1(1)c_3(1) $
      & 0 & 0 & 0
      & 0 & 0 & 1
      & 1 & 0 & 1
      & (7,2)   \\
  \hline
   58 & $a_2(1)b_3(1)c_1(1) $
      & 0 & 0 & 0
      & 1 & 0 & 0
      & 1 & 0 & 1
      & (7,2)   \\
  \hline
   59 & $a_3(1)b_1(1)c_2(1) $
      & 0 & 0 & 0
      & 0 & 1 & 0
      & 1 & 1 & 0
      & (3,6)   \\
  \hline
   60 & $a_3(1)b_2(1)c_1(1) $
      & 0 & 0 & 0
      & 1 & 0 & 0
      & 1 & 1 & 0
      & (3,6)   \\
  \hline
   61 & $a_1(1)a_2(1)a_3(1) $
      & 0 & 0 & 0
      & 0 & 0 & 0
      & 0 & 0 & 0
      & (9,0)   \\
  \hline
   62 & $b_1(1)b_2(1)b_3(1) $
      & 0 & 0 & 0
      & 0 & 0 & 0
      & 1 & 1 & 1
      & (6,3)   \\
  \hline
   63 & $c_1(1)c_2(1)c_3(1) $
      & 0 & 0 & 0
      & 1 & 1 & 1
      & 1 & 1 & 1
      & (3,6)   \\
  \hline \hline
   64 & $I_1 $
      & 1 & 1 & 1
      & 1 & 1 & 1
      & 1 & 1 & 1
      & (1,8)   \\
  \hline \hline
  \end{tabular}
  \vspace{2mm}
  \caption{\hspace{3mm}The characteristic function on the $H^+(1)$.
        The elements 10--64}
  \label{table_char_fun1}
 \end{table}
\normalsize

\subsection{A relation between the upper and the lower generators of $H^+(n)$}
 \label{sect_upper_lower}

 ~\vspace{5mm}\\
 {\bf Proposition \ref{two_set_gen}}
  {\it Table \ref{mult_add_forms} from \S\ref{another_gen_eq} gives
  relations between the upper and the lower generators.
  These relations are true$\mod\theta$}. \vspace{2mm}

\PerfProof
  For definition of $a_i(n), b_i(n), c_i(n)$
  (resp. $p_i(n), q_i(n), s_i(n)$), see
  (\ref{def_gen_abc}) (resp. (\ref{another_gen_eq})).

1) By (\ref{x_cap_a_0}) we have
           $x_0(n+2) \subseteq a_i(n)$ and
\begin{equation*}
\begin{split}
         &  s_j(n) + s_k(n) =  \\
         &  x_j(n) + y_j(n+1) + x_k(n) + y_k(n+1) + x_0(n+2) = \\
         &  a_i(n) + x_0(n+2) = a_i(n).
         \vspace{2mm}
\end{split}
\end{equation*}

2) Here, we have
\begin{equation*}
\begin{split}
     & s_j(n) + s_k(n) + p_i(n) = \\
     & a_i(n) + p_i(n) =    \\
     & x_j(n) + y_j(n+1) + x_k(n) + y_k(n+1) + x_0(n+2) + x_i(n+1) = \\
     & a_i(n) + x_i(n+1) = b_i(n).
     \vspace{2mm}
\end{split}
\end{equation*}

3) Further,
\begin{equation*}
\begin{split}
  & s_j(n) + s_k(n) + q_i(n) = a_i(n) + q_i(n) =   \\
  & x_j(n) + y_j(n+1) + x_k(n) + y_k(n+1) + x_0(n+2) + y_i(n+1) = \\
  & a_i(n) + y_i(n+1) = c_i(n).
         \vspace{2mm}
\end{split}
\end{equation*}

4) Consider the intersection $a_k(n)a_j(n)$:
\begin{equation*}
\begin{split}
    & a_k(n)a_j(n) = \\
    & [x_i(n) + y_i(n+1) + x_j(n) + y_j(n+1)]
      [x_i(n) + y_i(n+1) + x_k(n) + y_k(n+1)] = \\
    & x_i(n) + y_i(n+1) + (x_j(n) + y_j(n+1))a_j(n).
\end{split}
\end{equation*}
       Since  $x_j(n) + y_j(n+1) \subseteq a_i(n)a_k(n)$,
       we see that
\begin{equation*}
\begin{split}
    &   a_k(n)a_j(n) \subseteq
        x_i(n) + y_i(n+1) + a_i(n)a_j(n)a_k(n) = \\
    &   x_i(n) + y_i(n+1) + x_0(n+2) = s_i(n).
\end{split}
\end{equation*}
On the other hand,
$$
   s_i(n) = x_i(n) + y_i(n+1) + x_0(n+2)  \subseteq  a_j(n)a_k(n).
$$
Therefore,
$$
  s_i(n) = a_j(n)a_k(n). \vspace{2mm}
$$

5) Since $x_i(n+1) \subseteq a_j(n)$ and
       $x_j(n+1) \subseteq a_i(n)$, by heading 4) we have
\begin{equation*}
\begin{split}
     &  b_i(n)b_j(n) = (a_i(n) + x_i(n+1))(a_j(n) + x_j(n+1)) = \\
     & x_i(n+1) + x_j(n+1) + a_i(n)a_j(n) = \\
     & x_i(n+1) + x_j(n+1) + a_i(n)a_j(n) +
       x_0(n+2) = p_i(n) + p_j(n) + s_k(n).
         \vspace{2mm}
\end{split}
\end{equation*}

6) Since
     $y_i(n+1) \subseteq a_j(n)$ and
     $y_j(n+1) \subseteq a_i(n)$, by heading 4) we get
\begin{equation*}
\begin{split}
     &  c_i(n)c_j(n) =
       (a_i(n) + y_i(n+1))(a_j(n) + y_j(n+1)) = \\
     & y_i(n+1) + y_j(n+1) + a_i(n)a_j(n) = \\
     & y_i(n+1) + y_j(n+1) + a_i(n)a_j(n) + x_0(n+2) =
       q_i(n) + q_j(n) + s_k(n).
         \vspace{2mm}
\end{split}
\end{equation*}

7) Since $x_i(n+1) \subseteq b_j(n)$, by heading 4) we see that
\begin{equation*}
\begin{split}
     & a_i(n)b_j(n) =
       a_i(n)(a_j(n) + x_j(n+1)) = x_j(n+1) + a_i(n)a_j(n) = \\
     & x_j(n+1) + a_i(n)a_j(n) + x_0(n+2) =
       p_j(n) + s_k(n).
\end{split}
\end{equation*}

8) Since $y_i(n+1) \subseteq c_j(n)$, by heading 4) we have
\begin{equation*}
\begin{split}
     &  a_i(n)c_j(n) = a_i(n)(a_j(n) + y_j(n+1)) =
        y_j(n+1) + a_i(n)a_j(n) = \\
     &  y_j(n+1) + a_i(n)a_j(n) + x_0(n+2) =
        q_j(n) + s_k(n).
\end{split}
\end{equation*}

9) Since $x_i(n+1) \subseteq c_j(n)$ and
 $y_j(n+1) \subseteq b_i(n)$, by heading 4) we have
\begin{equation*}
\begin{split}
     & b_i(n)c_j(n) = (a_i(n) + x_i(n+1))(a_j(n) + y_j(n+1)) = \\
     & x_i(n+1) + y_j(n+1) + a_i(n)a_j(n) =  \\
     & x_i(n+1) + y_j(n+1) + a_i(n)a_j(n) + x_0(n+2) =
       p_j(n) + q_j(n) + s_k(n).
\end{split}
\end{equation*}

10) Since $x_k(n+1) \subseteq a_i(n)a_j(n)$, we see that
\begin{equation*}
\begin{split}
   & a_i(n)a_j(n)b_k(n) =
     a_i(n)a_j(n)(a_k(n) + x_k(n+1)) = \\
   & x_k(n+1) + a_i(n)a_j(n)a_k(n) = x_k(n+1) + x_0(n+2) =
     p_k(n).
\end{split}
\end{equation*}

11) Since $y_k(n+1) \subseteq a_i(n)a_j(n)$, we have
\begin{equation*}
\begin{split}
  & a_i(n)a_j(n)c_k(n) =
    a_i(n)a_j(n)(a_k(n) + y_k(n+1)) = \\
  & y_k(n+1) + a_i(n)a_j(n)a_k(n) = y_k(n+1) + x_0(n+2) =
    q_k(n).
\end{split}
\end{equation*}

12) Here we have
\begin{equation*}
\begin{split}
   & b_i(n)b_j(n)a_k(n) =
     (a_i(n) + x_i(n+1))(a_j(n) + x_j(n+1))a_k(n) = \\
   & x_i(n+1) + x_j(n+1) + a_i(n)a_j(n)a_k(n) = \\
   & x_i(n+1) + x_j(n+1) + x_0(n+2) = p_i(n) + p_j(n).
\end{split}
\end{equation*}

13) Further,
\begin{equation*}
\begin{split}
  & b_i(n)b_j(n)c_k(n) =
    (a_i(n) + x_i(n+1))(a_j(n) + x_j(n+1))(a_k(n) + y_k(n+1)) = \\
  & x_i(n+1) + x_j(n+1) + y_k(n+1) + a_i(n)a_j(n)a_k(n) = \\
  & x_i(n+1) + x_j(n+1) + y_k(n+1) + x_0(n+2) = p_i(n) + p_j(n) + q_k(n).
\end{split}
\end{equation*}

14) As above,
\begin{equation*}
\begin{split}
  & c_i(n)c_j(n)a_k(n) =
   (a_i(n) + y_i(n+1))(a_j(n) + y_j(n+1))a_k(n) = \\
  & y_i(n+1) + y_j(n+1) + a_i(n)a_j(n)a_k(n) =
   y_i(n+1) + y_j(n+1) + x_0(n+2) = \\
  & q_i(n) + q_j(n).
\end{split}
\end{equation*}

15) Similarly, we have
\begin{equation*}
\begin{split}
  & c_i(n)c_j(n)b_k(n) =
    (a_i(n) + y_i(n+1))(a_j(n) + y_j(n+1))(a_k(n) + x_k(n+1)) = \\
  & y_i(n+1) + y_j(n+1) + x_k(n+1) +
    a_i(n)a_j(n)a_k(n) = \\
  & y_i(n+1) + y_j(n+1) + x_k(n+1) + x_0(n+2) =
    q_i(n) + q_j(n) + p_k(n).
\end{split}
\end{equation*}

16) Similarly,
\begin{equation*}
\begin{split}
   & a_i(n)b_j(n)c_k(n) =
     a_i(n)(a_j(n) + x_j(n+1))(a_k(n) + y_k(n+1)) = \\
   & x_j(n+1) + y_k(n+1) + a_i(n)a_j(n)a_k(n) =  \\
   & x_j(n+1) + y_k(n+1) + x_0(n+2) = p_j(n) + q_k(n).
\end{split}
\end{equation*}

17) It is just (\ref{x_cap_a_0}):
\begin{equation*}
   a_i(n)a_j(n)a_k(n) = x_0(n+2).
\end{equation*}

18) Here we have
\begin{equation*}
\begin{split}
  & b_i(n)b_j(n)b_k(n) =
    (a_i(n) + x_i(n+1))(a_j(n) + x_j(n+1))(a_k(n) + x_k(n+1)) = \\
  & x_i(n+1) + x_j(n+1) + x_k(n+1) + a_i(n)a_j(n)a_k(n) =  \\
  & x_i(n+1) + x_j(n+1) + x_k(n+1) + x_0(n+2) = p_i(n) + p_j(n) + p_k(n).
\end{split}
\end{equation*}

19) Similarly,
\begin{equation*}
\begin{split}
  & c_i(n)c_j(n)c_k(n) =
    (a_i(n) + y_i(n+1))(a_j(n) + y_j(n+1))(a_k(n) + y_k(n+1)) = \\
  & y_i(n+1) + y_j(n+1) + y_k(n+1) + a_i(n)a_j(n)a_k(n) = \\
  & y_i(n+1) + y_j(n+1) + y_k(n+1) + x_0(n+2) =
    q_i(n) + q_j(n) + q_k(n).
\end{split}
\end{equation*}

20) As follows from Line 1) of this table,
\begin{equation*}
     a_i(n) + a_j(n) = (s_j(n) + s_k(n)) + (s_i(n) + s_k(n)) =
           s_i(n) + s_j(n) + s_k(n).
\end{equation*}

Proposition \ref{two_set_gen} is proven. \qedsymbol

\end{appendix}
\chapter*{Conjectures}
 \label{sect_list_of_conj}

1. {\bf Conjecture \ref{conj_1}} {\it
  The lattice $H^+ \bigcup H^-$ contains all perfect
  elements of $D^{2,2,2}\mod\theta$.} \vspace{2mm}

  For $D^4$, the similar conjecture (due to Gelfand-Ponomarev) was proved by
  Dlab and Ringel in \cite{DR80} and Cylke in \cite{Cyl82},
  see \S\ref{biblio_notes}.

  \vspace{3mm}

2. {\bf Conjecture \ref{conj_2}} {\it  The relation
\begin{equation*}
      x_0(n+2) = \bigcap_\text{$i=1,2,3$}a_i(n) \quad \mod\theta,
\end{equation*}
 see (\ref{x_cap_a_0}, Proposition \ref{prop_x_cap}) takes place
 without restriction$\mod\theta$.}

 Here, $x_0(n+2)$ is the
 cumulative polynomial defined in \S\ref{def_cumul}, and the $a_i(n)$
 are perfect elements defined in \S\ref{def_gen_abc}. The elements
$x_0(n+2)$ are also perfect, see Corollary \ref{x0_is_perfect}.

 \vspace{3mm}

3. {\bf Conjecture \ref{conj_3}} {\it The relation
\begin{equation*}
      a_i(n)x_i(n+1) \subseteq  y_i(n+2) \quad \mod\theta,
\end{equation*}
(see (\ref{long_incl_5})), takes place without
 restriction$\mod\theta$.}

 Here, $x_i(n+1), y_i(n+2)$ are the
 cumulative polynomial defined in \S\ref{def_cumul}, and the $a_i(n)$
 are perfect elements defined in \S\ref{def_gen_abc}.
 \vspace{3mm}

4. {\bf Conjecture \ref{conj_4}} {\it  For every admissible
  sequence $\alpha$,
  the elements $e_\alpha$ (resp. $f_{{\alpha}0}$) and
  $\tilde{e}_\alpha$ (resp. $\tilde{f}_{{\alpha}0}$) coincide
  without restriction$\mod\theta$
  (see Proposition \ref{coincidence_GP}).}

  Here, the elements $e_\alpha$ (resp. $f_{{\alpha}0}$) given by
  Table \ref{table_adm_elem_D4}, and the elements
  $\tilde{e}_\alpha$ (resp. $\tilde{f}_{{\alpha}0}$) determined by
  Gelfand and Ponomarev, see \S\ref{coinc_GP_D4}. For small
  admissible sequences, this coincidence is proved in
  Propositions \ref{coincidence_E}, \ref{coincidence_F}.

  \vspace{3mm}

\renewcommand{\appendixname}{}

\printindex

\end{document}